\documentclass{svmult}
\usepackage{amsmath,amssymb,amsfonts,verbatim}
\spnewtheorem{theorem*}{Theorem}[section]{\bf}{\it}
\spnewtheorem{lemma*}[theorem*]{Lemma}{\bf}{\it}
\spnewtheorem{proposition*}[theorem*]{Proposition}{\bf}{\it}
\spnewtheorem{corollary*}[theorem*]{Corollary}{\bf}{\it}
\spnewtheorem{definition*}[theorem*]{Definition}{\bf}{\it}
\spnewtheorem{example*}[theorem*]{Example}{\bf}{\it}
\spnewtheorem{remark*}[theorem*]{Remark}{\bf}{\it}

\def\Ab{\bar{\rm A}}
\def\ad{\text{ad}}

\def\al{\alpha}
\def\Ar{{\rm A}}

\def\be{\beta}
\def\bg{\tilde g}
\def\bi{\tilde\imath}
\def\bj{\tilde\jmath}

\def\bk{\tilde k}
\def\bl{\tilde l}
\def\Br{{\rm E}}

\def\CC{\mathbb{C}}
\def\com{\ts,\hskip-.5pt}

\def\d{\partial}
\def\de{\delta}
\def\De{\Delta}
\def\DA{\mathfrak{T}}
\def\Dwo{\mathrm{Q}}

\def\End{\operatorname{End}\ts}
\def\ep{\varepsilon}
\def\ES{\mathfrak{R}} %% Extended Weyl group

\def\Fr{{\rm F}}
\def\f{\mathfrak{f}}

\def\g{\mathfrak{g}}
\def\ga{\gamma}
\def\Gd{\G^{\tts\prime\tts}}
\def\gd{\g^{\tts\prime\tts}}
\def\ge{\geqslant}
\def\GL{\mathrm{GL}}
\def\gl{\mathfrak{gl}}

\def\h{\mathfrak{t}}
\def\Heist{\mathrm{H}\ts(U)}

\def\I{{\rm I}}
\def\Ib{\tts\overline{\nns\rm I\nns}\tts}

\def\io{\iota}
\def\IV{\rm K^{\tts\prime}}

\def\J{{\rm J}}
\def\Jb{\,\overline{\!\rm J\ns}\ts}
\def\Jp{{\rm J}^{\ts\prime\prime}}
\def\Jm{{\rm J}^{\ts\prime}}
\def\Jmb{\,\overline{\!\rm J\ns}\ts^{\ts\prime}}
\def\Jpb{\,\overline{\!\rm J\ns}\ts^{\ts\prime\prime}}
\def\JV{\rm K}

\def\G{{\mathrm G}}
\def\ka{\kappa}
\def\Ker{\mathrm{Ker}}

\def\la{\lambda}
\def\lcd{\ts,\ldots,}
\def\le{\leqslant}
\def\Lr{{\mathrm L}}

\def\Mr{{\mathrm M}}
\def\mult{\diamond}
\def\mut{\breve\mu}

\def\n{\mathfrak{n}}

\def\NN{\mathbb{N}}
\def\nns{\hskip-.5pt}
\def\Norm{\mathrm{Norm}}
\def\np{\n^{\ts\prime}}
\def\nplus{\n'}
\def\Nr{{\mathrm N}}
\def\ns{\hskip-1pt}
\def\nut{\breve\nu}

\def\om{\omega}
\def\omb{\ts\overline{\ns\om\ns}\ts}
\def\omh{\widehat{\om}}

\def\op{\oplus}
\def\ot{\otimes}

\def\p{P}
\def\P{\operatorname{P}}
\def\Ph{\Phi}
\def\ph{\varphi}
\def\psicon{\,\overline{\!\psi}}
\def\Psicon{\,\overline{\ns\Psi\ns}\ts}

\def\R{{\rm R}}
\def\Rb{\bar{\rm R}}
\def\rhot{\breve\rho}

\def\S{\operatorname{S}}
\def\SA{\mathrm{Z}}
\def\SAzero{\S}
\def\sgn{\operatorname{sign}}
\def\si{\sigma}
\def\sib{\bar{\si}}
\def\sih{\widehat{\si}}
\def\SO{\mathrm{O}}
\def\so{\mathfrak{so}}
\def\SP{\mathrm{Sp}}
\def\sp{\mathfrak{sp}}
\def\St{\widetilde{S}}
\def\SY{\operatorname{SY}}
\def\Sym{\mathfrak{S}}

\def\ts{\hskip1pt}
\def\tts{\hskip.5pt}

\def\th{\theta}

\def\TT{\mathrm{T}}
\def\Tt{\widetilde{T}}
\def\ttau{\widehat{\tau}}

\def\U{\operatorname{U}}
\def\Uh{\,\overline{\!\U(\h)\!\! }\ }
\def\Upde{\De}
\def\Upga{\mathrm{\Gamma}}
\def\ups{\upsilon}

\def\V{{\mathrm V}}

\def\X{\operatorname{X}}
\def\xib{\bar{\xi}} 
\def\xic{\check\xi}

\def\Y{\operatorname{Y}}

\def\ZA{\bar{\mathrm{Z}}}

\def\ZZ{{\mathbb Z}}

\begin{document}

\title*{Mickelsson algebras and representations of Yangians}

\author{Sergey Khoroshkin\inst{1,2}\and Maxim Nazarov\inst{3}}
\institute{Institute for Theoretical and
Experimental Physics, Moscow 117259, Russia
\and
Department of Mathematics,
Higher School of Economics,
Moscow 117312, Russia
\and 
Department of Mathematics, 
University of York, 
York YO10 5DD, England
}

\maketitle
\smartqed
\renewcommand{\theequation}{\thesection.\arabic{equation}}
\makeatletter 
\@addtoreset{equation}{section} 
\makeatother

%==============================================================================

\begin{abstract}
Let $\Y(\gl_n)$ be the Yangian of the general linear Lie algebra $\gl_n\ts$.
Let $\Y(\sp_n)$ and $\Y(\so_n)$ be the twisted Yangians corresponding to 
the symplectic and orthogonal subalgebras in the Lie algebra $\gl_n\ts$. 
These twisted Yangians are one-sided coideal
subalgebras in the Hopf algebra $\Y(\gl_n)\ts$.
We give realizations of irreducible
modules of the algebras $\Y(\sp_n)$ and $\Y(\so_n)\ts$,
as certain quotients of tensor products of symmetic and exterior powers
of the vector space $\CC^n$. For the Yangian $\Y(\gl_n)$
such realizations have been known, but we give new proofs of these results.
For the twisted Yangian $\Y(\sp_n)\ts$,
we realize all irreducible finite-dimensional modules.
For the twisted Yangian $\Y(\so_n)\ts$,
we realize all those irreducible finite-dimensional modules,
where the action of the Lie algebra $\so_n$ integrates to
an action of the complex special \text{orthogonal} 
Lie group $\mathrm{SO}_n\ts$.
Our results are based on the theory of reductive dual pairs due to Howe,
and on the representation theory of Mickelsson algebras.
\end{abstract}
 
%==============================================================================

%\newpage%%%%%%%%%%%%%%%%%%%%%%%%%%%%%%%%%%%%%%%%%%%%%%%%%%%%%%%%%%%%%%%%%%%%%%%
 
\section*{\normalsize 0.\ Preface}
\setcounter{section}{0}

%------------------------------------------------------------------------------

\subsection*{\it\normalsize 0.1.\  Brief introduction}

Mickelsson algebras are natural objects used for the study of
Harish-Chandra modules of real reductive Lie groups. They appeared
in \cite{M} first. Their detailed study was later undertaken by
Zhelobenko \cite{Z1,Z3} who employed the theory of
extremal projectors for reductive Lie algebras due to
Asherova, Smirnov and Tolstoy \cite{AST}.
The Mickelsson algebras were further studied by
Ogievetsky and the first named author of this article \cite{K,KO}.
In particular, these works dealt with the \textit{extremal cocycle\/} 
on the Weyl group of any reductive Lie algebra, introduced in \cite{Z1}. 

In \cite{Z1} a Mickelsson algebra was defined by
a finite-dimensional complex Lie algebra $\mathfrak{a}$ 
and its reductive subalgebra $\g\ts$. Following \cite{K,KO} 
we define Mickelsson algebras in a more general setting. 
Let us take any associative algebra $\Ar$ over the complex field $\CC\ts$,
containing the universal enveloping algebra $\U(\g)$ as a subalgebra.
Suppose there is a vector subspace $\V\subset\Ar\ts$, invariant 
and locally finite under the adjoint action of $\g\ts$, such that
the multiplication map $\U(\g)\ot\V\to\Ar$ is bijective.
Choose a Borel subalgebra of $\g$ with the Cartan subalgebra $\h$
and the nilpotent radical $\n\ts$.
Take the right ideal $\J=\n\ts\Ar$ 
of the algebra $\Ar$ and consider the normalizer 
$\Norm\,\J\subset\Ar$ of this ideal.
Our Mickelsson algebra
$\R$ is the quotient of the algebra $\Norm\,\J$ by
its two-sided ideal $\J\ts$. In the case when $\Ar=\U(\mathfrak{a})\ts$,
our $\R$ becomes the Mickelsson algebra considered in \cite{Z1}.

Consider the ring of fractions $\Ab$ of the algebra 
$\Ar$ relative to the set of denominators \eqref{M2}.
Take the right ideal $\Jb=\n\ts\Ab$ of the algebra $\Ab$ and its
normalizer $\Norm\,\Jb\subset\Ab\ts$. Then define
the algebra $\Rb$ as the quotient of $\Norm\,\Jb$ by $\Jb\ts$. 
Unlike $\R\ts$, the algebra $\Rb$ admits a rather nice description.
Let $\np$ be the nilpotent subalgebra of $\g$ opposite to $\n\ts$.
Let $\Jm=\Ar\ts\np$ and $\Jmb=\Ab\ts\np$ be the corresponding left ideals
of %the algebras 
$\Ar$ and $\Ab\ts$.
Consider the quotient vector space $\ZA=\Ab/(\ts\Jb+\Jmb)\ts$.
Taking the elements of $\Rb$ modulo $\Jmb$
defines a map $\Rb\to\ZA\ts$. %It turns out that 
This map is bijective. %Moreover, 
The multiplication in $\Rb$
can be described in terms of the vector space $\ZA\ts$,
by using the extremal projector for $\g$\ts;
%due to Asherova, Smirnov and Tolstoy \cite{AST}\ts;
see Subsection~1.2.
Furthermore, according to \cite{KO}
the extremal cocycle corresponding to the
Weyl group of $\g$ determines an action of the braid
group of $\g$ by automorphisms of the algebra $\ZA\ts$.
We call them the Zhelobenko automorphisms of $\ZA\ts$.
This braid group action is closely 
related to the \textit{dynamical Weyl group action\/}
of Etingof, Tarasov and Varchenko \cite{EV,TV1,TV2}.

The adjoint action of the Lie algebra $\g$ on $\Ab$ 
determines an action of the Cartan subalgebra $\h\subset\g$ on $\ZA\ts$. 
%Let $\ZA^{\ts0}$ be 
Consider the subalgebra of $\ZA$
consisting of all elements of zero weight relative to $\h\ts$.
The Zhelobenko automorphisms of $\ZA$
preserve this subalgebra, %$\ZA^{\ts0}$ 
and moreover
determine on it an action of the Weyl group of $\g\ts$,
%$\ZA^{\ts0}$,
see Subection 1.4 below.  

Let $\G$ be a reductive algebraic group over $\CC$ 
with the Lie algebra $\g\ts$. 
%We do not assume that $\G$ is connected. 
Suppose there is an action
of the group $\G$ by automorphisms on the algebra $\Ar\ts$,
extending the adjoint action of $\G$ on the subalgebra
$\U(\g)\subset\Ar\ts$. Suppose the action of $\G$
on $\Ar$ preseves the subspace $\V$, and that
the action of $\G$ on $\V$ is locally finite.
Suppose
the action of $\ts\g$ on $\Ar$ corresponding to that of $\G$ coincides
with the adjoint action.
Our joint work with Vinberg
\cite{KNV} describes the image $\Dwo$
of the subalgebra of $\G\ts$-invariant elements 
$\Ar^{\G}\subset\Ar$
under the projection $\Ar\to\Ar/(\ts\J+\Jm)\ts$.
It shows that %the image 
$\Dwo$ consists of all elements which 
have weight zero relative to $\h\ts$,
are invariant under the Zhelobenko automorphisms of $\Ab/(\ts\Jb+\Jmb)\ts$,
and are also invariant under certain other automorphisms
arising when the group 
$\G$ is not connected\ts; see Subsection 1.5.
Here $\Ar/(\ts\J+\Jm)$ is regarded
as a subspace of $\ZA=\Ab/(\ts\Jb+\Jmb)\ts$.
The Zhelobenko automorphisms do not preserve 
this subspace in general. 
Thus using the ring of fractions $\Ab$ is necessary
for this description~of~$\Ar^{\G}$.

The present work is a continuation of \cite{KNV}
and also of a series of our publications 
\cite{KN1,KN2,KN3,KN4}. 
The latter series
%was inspired by the works of Tarasov and Varchenko \cite{TV1,TV2}. It 
has established correspondences between the Zhelobenko automorphisms 
of certain algebras of the form $\ZA$ discussed above,
and the canonical interwining operators of tensor products
of representations of Yangians. These are
the Yangian $\Y(\gl_n)$ of the general linear Lie algebra 
$\gl_n\ts$, and its twisted analogues $\Y(\sp_n)$ and $\Y(\so_n)$
corresponding to the symplectic and orthogonal subalgebras
$\sp_n$ and $\so_n$ of $\gl_n\ts$.
For an introduction to %the theory of 
these Yangians see the recent book \cite{M3}. 
The Yangian $\Y(\gl_n)$ is a Hopf algebra while 
the twisted Yangians $\Y(\sp_n)$ and $\Y(\so_n)$ are 
one-sided coideal subalgebras of $\Y(\gl_n)\ts$.
%Hence the tensor products
%of representations of all these Yangians are defined.
These Yangians admit homomorphisms
to the rings of the form $\Ar^\G$ discussed above.
Moreover, together with the subalgebra $\U(\g)^{\ts\G}\subset\Ar^\G\ts$,
the image of each of these homomorphisms generates
the corresponding ring of invariants $\Ar^\G\ts$.
These homomorphisms first arose in the works of
Olshanski \cite{O1,O2} and are also discussed in our Subsection 4.3.
Another connection between the works \cite{O1,O2} and the 
theory of Mickelsson algebras was used by 
Molev to construct weight bases of irreducible finite-dimensional
modules of the Lie algebras $\sp_n$ and $\so_n\,$; see 
for instance \cite{M3}.

In the present work we combine
the results of \cite{KN1,KN2,KN3,KN4,KNV}
to give realizations of irreducible 
representations of the twisted Yangians.
In particular,
every irreducible finite-dimensional module
of the twisted Yangian $\Y(\sp_n)$ will be realized,
up to changing the action of the centre of this algebra, 
as the image
of an intertwining operator between
representations of $\Y(\sp_n)$ in tensor products
of exterior powers of $\CC^n\ts$.
For the Yangian $\Y(\gl_n)\ts$,
such realizations have been provided by the works of
Akasaka and Kashiwara \cite{AK}
and of Cherednik \cite{C2}.
%and the second author of this article~\cite{N}. 
%where some irreducible representations of Yangians and
%their twisted analogues were constructed explicitly.
Our work gives new proofs of these results
for $\Y(\gl_n)\ts$.
For the twisted Yangian $\Y(\so_n)\ts$, the images
of our intertwining operators realize,
up to changing the action of the centre of $\Y(\so_n)\ts$,
all those irreducible finite-dimensional modules, 
where the action of the subalgebra $\U(\so_n)\subset\Y(\so_n)$
integrates to an action of the complex special orthogonal Lie group 
$\textrm{SO}_n\ts$. %see Subsections 5.4~and~5.5.

Now let $\U(\g)\subset\Ar$ be an arbitrary pair
as described above. Let $\Ar^0$ be the zero weight
component of the algebra $\Ar$ relative to the adjoint action of
$\h\ts$. We denote by $\SAzero$ the subalgebra in $\ZA$
generated by the image of $\Ar^0\subset\Ab$ under the
projection $\Ab\to\ZA\ts$. Then $\Dwo\subset\SAzero$ by definition.
For any weight $\la\in\h^*$ let $N$ be an
irreducible $\SAzero\ts$-module whereon
the subalgebra $\U(\h)\subset\SAzero$
acts via the mapping $\la:\h\to\CC\ts$,
extended to a homomorphism $\U(\h)\to\CC\ts$.
Then our Proposition~\ref{proposition3.9}.
gives sufficient conditions for 
irreducibility of the restriction
of $N$ to the subalgebra $\Dwo\subset\SAzero\ts$.
One of the conditions is that the weight $\la+\rho$
is nonsingular. Here $\rho$ denotes the half-sum
of positive roots, and nonsingularity means
that the value of %the weight 
$\la+\rho$ on any positive coroot vector
of $\h$ is not a negative integer. 

Our results on irreducibility of representations of the Yangians 
$\Y(\gl_n)\ts$, $\Y(\sp_n)$ and $\Y(\so_n)$
are based on Proposition~\ref{proposition3.9}.
By using the homomorphism $\Ar^\G\to\Dwo\ts$
determined by the projection $\Ar\to\Ar/(\ts\J+\Jm)\ts$,
any irreducible $\Dwo\ts$-module becomes
an irreducible module over the algebra $\Ar^\G\ts$.
By using the homomorphisms from the Yangians to the algebras of the
form $\Ar^\G$, we then obtain irreducible representations of these Yangians.

Proposition~\ref{proposition3.9} is the main result
of Section 1. In that section we also collect several other
results on the algebra $\ZA$ corresponding to an arbitrary
pair $\U(\g)\subset\Ar\ts$. 
Below we describe the contents of all other sections of 
our article in more detail.

%------------------------------------------------------------------------------

\subsection*{\it\normalsize 0.2.\ Overview of Section 2}

In this article, we employ the theory of \textit{reductive dual pairs}
due to Howe \cite{H1}. Namely, we use the pairs $(\G\com\Gd)$
of complex reductive algebraic groups from the following list:
\begin{equation}
\label{holist}
(\GL_m\com\GL_n)\ts,\ 
(\SO_{\tts2m}\com\SP_n)\ts,\ 
(\SP_{\tts2m}\com\SO_n)\ts,\ 
(\SP_{\tts2m}\com\SP_n)\ts,\ 
(\SO_{\tts2m}\com\SO_n)\ts. 
\end{equation}
It is the Lie algebra $\gd$ of the group $\Gd$ dual to $\G$
that will determine our Yangian $\Y(\gd)\ts$.
We will employ both the symmetric and skew-symmetric versions
of the Howe duality. Let $U$ be the tensor product $\CC^m\ot\CC^n$
of vector spaces. The group $\Gd=\GL_n\com\SP_n\com\SO_n$ acts on 
$U$ naturally, via its defining action on $\CC^n\ts$.
In the symmetric case, this action of $\Gd$ is extended 
from $U$ to the Weyl algebra of $U$. The subalgebra of
$\Gd$-invariant elements in the Weyl algebra then
becomes a homomorphic image of the universal enveloping algebra
$\U(\g)$ where $\g$ is the Lie algebra of 
$\G=\GL_m\com\SO_{\tts2m}\com\SP_{\tts2m}$ respectively.
In the skew-symmetric case, the action of $\Gd$ on $U$
is extended from $U$ to the Clifford algebra of $U$. 
Then the subalgebra of
$\Gd$-invariant elements in the Clifford algebra
becomes a homomorphic image of
$\U(\g)$ where $\g$ is the Lie algebra of 
$\G=\GL_m\com\SP_{\tts2m}\com\SO_{\tts2m}$ respectively.
Thus the first pair in the list \eqref{holist} appears
in both the symmetric and skew-symmetric cases,
the next two pairs appear only in the symmetric case, 
while the last two pairs appear only in the skew-symmetric case.
 
We will denote by $\Heist$ either the Weyl or the Clifford
algebra, and will distinguish the symmetric and the
skew-symmetric cases by using a parameter $\th$ which
equals $1$ or $-1$ respectively.
The homomorphism $\U(\g)\to\Heist$
will be denoted by $\zeta\ts$. 
Our algebra $\Ar$ will be the tensor product
$\U(\g)\ot\Heist\ts$. The algebra $\U(\g)$ 
will be embedded into this tensor product diagonally\tts: 
any element $X\in\g$ will be identified with  
the element \eqref{xdiag} of $\Ar\ts$.
The subspace $\V\subset\Ar$ will be then $1\ot\Heist\ts$. 
The group $\G$ acts by automorphisms of
the algebra $\Heist\ts$, and
the corresponding action of $\g$ on $\Heist$ coincides
with the adjoint to $\zeta\ts$. Thus we get a (diagonal) action
of the group $\G$ by automorphisms of the algebra $\Ar\ts$,
obeying the assumptions from the beginning of this section.
In Subsection~2.1 we summarize the properties
of the homomorphism $\zeta:\U(\g)\to\Heist$ that we will use.

In Subsection 2.2 we introduce our main tool,
the Shapovalov form on the %double coset 
algebra $\ZA$ corresponding to $\Ar=\U(\g)\ot\Heist\ts$.
It takes values in the ring of fractions
$\Uh$ of the algebra $\U(\h)$ relative to the
set of denominators \eqref{M2}. %as the ring of fractions $\Ab$ has.
The algebra $\Heist$ is generated by its two subspaces,
$U$ and the dual $U^*\ts$. 
Let $\P\ts(U)$ be the subalgebra of $\Heist$ generated by the subspace $U$.
Then $\P\ts(U)$ is the symmetric algebra of $U$ when $\th=1$,
or the exterior algebra of $U$ when $\th=-1\ts$.
For any $\mu\in\h^*$ let $M_\mu$ be the corresponding
Verma module of the Lie algebra 
$\g\ts$. There is a unique weight $\ka\in\h^*$
satisfying the condition \eqref{deka}.
This weight is zero if $\g=\gl_m\ts$,
but differs from zero if $\g=\sp_{2m}$ or $\g=\so_{2m}\ts$.

In Subsection 2.3 %for any $\mu\in\h^*$ 
we take the space of
$\n\ts$-coinvariants of the $\g$-module $M_\mu\ot\P\ts(U)\ts$.
This space is denoted by $\Mr_\mu\ts$. It is identified
with a quotient space of the algebra $\Ar$
by the right ideal $\J$ and by a certain left ideal depending on $\mu\ts$, 
see the definition \eqref{nmummu}. If the weight $\mu+\ka$
is generic, that is if the value
of $\mu+\ka$ on any coroot vector of $\h$ is not an integer,
then the $\Uh\ts$-valued Shapovalov form on $\ZA$ defines
a $\CC\ts$-valued bilinear form on $\Mr_\mu\ts$. 
The latter form is denoted by $S_\mu\ts$.
Proposition~\ref{propositionSh14} relates the
form $S_\mu$ to the extremal projector for the Lie algebra $\g\ts$.
%as defined in Subsection 1.2.

In Subsection 2.4 for any $\la\in\h^*$
we take the subspace of weight $\la$ in the space $\Mr_\mu\ts$. 
This subspace is denoted by $\Mr_\mu^\la\ts$.
If $\mu+\ka$ is generic, 
then we define a $\CC\ts$-valued bilinear form $S_\mu^\la$ on 
the subspace $\Mr_\mu^\la\subset\Mr_\mu$ as the restriction of 
the form $S_\mu\ts$.
Letting the weights $\la$ and $\mu$ vary while
the difference $\la-\mu$ is fixed,
we extend the definition of $S_\mu^\la$
to all $\mu$ such that $\la+\rho$ is nonsingular,
by continuity. Thus here the weight $\mu+\ka$
needs not to be generic. This extension of $S_\mu^\la$
corresponds to the \textit{fusion procedure\/}
of Cherednik \cite{C2}, see below for an 
explanation of this correspondence.
Our Proposition \ref{proposition2}
states that the quotient of $\Mr_\mu^\la$
by the kernel of the bilinear form $S_\mu^\la$
is an irreducible $\SAzero\ts$-module. 
This is the main result of Section 2.
It shows that under the extra conditions 
of Proposition~\ref{proposition3.9}
the quotient of $\Mr_\mu^\la$
by the kernel of $S_\mu^\la$
is an irreducible $\Dwo\ts$-module, 
see Corollary~\ref{corollary3.11}.

In Subsection 2.5 we give another construction of 
irreducible $\Dwo\ts$-modules, by using the general results
on the Mickelsson algebras from our Subsection 1.6.
Let $L_\mu$ be the quotient of
the Verma module $M_\mu$ by its maximal proper submodule $N_\mu\ts$. 
Take the space of $\n\ts$-coinvariants of weight $\la$ of the
$\g\ts$-module $L_\mu\ot\P\ts(U)\ts$.
We denote this space by $\Lr_\mu^\la\ts\tts$. If
$\la+\rho$ is nonsingular, then $\Lr_\mu^\la$
is an irreducible $\SAzero\ts$-module.
%by our Proposition \ref{proplus}. 
Moreover, then under the extra conditions 
of Proposition~\ref{proposition3.9}
the space $\Lr_\mu^\la$ is an irreducible $\Dwo\ts$-module, 
see Corollary~\ref{corollary3.11plus}.
Note that then $\Lr_\mu^\la$ is an irreducible
$\Ar^\G\ts$-module by Proposition \ref{knv}. Thus we extend
the results of Harish-Chandra \cite{HC} and of
Lepowski and McCollum \cite{LM}.

Now take the space of $\n\ts$-coinvariants of weight $\la$ of the
$\g\ts$-module $N_\mu\ot\P\ts(U)\ts$.
We denote this space by $\Nr_\mu^\la\ts\tts$.
Then $\Lr_\mu^\la$ can be regarded
as the quotient of the space of $\Mr_\mu^\la$ by its subspace 
$\Nr_\mu^\la\ts$. In the case when $\la+\rho$ is nonsingular
while $\mu$ satisfies the conditions of 
Proposition \ref{lmeq}, the subspace 
$\Nr_\mu^\la\subset\Mr_\mu^\la$ coincides with
the kernel of the form $S_\mu^\la\ts$.
Hence the quotient of $\Mr_\mu^\la$ by the kernel
of $S_\mu^\la$ coincides with $\Lr_\mu^\la$ in this case. 
For each pair $(\G\com\Gd)$ from the list \eqref{holist},
weights $\mu$ satisfying the conditions of Proposition~\ref{lmeq} 
will be produced later on, in Section 4.

%------------------------------------------------------------------------------

\vspace{-6pt}%%%%%%%%%%%%%%%%%%%%%%%%%%%%%%%%%%%%%%%%%%%%%%%%%%%%%%%%%%%%%%%%%%

\subsection*{\it\normalsize 0.3.\ Overview of Section 3}

In Section 3 we give another interpretation
of the quotient space of $\Mr_\mu^\la$ by the kernel of %the form 
$S_\mu^\la\ts$, for any $\mu$ and nonsingular $\la+\rho\ts$.
The extremal cocycle corresponding to the
Weyl group $\Sym$ of $\g$ determines not only
an action of the braid group of $\g$
by automorphisms of the algebra $\ZA\ts$,
but also an action of the same braid group
by linear operators on the vector space $\Ab/\ts\Jb\ts$.
We call them the Zhelobenko operators on
$\Ab/\ts\Jb\ts$, see Subection~1.4.
Throughout this article
the symbol $\circ$ indicates the shifted
action of the group $\Sym$ on $\h^*$, 
see \eqref{shiftact}. 
Again suppose that
the weight $\mu+\ka$ is generic. Then
the quotient space $\Mr_\mu$ of $\Ar$ can be identified
with the quotient space of $\Ab$ by the right ideal $\Jb$
and by a certain left ideal. In Subsection 3.1 we describe the 
image of the quotient $\Mr_\mu$ of $\Ab$
under the action of the Zhelobenko operator $\xic_{\ts\si}$
on $\Ab\ts/\ts\Jb\ts$, corresponding to any element
$\si\in\Sym\ts$. If $\g=\gl_m$ then this image
can be identified with $\Mr_{\ts\si\circ\mu}\ts$.
But if $\g=\sp_{2m}$ or $\g=\so_{2m}$ then
our description of the image is more involved,
see Corollary \ref{corolary2.13}. 
For $\g=\gl_m$ this result 
was obtained in \cite{KN1,KN2}.
For $\g=\sp_{2m}$ or $\g=\so_{2m}$ 
it was obtained in \cite{KN3,KN4}.
In Subsection 3.1 we give a new proof of this result,
uniform for all these %Lie algebras~
$\g\ts$.

For generic $\mu+\ka$ the operator $\xic_{\ts\si}$ on $\Ab\ts/\ts\Jb$
maps the subspace $\Mr_\mu^\la\subset\Mr_\mu$ of weight $\la$
to the subspace of weight $\si\circ\la$ in the image of $\Mr_\mu\ts$.
Our Proposition \ref{proposition3.7} relates this map to
the extremal projector for $\g\ts$.
Letting the weights $\la$ and $\mu$ vary while
their difference $\la-\mu$ is fixed,
we extend the definition of this map %on $\Mr_\mu^\la$
to all $\mu$ such that $\la+\rho$ is nonsingular,
again by continuity\tts; see Corollary \ref{corollary3.7}.
Thus the weight $\mu+\ka$ needs not to be generic anymore.
Yet if $\mu+\ka$ is not generic then
the operator $\xic_{\ts\si}$ on $\Ab\ts/\ts\Jb$
does not necessarily define any map on the whole space $\Mr_{\mu}\ts$. 
In Subsection 3.3 we consider the operator %$\xic_{\tts0}=
$\xic_{\ts\si_0}$ where $\si_0$ is the longest
element of %the Weyl group 
$\Sym\ts$. Our Proposition \ref{corollary3.8}
states that the kernel of the corresponding map of
$\Mr_\mu^\la$ coincides with the kernel of the bilinear form 
$S_\mu^\la\ts$, for any $\mu$ and nonsingular $\la+\rho\ts$. 
The proof is based on 
Propositions \ref{propositionSh14} and \ref{proposition3.7}.
Corollary~\ref{corollary3.11} then implies
that the quotient space of $\Mr_\mu^\la$ 
by the kernel of the map defined by $\xic_{\tts\si_0}$
is an irreducible $\Dwo\ts$-module, 
under the extra conditions 
of Proposition~\ref{proposition3.9}.
Moreover, then by Proposition \ref{knv}
this quotient space of $\Mr_\mu^\la$ is
an irreducible $\Ar^\G\ts$-module.
These are the main results of Section 3, they
are stated as Theorem \ref{theorem2}
and Corollary \ref{ircor}. %respectively.

%------------------------------------------------------------------------------

\vspace{-6pt}%%%%%%%%%%%%%%%%%%%%%%%%%%%%%%%%%%%%%%%%%%%%%%%%%%%%%%%%%%%%%%%%%%

\subsection*{\it\normalsize 0.4.\ Overview of Section 4}

In Section 4 we apply the results of Section 3
to the representation theory of
the Yangians $\Y(\gd)$ where $\gd=\gl_n\com\sp_n\com\so_n$ 
according to %the list 
\eqref{holist}. In Subsection~4.1
we fix realizations of the groups $\G$ and $\Gd$
appearing in \eqref{holist}, and define the corresponding
homomorphisms $\zeta:\U(\g)\to\Heist\ts$. In Subsection 4.2
we recall the definitions of the Yangian $\Y(\gl_n)$ and the
twisted Yangians $\Y(\sp_n),\Y(\so_n)\ts$.
The latter two are defined as certain subalgebras of $\Y(\gl_n)\ts$.
Their definition implies that they are also
right coideals of $\Y(\gl_n)\ts$, see \eqref{rci}.

The twisted Yangians $\Y(\sp_n),\Y(\so_n)$
can also be defined in terms of generators and relations,
see \eqref{xrel} and \eqref{srel}.
The collection of all relations \eqref{xrel}
can be written as the \textit{reflection equation\/},
introduced by Cherednik \cite{C1} and Sklyanin \cite{S}.
The algebras defined by the relations \eqref{xrel} only,
that is without imposing the relations \eqref{srel},
are called the extended twisted Yangians. They are denoted
by $\X(\sp_n),\X(\so_n)$ respectively.
Thus for $\gd=\sp_n\com\so_n$ we have a surjective
homomorphism $\X(\gd)\to\Y(\gd)\ts$. Its kernel 
%of this homomorphism 
is generated by certain central elements of %the algebra
$\X(\gd)\ts$, see Subsection 4.2 for their definition. 
%We will prefer to work with the algebras $\X(\sp_n),\X(\so_n)$
%rather then with $\Y(\sp_n),\Y(\so_n)\ts$.

In Proposition \ref{propo1} 
for $\Ar=\U(\gl_m)\ot\Heist$ and $\th=1,-1$ 
we define a homomorphism $\Y(\gl_n)\to\Ar^{\GL_m}$.
Its image and the subalgebra $\U(\gl_m)^{\tts\GL_m}\ot1\subset\Ar$
together generate $\Ar^{\GL_m}$. %This homomorphism 
It is related to the homomorphism
$\Y(\gl_n)\to\U(\gl_{\ts n+l})^{\tts\GL_{\tts l}}$ defined for 
$l=1,2,\ldots$ by Olshanski \cite{O1}. It is also related
to the \textit{Cherednik functor\/} 
studied by Arakawa and Suzuki \cite{AS}
and to the \textit{Drinfeld functor\/} \cite{D1}.   
We explained the relations in \cite{KN1,KN2}.

Now consider any pair $(\G\com\Gd)$ from the list \eqref{holist}
other than the pair $(\GL_m\com\GL_n)\ts$.
In Proposition \ref{propo2} 
for $\Ar=\U(\g)\ot\Heist$ 
we define a homomorphism $\X(\gd)\to\Ar^{\G}$.
Its image together with the subalgebra $\U(\g)^{\tts\G}\ot1\subset\Ar$
generate $\Ar^{\G}$.
This homomorphism is related to the homomorphisms
$\Y(\sp_n)\to\U(\sp_{\ts n+2l})^{\ts\SP_{\tts2l}}$ 
and
$\Y(\so_n)\to\U(\so_{\ts n+l})^{\tts\SO_{\tts l}}$ 
defined for any $l=1,2,\ldots$ by Olshanski \cite{O2}. 
We explained the relation in \cite{KN3,KN4}.
To prove our  Propositions \ref{propo1} and \ref{propo2}
we use the classical invariant theory, like
Molev and Olshanski \cite{MO} did
when studying the homomorphisms defined in \cite{O1,O2}. 
%We also follow the approach of Howe \cite{H1} to that theory.

For any $\gd=\gl_n\com\sp_n\com\so_n$ the algebra
$\Y(\gd)$ contains the universal enveloping algebra $\U(\gd)$
as a subalgebra, and admits a homomorphism $\Y(\gd)\to\U(\gd)$
identical on that subalgebra. For $\gd=\gl_n$ this fact
has been well known, see for instance the seminal work of
Kulish, Reshetikhin and Sklyanin \cite{KRS}.
For $\gd=\sp_n\com\so_n$ this is another result from \cite{O2}. 
If $\th=1$ then denote by $\Ph^{\ts k}$
the $k\,$th symmetric power of the defining $\gl_n$-module $\CC^n\ts$.
If $\th=-1$ then denote by the same symbol $\Ph^{\ts k}$
the $k\,$th exterior power of the $\gl_n$-module $\CC^n\ts$.
%The group $\Gd=\GL_n$ also acts on $\Ph^{\ts k}\ts$. 
Using the homomorphism $\Y(\gl_n)\to\U(\gl_n)\ts$,
regard $\Ph^{\ts k}$ as a module over the Yangian $\Y(\gl_n)\ts$. 
For any $t\in\CC$
denote by $\Ph^{\ts k}_{\tts t}$ the $\Y(\gl_n)\ts$-module obtained by
pulling the $\Y(\gl_n)\ts$-module $\Ph^{\ts k}$ back through the
automorphism \eqref{tauz} where $z=\th\,t\ts$.

The definition of the Hopf algebra $\Y(\gl_n)$ employs a certain
$n\times n$ matrix $T(x)\ts$. The entries of this matrix
are formal power series in $x^{-1}$ where $x$ is the
\textit{spectral parameter\/}.  The coefficients of
these series are generators of %the algebra 
$\Y(\gl_n)\ts$. For $\Gd=\SP_n$ or $\Gd=\SO_n$ we denote by 
$\Tt(x)$ the transpose to the matrix $T(x)$
relative to the bilinear form on $\CC^{\ts n}$
preserved by the subgroup $\Gd\subset\GL_n\ts$.
Then the assignment \eqref{transauto}
defines an involutive automorphism of the
algebra $\Y(\gl_n)\ts$. 
For any non-negative integer $k$
we denote by $\Ph^{\ts-k}_{\ts t}$
the $\Y(\gl_n)\ts$-module obtained by
pulling the % $\Y(\gl_n)\ts$-module 
$\Ph^{\ts k}_{\tts t}$ back through the automorphism \eqref{transauto}.

By definition, for $\gd=\sp_n\com\so_n$ the subalgebra
$\Y(\gd)\subset\Y(\gl_n)$ is generated by the coefficients
of the series in $x$ which arise as the entries of the $n\times n$ 
matrix $\Tt(-x)\,T(x)\ts$.
It turns out that if $\gd=\sp_n$ and $\th=1\ts$, or if
$\gd=\so_n$ and $\th=-1\ts$, then the restriction of 
the above defined $\Y(\gl_n)\ts$-module $\Ph^{\ts-k}_{\ts t}$
to the subalgebra $\Y(\gd)\subset\Y(\gl_n)$
coincides with the restriction of the
$\Y(\gl_n)\ts$-module $\Ph^{\ts k}_{\tts t}\ts$.
Both cases correspond to $\g=\so_{\tts2m}\ts$.
For further explanation of this phenomenon, 
see the end of Subsection 4.4.
 
For any pair $(\G\com\Gd)$ from the list 
\eqref{holist} and for any weights
$\la\com\mu\in\h^*$ the subspace $\Mr_\mu^\la$ of the
quotient $\Mr_\mu$ of the algebra $\Ar$ is a module over 
the subalgebra $\Ar^\G$ by definition. For
$(\G\com\Gd)=(\GL_m\com\GL_n)$ 
by using the homomorphism $\Y(\gl_n)\to\Ar^{\GL_m}$
we can regard $\Mr_\mu^\la$ as a module over the algebra $\Y(\gl_n)\ts$.
If non-zero, this $\Y(\gl_n)\ts$-module
is equivalent to a tensor product of certain
modules of the form $\Ph^{\ts k}_{\tts t}\ts$.
This is stated as Proposition~\ref{iso1}, see
\cite{KN1,KN2} for the proof. Here we use the comultiplication 
\eqref{1.33} on %the algebra 
$\Y(\gl_n)\ts$. 

For any other pair $(\G\com\Gd)$ from %the list 
\eqref{holist},
we can regard $\Mr_\mu^\la$ as a module of the algebra 
$\X(\gd)$ by using the homomorphism $\X(\gd)\to\Ar^{\G}$.
If non-zero, this $\X(\gd)\ts$-module
is also equivalent to a tensor product of certain
modules of the form $\Ph^{\ts k}_{\tts t}\ts$,  %but
pulled back through an automorphism \eqref{fus} of $\X(\gd)\ts$.
This is a particular case of our Proposition~\ref{iso2},
see \cite{KN3,KN4} for the proof of that proposition. 
Here we first take the tensor product of the
$\Y(\gl_n)\ts$-modules $\Ph^{\ts k}_{\tts t}\ts$,
then restrict the tensor product to the subalgebra
$\Y(\gd)\subset\Y(\gl_n)\ts$, and then use the
homomorphism $\X(\gd)\to\Y(\gd)$ mentioned above.
Equivalently, we can use the right $\Y(\gl_n)\ts$-comodule 
structure on the algebra $\X(\gd)\ts$; see %the end of 
Subsection~4.2.
 
Now suppose $\la+\rho$ is nonsingular, while
$\mu$ is arbitrary. Consider again the 
mapping of the weight subspace $\Mr_\mu^\la$
defined by the Zhelobenko operator $\xic_{\ts\si_0}$
on $\Ab\ts/\ts\Jb\ts$. In the case $(\G\com\Gd)=(\GL_m\com\GL_n)$ 
by replacing $\Mr_\mu^\la$ by an equivalent $\Y(\gl_n)\ts$-module
we obtain an intertwining operator from a tensor product
of $\Y(\gl_n)\ts$-modules of the form $\Ph^{\ts k}_{\tts t}\ts$.
The target module of this operator can be identified with
the tensor product of the same $\Y(\gl_n)\ts$-modules
$\Ph^{\ts k}_{\tts t}$ as for the source module, 
but taken in the reversed order\tts; see \eqref{ak1}.
The quotient by the kernel, or equivalently the image 
of this intertwining operator is an irreducible $\Y(\gl_n)\ts$-module
due to Corollary~\ref{ircor} and Proposition~\ref{propo1}.

In the case $(\G\com\Gd)=(\GL_m\com\GL_n)$ we can regard 
the quotient space $\Lr_\mu^\la$ of $\Mr_\mu^\la$ as an 
$\Y(\gl_n)\ts$-module. %Note that 
The Yangian
$\Y(\gl_n)$ acts on $\Mr_\mu^\la$ and hence on $\Lr_\mu^\la$
via the homomorphism $\Y(\gl_n)\to\Ar^{\GL_m}\ts$.
If $\la+\rho$ is nonsingular, then the
$\Y(\gl_n)\ts$-module $\Lr_\mu^\la$ is irreducible for any $\mu\ts$, 
see Subsection 4.6. 
Moreover, then the conditions of Proposition~\ref{lmeq}
are satisfied, if for each positive root $\al$ of $\gl_m$ 
the number $z_\al$ defined just before stating 
Proposition~\ref{norm1}, is not zero. Under these conditions,
the image of our intertwining operator \eqref{ak1} is not zero,
and is equivalent to $\Lr_\mu^\la$ as $\Y(\gl_n)\ts$-module.
Thus we extend the results of \cite{AS}.

Now take a pair $(\G\com\Gd)$ from the list 
\eqref{holist} other than $(\GL_m\com\GL_n)\ts$.
By replacing $\Mr_\mu^\la$ by an equivalent $\X(\gd)\ts$-module
we again obtain an intertwining operator from a certain tensor product
of modules of the form $\Ph^{\ts k}_{\tts t}\ts$.
%If $\g=\sp_{\tts2m}\ts$, or if $\g=\so_{\tts2m}$ where $m$ is even,
The target module of this %intertwining 
operator can be identified with the tensor product of the modules
$\Ph^{\ts-k}_{\ts t}$ corresponding to the tensor factors
$\Ph^{\ts k}_{\tts t}$ of the source module, and
taken in the same order\tts; 
see \eqref{ak2}.
Here we need not to apply to the source and the target
$\X(\gd)\ts$-modules the automorphisms \eqref{fus}, 
because by Proposition~\ref{iso2}
these automorphisms are the same
for the source and for the target. Hence
we may also regard both the source and the target
as $\Y(\gd)\ts$-modules.
Note that if $\g=\so_{\tts2m}\ts$, then 
instead of the first tensor factor 
$\Ph^{\ts-k}_{\ts t}$ of the target
we can use $\Ph^{\ts k}_{\tts t}$ as well\tts;
see \eqref{4.56}.
This is because for $\g=\so_{\tts2m}$
the restrictions of %the $\Y(\gl_n)\ts$-modules 
$\Ph^{\ts k}_{\tts t}$
and $\Ph^{\ts-k}_{\tts t}$ to the subalgebra 
$\Y(\gd)\subset\Y(\gl_n)$ coincide, and because 
we regard this subalgebra as a right coideal
to define the tensor products.

If $\G=\SP_{\tts2m}$ then 
by Corollary~\ref{ircor} and Proposition~\ref{propo2}
the image %quotient by the kernel
of our intertwining 
operator \eqref{ak2}
%of the tensor products
is an irreducible $\Y(\gd)\ts$-module
for both $\th=1$ and $\th=-1$, that is for both 
$\gd=\so_n$ and $\gd=\sp_n\ts$.
But if $\G=\SO_{\tts2m}$ then our
Corollary~\ref{ircor} and Proposition~\ref{propo2}
imply $\Y(\gd)\ts$-irreducibility of the image
%by the kernel of our intertwining operator of the tensor products
%as an $\X(\gd)\ts$-module
only under an extra condition, that
the stabilizer of the weight $\la+\rho$ in
the extended Weyl group of $\g=\so_{\tts2m}$ is 
contained in the proper Weyl group.
In Subsection 4.5 we remove this extra condition,
but only for $\th=-1\ts$, that is only for $\gd=\so_n\ts$.
Namely, for $\G=\SO_{\tts2m}$ and $\th=-1$ we prove that the 
image %quotient by the kernel 
of our intertwining operator
is irreducible under the joint action of the algebra $\Y(\so_n)$ and 
the group $\SO_n\ts$. Moreover, if $n$ is odd then
the image is irreducible under the action of %the algebra 
$\Y(\so_n)$ alone. But if $n$ is even then the image is either 
an irreducible $\Y(\so_n)\ts$-module, or is a direct sum of two
non-equivalent irreducible $\Y(\so_n)\ts$-modules\ts; 
see Corollary \ref{iror}.

Now let $(\G\com\Gd)$ again be any pair from the list 
\eqref{holist} other than $(\GL_m\com\GL_n)\ts$.
Regard the quotient space $\Lr_\mu^\la$ of $\Mr_\mu^\la$ as
$\X(\gd)\ts$-module. The extended twisted Yangian
$\X(\gd)$ acts on $\Mr_\mu^\la$ and hence on $\Lr_\mu^\la$
via the homomorphism $\X(\gd)\to\Ar^{\G}\ts$.
In Subsection~4.6 we derive
the same irreducibilty properties
of the $\X(\gd)\ts$-module $\Lr_\mu^\la\ts$,
as described above for the image of our %intertwining 
operator \eqref{ak2}. Suppose that
$\la+\rho$ is nonsingular. If $\G=\SP_{\tts2m}$ then 
$\Lr_\mu^\la$ is an irreducible $\X(\gd)\ts$-module
for both $\gd=\so_n$ and $\gd=\sp_n\ts$.
If $\G=\SO_{\tts2m}$ then 
the $\X(\gd)\ts$-module $\Lr_\mu^\la$
is irreducible under an extra condition, that
the stabilizer of the weight $\la+\rho$ in
the extended Weyl group of $\g=\so_{\tts2m}$ is 
contained in the proper Weyl group.
If $\G=\SO_{\tts2m}$ and $\th=-1$ then $\Lr_\mu^\la$
is irreducible under the joint action of %the algebra 
$\X(\so_n)$ and %the group 
$\SO_n\ts$. Moreover, if $n$ is odd then $\Lr_\mu^\la$
the image is irreducible under the action of the algebra 
$\X(\so_n)$ alone. But if $n$ is even then $\Lr_\mu^\la$ is either 
an irreducible $\X(\so_n)\ts$-module, or is a direct sum of two
non-equivalent irreducible $\X(\so_n)\ts$-modules\ts.

Let $(\G\com\Gd)$ once again be any pair from
\eqref{holist} other than $(\GL_m\com\GL_n)\ts$.
Let $\la+\rho$ be nonsingular.
When $\th=-1$ or $n=1\ts$,
the conditions of Proposition~\ref{lmeq}
are satisfied if for each positive root $\al$ of $\g$ 
the number $z_\al$ defined just before stating 
Proposition \ref{norm2}, is not zero. 
When $\th=1$ and $n>1\ts$,
the conditions of Proposition~\ref{lmeq}
are satisfied if $z_\al$ is not zero  
for each compact positive root $\al$ of $\g\ts.$ 
Under these conditions,
the image of our intertwining operator \eqref{ak2} is not zero,
and is equivalent to $\Lr_\mu^\la$ as $\X(\gd)\ts$-module.
Moreover, then the image is equivalent to $\Lr_\mu^\la$
under the joint actions of $\X(\gd)\ts$ and~$\Gd$.

In Subsection 4.7 we once again consider the
$\Y(\gd)\ts$-intertwining operators \eqref{ak1} and \eqref{ak2},
for $\gd=\gl_n$ and $\gd=\sp_n\com\so_n$ respectively.
We show that the target $\Y(\gd)\ts$-modules of these operators are 
dual to the source modules. We also show
how this duality arises from the theory of Mickelsson algebras,
when the weight $\la+\rho$ is nonsingular.

%------------------------------------------------------------------------------

\subsection*{\it\normalsize 0.5.\ Overview of Section 5}

In Section 5 we use our intertwining operators \eqref{ak1} and
\eqref{ak2} to give realizations of irreducible 
representations of the Yangian $\Y(\gl_n)$ and of the twisted Yangians 
$\Y(\sp_n)$, $\Y(\so_n)\ts$. We call
two $\Y(\gl_n)\ts$-modules similar if they differ
only by an automorphism \eqref{fut} of %the algebra 
$\Y(\gl_n)\ts$, where $g(x)$ is any 
formal power series in $x^{-1}$ with coefficients from
$\CC$ and leading term $1\ts$.
Up to equivalence and similarity, the irreducible
finite-dimensional $\Y(\gl_n)\ts$-modules
were classified by Drinfeld \cite{D2}, 
who generalized the classification 
given for $n=2$ by Tarasov \cite{T1,T2}.
It was then proved by Akasaka and Kashiwara \cite{AK}
that any of these modules can be realized as
a quotient of a tensor product of
$\Y(\gl_n)\ts$-modules of the form 
$\Ph^{\ts k}_{\tts t}$ with $\th=-1\ts$; 
see also the work of Chari and Pressley \cite{CP}.
Further results %in this direction 
were obtained by Chari \cite{C} and by Brundan and Kleshchev \cite{BK}.
Note that the works \cite{AK} and \cite{C} deal with
representations of quantum affine algebras. For
a connection to the representation theory of Yangians
see the work of Molev, Tolstoy and Zhang \cite{MTZ}.

In Subsection~5.1 we give new proofs of the results from \cite{AK}
for $\Y(\gl_n)\ts$. Namely, we prove that up to equivalence and similarity, 
any irreducible finite-dimensional $\Y(\gl_n)\ts$-module
arises as the quotient by the kernel of intertwining
operator \eqref{ak1} for $\th=-1\ts$, some $m$ and certain weights
$\la\com\mu$ of $\gl_m\ts$. Here the weight $\la+\rho$ 
is nonsingular, that is satisfies the conditions \eqref{domcon1}.
The difference $\nu=\la-\mu$ satisfies the conditions
\eqref{chercon}, which %in our context 
come from Proposition \ref{norm1}. This proposition 
was obtained in \cite{KN1,KN2} and
gives an explicit formula for the image 
under the operator \eqref{ak1} of a certain distinguished
vector, called highest. 
For $\th=-1$ 
the conditions \eqref{chercon} guarantee that the image is not zero\tts;
see Theorem \ref{drin1}.

Now consider the twisted Yangians $\Y(\gd)$ where
$\gd=\sp_n\com\so_n\ts$. Any automorphism \eqref{fut}
of $\Y(\gl_n)$ determines an automorphism of the
subalgebra $\Y(\gd)\subset\Y(\gl_n)\ts$.
We call two $\Y(\gd)\ts$-modules similar if they differ only
by such an automorphism. %of $\Y(\gd)\ts$.
Up to equivalence and similarity, the irreducible
finite-dimensional $\Y(\gd)\ts$-modules
have been classified by Molev\tts;
see \cite{M3} for an exposition
of the classification. No explicit realization of these 
$\Y(\gd)\ts$-modules had been known so far in general.
It is provided by our results. 
Thus we extend the works of Molev \cite{M2} and
the second named author of this article~\cite{N}
which give explicit realizations of irreducible 
$\Y(\gd)\ts$-modules from a particular class. 

In Subsection 5.3 we prove that up to equivalence and similarity, 
any irreducible finite-dimensional $\Y(\sp_n)\ts$-module
arises as the quotient by the kernel of the intertwining
operator \eqref{ak2} for $\th=-1\ts$, some $m$ and certain weights
$\la\com\mu$ of $\sp_{\tts2m}\ts$. Here the weight $\la+\rho$ 
is nonsingular, that is satisfies the conditions 
\eqref{domcon31},\eqref{domcon32} and \eqref{domcon33}.  
The weight $\nu=\la-\mu-\ka$ of $\sp_{\tts2m}$
satisfies the conditions
\eqref{chercon31},\eqref{chercon32} and \eqref{chercon33}\tts;
see Theorem \ref{drin3}.

Recall that the twisted Yangian $\Y(\so_n)$ 
contains $\U(\so_n)$ as a subalgebra.
We shall call a finite-dimensional module over the algebra $\Y(\so_n)$
an $(\ts\Y(\so_n),\SO_n)\ts$-module, if
the group $\SO_n$ also acts on this module,
and the corresponding action of the Lie algebra 
$\so_n$ of $\SO_n$ on this module coincides with the 
action obtained by restricting the action of $\Y(\so_n)$
to the subalgebra $\U(\so_n)\ts$. If $n$ is odd then
any irreducible $(\ts\Y(\so_n),\SO_n)\ts$-module
is also irreducible over $\Y(\so_n)\ts$. But if $n$ is even
then any $(\ts\Y(\so_n),\SO_n)\ts$-module
is either irreducible over $\Y(\so_n)\ts$,
or splits into a direct sum of two irreducible
$\Y(\so_n)\ts$-modules, not equivalent to each other.
The irreducible $\Y(\so_n)\ts$-modules occuring in this way
for any $n$ are all those
whose restriction to the subalgebra $\U(\so_n)\subset\Y(\so_n)$
integrates to a module of the special orthogonal group 
$\mathrm{SO}_n\subset\SO_n\ts$.

In Subsections 5.4 and 5.5 we prove
that up to equivalence and similarity, 
any finite-dimensional irreducible
$(\ts\Y(\so_n),\SO_n)\ts$-module
arises as the quotient by the kernel of the intertwining
operator \eqref{ak2} for $\th=-1\ts$, some $m$ and certain weights
$\la\com\mu$ of $\so_{\tts2m}\ts$. Here $\la+\rho$ 
is nonsingular, that is satisfies the conditions 
\eqref{domcon41},\eqref{domcon42}.  
The weight $\nu=\la-\mu-\ka$ of $\so_{\tts2m}$
satisfies the conditions \eqref{chercon41},\eqref{chercon42}\tts;
see Theorems \ref{drin4} and \ref{drin5}.

Here for both $\gd=\sp_n\com\so_n$ the weight $\ka$ of 
$\g=\sp_{\tts2m}\com\so_{\tts2m}$ respectively 
is defined by \eqref{deka}.
The above mentioned conditions on the weight 
$\nu$ of $\g$ come from Proposition~\ref{norm2}. 
This proposition was obtained in \cite{KN3,KN4} and
gives an explicit formula for the image 
under the operator \eqref{ak2} of a certain distinguished
vector, again called highest. 
For $\th=-1$ 
our conditions on %the weight 
$\nu$ guarantee that the image is not zero.
The notion of a highest vector relative to
$\Y(\gd)$ for $\gd=\sp_n\com\so_n$ is discussed in Subsection 5.2.

%==============================================================================

\section*{\bf\normalsize 1.\ Mickelsson algebras}
\setcounter{section}{1}
\setcounter{theorem*}{0}
\setcounter{equation}{0}

%------------------------------------------------------------------------------

\subsection*{\it\normalsize 1.1.\ Generalities}

Let $\g$ be any reductive complex Lie algebra of semisimple rank $r\ts$.
Choose a \textit{triangular decomposition\/}
\begin{equation}
\label{hc0}
\g=\n\op\h\op\np
\end{equation}
where $\h$ is a Cartan subalgebra, while
$\n$ and $\np$ are the nilpotent radicals of two opposite Borel
subalgebras of $\g$ containing $\h\ts$. Consider 
the root system of $\g$ in $\h^*$. The set of the positive
roots of $\g$ will be denoted by $\De^+\ts$.

Let $\al_1\lcd\al_r\in\De^+$ be simple roots.
For each $c=1\lcd r$ let $H_c=\al_c^\vee\in\h$ be the coroot
corresponding to the simple root $\al_c\ts$, 
Let $E_c\in\np$ and $F_c\in\n$ be the root vectors
corresponding to the roots $\al_c$ and $-\al_c\ts$. 
We suppose that $[E_c\com F_c]=H_c\ts$. 
Let $\ep$ be a Chevalley anti-involution on $\g\ts$.
This is an involutive anti-automorphism of $\g$
identical on $\h\ts$, such that for every $c=1\lcd r$ one has
$\ep(E_c)=\ep_c\,F_c$ and $\ep(F_c)=\ep_c^{-1}E_c$
for some non-zero $\ep_c\in\CC\ts$.

Let $\Sym$ be the Weyl group of the root system of $\g\ts$. Let
$\si_1\lcd\si_r\in\Sym$ be the reflections in $\h^*$ corresponding
to the simple roots $\al_1\lcd\al_r\ts$. We also use the
induced action of the Weyl group $\Sym$ on the vector space $\h\ts$.
It is defined by setting $\la\ts(\si(H))=\si^{-1}(\la)(H)$ for all
$\si\in\Sym$, $H\in\h$ and $\la\in\h^*$. Let $\rho\in\h^*$ be the
half-sum of the positive roots. Then the \textit{shifted action\/}
$\circ$ of the group $\Sym$ on the vector space $\h^*$ is defined by
setting
\begin{equation}
\label{shiftact}
\si\circ\la=\si(\la+\rho)-\rho\ts.
\end{equation}
It induces an action $\circ$ of $\Sym$ on the commutative
algebra $\U(\h)\ts$,
by regarding the elements of this algebra as
polynomial functions on $\h^\ast$. In particular, then
$(\si\circ H)(\la)=H(\si^{-1}\circ\lambda)$ for $H\in\h\ts$. 
For any left
$\g\ts$-module $K$ and any $\la\in\h^*$ we will denote by
$K^\la$ the subspace \text{of weight\/} $\la$ in $K\ts$; we have
$v\in K^\la$ if and only if $H\ts v=\la(H)\,v$ for every $H\in\h\ts$.

Let $\Ar$ be any complex associative algebra, 
containing as a subalgebra the universal
enveloping algebra $\U(\g)\ts$. Suppose there is also
a vector subspace $\V\subset\Ar\ts$,
invariant under the adjoint action of the Lie algebra 
$\g$ on $\Ar\ts$, such that:
\\
\vbox{
\begin{enumerate}
\item[(a)]
the multiplication map $\U(\g)\ot\V\to\Ar:X\ot Y\mapsto X\,Y$ is
bijective\ts;
%and equivariant relative to the adjoint action of $\g\ts$
%on $\U(\g)$, on $\V$ and on $\Ar\ts;
\item[(b)]
the multiplication map $\V\ot\U(\g)\to\Ar:Y\ot X\mapsto YX$ is
bijective\ts;
\item[(c)]
the space $\V$ is locally finite relative to the adjoint action of $\g$\ts.
\end{enumerate}
}
Since the subspace 
$\V\subset\Ar$ is invariant under the adjoint
action of $\g$, the conditions (a) and (b) are equivalent. The
condition (c) means that $\V$ can be decomposed into direct sum of
irreducible finite-dimensional $\g\ts$-modules. Since the
adjoint action of $\g$ on $\U(\g)$ is locally finite, due to (c)
the same is true for the adjoint action of $\g$ on $\Ar\ts$.

Let $\J$ be the right ideal of the algebra $\Ar$ generated by the
elements of the subalgebra $\n\subset\g\ts$.
Let ${\rm Norm}\,\ts\J\ts\subset\ts\Ar$ be the \text{normalizer\/} of this
right ideal, so that $Y\in{\rm Norm}\ts\,\J$ if and only if
$\,Y\ts\J\subset\J\ts$. Then $\J$ is a two-sided ideal of 
${\rm Norm}\ts\,\J\ts$. The quotient algebra
\begin{equation*}
%\label{malg}
\R\,=\,{\rm Norm}\,\J\ts/\ts\J
\end{equation*}
is called a \textit{Mickelsson algebra}. 
%The backslash here
%indicates that $\J$ is a right rather than a left ideal of $\Ar\ts$.
By definition, the algebra $\R$ acts on the space $K/\n\ts K$ of
$\n\ts$-coinvariants of any left $\Ar$-module $K\ts$. This space
itself will be denoted by $K_\n\ts$.

Recall that $\U(\g)\subset\Ar\ts$. Hence we have
$\U(\h)\ts\subset\ts{\rm Norm}\ts\,\J\ts$. It follows that the
Mickelsson algebra $\R$ also contains $\U(\h)$ as subalgebra. Due to
the condition (a) above, $\R$ is a torsion free module over $\U(\h)$
relative to the left multiplication.
%that is the equality $X\,Y=0$ for $X\in\U(\h)$ and $Y\in\R$ implies that
%$X=0\ts$.
The condition (b) implies that $\R$ is a torsion free module over
$\U(\h)$ relative to the right multiplication,

For each root $\al\in\h^*$ let $H_\al=\al^\vee\in\h$ be the corresponding
coroot. The Weyl group $\Sym$ acts on the vector space $\h$ so that
$\si:\,H_\al\mapsto H_{\si(\al)}$ for any $\si\in\Sym\ts$.
Denote by $\Uh$ the ring of fractions of the
commutative algebra $\U(\h)$ relative to the set of denominators
\begin{equation}
\label{M2}
\{\,H_\al+z\ |\ \al\in\De^+\ts,\ z\in\ZZ\,\ts\}\,.
\end{equation}
The elements of this ring can be regarded as rational functions
on the vector space $\h^\ast\ts$. The elements of
$\U(\h)\subset\,\overline{\!\U(\h)\!\!\!}\,\,\,$ are then regarded
as polynomial functions on $\h^\ast\ts$. Let
$\overline{\!\U(\g)\!\!\!}\,\,\,$ be the ring of fractions of the
algebra $\U(\g)$ relative to the same set of denominators \eqref{M2}.
This ring is well defined, because $\U(\g)$ satisfies the Ore
condition relative to \eqref{M2}.

Now regard \eqref{M2} as a subset of $\Ar$ using the embedding of
$\h\subset\g\ts$. Due to the conditions (a,b,c) the algebra $\Ar$
also satisfies the Ore condition relative to its subset \eqref{M2}.
Let $\Ab$ be the ring of fractions of $\Ar$ relative to the same set
of denominators \eqref{M2}. Then $\Ab$ is a free left and free right
$\Uh\,$-module, which contains $\Ar$ as a subalgebra.

Denote by $\Jb$ the right ideal of the algebra $\Ab$ generated by
the elements of $\n\ts$. Consider the normalizer 
${\rm Norm}\ts\,\Jb\subset\Ar$ of this right ideal. This normalizer
coincides with the ring of fractions of ${\rm Norm}\ts\,\J$ relative
to the same set of denominators \eqref{M2} as before. Now $\Jb$ is a
two-sided ideal of ${\rm Norm}\ts\,\Jb\ts$. The quotient algebra
$$
\Rb\,=\,{\rm Norm}\,\Jb\ts/\ts\Jb
$$
bears the same name of \textrm{Mickelsson algebra\/}, as the
quotient algebra $\R$ does. Note that $\Jb\cap\Ar=\J\ts$, because
the subalgebra $\h\subset\g$ normalizes the subalgebra
$\n\subset\g\ts$. Hence the natural embedding
$\Ar\to\Ab$ determines an embedding %of Mickelsoon algebras
$\R\to\Rb\ts$.

%------------------------------------------------------------------------------

\vspace{-6pt}%%%%%%%%%%%%%%%%%%%%%%%%%%%%%%%%%%%%%%%%%%%%%%%%%%%%%%%%%%%%%%%%%%

\subsection*{\it\normalsize 1.2.\ Double coset algebra}

For each positive root $\al\in\De^+$ let $E_\al\in\np$ and
$F_\al\in\n$ be the corresponding root vectors.
In paricular, for any simple root
$\al=\al_c$ we have $E_\al=E_c$ and $F_\al=F_c\ts$. For any weight
$\la\in\h^*$ consider the infinite sum
\begin{equation}
\label{pal} 
P_\al[\la]\,=\,1+ \sum_{s=1}^\infty\,(-1)^s\,\bigl(\ts
s!\,(H_\al+\la\ts(H_\al)+1)\ldots(H_\al+\la\ts(H_\al)+s) \bigr)^{-1}
F_\al^s\,E_\al^s
\end{equation}
where the denominators do not belong to the set \eqref{M2}
in general. The sum belongs to certain algebra which contains
$\,\overline{\!\U(\g)\!\!\!}\,\,\,\ts$ and needs not to be defined here; 
see \cite[Section~1]{Z3}. %for its definition.
Take $\si\in\Sym$ %be an element of the Weyl group $\Sym\ts$. 
and any reduced decomposition $\si=\si_{d_l}\ldots\si_{d_1}$. 
Here $l$ is the
length of $\si$. Any such a decomposition defines a 
\textit{normally ordered\/} sequence of positive roots\ts:
$$
\be_1=\al_{d_1}\,, \quad \be_2=\si_{d_1}(\alpha_{d_{2}})\,, \quad
\ldots\,, \quad
\be_l=\si_{d_1}\si_{d_2}\ldots\si_{d_{l-1}}(\al_{d_l})\,.
$$
The product
$$
P_\si\ts[\la]=P_{\ts\be_l}[\la]\ldots P_{\ts\be_1}[\la]
$$
in the above mentioned algebra does not depend on the choice
of a decomposition; see \cite[Section 2]{Z3}.
Let $\si_0$ be the longest element of the Weyl group $\Sym\ts$. Put
$P\ts[\la]=P_{\si_0}[\la]\ts$.

The element $P=P\ts[\la]$ with $\la=\rho$ is called the \textit{extremal
projector\/} for the reductive Lie algebra~$\g\ts$; its definition
is due to Asherova, Smirnov and Tolstoy \cite{AST}.
The element $P$ can be presented as an infinite sum 
of elements of the algebra $\,\,\overline{\!\U(\g)\!\!\!}\,\,\,\ts$.
Thus it belongs to a certain completion 
$\,\widetilde{\!\U(\g)\!\!\!}\,\,\,$ of the latter algebra;
the completion needs not to be defined here.
Below are the basic properties of the extremal projector:
\begin{gather}
\label{M35} 
P^{\ts2}=P\ts;
\\[3pt]
\label{M4} 
E_\al P=0 
\quad\text{and}\quad 
P\ts F_\al=0
\quad\text{for}\quad 
\al\in\De^+\ts;
\\
\label{M5} 
P\in1+\n\,\,\widetilde{\!\U(\g)\!\!\!}\,\,\,
\quad\text{and}\quad P\in1+\,\widetilde{\!\U(\g)\!\!\!}\,\,\,\,\np.
\end{gather}
For the proofs of \eqref{M35} and \eqref{M4} see
\cite{AST} and \cite[Section 3]{Z3}. The~two properties
\eqref{M5} follow directly from the definition \eqref{pal},
which also implies that $\ep(P)=P\ts$. Moreover,
\begin{equation}
\label{M55}
P\,\in\,1+\sum_{s=1}^\infty\,
\n^{\ts s}\,\,\overline{\!\U(\g)\!\!\!}\,\,\,\ts
\quad\text{and}\quad
P\,\in\,1+\sum_{s=1}^\infty\,
\,\,\overline{\!\U(\g)\!\!\!}\,\,\,\ts(\ts\np\ts)^s\,.
\end{equation}

Now let $\Jm$ be the left ideal of $\Ar$ generated by elements of
$\np$. Similarly, let $\Jmb$ be the left ideal of $\Ab$ generated by
elements of $\np$. First consider the double coset vector space
$$
\SA
%=\J\backslash\Ar/\Jm
=\Ar/(\J+\Jm)\,,
$$
where the quotient is taken by a left and a right ideal.
Since the Cartan subalgebra $\h\subset\g$ normalizes the subalgebras
$\n$ and $\np$, the left and right multiplications in $\Ar$ by the
elements of the algebra $\U(\h)$ make the vector space $\SA$ a
$\U(\h)\ts$-bimodule. 
Relative to the adjoint action of $\h$ on $\SA\ts$,
we have a weight~decomposition
\begin{equation}
\label{sawede}
\SA\,=\underset{\lambda\in\h^*}\op\,\SA^\la\ts.
\end{equation}

Next consider the double coset vector space
$$
\ZA
%=\Jb\backslash\Ab/\Jmb
=\Ab/(\Jb+\Jmb)\,.
$$
This vector space is a bimodule over the algebra
$\overline{\!\U(\h)\!\!\!}\,\,\,$. Since
$(\Jb+\Jmb)\cap\Ar=\J+\Jm\ts$, the natural embedding $\Ar\to\Ab$
determines an embedding of $\U(\h)\ts$-bimodules $\SA\to\ZA\ts$.
Let us now equip the vector space $\ZA$ with a binary operation
$\mult$ defined by setting
$$
X\mult Y=XP\,Y \quad\text{for}\quad X\com Y\in\ZA\,.
$$
We first define the coset $XP\,Y\in\ZA$ in the case when
$X\com Y\in\SA$ by using the relations \eqref{M4}. Although $P$ is an
infinite sum of elements of $\,\overline{\!\U(\g)\!\!\!}\,\,\,$,
the condition (c) implies that only finitely many summands
of $XP\,Y\ts$ are non-zero cosets in $\ZA$ in this case.
Here we use %the first property in~
\eqref{M55}.
This definition of $XP\,Y$ extends to any cosets 
$X\com Y\in\ZA$ by $\,\overline{\!\U(\h)\!\!\!}\,\,\,\ts$-linearity. Note
that each summand of $P$ commutes with every element of the
subalgebra $ \,\overline{\!\U(\h)\!\!\!}\,\,\,\subset
\,\overline{\!\U(\g)\!\!\!}\,\,\,\ts$.

So $\ZA$ becomes an associative algebra over $\CC\ts$. We
call it the \textit{double coset algebra}. The quotient space
$\Ab\ts/\ts\Jb$ is a left $\ZA\ts$-module relative to an action,
denoted by $\mult$ too, defined~by
\begin{equation}
\label{pact}
X\mult Z=XP\,Z \quad\text{for}\quad X\in\ZA
\quad\text{and}\quad Z\in
%\Jb\backslash
\Ab\ts/\ts\Jb\,.
\end{equation}

Recall that 
$\R\subset\Rb\subset
%\Jb\backslash
\Ab\ts/\ts\Jb\ts$. Now let
$\pi:
%\Jb\backslash
\Ab\ts/\ts\Jb\to\ZA$ be the canonical projection, so that
$$
\pi(A+\Jb)=A+\Jb+\Jmb \quad\text{for}\quad A\in\Ab\ts.
$$

\begin{proposition*}
\label{prop1.1}
\begin{enumerate}
\item[\rm(i)]
The restriction of $\pi$ to $\R$ is a monomorphism of algebras.
\item[\rm(ii)]
The restriction of $\pi$ to $\Rb$ is an isomorphism of algebras
$\Rb$ and\/ $\ZA$.
\end{enumerate}
\end{proposition*}

\begin{proof}
The properties \eqref{M4} of $P$
imply that the assignment $X\to X\ts P$ defines a map 
$\ZA\to\Rb\ts$. Namely, the first of the properties \eqref{M4}
implies that the double coset $X\in\ZA$ gets mapped to
a single coset in 
$%\Jb\,\backslash
\Ab\ts/\ts\Jb\ts$. 
The second property 
implies that the latter coset belongs to $\Rb\ts$. 
Although $P$ is an
infinite sum of elements of $\,\overline{\!\U(\g)\!\!\!}\,\,\,$,
by \eqref{M55} the condition (c) implies that only finitely many summands
of $X\ts P\ts$ are non-zero cosets in 
$%\Jb\backslash
\Ab\ts/\ts\Jb\ts$.

Let us denote by $\pi^*$ this map from $\ZA$ to $\Rb\ts$.
Further, denote by $\pi_*$ the restriction to $\Rb$ of the 
projection map $\pi:%\Jb\backslash
\Ab\ts/\ts\Jb\to\ZA\ts$.
The first property in \eqref{M5} implies the equality 
\begin{equation}
\label{xpx}
X\ts P=X
\quad\text{for any}\quad
X\in\Rb\ts.
\end{equation} 
Hence the composition $\pi^*\ts\pi_*$ 
is the identity map on $\Rb\ts$. The equality \eqref{xpx} 
also implies that $\pi_*$ is a homomorphism of algebras.
Finally, the second property in \eqref{M5} implies that
the composition $\pi_*\ts\pi^*$ 
is the identity map on $\ZA\ts$.
Thus we get (ii), and hence (i).
\qed
\end{proof}

Thus the algebra $\ZA$ contains as a subalgebra a copy $\pi(\R)$ of
the algebra $\R\ts$. It also contains as a subalgebra a copy of
another Mickelsson algebra, defined as a quotient of normalizer of
the left ideal $\Jm\subset\Ar$ relative to this ideal. However, the
latter quotient algebra will not be used in this article. 
Note that both algebras $\Rb$ and $\ZA$ act on the 
quotient vector space 
$%\Jb\backslash
\Ab\ts/\ts\Jb\ts$.
Here $\Rb$ acts via left multiplication from the algebra $\Ab$,
while the action of $\ZA$ is defined by \eqref{pact}.
The proof of Proposition \ref{prop1.1}
also demonstrates that the isomorphism
$\pi:\Rb\to\ZA$ identifies these two actions.

The algebra $\ZA$ contains $\overline{\!\U(\h)\!\!\!}\,\,\,$ as a
subalgebra and is free both as a left and as a right module over it.
Further, there is a weight decomposition of $\ZA$
%$$\ZA\,=\underset{\lambda\in\h^*}\op\,\ZA^\la$$
relative to the adjoint action of $\h$ on $\ZA\ts$,
similar to \eqref{sawede}.
The zero weight component $\ZA^{\ts0}$ is a subalgebra of $\ZA\ts$.
The algebra $\ZA^{\ts0}$ has several natural subalgebras. One
of them is the image $\pi(\R^{0})$ of the zero weight component of
the algebra $\R\ts$. Further, consider the image of the zero weight
component $\Ar^0$ of $\Ar\subset\Ab$ relative to the canonical
projection $\Ab\to\ZA\ts$. This image is a vector subspace in
$\ZA\ts$. Denote by $\SAzero$ the subalgebra in $\ZA$
generated by the elements of this subspace. Then
\begin{equation}
\label{shortchain} 
\pi(\R^0)\subset\SAzero\subset\ZA^{\ts0}\ts.
\end{equation}

%------------------------------------------------------------------------------

\subsection*{\it\normalsize 1.3.\ Algebraic group action}

Let $\G$ be a reductive  algebraic group over $\CC$ with the Lie algebra
$\g\ts$. Suppose there is an action
of the group $\G$ by automorphisms on the algebra $\Ar\ts$,
extending the adjoint action of $\G$ on the subalgebra
$\U(\g)\subset\Ar\ts$. Suppose that the action of $\G$
on $\Ar$ preseves the subspace $\V$, and that
the action of $\G$ on $\V$ is locally finite.
Moreover, suppose that
the action of $\ts\g$ on $\Ar$ corresponding to that of $\G$ coincides
with the adjoint action:
$$
X:\,A\mapsto\,\ad_X(A)=XA-A\ts X
\quad\text{for}\quad
X\in\g
\quad\text{and}\quad
A\in\Ar\ts.
$$

Let $\G_{\ts0}$ be the connected component of $\G$ and  
$\TT\subset\G_{\ts0}$ be the maximal torus of $\G$ with
the Lie algebra $\h\ts$. Let $\Norm\,\TT$ and $\Norm_{\,0}\TT$ 
be the normalizers of $\TT$ in $\G$ and $\G_{\ts0}$ respectively.
The adjoint action of the group $\Norm_{\,0}\TT$ on $\h$ 
identifies the quotient group $\,\Norm_{\,0}\TT\ts/\ts\TT$
with the Weyl group $\Sym$ of $\g\ts$.
Choose a representative $\sih_c\in\Norm_{\,0}\TT$ of $\si_c\ts$.
The elements $\sih_1\lcd\sih_r\in\G$
can be chosen to satisfy the braid relations
\begin{equation}
\label{M8a}
\underbrace{\sih_c\,\sih_d\,\sih_c\,\ldots }_{m_{cd}}=
\underbrace{\sih_d\,\sih_c\,\sih_d\,\ldots }_{m_{cd}}\,
\quad\text{for}\quad c\neq d
\end{equation}
where $m_{cd}$ is the order of the element
$\si_c\ts\si_d\in\Sym\ts$, see \cite{T}.
For any
reduced decomposition $\si=\si_{c_1}\ldots\si_{c_k}$ in $\Sym$ put
\begin{equation*}
%\label{sih}
\sih=\sih_{c_1}\ldots\,\sih_{c_k}.
\end{equation*}
Due to \eqref{M8a} this definition of $\,\sih$ is independent of the
choice of the decomposition~of~$\si\ts$.
Thus we get an action of the braid group of $\g$
by automorphisms $\sih$ of the algebra $\Ar\ts$.

The adjoint action of the group $\Norm\,\TT$ on $\h$ 
identifies $\,\Norm\,\TT\ts/\ts\TT$
with a semidirect product $\ES=\DA\ltimes\Sym$
where $\DA$ is a certain finite subgroup
of $\,\Norm\,\TT\ts/\ts\TT\ts$.
Each element $\tau\in\DA$ acts on $\h$
by permuting $H_1\lcd H_r\ts$. 
%and leaving all central elements of $\g$ fixed.
Hence $\tau$ determines an automorphism 
of the semisimple part 
%(the derived subalgebra) 
%$\g^{\ts\prime}$ 
of $\g\ts$, which permutes $E_1\lcd E_r$ and in the same way
permutes $F_1\lcd F_r\ts$. For each $\tau\in\DA$
choose a representative $\ttau\in\Norm\,\TT$
whose adjoint action on the semisimple part of $\g$ 
yields this automorphism.
Then $\ttau\mapsto\tau$ and $\sih\mapsto\si$
under the canonical map
$$
\Norm\,\TT\to\Norm\,\TT\ts/\ts\TT=\ES\ts.
$$
Moreover, the representatives $\sih_1\lcd\sih_r$ above
can be chosen so that the adjoint action of $\ttau$ on $\G$
permutes them in the same way as $\tau$ permutes $H_1\lcd H_r\ts$.

Regard $\g$ as a subspace
of $\Ar\ts$. Note that because the action of $\sih\in\G$
on the subspace $\h\subset\Ar$ factors through that of $\si\in\Sym\ts$,
the action of $\sih$ on $\Ar$ maps any vector of weight $\la$
to a vector of weight $\si(\la)\ts$.
Let $\Ar^0\subset\Ar$ be the zero weight subspace
relative to the adjoint action of $\h\ts$. This
subspace is preserved by the action of all the 
elements $\sih\com\ttau\in\G\ts$. Since the subgroup
$\TT\subset\G$ acts on $\Ar^0$ trivially, we get
an action of the group $\ES$ on %the subspace 
$\Ar^0\ts$. 

Since the action of the elements $\sih\com\ttau\in\G$ on $\Ar$ 
preserves the set of denominators~\eqref{M2}, it extends
from $\Ar$ to $\Ab\ts$. Further, the action of $\ttau$ on $\Ar$
preserves the ideals $\J$ and $\Jm$. Hence $\ttau$ acts on 
$%\J\ts\backslash
\Ar\ts/\ts\J$ and on $\SA\ts$.
Similarly, $\ttau$ acts on 
$%\Jb\ts\backslash
\Ab$ and on $\ZA\ts$.

%------------------------------------------------------------------------------

\vspace{-8pt}%%%%%%%%%%%%%%%%%%%%%%%%%%%%%%%%%%%%%%%%%%%%%%%%%%%%%%%%%%%%%%%%%%

\subsection*{\it\normalsize 1.4.\ Zhelobenko operators}

For each $c=1\lcd r$ define a linear map $\xi_{\ts c}:\Ar\to\Ab$
by setting for any $A\in\Ar\ts$
\begin{equation}
\label{M10} 
\xi_{\tts c}\ts(A)=A+\,\sum_{s=1}^\infty\,\, 
\bigl(\ts
s\ts!\, H_c(H_c-1)\ldots(H_c-s+1) \bigr)^{-1}\ts E_c^{\ts s}\,
\ad^{\ts s}_{\ts F_c}(A\ts)
\end{equation}
By (c) for any given $A\in\Ar$ only finitely many terms of the sum
\eqref{M10} differ from zero, hence the map $\xi_{\ts c}$ is well
defined. The definition \eqref{M10} and the next two propositions go
back to \cite[Section 2\ts]{Z1}. See \cite[Section 3]{KN1} for
detailed proofs of these propositions.

\begin{proposition*}
\label{prop1} For any\/ $H\in\h$ and\/ $A\in\Ar$ we have
\begin{align*}
\xi_{\ts c}(HA)
&\in (\ts H+\al_c(H))\,\ts\xi_{\ts c}(A)\ts+\ts\Jb\ts,
\\
\xi_{\ts c}(A\ts H)
&\in\,\xi_{\ts c}(A)\ts(\ts H+\al_c(H))\ts+\ts\Jb\ts.
\end{align*}
\end{proposition*}

Proposition \ref{prop1} allows us to define a linear map
$\bar\xi_{\ts c}:\Ab\to
%\Jb\,\ts\backslash\,
\Ab\ts/\ts\Jb$ by setting
$$
\bar\xi_{\ts c}(XA)=Z\,\xi_{\ts c}(\ts A)\ts+\ts\Jb\quad\
\text{for}\quad X\in\,\overline{\!\U(\h)\!\!\!}\,\,\,
\quad\text{and}\quad A\in\Ar\ts,
$$
where the element $Z\in\,\overline{\!\U(\h)\!\!\!}\,\,\,$ is defined
by the equality
$Z(\mu)=X(\ts\mu+\al_c)$
for %\quad\ \text{for}\quad
$\mu\in\h^\ast$
when $X$ and $Z$ are regarded as rational functions on $\h^\ast\ts$.

The action of the Weyl group $\Sym$ on the algebra $\U(\h)$ extends
to an action on $\,\overline{\!\U(\h)\!\!\!}\,\,\,\ts$, so that for
any $\si\in\Sym$
$$
(\ts\si\ts X)(\mu)=X(\ts\si^{-1}(\mu))
$$
if %the element 
$X\in\,\overline{\!\U(\h)\!\!\!}\,\,\,$ is regarded
as a rational function on $\h^\ast\ts$. 
Take the images
$\sih_c(\ts\Jb\ts)$~and~$\sih_c(\ts\Jmb\ts)$
relative to the action of $\sih_c\in\G$ on $\Ab\ts$;
they are respectively right and left ideals of $\Ab\ts$.

\begin{proposition*}
\label{prop2} We have\/ $\sih_c(\ts\Jb\ts)\subset\Ker\ts\,\bar\xi_{\ts
c}$ and\/ $\,\bar\xi_{\ts
c}(\,\sih_c(\Jmb\ts))\ts\subset\ts\Jmb\ns+\ts\Jb\ts$.
\end{proposition*}

This proposition allows us to define the linear maps
$$
\xic_{\ts c}:%\,\Jb\ts\backslash
\Ab\ts/\ts\Jb\to%\Jb\ts\backslash
\Ab\ts/\ts\Jb
\quad\text{and}\quad \xic_{\ts c}:\,\ZA\to\ZA
$$
as the composition $\xib_{\ts c}\,\sih_c$ applied to elements of
$\Ab$ which are taken modulo $\Jb\ts$ and $\Jb+\Jmb$ respectively. 
%Denoting both maps by the same symbol $\xic_{\ts c}$
%should cause no confusion.
In their present form, the maps $\xic_{\ts c}$ have been
introduced in \cite{KO}. 
%We will call them \textit{Zhelobenko operators}. 
The next proposition states their key properties.
For the proof of Part (i) see \cite[Section~6]{Z1}; 
for the proofs of (ii) and (iii) see \cite[Section~5]{KO}.

\begin{proposition*}
\label{proposition1.3}
{\rm(i)}
The maps\/ $\xic_1\lcd\xic_r$ of\/ 
$%\Jb\ts\backslash
\Ab\ts/\ts\Jb$ and\/
$\ZA$ satisfy the braid relations
$$
\underbrace{\xic_c\,\xic_d\,\xic_c\,\ldots}_{m_{cd}}\,=\,
\underbrace{\xic_d\,\xic_c\,\xic_d\,\ldots}_{m_{cd}}\,
\quad\text{for}\quad c\neq d\ts.
$$
{\rm(ii)}
The map\/ $\xic_{c}$ of\/ $\ZA$ is an automorphism of the
double coset algebra:
$$
\xic_c(X\mult Y)\,=\,\xic_c(X)\,\mult\,\xic_c(Y)
\quad\text{for}\quad X\com Y\in\ZA\,.
$$
{\rm(iii)}
The map\/ $\xic_{c}$ of\/ 
$%\Jb\ts\backslash
\Ab\ts/\ts\Jb$ is $\ZA\ts$-equivariant 
%$\ZA\ts$-module map 
in the following sense:
$$
\xic_c(X\mult Z)=\xic_c(X)\,\mult\,\xic_c(Z) \quad\text{for}\quad
X\in\ZA \quad\text{and}\quad Z\in%\Jb\ts\backslash
\Ab\ts/\ts\Jb\,.
$$
\end{proposition*}

Part (i) of Proposition \ref{proposition1.3} implies that
for any reduced decomposition of an element
$\si=\si_{c_1}\ldots\ts\si_{c_k}$ in the Weyl group $\Sym$ the map
\begin{equation}
\label{xisi}
\xic_{\ts\si}=\xic_{\ts c_1}\ldots\,\xic_{\ts c_k}
\end{equation}
of $%\Jb\ts\backslash
\Ab\ts/\ts\Jb$ or $\ZA$ does not depend on the choice of the
decomposition. In view of Part (ii), any map $\xic_{\ts\si}$ of $\ZA$
will be called a \textit{Zhelobenko automorphism\/} of $\ZA\ts$.
Any map $\xic_{\ts\si}$ of $\Ab\ts/\ts\Jb$
will be called a \textit{Zhelobenko operator\/} on
$\Ab\ts/\ts\Jb\ts$. Then using the same symbol
$\xic_{\ts\si}$ for the maps of both $\Ab\ts/\ts\Jb$
and $\ZA$ should cause no confusion.
By Proposition~\ref{prop1}
\begin{align}
\label{M12}
&\xic_{\ts\si}(H\ts Z)=(\si\circ H)\,\xic_{\ts\si}(Z)\ts,
\\
\label{M21}
&\xic_{\ts\si}(Z\ts H)=\xic_{\ts\si}(Z)\,(\si\circ H)
\end{align}
for any $\si\in\Sym$, $H\in\h$ and $Z\in%\Jb\backslash
\Ab\ts/\ts\Jb$ or
$Z\in\ZA\ts$. Here for $c=1\lcd r$ we use the equality
$$
\si_c\ts(\mu+\al_c)=\si_c\circ\mu \quad\text{for}\quad \mu\in\h^*.
$$
The relations \eqref{M12} imply that all the operators $\xic_{\ts\si}$
on $\ZA$  preserve the %zero weight 
subspace $\ZA^{\ts0}$.

We will complete this subsection with an explicit formula
for the operator $\xic_c$ on $\ZA\ts$. 
Let us denote by $\g_{\tts c}$ the 
$\mathfrak{sl}_2\ts$-subalgebra of 
$\g$ spanned by the elements $E_c$, $F_c$ and $H_c\ts$. 
It acts on the vector space $\Ar$,
by restricting to $\g_{\tts c}$ the adoint action of $\g$ on $\Ar\ts$.
As an $\g_{\tts c}\ts$-module, $\Ar$
splits into a direct sum of irreducible finite-dimensional
submodules. Choose $j\in\{\ts0,{\textstyle\frac12},1,\ldots\ts\}$
and take any element $A\in\Ar$ from an irreducible $\g_{\tts c}\ts$-submodule
of dimension $2j+1\ts$. Suppose that $A$ is of weight $2\ts h$
relative to the action of $H_c\ts$,
that is $[\ts H_c\com A\ts]=2\ts h\ts A$ for
$h\in\{-j\com -j+1\lcd j-1\com j\}\ts$.
Since $\sih_c$ is an automorphism of the algebra $\Ar\ts$,
the element $\sih_c(A)\in\Ar$ belongs to
an irreducible $\g_{\tts c}\ts$-submodule of the same dimension $2j+1\ts$,
and is of weight $-\ts2h$ relative to the action of $H_c\ts$. 
%Here we also use the condition (f).
Therefore
the operator $\xic_c$ on $\ZA$ maps the coset of~$A$ to the coset of
%the following element of $\Ab\ts$,
\begin{gather*}
\sum_{s=0}^\infty
\,\,\bigl(\ts
s\ts!\, H_c(H_c-1)\ldots(H_c-s+1) \bigr)^{-1}\ts E_c^{\ts s}\,
\ad^{\,s}_{\,F_c}(\sih_c(A))\,=
\\
\sum_{s=0}^{j-h}
\,\bigl(\ts
s\ts!\, H_c(H_c-1)\ldots(H_c-s+1) \bigr)^{-1}\ts E_c^{\ts s}\,
\ad^{\,s}_{\,F_c}(\sih_c(A))\,.
\end{gather*}
By the definition of subspace $\Jb\subset\Ab\ts$,
the coset of this sum coincides with the coset of
$$
\sum_{s=0}^{j-h}
\,\bigl(\ts
s\ts!\, H_c(H_c-1)\ldots(H_c-s+1) \bigr)^{-1}\ad^{\,s}_{\,E_c}\,
\ad^{\,s}_{\,F_c}(\sih_c(A))\ts.
$$
But the coset of the last sum in $\ZA$ coincides with the coset of
\begin{equation}
\label{sumprod}
\sum_{s=0}^{j-h}\,\,
\prod_{i=0}^{s-1}\,\,\frac{(\ts j-h-i)\,(\ts j+h+i+1)}{(\ts i+1)\,(H_c-i)}
\,\,\,\sih_c(A)\,;
\end{equation}
see for instance the proof of \cite[Proposition 3.7]{KN1}.
The sum over $s=0\lcd j-h$
of the products of fractions in \eqref{sumprod}
is a particular value $\mathrm{F}\ts(\ts h-j\com\ts j+h+1\com-H_c\,;1)$
of the hypergeometric function. %\cite[Section I.2.1]{BE}.
Using the well known
formula % (14) Subsection 2.1.3
$$
\label{gaga} 
\mathrm{F}\ts(u,v,w\ts;1)=
\frac{\,\Upga(w)\,\Upga(w-u-v)}
{\,\Upga(w-u)\,\Upga(w-v)}
$$
valid for $u\com v\com w\in\CC$ with $w\neq0,-1,\ldots$
and $\mathrm{Re}\,(w-u-v)>0\ts$, the sum \eqref{sumprod}~equals
\begin{gather*}
\frac{\Upga(-H_c)\,\Upga(-H_c-2h-1)}
{\Upga(-H_c-h+j)\,\Upga(-H_c-h-j-1)}
\,\,\,\sih_c(A)\,=\,
\prod_{i=1}^{j-h}\,\frac{H_c+2h+i+1}{H_c-i+1}\,\,\,\sih_c(A)
\\[2pt]
=\ \prod_{i=1}^{j-h}\,(H_c-i+1)^{-1}\ts\cdot\ts\sih_c(A)\cdot
\prod_{i=1}^{j-h}\,(H_c+i+1)\ts.
\end{gather*}

%------------------------------------------------------------------------------

\subsection*{\it\normalsize 1.5.\ Invariants of Zhelobenko automorphisms}

Consider the subspace $\Ar^\G$ of $\G\ts$-invariants in $\Ar\ts$.
Define the linear map $\ga:\Ar^\G\to\SA$ as the restriction to 
$\Ar^\G$ of the canonical projection $\Ar\to\SA\ts$.
The natural embedding of algebras $\Ar\to\Ab$
determines an embedding of $\U(\h)\ts$-bimodules $\SA\to\ZA.$
Consider the zero weight component $\SA^{\ts0}$ of the vector space
$\SA$. Here we again refer to the adjoint action of $\h$ on
$\SA\ts$. The definition of the action of $\xic_{\ts\si}$ and $\ttau$ on $\ZA$
implies immediately
that the image of the map $\ga$ is contained in the subspace of $\SA\ts$,
$$
\Dwo=\{\,Z\in\SA^{\ts0}
\ |\ 
\xic_{\ts\si}(Z)=Z
\ \,\text{and}\ \,
\ttau\ts(Z)=Z
\ \,\text{for all}\ \,
\si\in\Sym
\ \,\text{and}\ \,
\tau\in\DA \,\}\,.
$$
By \cite[Remark 3]{KNV} this $\Dwo$ is
a subalgebra of $\ZA\ts$, and is included in the chain \eqref{shortchain}:
\begin{equation}
\label{chain}
\Dwo\subset\pi(\R^0)\subset\SAzero%\subset\ZAstar
\subset\ZA^{\ts0}\ts.
\end{equation}
The next proposition has been also proved in \cite{KNV}. 

\begin{proposition*}
\label{knv}
The map $\ga:\Ar^\G\to\SA$ is injective,
and its image is equal to $\Dwo\ts$. 
\end{proposition*}

%In this subsection,
%we will consider the action of $\Dwo$ on $\SAzero\ts$-modules
%by restriction.

The squares of the Zhelobenko automorphisms $\xic_1\lcd\xic_r$ 
of $\ZA$ are given by the formula
$$
%\label{koceq}
\xic_c^{\ts2}(Z)=(H_c+1)\,\sih_c^{\ts2}(Z)\,(H_c+1)^{-1}
\quad\text{for all}\quad
Z\in\ZA\ts,
$$
see \cite[Corollary 7.5]{KO}. Here $\sih_c^{\ts2}\in\TT$
by definition, so that the squares $\sih_c^{\ts2}$ 
and hence $\xic_c^{\ts2}$ 
act trivially on the zero weight subspace
$\ZA{}^{\,0}\subset\ZA\ts$. 
This means that the restrictions of all the operators $\xic_{\ts\si}$ 
to $\ZA{}^{\,0}$ determine an action of the Weyl group $\Sym\ts$.
Furthermore, for any $\tau\in\DA$
the adjoint action of $\ttau\in\G$ on $\g$
permutes the operators $\xic_1\lcd\xic_r$ on $\ZA$
in the same way, as it permutes the elements $H_1\lcd H_r$ of $\h\ts$.
This implies that the latter action of the group $\Sym$ on $\ZA{}^{\,0}$
extens to that of the semidirect product $\ES=\DA\ltimes\Sym\ts$.

We call a weight $\mu\in\h^*$ 
\textit{nonsingular\/} if
$\mu(H_\al)\neq-1,-2,\,\ldots\,$ for all $\al\in\De^+$.
Fix any $\lambda\in\h^*$ such that
$\lambda+\rho$ is nonsingular. Let $N$ be any
irreducible $\SAzero\ts$-module of weight $\lambda\ts$. 
The latter condition means that the subalgebra
$\U(\h)\subset\SAzero$ acts on $N$ via the mapping $\la:\h\to\CC\ts$.
This mapping extends to a homomorphism $\U(\h)\to\CC\ts$,
also denoted by $\la\ts$.

Due to \eqref{chain}, $N$ is a module over the algebra $\Dwo$
by restriction. Let $\Sym_\la$ and $\ES_\la$ denote
the stabilizers of $\lambda\in\h^*$
in $\ES$ and $\Sym$ relative to 
the shifted actions of these groups on $\h^*$; 
$$
\Sym_\la=\{\ts\si\in\Sym\ |\ \si\circ\la=\la\ts\}\ts.
$$ 
The shifted action of $\Sym$ on $\h^*$ extends to an action
of $\ES$ since $\tau(\rho)=\rho$ for every $\tau\in\DA\ts$.

\begin{proposition*}
\label{proposition3.9} 
Suppose that $\la+\rho\in\h^*$ is nonsingular, and that\/ 
$\ES_\la=\Sym_\la\ts$. 
Then $N$ is an irreducible $\Dwo$-module.
\end{proposition*}

We shall prove Proposition~\ref{proposition3.9} 
in the remainder of this subsection. Let
$\chi:\SAzero\to\End N$ the defining homomorphism
of the $\SAzero\ts$-module $N$.
For the proof of the proposition it is sufficient to find for each
$Y\in\SAzero$ an element $Z\in\Dwo$
such that $\chi\ts(Y)=\chi\ts(Z)$.
Here we will assume that $Y$ is the image of some element
$A\in\Ar^0$ under the canonical projection
$\Ab\to\ZA\,$; see the definition of the subalgebra
$\SAzero$ given at the end of Subsection~1.2.

\begin{lemma*}
\label{19-5}
There exists $X\in\U(\h)$ such that $\la(X)\neq0$
and for any $\si\in\Sym$ the element 
$\xic_{\ts\si}(X\ts Y)\in\ZA$ belongs to the image of $\Ar^0$
under the canonical projection
$\Ab\to\ZA\ts$.
\end{lemma*}

\begin{proof}
All elements $\sih\ts(A)$
with $\si\in\Sym$ belong to some finite-dimensional submodule $M$
of $\Ar^0$ under the adjoint action of $\g\ts$.
For every positive root 
$\al\in\De^+$ 
we can choose a non-negative integer $n_\al$ such that
\begin{equation}
\label{ourass}
\ad_{F_\al}^{\,n_\al+1}(B)=0
\quad\text{for all}\quad
B\in M\ts.
\end{equation}
Put
\begin{equation}
\label{ah}
X=\prod_{\al\in\De^+}\prod_{s=1}^{n_\al}\,
(H_\al+\rho\ts(H_\al)+s\ts)\,.
\end{equation}
Then we have $\la(X)\neq0\ts$, because the weight $\la+\rho$ is 
nonsingular by our assumption. 

Now take any element $\si\in\Sym$ and a reduced decomposition
$\si=\si_{c_1}\ldots\ts\si_{c_k}$. Denote 
$$
\om_{\ts l}=\si_{c_1}\ldots\si_{c_l}%\in\Sym
\quad\text{and}\quad
\be_{\ts l}=\om_{\ts l-1}(\alpha_{\ts c_l})
\quad\text{for}\quad
l=1\lcd k\ts.
$$
Note that $\be_1\lcd\be_{\ts k}$ are all those roots $\al\in\De^+$
for which $\si^{-1}(\al)\notin\De^+$.
By \eqref{M12},
$$
\xic_{\ts\si}(X\ts Y)=(\si\circ X)\,\xic_{\ts\si}(Y)\ts.
$$
In each factor $\xic_{\ts c_l}$ of $\xic_{\ts\si}$ in the product \eqref{xisi}
we can replace, modulo the left ideal $\Jmb$ of $\Ab$, every entry of
$E_{\ts c_l}$ by the corresponding adjoint operator. Using the property
\eqref{M12} repeatedly, in $\xic_{\ts\si}(Y)$ we can also move to the left
all denominators from $\U(\h)$. Hence $\xic_{\ts\si}(Y)$
equals to the coset in $\ZA$ of the image of $\sih\ts(A\ts)$
under the operator
$$
\sum_{s_1,\ldots,s_k=0}^\infty\ \prod_{l=1}^k\ 
\bigl(s_l\ts!\,
(\ts\om_{\ts l-1}\circ H_{\ts c_l})
\ldots
(\ts\om_{\ts l-1}\circ H_{\ts c_l}-s_l+1)
\bigr)^{-1}
\cdot
\ad^{\,s_1}_{\ts E_{\be_1}}\!\ad^{\,s_1}_{\ts F_{\be_1}}
\!\ldots\,
\ad^{\,s_k}_{\ts E_{\be_k}}\!\ad^{\,s_k}_{\ts F_{\be_k}}.
$$
By the assumption \eqref{ourass}, all summands above
with at least one index $s_l>n_{\be_{\ts l}}$ vanish.
Further, for all $l=1\lcd k$ and $s=1\com2\com\ldots$ we have
\begin{align}
\notag
\om_{\ts l-1}\circ H_{\ts c_l}-s+1
&=H_{\ts\be_{\ts l}}+(\ts\om_{\ts l-1}^{\,-1}(\rho)-\rho\ts)(H_{\ts c_l})-s+1
\\
\notag
&=H_{\ts\be_{\ts l}}+\rho\ts(H_{\ts\be_{\ts l}})-1-s+1
\\
\label{cancel}
&=H_{\ts\be_{\ts l}}+\rho\ts(H_{\ts\be_{\ts l}})-s\ts.
\end{align}

Now for $l=1\lcd k$ consider the factor in \eqref{ah}
corresponding to the positive root 
$\al=-\ts\si^{-1}(\ts\be_{\ts l})$
and an index $s=1\lcd n_\al\ts$.
The shifted action of $\si$ on this factor~yields  
$$
\si\circ H_{\al}+\rho\ts(H_{\al})+s=
-H_{\be_{\ts l}}+(\si^{\,-1}(\rho)-\rho\ts)(H_\al)+\rho\ts(H_{\al})+s=
-H_{\be_{\ts l}}-\rho\ts(H_{\be_{\ts l}})+s\ts.
$$
The yielded factor cancels, up to the minus sign, 
the factor \eqref{cancel} in the denominator of $\xic_{\ts\si}(Y)\ts$.
Hence the product $(\si\circ X)\,\xic_{\ts\si}(Y)\in\ZA$
belongs to the image of $\Ar^0\ts$.
\qed
\end{proof}

Take the element $X\in\U(\h)$ from Lemma \ref{19-5}.
We have 
$\chi\ts(X\ts Y)=\la(X)\,\chi\ts(Y)$ because the $\SAzero\ts$-module
is of weight $\la\ts$ by our assumption. Hence
$$
\chi\ts(Y)=\chi\ts(X\,Y/\la(X))\,.
$$
Replacing the given element $Y\in\SAzero$
by $X\,Y/\la(X)\ts$, we may from now
assume that for any $\si\in\Sym$ the element
$\xic_{\ts\si}(Y)\in\ZA$ belongs to the image of $\Ar^0$
under the projection~$\Ab\to\ZA\ts$.
  
Consider the orbit of the weight $\lambda$ relative to the shifted action of
the group $\ES$ on $\h^*$,
$$
{\mathcal O}_\la=\{\ts\om\circ \lambda\ |\ \om\in\ES\ts\}\,.
$$
This is a finite subset of $\h^*$. Hence there exists an element
$X'\in\U(\h)$, such that $\la(X')=1$ and
$\mu(X')=0$ for any weight $\mu\in{\mathcal O}_\la$
with $\mu\not=\lambda\ts$. Put
\begin{equation}
\label{essym}
Z=|\Sym_\lambda|^{-1}\sum_{\om\in\ES}\xic_{\ts\om}(X'\,Y)\,.
\end{equation}
Here $\xic_{\ts\om}$ denotes the composition $\ttau\,\ts\xic_{\ts\si}$
of operators on $\ZA\ts$, if $\om=\tau\,\si$ for some 
$\si\in\Sym$ and $\tau\in\DA\,$.
Due to our assumptions on the given element $Y$, we have $Z\in\SAzero\ts$.
Moreover, for every $\om\in\ES$ we have $\xic_{\ts\om}(Z)=Z\ts$. 
%due to \eqref{koceq}. 
Thus $Z\in\Dwo\ts$. Let us show that $\chi\ts(Y)=\chi\ts(Z)$.

Until the end of this section, the symbol $\equiv$ will indicate
equalities in the algebra $\ZA$ modulo the right ideal
generated by all the elements $H-\lambda(H)$ with $H\in\h\ts$.
Firstly, take any $\om\in\ES$ such that $\om\circ\la\not=\la\ts$.
By our choice of the element $X'$, then we have
$$
\xic_{\ts\om}(X'\,Y)=(\ts\om\circ X')\,\xic_{\ts\om}(Y)
\equiv
\la(\ts\om\circ X')\,\xic_{\ts\om}(Y)=
(\om^{-1}\circ\la)(X')\,\xic_{\ts\om}(Y)=0\ts.
$$
Hence $\chi\ts(\xic_{\ts\om}(X'\,Y))=0$ unless
$\om\in\ES_\lambda=\Sym_\lambda\ts$. Let us now prove another lemma. 

\begin{lemma*}
\label{lemma3.10}
For any $\si\in\Sym_\la$ we have $\xic_{\ts\si}(Y)\equiv Y\ts$.
\end{lemma*}

\begin{proof}
By \cite[Proposition V.3.2]{B}
the subgroup $\Sym_\lambda\subset\Sym$ is generated by the reflections 
$\si_\al$ corresponding to the positive roots $\al\in\De^+$ 
such that $(\lambda+\rho)(H_\al)=0\ts$. Fix such~an~$\al$
and write $\al=\si(\al_c)$ where $\al_c$ is a simple root and $\si\in\Sym$.
Since $Y\in\ZA{}^{\ts0}$, then we have
$$
\xic_{\ts\si_\al}(Y)=(\ts\xic_{\ts\si}\,\xic_{c}\,\xic_{\ts\si^{-1}})(Y)\ts.
$$ 

Consider the element $Y'=\xic_{\ts\si^{-1}}(Y)$ of $\ZA\ts$.
By our assumptions on $Y$, the element $Y'$
is the image of some element $A^{\ts\prime}\in\Ar^0$ under the projection
$\Ab\to\ZA\ts$. Consider the $\mathfrak{sl}_2\ts$-subalgebra
$\g_{\tts c}\subset\g\ts$. 
First suppose that $A^{\ts\prime}$ belongs
to an irreducible $\g_{\tts c}\ts$-submodule of $\Ar\ts$.
If $A^{\ts\prime}\neq0$, then this submodule has an odd dimension,
say $2j+1\ts$. By the calculation made in the end of
Subsection 1.4, the operator $\xic_c$
maps the coset of $A^{\ts\prime}$ in $\ZA$ to
$$
(-1)^{\ts j}\,\prod_{i=1}^j\,\frac{H_c+i+1}{H_c-i+1}\ts\cdot\ts Y'\ts.
$$
Here we also used the following observation: 
because $A^{\ts\prime}\in\Ar^0$ belongs to an
irreducible $\g_{\tts c}\ts$-submodule of dimension $2j+1$, we have
$\sih_c(A^{\ts\prime}) = (-1)^j\,A^{\ts\prime}\ts$. Thus
$$
\xic_{\ts\si_\al}(Y)=
(-1)^{\ts j}\,\prod_{i=1}^j\,\frac{\si\circ H_c+i+1}{\si\circ H_c-i+1}\ts
\cdot\ts\xic_{\ts\si}(Y')=
(-1)^{\ts j}\,\prod_{i=1}^j\,\frac{\si\circ H_c+i+1}{\si\circ H_c-i+1}\ts
\cdot\ts Y\ts.\hspace{-20pt}
$$  
But
\begin{align*}
\la(\si\circ H_c)
&=(\si^{-1}\circ\lambda)(H_c)=
(\si^{-1}(\lambda+\rho)-\rho)(H_c)
\\
&=(\lambda+\rho)(\si(H_c))-
\rho\ts(\si(H_c))=
(\lambda+\rho)(H_\al)-1=-1\,,
\end{align*}
so that
$$
\xic_{\ts\si_\al}(Y)\equiv
(-1)^{\ts j}\,\prod_{i=1}^j\,\frac{-1+i+1}{-1-i+1}
\cdot\ts Y=Y\ts.
$$
The assumption that $A^{\ts\prime}$ belongs
to an irreducible $\g_{\tts c}\ts$-submodule
can now be~removed without any loss of generality.
Lemma \ref{lemma3.10} is thus proved 
for any $\si=\si_\al\ts$.
\qed
\end{proof}

Using Lemma \ref{lemma3.10}, for any $\si\in\Sym_\la$ we now get
$$
\xic_{\ts\si}(X'\,Y)
=(\si\circ X')\,\xic_{\ts\si}(Y)\equiv
\la\ts(\ts\si\circ X')\,Y
=(\si^{-1}\circ\la)(X')\,Y=\la(X')\,Y\equiv X'Y\ts.
$$
This %argument 
completes the proof of the equality
$\chi\ts(Y)=\chi\ts(Z)$ and hence
that of Proposition~\ref{proposition3.9}.

%------------------------------------------------------------------------------

\vspace{-6pt}%%%%%%%%%%%%%%%%%%%%%%%%%%%%%%%%%%%%%%%%%%%%%%%%%%%%%%%%%%%%%%%%%%

\subsection*{\it\normalsize 1.6.\ Irreducible\/ {\rm S}-modules}

In this subsection we will introduce a class of
irreducible $\SAzero\ts$-modules, to which we can then 
apply Proposition~\ref{proposition3.9}. Let $K$ be any
left $\Ar\ts$-module. Since $\U(\g)$ is a subalgebra of $\Ar\ts$,
we can regard $K$ as a $\g\ts$-module by restriction.
Suppose the action of the subalgebra $\h\subset\g$ on $K$
is semisimple, so that $K$ splits into direct sum
of the weight subspaces $K^\la$ where $\la$ ranges over $\h^*$.
Also suppose the action of the subalgebra $\np\subset\g$
on $K$ is locally nilpotent, so that for any
$v\in K$ there is a positive integer $s$ such that
$(\np)^{\ts s}\,v=\{0\}\ts$.
%belongs to the category $\mathcal{O}$ of
%Bernstein, Gelfand and Gelfand \cite{BGG}. By definition,
%$K$ is finitely generated over $\U(\g)$, is locally finite over $\U(\np)$, 
%and splits into direct sum
%of the weight subspaces $K^\la$ where $\la$ ranges over $\h^*$.

Consider the space $K_\n$ of $\n\ts$-coinvariants of $K$.
Take any $\la\in\h^*$ and consider the weight subspace
$K_\n^\la\subset K_\n\ts$. Denote by $N$ the subspace
of $K$ consisting of all $\np$-invariants of weight $\la\ts$.
By restricting the canonical projection $K\to K_\n$ to
the subspace $N\subset K\ts$, we get a natural 
linear map $N\to K_\n^\la\ts$. In general, the restriction map
may be not bijective. 

Now suppose that the weight
$\la+\rho$ is nonsingular. It turns out that then
the map $N\to K_\n^\la$ is bijective. To define the inverse map,
take any coset $f\in K_\n^\la$
and choose its representative $u\in K$. We may assume
that $u\in K^\la$. Then we can define
a linear map $K_\n^\la\to N$ by mapping $\ts f\mapsto P\,u\ts$.
Here we use the properties \eqref{M4}, the second property in \eqref{M55},
the local nilpotency of $K$ relative to $\np$, and the nonsingularity
of $\la+\rho\ts$. This map is the right and left inverse
to the projection $N\to K_\n^\la\,$, due to
the first and the second properties in \eqref{M5} respectively.
See also the proofs of the properties \eqref{M4} we~referred~to. 

Since $\la+\rho$ is nonsingular,
the weight subspace $K_\n^\la\subset K_\n$ has 
a structure of a module over the algebra $\SAzero\ts$. 
The action of $\SAzero$ on $K_\n^\la$ will be
denoted by $\mult$ like in Subsection 1.2, and defined by
\begin{equation*}
X\mult f=XP\ts f \quad\text{for}\quad X\in\SAzero
\quad\text{and}\quad f\in K_\n^\la\,.
\end{equation*}
At the right hand side of the above equality,
we use the action of the algebra $\Ar$ on $K\ts$.

\begin{proposition*}
\label{proplus}
Suppose that $K$ is an irreducible $\Ar$-module,
while the weight $\la+\rho$ is nonsingular.
Then $K_\n^\la$ is an irreducible $\SAzero\tts$-module.
\end{proposition*}

\begin{proof}
Take any two cosets $f,g\in K_\n^\la$
and choose their representatives $u\com v\in K\ts$. 
We may assume that both $u\com v\in K^\la\ts$.
Determine the element $P\,u\in K^\la$ as above. 
Because $K$ is an irreducible $\Ar\ts$-module, we can find
an element $A\in\Ar^0$ such that $A\,P\,u=v\ts$.
Let $X\in\SA$ be the coset of $A\ts$. Then we have
$X\in\SAzero$ and $X\mult f=g\ts$.
\qed
\end{proof}

%==============================================================================

\section*{\bf\normalsize 2.\ Shapovalov forms}
\setcounter{section}{2}
\setcounter{theorem*}{0}
\setcounter{equation}{0}

%------------------------------------------------------------------------------

\subsection*{\it\normalsize 2.1.\ Howe systems}

Suppose that the reductive Lie algebra $\g$ is {\it symmetric}.
That is, $\g$ is equipped with an involutive automorphism
identical on the Cartan subalgebra $\h$ and preserving
each of the nilpotent subalgebras $\n$ and $\np$.
%preserving the triangular decomposition \eqref{hc0}  
We have a \textit{Cartan decomposition\/}
$\g=\g_{\tts+}\op\g_{\tts-}$ where 
$\g_{\tts+}$ and $\g_{\tts-}$ are
the eigenspaces of this automorphism with the eigenvalues $1$ and $-1$.
Then $\g_{\tts+}$ is a Lie subalgebra of $\g$ containing $\h\ts$, while
$$
[\ts\g_{\tts+},\g_{\tts-}]\subset\g_{\tts-}
\quad\text{and}\quad 
[\ts\g_{\tts-},\g_{\tts-}]\subset \g_{\tts+}\,.
$$
Note that then the involutive anti-automorphism $\ep$ of $\g$
preserves the decomposition $\g=\g_{\tts+}\op\g_{\tts-}\ts$.
This decomposition induces decompositions of the nilpotent
subalgebras: %of $\g\ts$: 
we have 
$\n=\n_{\tts+}\op\n_{\tts-}$ and $\np=\np_{\tts+}\op\np_{\tts-}$
where 
$$
\n_{\tts+}=\n\cap\g_{\tts+}\com\ 
\n_{\tts-}=\n\cap\g_{\tts-}\com\  
\np_{\tts+}=\np\cap\g_{\tts+}\com\ 
\np_{\tts-}=\np\cap\g_{\tts-}\ts.
$$

Let $U$ be any finite-dimensional complex vector space
and $\th\in\{1,-1\}\ts$. Depending on whether
$\th=1$ or $\th=-1$,
denote by $\Heist$ the {\it Weyl algebra\/} or the {\it Clifford algebra\/}
of $U\op U^*$.
This is the complex associative unital algebra generated by the elements 
of $U$ and of the dual space $U^*$ subject to the relations
for $u\com v\in\U$ and $u'\com v'\in U^*$
\begin{equation}
\label{H0}
v\,u=\th\,u\,v\com\quad 
v'u'=\th\,u'v'\com\quad
u'u-\th\,u\,u'=u'\ts(u)\,. 
\end{equation}

Equip the vector space $W=U\op U^*$ with the bilinear form $B\ts$,
symmetric if $\th=-1$ and alternating if $\th=1\ts$, such that
$B(u,v)=B(u',v')=0$ and $B(u,u')=u'\ts(u)\ts$.
The relations \eqref{H0} are then equivalent to the relations
for all $w\com w'\in W$
\begin{equation}
\label{H00}
w'w-\th\,w\,w'=B(w,w')\ts.
\end{equation}
Then any isotropic subspace $V\subset W$
generates a subalgebra of $\Heist$, which is
a free commutative algebra if $\th=1$, or free
skew-commutative if $\th=-1$. This subalgebra
will be denoted by $\P\ts(V)$.
In particular, the algebra $\Heist$ contains two distinguished
subalgebras $\P\ts(U)$ and $\P\ts(U^*)\ts$,
generated by the elements of $U$ and of $U^*$ respectively.

Further, equip the vector space $W$
with a grading so that the direct summands 
$U$ and $U^*$ of $W$ have degrees $1$ and $-1$ respectively.
This grading naturally extends to a $\ZZ\ts$-grading
on the algebra $\Heist\ts$. 
In particular, the subalgebras $\P\ts(U)$ and $\P\ts(U^*)$ of $\Heist$
are graded respectively by non-negative and non-positive integers.

Now suppose there is
a homomorphism of associative algebras $\zeta:\U(\g)\to\Heist\ts$.
Further suppose there is an action of the group $\G$ on
the vector space $W$ preserving the bilinear form $B$.
Since the form is preserved, this action
extends uniquely to a action of $\G$ by automorphisms of
the algebra $\Heist\ts$. In particular, we get
an action of the Lie algebra $\g$ on $\Heist\ts$, 
and the weight decomposition
\begin{equation}
\label{wede}
\Heist\,=\underset{\lambda\in\h^*}\op\,\Heist{}^\la\ts.
\end{equation}

Choose a non-degenerate symmetric bilinear form $\langle\ ,\,\rangle$ on
the vector space $U\ts$. %which defines the euclidian structure
Using this form, define a linear map $\ep:W\to W$ such that
$\ep:U\to U^*$, $\ep:U^*\to U$ while
\begin{equation}
\label{epsdef}
\ep(u)(v)=\langle\ts u\com v\rangle
\quad\text{and}\quad
\langle\ts u\com\ep(u')\ts\rangle=u'\ts(u)
\end{equation}
for any $u\com v\in U$ and $u'\in U^*$.
The map $\ep$ is involutive and for all $w\com w'\in W$ satisfies
$$
B(\ep(w),\ep(w'))=-\ts\th\,B(w,w')\ts.
$$
Therefore this map admits a unique extension to
an involutive anti-automorphism of the algebra $\Heist$.
We denote the extension by $\ep$ again.

We shall say that that the homomorphism 
$\zeta:\U(\g)\to\Heist$
and the action of $\G$ on $W$ 
form a {\it Howe system\/} on $U$
if the following six conditions are satisfied\ts:

\begin{enumerate}
\item[(1)]
the map $\zeta:\U(\g)\to\Heist$ is $\G\ts$-equivariant\ts;
\item[(2)]
the  action of $\g$ on $\Heist$ 
corresponding to that of $\G$ is adjoint to $\zeta\ts$;
\item[(3)]
$\ep(\zeta(X))=\zeta({\ep}(X))$ for all $X\in\g\ts$, where
${\ep}:\g\to\g$ is the Chevalley anti-involution\ts;
\item[(4)]
%we have the inclusions
$\zeta(\n)\subset U\cdot\Heist$ 
and
$\zeta(\h)\subset\CC\cdot 1+U\cdot\Heist\ts$;
\item[(5)]
%we have the inclusions
$[\ts\zeta(\g_{\tts+})\com U]\subset U$
and
$[\ts\zeta(\n_{\tts-})\com U^*]\subset U$
while
$[\ts\zeta(\n_{\tts-})\com U]=\{0\}\ts$;
\item[(6)]
any weight element of $\P\ts(U)$ 
has a $\ZZ\ts$-degree uniquely determined by the weight.
\end{enumerate}

The property (2) of a Howe system means that
the action of Lie algebra $\g$ on $\Heist$ 
corresponding to that of $\G$ is given~by
$$
X:Y\mapsto[\zeta(X),Y]
\quad\text{for}\quad 
X\in\g
\quad\text{and}\quad
Y\in\Heist\ts.
$$
The properties (2) and (3) imply that the automorphism $\ep$ of $\Heist$
changes the signs of the degree and of the weight. Here and in (6) 
we refer to the $\ZZ\ts$-grading on $\Heist$ and the weight decomposition
\eqref{wede}.
By combining (3) and (4),(5) we get
two more properties: %of a Howe system:
\begin{enumerate}
\item[(7)]
%we have the inclusions 
$\zeta(\np)\subset \Heist\cdot U^*$ and
$\zeta(\h)\subset\CC\cdot 1+\Heist\cdot U^*$;
\item[(8)]
%we have the inclusions
$[\ts\zeta(\g_{\tts+})\com U^*]\subset U^*$
and
$[\ts\zeta(\np_{\tts-})\com U]\subset U^*$
while
$[\ts\zeta(\np_{\tts-})\com U^*]=\{0\}\ts$;
%\item[(0)]
%any weight element of $\P\ts(U^*)$ 
%has a degree uniquely determined by its weight.
\end{enumerate}
The properties (5) and (8) imply that for 
$X\in\g_{\tts+}\com\n_{\tts-}\com\np_{\tts-}$ 
the adjoint operators $\ad_{\ts\zeta(X)}$ on $\Heist$ have the degrees 
$0\com2\com-2$ respectively.

Later on we will work with known examples \cite{H1}
of Howe systems. The list of (1)~to~(8)
summarizes the common properties of these examples that we shall use.   
Now fix any Howe system on $U$ 
and set 
\begin{equation}
\label{H2}
\Ar=\U(\g)\ot\Heist\ts.
\end{equation}
The group $\G$ acts (diagonally) by automorphisms of the algebra $\Ar\ts$.
Let us identify $\U(\g)$ with subalgebra of $\Ar$ generated
by the elements 
\begin{equation}
\label{xdiag}
X\ot 1+1\ot\zeta(X)
\quad\text{where}\quad
X\in\g\ts.
\end{equation}
Then set 
$$
\V=1\ot\Heist\ts.
$$ 
The condition (2) on a Howe system imply that
the condition (c) on the algebra $\Ar$ is satisfied.
Note that in this case $\V$ is a subalgebra of $\Ar\ts$,
not only a vector subspace.
In %the rest of 
this section
we will investigate the Mickelsson algebras
$\R$ and $\ZA$ corresponding
to the associative algebra \eqref{H2} determined by any Howe system.
%For simplicity of notation 
We will usually
identify any element $X\in\g$ with the element \eqref{xdiag}
of $\Ar$, and any element $Y\in\Heist$ with $1\ot Y\in\Ar\ts$.
Note that for $X\in\g$ the commutator $[X,Y]$ in the algebra $\Ar$
equals $[\ts\zeta(X),Y]\in\Heist\ts$.

%------------------------------------------------------------------------------

\subsection*{\it\normalsize 2.2.\ Shapovalov form on double coset algebra}

Let us extend the map $\ep$ from $\Heist$ to the algebra $\Ar$
defined by \eqref{H2}, so that under the extension
$X\ot Y\mapsto\ep(X)\ot\ep(Y)$ for $X\in\U(\g)$ and
$Y\in\Heist\ts$. Here the symbol $\ep$ in the first tensor factor
denotes the Chevalley anti-involution on $\U(\g)\ts$.
The same symbol $\ep$ in the second tensor factor denotes 
the involutive anti-automorphism of $\Heist$ defined by \eqref{epsdef}.
The extended map 
is an involutive anti-automorphism of the algebra $\Ar\ts$.
Its restriction to the subalgebra $\U(\g)$ generated 
by the elements \eqref{xdiag} coincides
with the Chevalley anti-involution on $\U(\g)\ts$,
due to the property (3) of a Howe system. 
This anti-automorphism further extends from $\Ar$ to $\Ab\ts$.

Since the latter extension preserves the subspace $\Jb+\Jmb\subset\Ab\ts$,
it defines an involutive linear map 
$\varepsilon:\ZA\to\ZA\ts$. Because
$\varepsilon(\p)=\p$ where $\p$ is the extremal projector for $\g\ts$, 
this map is an anti-automorphism of the double coset algebra $\ZA\ts$. 
We will denote it by the same symbol $\ep\ts$.
Clearly,
$\ep$ maps the subalgebra $\SAzero\subset\ZA$ to itself.
Denote by $\JV$ the right ideal of the algebra $\ZA\ts$,
generated by the cosets of all elements  $u\in U$.
Similarly, denote by $\IV$ the left ideal of the algebra $\ZA\ts$,
generated by the cosets of all elements  $u'\in U^*$.
We identify the elements $u\com u'\in\Heist$ with
the elements $1\ot u\com1\ot u'\in\Ar$ respectively. 

\begin{lemma*}
\label{lemma2}
{\rm(i)}
$\JV$  is spanned by the cosets of products 
$u\,X$ in $\Ab$ where $u\in U$, $X\in\Ab\ts$.
\\
{\rm(ii)}
$\IV$ is spanned by the cosets of products 
$X\ts u'$ in $\Ab$ where $u'\in U^*$, $X\in\Ab\ts$.
\end{lemma*}

\begin{proof}
Take any $u\in U$ and $X\in\Ab\ts$.
The property (5) of a Howe system implies that 
for any $A\in\n$ the commutator $[\ts u\ts,\ns A]\in U\ts$.
We may assume that $u$ is a weight element, say $u\in U^\la$
for some $\la\in\h^*$. If
$A\in\n$ is also a weight element then 
$[u\ts,\ns A]$ has a weight less than $\la\ts$.
The inequality
$\mu<\la$ for some $\mu\in\h^*$ 
means that the difference $\la-\mu$ is
a non-zero sum of simple roots in $\De^+$ 
with non-negative integral coefficients.
By the first relation~in~\eqref{M5},
\begin{equation}
\label{Sh3}
u\,\p\ts X\,\in\,u\,X+\sum_{v}\,v\ts Y+\,\Jb
\end{equation}
where $v$ ranges over a certain finite subset of $U$ and
has a weight $\mu<\la\ts$, while
$Y\in\Ab$ corresponds to $v$.
But \eqref{Sh3} shows that the product of the cosets
of $u$ and of $X$ in the algebra $\ZA$ is the coset of %the sum
$$
u\,X+\sum_{v}v\,Y\ts.
$$

Furthermore, we can recursively invert \eqref{Sh3} and 
get the relation in the algebra $\Ab\ts$,
$$
u\,X\,\in\,u\,\p\ts X+\sum_{w}\,w\ts\p\ts Z+\,\Jb
$$
where $w$ ranges over some finite subset of $U$ and
has a weight less than $\la\ts$, while
$Z\in\Ab$ corresponds to $w\ts$. So
the coset in $\ZA$ of the product $u\,X$ in the algebra $\Ab$
is the coset of
$$
\,u\,\p\ts X+\sum_{w}\,w\ts\p\ts Z
$$
and therefore belongs to the ideal $\JV$ of the algebra $\ZA\ts$.
This completes the proof of the part (i) of the lemma.
The proof of the part (ii) is very similar, and is omitted here.
\qed
\end{proof}

\begin{corollary*}
For any $Z\in\ZA$ there is a unique element $(Z)_{\tts0}\in\Uh$ such that
$$
Z-(Z)_{\tts0}\in\JV+\IV\ts.
$$
\end{corollary*}

\begin{proof}
Choose any ordered basis in each of the vector spaces 
$\n\com\n'\com U\com U'$ and $\h\ts$. By the classical
Poincar\'e\ts-Birkhoff\ts-Witt theorem, each of
the vector spaces $\U(\n)\com\U(\np)\com\P\ts(U)$, $\P\ts(U^*)$
and $\U(\h)$ has a basis made of all monomials in the corresponding
basis vectors. If $A\com A'\com B\com B^{\ts\prime}$ and $C$
run through these monomials, then
the products $A\ts B \ts C\ts B^{\ts\prime}A'$ form a basis in $\Ar\ts$.
So the vector space $\ZA$ is spanned by the cosets
of products in $\Ab$ of the form $B \ts D\ts B^{\ts\prime}$
where $B$ and $B^{\ts\prime}$ are as above whereas $D\in\Uh\ts$.
By Lemma \ref{lemma2}, any such a coset~belongs to $\JV+\IV$ unless 
$B=B^{\ts\prime}=1$. Thus $(Z)_{\tts0}\in\Uh$ exists for any $Z\in\ZA\ts$.

Now suppose that $(Z)_{\tts0}\in\Uh$ is not unique. Then
Lemma \ref{lemma2} implies the existence of a nonzero element of $\Uh$ 
which belongs to the sum of the right ideal of $\Ab$ generated 
by $\n$ and $U\ts$, and of the left ideal
of $\Ab$ generated by $\n'$ and $U^*$. 
But this is not possible due to existence of
the basis in $\Ar$ described above. 
\qed
\end{proof}

Now regard $\ZA$ as a left $\Uh$-module.
Define the {\it Shapovalov form} $S:\ZA\times\ZA\to\Uh$ by %setting
$$
S(X,Y)=(\varepsilon(X)\mult Y)_{\tts0}
\quad\text{for}\quad
X,Y\in\ZA\ts.
$$
This form is symmetric, contravariant and $\Uh$-linear by definition.
That is, we have
\begin{align}
\label{SF1}
S(X,Y)&=S(Y,X)\ts;
\\[2pt]
\label{SF2}
S(Z\mult X,Y)&=S(X,\varepsilon(Z)\mult Y)
\quad\text{for}\quad
Z\in\ZA\,;
\\
\label{SF3}
S(Z\ts X,Y)&=Z\ts S(X,Y)
\quad\text{for}\quad
Z\in\!\Uh\,.
\end{align}
Note that
\begin{equation}
\label{Sh8}
S(X,Y)=0
\quad\text{when}\quad
Y\in\IV\,.
\hspace{26pt}
\end{equation}

%------------------------------------------------------------------------------

\subsection*{\it\normalsize 2.3.\ Shapovalov form on coinvariants space}

Let $\mu\in\h^*$ be any weight. 
Denote by $\Jm_\mu$ the left ideal of $\Ar$, generated
by the elements $X$ of $\np$, and by the elements 
$H-\mu(H)$ where $H\in\h\ts$. Here we identify any element
$X$ of $\g$ with the element \eqref{xdiag} of the algebra $\Ar\ts$.
Further, denote by $\Jp_\mu$ the left ideal of $\Ar$, generated
by the elements $X\ot 1$ where $X\in\np$, and the elements
$H\ot 1-\mu(H)$ where $H\in\h$.
Let $\I$ be the left ideal of $\Ar$ generated by the elements of
$U^*\subset\Heist\ts$. 
The first inclusion in the property (7)
of a Howe system implies that 
$X-X\ot1\in\I$
for 
$X\in\np\ts$.
The second inclusion in the property (7)
implies that there is a weight $\ka\in\h^*$ such that
\begin{equation}
\label{deka}
\zeta(H)-\ka(H)\in\Heist\cdot U^*
\quad\text{for all}\quad 
H\in\h\ts.
\end{equation}
Hence 
$H-H\ot1-\ka(H)\in\I$
for
$H\in\h\ts$. Therefore
\begin{equation}
\label{ijp}
\Jp_\mu+\I=\Jm_{\mu+\ka\ts}+\I\ts.
\end{equation}

Let $\Ar_{\ts\mu}$ be the quotient space of $\Ar$ 
relative to the left ideal $\Jp_\mu+\I\ts$. Let
$\Mr_\mu$ be the vector space of double cosets of
$\Ar$ relative to the left ideal $\Jp_\mu+\I$ and 
the right ideal $\J\ts$:
\begin{equation}
\label{nmummu}
\Ar_{\ts\mu}=\Ar\ts/\ts(\Jp_\mu+\I\ts)\ts,
\qquad 
\Mr_\mu=%\J\,\backslash
\Ar_{\ts\mu}\ts/\ts\J=
\Ar\ts/\ts(\Jp_\mu+\I+\J\ts)\,.
\end{equation}
Regard the ring $\P\ts(U)$ as a left $\Heist$-module, 
by identifying this ring with the quotient
of $\Heist$ over the left ideal generated by $U^*$. 
Thus for
$X\in\Heist$ and $u\com v\in\P\ts(U)$
\begin{equation}
\label{Sh8a} 
X(u)=v
\quad\text{if}\quad 
X\ts u-v\in\Heist\cdot U^*\ts.
\end{equation}
Then %by using the definition \eqref{H2} and the equality \eqref{ijp},
we can identify the $\Ar\ts$-module $\Ar_{\ts\mu}$ with the tensor
product $M_\mu\ot\P\ts(U)$ where $M_\mu$ is the Verma module
of the algebra $\U(\g)$ appearing as the first tensor factor in
\eqref{H2}. The vector space $\Mr_\mu$ gets identified with the space
of coinvariants
\begin{equation}
\label{coinv}
(M_\mu\ot\P\ts(U))_{\ts\n}\ts.
\end{equation}
The Cartan subalgebra $\h\subset\g$ acts on this space
via left multiplication in %the algebra 
$\Ar\ts$. Moreover, this
space is a left module over the Mickelsson algebra $\R\ts$.
It is generated by the image~of %the vector 
\begin{equation}
\label{1mu1}
1_\mu\ot 1\in M_\mu\ot\P\ts(U)
\end{equation}
where $1_\mu$ is the highest weight
vector of the Verma module $M_\mu\ts$. By using the equality \eqref{ijp},
any~element $H\in\h$ acts on the vector \eqref{1mu1} 
as multiplication by $(\mu+\ka)(H)\ts$.
Hence the image of the vector \eqref{1mu1} in $\Mr_\mu$ belongs to the
weight subspace $\Mr_\mu^{\tts\mu+\ka}\ts$. Due to the condition (6)
on a Howe system, the subspace $\Mr_\mu^{\tts\mu+\ka}\subset\Mr_\mu$
is one-dimensional, and is spanned by the image of \eqref{1mu1}.
Note that there is an isomorphism of vector spaces
$\P\ts(U)\to\Mr_\mu\ts$. It can be defined
by assigning to any $u\in\P\ts(U)$ the coset
of $1\ot u\in\Ar$ in the quotient $\Mr_\mu=\Ar/(\Jm_\mu+\I+\J\ts)\ts$.
This isomorphism will be denoted by $\io_\mu\ts$.

Now consider the left ideals of the algebra $\Ab\ts$, 
$$
\Jpb_{\mu}=\Uh\,\Jp_{\mu}
\quad\text{and}\quad
\Ib=\Uh\,\I\ts.
$$ 
Suppose the weight $\mu+\ka\in\h^\ast$ is \textit{generic},
that is $(\mu+\ka)(H_\al)\notin\ZZ$ for all $\al\in\De^+\ts$.
Then the spaces of double cosets 
$%\Jb\ts\backslash\ts
\Ab\ts/\ts(\Jpb_\mu+\Ib+\Jb\ts)$ and
$%\J\ts\backslash\ts
\Ar\ts/\ts(\Jp_\mu+\I+\J\ts)=\Mr_\mu$ 
are naturally isomorphic. Here we once again use the equality \eqref{ijp}.
Further, $\Mr_\mu$ is isomorphic to
the left $\ZA\ts$-module, defined as
the quotient of the algebra $\ZA$ by the left ideal generated by 
the cosets of the elements of $1\ot U^*$ and the cosets
of the elements $H-(\mu+\ka)(H)$ where $H\in \h\ts$.
%\begin{equation}
%\label{mmmu}
%\Mr_\mu\ts\cong\,\ZA\,/\,\ZA\,
%\{\ts v\com H\ot 1-\mu(H)\,|\,v\in U^*,H\in \h\ts\}\ts.
%\end{equation}
Here we use the isomorphism of the algebras $\Rb$ and $\ZA\ts$;
see Proposition~\ref{prop1.1} and remarks made immediately
after stating it. Accordingly, we will use the symbol
$\mult$ to denote the action of the algebra $\ZA$ on $\Mr_\mu\ts$.
Then any element of $\Mr_\mu$ can be presented as
$Z\mult\io_\mu(1)$ for some $Z\in\ZA\ts$.
Since $\io_\mu(1)\in\Mr_\mu^{\ts\mu+\ka}$, 
we get $Z\mult\io_\mu(1)\in\Mr_\mu^{\la}$
if $Z\in\ZA^\nu$ and
\begin{equation}
\label{lamuka}
\nu=\lambda-\mu-\ka\,.
\end{equation}
We may also choose $Z$ from the
image in $\ZA$ of the subalgebra $1\ot\P\ts(U)\subset\Ab\ts$.
Indeed, due to %the second relation in 
\eqref{M5} the cosets of $(1\ot u)\ts P$ and $1\ot u$ in $\Mr_\mu$ coincide
for any element $u\in\P\ts(U)\ts$. 

Due to \eqref{SF1}-\eqref{Sh8},
% and to the latter description of the left ideal of $\ZA\ts$,
for a generic weight $\mu+\ka$
the Shapovalov form $S:\ZA\ot\ZA\to\Uh$
defines a symmetric contravariant form 
$S_\mu:\Mr_\mu\ot\Mr_\mu\to\CC$ by setting
\begin{equation}
\label{S8}
S_\mu(X\mult\io_\mu(1)\com Y\mult\io_\mu(1))= S(X,Y)(\mu+\ka)
\quad\text{for}\quad 
X\com Y\in\ZA\ts.
\end{equation}
%Here $X,Y\in\ZA$ while $1$ is the identity in $\ZA$ 
%regarded as a generator of $\ZA$-module $\Mr_\mu\ts$.
Here $S(X,Y)(\mu+\ka)$ is the evaluation at $\mu+\ka$ 
of an element of $\Uh$, 
regarded as a rational function on $\h^*$.
The contravariance of the form $S_\mu$ means the equality
\begin{equation}
\label{S8a}
S_\mu(Z\mult f,g)=S_\mu(f,\varepsilon(Z)\mult g)
\quad\text{for}\quad
f,g\in\Mr_\mu
\quad\text{and}\quad Z\in\ZA\,.
\end{equation}
Note that if $f=\io_\mu(u)$ and $g=\io_\mu(v)$ for some
$u$ and $v$ from $\P\ts(U)\ts$, then
\begin{equation}
\label{better}
S_\mu(f,g)= S(X,Y)(\mu+\ka)
\end{equation}
where $X$ and $Y$ are respectively the images of $1\ot u$ and $1\ot v$
under the projection map $\Ar\to\SA\ts$. This follows from the remark
made at the end of the previous paragraph.

In Subsection 2.1 we selected
a non-degenerate symmetric bilinear form $\langle\ ,\,\rangle$ on
the vector space $U\ts$. Let us now extend this form
from $U$ to $\P\ts(U)$ in a natural way, as follows. 
Choose an orthonormal basis 
$u_1\lcd u_n$ of $U$. Thus $\langle u_i,u_j\rangle=\de_{ij}$
for $i\com j=1\lcd n\ts$.~Put
$$
\langle\ts u_1^{p_1}\ldots u_n^{p_n},u_1^{q_1}\ldots u_n^{q_n}\ts\rangle
\,=\,\prod_{k=1}^n\de_{\ts p_{k}\ts q_{k}}\,p_{k}\ts!
$$
where for every $k=1\lcd n$ we take
$p_k\com q_k\in\{0\com1\com2\com\ldots\}$ in the case $\th=1\ts$,
whereas $p_k\com q_k\in\{0\com1\}$ in the case $\th=-1\ts$.
This form on $\P\ts(U)$
is uniquely determined by setting $\langle\ts1\com1\ts\rangle=1\ts$
and by declaring
that the transpose to the operator of left multiplication by any $u_k$  
is given by the action of $u'_k\in U^*$. Here $u'_1\lcd u'_n$ 
is the basis in $U^*$ dual to the basis $u_1\lcd u_n$ in $U\ts$.
Thus for each $k=1\lcd d$ we have
$$
\langle\ts u_{k}\ts u\com v\ts\rangle=\langle\ts u\com u'_{k}\ts v\ts\rangle
\quad\text{for}\quad
u\com v\in\P\ts(U)\ts.
$$
Here we regard $\P\ts(U)$ as a left $\Heist$-module, see \eqref{Sh8a}.
Since $u'_k=\ep(u_k)\ts$, the above displayed 
equality implies a more general equality,
\begin{equation}
\label{Sh12}
\langle\ts X\ts u\com v\ts\rangle=\langle\ts u\com\ep(X)\ts v\ts\rangle
\quad\text{for}\quad
u\com v\in\P\ts(U) 
\quad\text{and}\quad
X\in\Heist\ts.
\end{equation}
%On the left of the equality displayed in the next lemma
%we will regard $u\com u'$ as elements of $\Mr_\mu$ by using the
%identification $\Mr_\mu=\P\ts(U)$ adopted above.

\begin{proposition*}
\label{propositionSh14}
Let $\mu+\ka\in\h^\ast$ be generic. 
Then for any\/ $u\com v\in\P\ts(U)$
$$
S_\mu(\ts f,g\ts)=
\langle\ts u\,,\zeta(\p[\mu+\rho])\ts v\ts\rangle
\quad\text{if}\quad
f=\io_\mu(u)
\quad\text{and}\quad
g=\io_\mu(v)\ts.
$$
\end{proposition*}

\begin{proof}
By definition, the extremal projector $P$ is a product in
$\widetilde{\!\U(\g)\!\!\!}\,\,\,$ over the set $\De^+$ 
of positive roots equipped with any normal ordering.
Any such an ordering on $\De^+$
has the following basic property \cite{AST}. 
Take any $\al\com\be\in\De^+$ such that $[E_\al\com F_\be]\neq0$.
Then $\al-\be$ is a root. The commutator $[E_\al\com F_\be]$
is proportional to $E_{\ts\al-\be}$ or to $F_{\ts\be-\al}$
respectively if $\al-\be\in\De^+$ or $\be-\al\in\De^+$.
Now suppose that
the root $\al$ precedes $\be$ in the ordering.
Then the basic property is that either the positive
root $\al-\be$ precedes $\al\ts$, or the positive root 
$\be-\al$ is preceded by $\be\ts$. 
Using repeatedly this property and the definition \eqref{pal}, 
we can present the extremal projector $\p=\p[\rho]$ as the sum in 
$\widetilde{\!\U(\g)\!\!\!}\,\,\,$
of the products of the form $A\,A'\,C$ where
$A$ and $A'$ are normally ordered monomials in the  
generators $F_\al$ and $E_\al$ respectively.
Here $\al$ ranges over the set $\De^+\ts$.
%of all positive roots equipped with any normal order. 
These monomials
form bases respectively in the vector spaces $\U(\np)$ and $\U(\n)\ts$.
Further, any monomial $A\,A'$ appears in the sum $P$ with
a unique factor $C\in\Uh$ on the right of~it.
This $C$ is a ratio of a complex number to a finite
product of factors %of the~form 
\begin{equation}
\label{everyfactor}
\ts H_\al+\rho\ts(H_\al)+s
\end{equation}
where $\al$ is a positive root and $s$ is a positive integer.
We will also regard the elements of $\Uh$ 
as rational functions on $\h^*\ts$. By %the definition 
\eqref{S8}
the value $S_\mu(\ts f,g\ts)$ of the Shapovalov form %on $\Mr_\mu$
equals the sum of the values
$(X)_{\tts0}\ts(\mu+\ka)$ where $X$ is the coset
in $\ZA$ of the product in~$\Ab\ts$,
\begin{equation}
\label{aac}
1\ot\ep(u)\cdot A\,A'\ts C\cdot1\ot v\ts.
\end{equation}

Without loss of generality we may assume that the element 
$v\in\P\ts(U)$ has a weight, say $\la\in\h^*$,
relative to the adjoint action of $\h\ts$:

\begin{equation}
\label{nuka}
[\,\zeta(H)\com v\,]=\la(H)\,v
\quad\text{for all}\quad
H\in\h\ts.
\end{equation}
Then 
\begin{equation}
\label{cde}
C\cdot1\ot v=1\ot v\cdot D
\end{equation}
where $D\in\Uh$ is obtained from $C$ by replacing every 
factor \eqref{everyfactor} respectively by 
$$
H_\al+(\la+\rho)(H_\al)+s\,.
$$

Consider the product
\begin{equation}
\label{compro}
1\ot\ep(u)\cdot A\,A'\cdot1\ot v
\end{equation}
in the algebra $\Ab$.
Here $A=F_\al\ldots F_\be$ and 
$A=E_{\ts\al'}\ldots E_{\ts\be'}$ for some positive roots
$\al\lcd\be$ and $\al'\lcd\be'$.
By the definition of the ideals $\Jb$ and $\Jmb$ of the algebra $\Ab$,
the coset in $\ZA$ of the product \eqref{compro}
coincides with that of the product of commutators
\begin{align*}
[F_\al\lcd [F_\be\com1\ot\ep(u)]\ldots\ts]
&\cdot
[E_{\ts\al'}\lcd[E_{\ts\be'}\com1\ot v]\ldots\ts]=
\\[4pt]
1\ot[\ts\zeta(F_\al)\lcd [\ts\zeta(F_\be)\com\ep(u)]\ldots\ts]
&\cdot
1\ot[\ts\zeta(E_{\ts\al'})\lcd[\ts\zeta(E_{\ts\be'})\com v]\ldots\ts]\,.
\end{align*}
Denote by $Y$ the coset of the last product. Due to the properties
(4) and (7) of a Howe system, then we have
$(Y)_{\tts0}=(Z)_{\tts0}$ where
$Z$ denotes the coset in $\ZA$ of the product in $\Ab\ts$,
$$
1\ot(\ts\ep(u)
\cdot
\zeta(F_\al)\ldots\zeta(F_\be)
\cdot
\zeta(E_{\ts\al'})\ldots\zeta(E_{\ts\be'})
\cdot
v\ts)=
1\ot(\ts\ep(u)\,\zeta(A\,A'\ts)\,v\ts)\,.
$$

Note that here $(Z)_{\tts0}\in\CC\subset\Uh\ts$. 
Further, the adjoint action of the subalgebra $\h\subset\g$ normalizes the
ideals $\Jb$ and $\Jmb$ of $\Ab\ts$.
It also normalizes the ideals $\JV\com\IV$ of $\ZA\ts$.
Hence the above argument together with \eqref{cde}
implies that $(X)_{\tts0}=(ZD)_{\tts0}=(Z)_{\tts0}\ts D\ts$. Therefore
\begin{equation}
\label{lare}
(X)_{\tts0}\ts(\mu+\ka)=
(Z)_{\tts0}\cdot D(\mu+\ka)=
(Z)_{\tts0}
\cdot
C(\mu+\ka+\la)\,.
\end{equation}
But by the definitions \eqref{deka} and \eqref{nuka} of the 
weights $\ka$ and $\la$ respectively,
for every element $H\in\h$ we have the relations in the algebra $\Heist\ts$,
$$
v\cdot(\ka+\la)(H)\in
v\cdot(\zeta(H)+\la(H))+\Heist\cdot U^*=
\zeta(H)\cdot v+\Heist\cdot U^*\,.
$$
They imply the equality in the left $\Heist\ts$-module $\P\ts(U)\ts$,
$$
\ep(u)\,\zeta(A\,A'\ts)\,v
\cdot 
C(\mu+\ka+\la)
=
\ep(u)\,\zeta(A\,A'B)\,v
$$
where $B\in\Uh$ is obtained from $C$ by replacing every 
factor \eqref{everyfactor} respectively by 
$$
H_\al+(\mu+\rho)(H_\al)+s\,.
$$
Since $\langle\ts 1\com U\cdot\P\ts(U)\ts\rangle=0\ts$,
the right hand side of \eqref{lare} can be now written as
$$
\langle\ts 1\ts,\ep(u)\,\zeta(A\,A'B)\,v\ts\rangle=
\langle\ts u\com\zeta(A\,A'B)\,v\ts\rangle\,.
$$
By using the definition of the element $B\in\Uh$,
we now get Proposition \ref{propositionSh14}.
\qed
\end{proof}

Consider again the symmetric
bilinear form $\langle\ ,\,\rangle$ on the vector space $\P\ts(U)\ts$.
Take the action of the Lie algebra $\g$ on the vector space $\Heist$ 
corresponding to that of the group $\G\ts$.
By the property (2) of a Howe system, 
this action of $\g$ is adjoint to the homomorphism
$\zeta:\U(\g)\to\Heist\ts$. 
By the properties (5) and (8), the action %on $\Heist$
of the Cartan subalgebra $\h\subset\g_+\subset\g$ 
preserves the subspaces $U$ and $U^*$ of $\Heist$.
By the definition \eqref{epsdef} for any
$u\in U$ and $u'\in U^*$ %we have
$$
\langle\ts u\com\ep(u')\ts\rangle=u'\ts(u)=B(u\com u')\ts,
$$
while the action of the group $\G$ on the vector space 
$W=U\oplus U^*$ preserves the bilinear form $B$ by our assumption.
Therefore for any $H\in\h$ and for any $u\com u'$ as above we have
$$
\langle\,[\ts\zeta(H),u\ts]\com\ep(u')\ts\rangle+
\langle\ts u\com\ep(\ts[\ts\zeta(H),u'\ts]\ts)\ts\rangle=0\ts.
$$
Here $\ep$ is an anti-automorphism of the algebra $\Heist$
obeying the property (3) of a Howe system.
Moreover, we have $\ep(H)=H$ for the Chevalley 
anti-involution $\ep$ on $\g\ts$. Therefore
$$
\langle\,[\ts\zeta(H),u\ts]\com\ep(u')\ts\rangle=
\langle\ts u\com[\ts\zeta(H),\ep(u')\ts]\ts)\ts\rangle\ts.
$$
Since any vector of $U$ can be written as $\ep(u')$
for some $u'\in U^*$, the last equality implies
that the action of the Cartan subalgebra
$\h$ on $\P\ts(U)$ is self-conjugate with respect to 
the bilinear form $\langle\ ,\,\rangle\ts$.
Therefore, because the form $\langle\ ,\,\rangle$
on $\P\ts(U)$ is non-degenerate, 
its restriction 
to any weight subspace of $\P\ts(U)$ is also non-degenerate.
%The latter observation will be used later on.

%------------------------------------------------------------------------------

\vspace{-6pt}%%%%%%%%%%%%%%%%%%%%%%%%%%%%%%%%%%%%%%%%%%%%%%%%%%%%%%%%%%%%%%%%%%

\subsection*{\it\normalsize 2.4.\ Shapovalov form on a weight subspace}

Now let $\mu\in\h^*$ be an arbitrary weight.
Take the subspace $\Mr_\mu^{\la}$ 
of $\Mr_\mu=\Ar/(\J+\Jp_\mu+\I)$ 
consisting of the elements of weight $\la\ts$. 
The Cartan subalgebra $\h\subset\g$ 
acts on %the quotient vector space 
$\Mr_\mu$
via the left multiplication in $\Ar$ by the elements of 
$\U(\g)\subset\Ar\ts$.
In this subsection, we will assume that
$\la+\rho$ is nonsingular. Then for any $f\in\Mr_\mu^\la$ and
$Z\in\SA$ the element $Z\mult f\in\Mr_\mu$ is well defined.
Indeed,
the denominators of $P=P[\rho]$ do not vanish on $\Mr_\mu^{\la}$
if $\la+\rho$ is nonsingular.
In partucular, $\Mr_\mu^{\la}$ is a
module over the subalgebra $\SAzero\subset\ZA\ts$. 
Let the weights $\la$ and $\mu$
vary under the constraint that their difference $\la-\mu$ is fixed.
Then for any given elements $X\com Y\in\ZA$ we can regard
$S(X\com Y)(\mu+\ka)$ as a rational function of $\mu\in\h^*$. Now
recall the definition \eqref{lamuka} of the weight $\nu\ts$.

\begin{proposition*}
\label{proposition3.6}
Suppose that the weight $\lambda+\rho\in\h^*$ is nonsingular.
Then for any given elements $X,Y\in\SA^{\tts\nu}$
the function $S(X\com Y)(\mu+\ka)$ of\/ $\mu$ has only finite values.
\end{proposition*}

\begin{proof}
Let $A$ and $B$ be any representatives in $\Ar$
of the cosets $X$ and $Y$ respectively. Then
the $\ep(X)\mult Y$ is the coset of %the product
\begin{equation}
\label{ezpz}
\ep(A)\,\p\,B
\end{equation}
where $\p=\p[\rho]\in\widetilde{\!\U(\g)\!\!\!}\,\,\,$ is the extremal
projector for $\g\ts$. We assume that the elements
$A$ and $B$ of $\Ar$ also have the weight $\nu$ 
relative to the adjoint action of $\h\ts$. 

Any poles of $(\ep(X)\mult Y)_{\tts0}$
as of a rational function on $\h^*$ 
may arise only from 
the denominators of the summands of $\p[\rho]\ts$. 
By the definition \eqref{pal}, up to non-zero scalar multipliers,
these denominators
are products of linear factors of the form
$H_\al+\rho\ts(H_\al)+s$ with $s=1,2,\ldots$ and $\al\in\De^+\ts$. 
Let us move these denominators to the right in \eqref{ezpz}
through the factor $B\ts$, and evaluate the resulting 
denominators at $\mu+\ka\ts$. Then the linear factors become %respectively
$$
(\ts\mu+\ka+\nu+\rho\ts)(H_\al)+s=(\lambda+\rho)(H_\al)+s
$$
which do not depend on $\mu\ts$, and
are nonzero numbers when $\lambda+\rho$ is nonsingular.
\qed
\end{proof}

When the weight 
$\mu+\ka$ is generic,
the Shapovalov form $S_\mu$ on $\Mr_\mu$
can be defined by the equation 
\eqref{better}
where $f=\io_\mu(u)$ and $g=\io_\mu(v)$ for some
$u$ and $v$ from $\P\ts(U)\ts$, while
$X$ and $Y$ are the images of $1\ot u$ and $1\ot v$
under the projection $\Ar\to\SA\ts$. The same
equation \eqref{better} can now be used to
define an $\SAzero\ts$-contravariant form on the subspace 
$\Mr_\mu^{\la}\subset\Mr_\mu$ for nonsingular $\la+\rho$ and 
any $\mu\ts$. Indeed, if here
$f,g\in\Mr_\mu^{\la}$ then $X,Y\in\SA^{\tts\nu}$
so that Proposition~\ref{proposition3.6} applies.
Denote by $S_\mu^\la$ the bilinear form
on the subspace $\Mr_\mu^{\la}\subset\Mr_\mu$ defined~by~\eqref{better},
$$
S_\mu^\la:\,
\Mr_\mu^\la
\times
\Mr_\mu^\la
\to\CC\,.
$$
The contravariance \eqref{SF2} of the Shapovalov form on $\ZA$ 
implies that $S_\mu^\la$ is a contravariant
form relative to the action of the subalgebra
$\SAzero\subset\ZA\ts$. %on $\Mr_\mu^{\la}\ts$.
Indeed, if $\mu+\ka$ is generic then by \eqref{S8a}
for any $f,g\in\Mr_\mu^\la$ and $Z\in\SAzero$ we have
\begin{equation}
\label{shalamu}	
S_\mu^\la(\ts Z\mult f\com g\ts)=
S_\mu^\la(\ts f\com\ep(Z)\mult g\ts)\,.
\end{equation}
When the weight $\la-\mu$ is fixed,
both sides of this equality are rational functions of $\mu\ts$.
Hence these two rational functions are the same. 
So the equality \eqref{shalamu} holds for any
$\mu\ts$, provided $\la+\rho$ is nonsingular,
which is assumed in this subsection.

In the next lemma, 
we regard $\P\ts(U)\subset\Heist$ as subalgebras
of $\Ar$ and use the adjoint action of $\h$ on~$\Ar\ts$;
see the condition (5) on a Howe system. 
We also identify the elements of $\Heist$ with their
images in the double coset algebra $\ZA\ts$.
The weights of these elements in $\ZA$ are taken
relative to the adjoint action of $\h\ts$. By the condition (2)
on a Howe system, these
weights are the same as relative to the adjoint action of $\h$ on
$\Heist\ts$. Recall \eqref{lamuka}.

\begin{lemma*}
\label{lemma3.7}
Suppose that for some
$f\in\Mr_\mu^{\la}$ and\/ $Y\in\P\ts(U^\ast)^{\ts-\nu}$
we have the equality $Y\mult f=\io_\mu(1)$ in\/ $\Mr_\mu\ts$.
Then for any $X\in\P\ts(U)^{\ts\nu}$
we also have $(X\ts Y)\mult f=\io_\mu(X)$ in\/ $\Mr_\mu\ts$.
\end{lemma*}

\begin{proof}
We have $f=\io_\mu(u)$ for a certain element $u\in\P\ts(U)\ts$,
which has weight 
$\nu$ under the adjoint action of $\h$ on $\P\ts(U)\ts$.
In Subsection 2.1 we equipped 
the algebra $\Heist$ with $\ZZ\ts$-grading so that
the elements of $U$ and $U^*$ have degrees $1$ and $-1$ respectively.
Since the weights of the elements $X,\ep(Y)$ and $u$ of $\Heist$
are the same, by the condition (6) on a Howe system
these elements also have the same (non-negative) degree. Thus
$$
\deg X=-\deg Y=\deg u\ts.
$$
%Further, for any positive root 
%$\al\in\De^+$
%the condition (5) on a Howe system implies that 
%$[F_\al\com U]\subset U$
%while the condition (8) implies that 
%$[E_\al\com U^*]\subset U^*$.

Consider the element $(X\ts Y)\mult f$ of $\Mr_\mu\ts$.
It corresponds to the product $X\ts Y\ts\p\,u$ in $\Ab\ts$.
Write the extremal projector $P$ as a sum, like we did
in the beginning of the proof of Proposition~\ref{propositionSh14}. 
Then move the generators
$F_\al$ and $E_\al$ of $\g$ respectively to the left of $X\ts Y$ and 
to the right of $u$ in the resulting summands of $X\ts Y\p\,u\ts$.
This procedure shows that modulo $\Jb+\Jmb$,
the product $X\ts Y\ts\p\,u\ts$ equals 
%$X\ts Y u$ plus summands of
the sum of products of the form $X'\ts Y' u'$ where
$$
X'\in\P\ts(U)^{\ts\nu-\al}\,,
\quad
Y'\in\P\ts(U^\ast)^{-\nu-\be}
\quad\text{and}\quad
u'\in\P\ts(U)^{\ts\nu+\al+\be}
$$
while $\al,\be$ are certain sums of positive roots. 

Suppose that the coset of $X'\ts Y' u'$
makes a non-zero contribution to
$(X\ts Y)\mult f\in\Mr_\mu\ts$. 
The condition (5) on a Howe system implies that
the adjoint operator $\ad_{\ts\zeta(F_\al)}$ on $\Heist$
either has zero degree, or vanishes on the subspace $\P\ts(U)\ts$. Hence 
$$
\deg X'=\deg X\ts.
$$
Thus the element $Y'\ts u'\in\Heist$ has the degree zero.
By our assumption, this element does not belong to the left ideal
$\Heist\,U^*\subset\Heist\ts$. It also has a weight
relative to the adjoint action of $\h$ on $\Heist\ts$.
But then the weight must be zero. Indeed, 
by the condition (8) on a Howe system,
the adjoint action of $\h$ on $\Heist$ preserves the left ideal
$\Heist\,U^*$. Modulo this ideal, the element $Y' u'$ of $\Heist$
equals a non-zero scalar, which has the weight zero.
Thus the weight of $Y' u'$ is zero as well.
Therefore $\al=0$ and $X'=X\ts$. 

The element $Y\mult f$ of $\Mr_\mu$
corresponds to the product $Y\ts\p\,u$ in $\Ab\ts$. 
Modulo $\Jb+\Jmb$, the latter product equals 
%$Y u$ plus summands
the sum of the products
$Y' u'\in\Heist$ where $Y'$ and $u'$ are 
the same as above in the particular case when $\al=0\ts$:
$$
Y'\in\P\ts(U^\ast)^{-\nu-\be}
\quad\text{and}\quad
u'\in\P\ts(U)^{\ts\nu+\be}
$$
where $\be$ ranges over sums of positive roots. 
Modulo the left ideal $\Heist\,U^*$ of $\Heist$,
the sum of the products $Y' u'$
equals a scalar. This scalar must be $1\ts$, 
because $Y\mult f=\io_\mu(1)$ in $\Mr_\mu\ts$.
Hence the sum of the products 
$X\ts Y'\ts u'$ equals $X$ modulo $\Heist\,U^*\subset\Heist\ts$. 
\qed
\end{proof}

The contravariance of the form $S_\mu^\la$ 
implies that the kernel $\Ker\,S_\mu^{\la}$ 
is an $\SAzero\ts$-submodule of $\Mr_\mu^{\la}\ts$.
Indeed, let $f\in\Ker\,S_\mu^\la$ and $Z\in\SAzero$. 
Then $Z\mult f\in\Ker\,S_\mu^\la\,$, because
for any $g\in\Mr_\mu^{\la}$ 
$$
S_\mu^\la(\ts Z\mult f\com g\ts)=
S_\mu^\la(\ts f\com\ep(Z)\mult g\ts)=0\,.
$$

\begin{proposition*}
\label{proposition2}
Let $\lambda+\rho\in\h^*$ be nonsingular. 
Then the quotient space
$\Mr_\mu^\lambda\,/\ts\Ker\,S^\lambda_\mu$ 
is an irreducible $\SAzero$-module.
\end{proposition*}

\begin{proof}
Take any $f\in\Mr_\mu^\lambda$ such that 
$f\notin\Ker\,S^\lambda_\mu\ts$, so that 
the image of $f$ in the quotient space
$\Mr_\mu^\lambda\,/\ts\Ker\,S^\lambda_\mu$ is not zero.
For any $Y\in\P\ts(U^\ast)^{\ts-\nu}$ consider the vector
$Y\mult f\in\Mr_\mu\ts$. This vector
is well defined, because
the denominators of $P=P[\rho]$ do not vanish on $f\ts$.
The weight of this vector is $\la-\nu=\mu+\ka\ts$,
while the subspace $\Mr_\mu^{\tts\mu+\ka}\subset\Mr_\mu$
is spanned by the vector $\io_\mu(1)\ts$.
Hence the vector $Y\mult f$ is a multiple of $\io_\mu(1)\ts$.

Suppose that $Y\mult f=0$ for every $Y\in\P\ts(U^\ast)^{\ts-\nu}\ts$. 
The map $\ep$ provides a bijection
$\P\ts(U^\ast)^{\ts-\nu}\to\P\ts(U)^{\ts\nu}$
of subspaces of $\Heist\ts$. Hence for
any $g\in\Mr_\mu^\lambda$ there is an element
$Y\in\P\ts(U^\ast)^{\ts-\nu}$ such that $g=\io_\mu(\ep(Y))\ts$.
If the weight $\mu+\ka$ is generic, then
$$
S_\mu^\la\ts(f,g)=
S_\mu\ts(f,\ep(Y)\mult\io_\mu(1))=
S_\mu\ts(Y\mult f,\io_\mu(1))=0\ts.
$$
By our definition of the form $S_\mu^\la$ then we must have 
$S_\mu^\la\ts(f,g)=0$ for all weights $\mu\ts$,
not only those where $\mu+\ka$ is generic.
Hence $f\in\Ker\,S^\lambda_\mu\ts$, 
a contradiction.
It shows that there is %an element 
$Y\in\P\ts(U^\ast)^{\ts-\nu}$
such that $Y\mult f\neq0\ts$. Moreover, we can choose $Y$ with 
$Y\mult f=\io_\mu(1)\ts$.

Now take any $h\in\Mr_\mu^\lambda$ and choose 
$X\in\P\ts(U)^{\ts\nu}\ts$
such that $\io_\mu(X)=h\ts$.
By Lemma \ref{lemma3.7}, then we have
$$
(X\ts Y)\mult f=h\ts.
$$
Here $X\ts Y$ stands for the coset in $\ZA$ of
the product of $X\in\P\ts(U)^{\ts\nu}$ and
$Y\in\P\ts(U^*)^{\ts-\nu}$ in the algebra $\Ar\ts$. 
This coset belongs to $\SAzero$ by definition,
because the product of $X$ and $Y$ in $\Ar$ has weight zero.
\qed
\end{proof}

Recall that $\Dwo$ is a subalgebra of $\SAzero\ts$.
By using Proposition \ref{proposition3.9}, we now get
the following 

\begin{corollary*}
\label{corollary3.11}
Suppose that $\la+\rho\in\h^*$ is nonsingular, and that\/ 
$\ES_\la=\Sym_\la\ts$. Then the quotient space
$\Mr_\mu^\lambda\,/\ts\Ker\,S^\lambda_\mu$ 
is an irreducible $\Dwo$-module.
\end{corollary*}

%------------------------------------------------------------------------------

\vspace{-8pt}%%%%%%%%%%%%%%%%%%%%%%%%%%%%%%%%%%%%%%%%%%%%%%%%%%%%%%%%%%%%%%%%%%

\subsection*{\it\normalsize 2.5.\ Irreducible\/ {\rm Q}-modules}

In this subsection we give another construction
of irreducible $\Dwo\ts$-modules, which employs Proposition \ref{proplus}.
For any weight $\mu\in\h^*$ denote by $L_\mu$ the quotient
of the Verma module $M_\mu$ by its maximal proper submodule $N_\mu\ts$.
This quotient is an irreducible $\g\ts$-module.
Note that $\P\ts(U)$ is irreducible as a module over the algebra $\Heist\ts$.
Hence the tensor product $L_\mu\ot\P\ts(U)$ is an irreducible module
over the algebra $\Ar=\U(\g)\ot\Heist\ts$.
Recall that the action of the algebra $\Heist$ on the vector space
$\P\ts(U)$ is defined by \eqref{Sh8a}. 

We can also regard $L_\mu\ot\P\ts(U)$ as a $\g\ts$-module,
by restriction from $\Ar$ to its subalgebra $\U(\g)\ts$. 
Here $\g$ acts on the tensor factor $\P\ts(U)$ via the homomorphism
$\zeta:\U(\g)\to\Heist\ts$. 
%Recall that the action of the algebra $\Heist$ on %the space 
%$\P\ts(U)$ defined by \eqref{Sh8a}. 
Let us show that then the actions on $\P\ts(U)$ of
the subalgebras $\h$ and $\np$ of $\g$ are
semisimple and locally nilpotent respectively.
Indeed, the action of the group $\G$ on $\Heist$ 
is locally finite, see Subsection 2.1.
Hence the action of $\g$ on $\Heist$ adjoint to %the homomorphism 
$\zeta$ is also locally finite,
by the property (2) of a Howe system.
In particular, the actions of $\h$ and $\np$
on $\Heist$ adjoint to $\zeta$ are 
respectively semisimple and locally nilpotent.
So are the actions of $\h$ and $\np$
on $\P\ts(U)\ts$, by \eqref{Sh8a}
and the property (7) of a Howe system.
It follows that the actions of $\h$ and $\np$
on $L_\mu\ot\P\ts(U)$ are also semisimple and locally nilpotent
respectively.

Take any $\la\in\h^*$ such that $\la+\rho$ is nonsingular.
Consider the space of $\n\ts$-coinvariants of
the $\g\ts$-module $L_\mu\ot\P\ts(U)\ts$, and 
denote by $\Lr_\mu^\la$ its weight subspace 
\begin{equation}
\label{2.4444}
(L_\mu\ot\P\ts(U))_{\tts\n}^\la\subset(L_\mu\ot\P\ts(U))_{\tts\n}\,.
\end{equation}
Then $\Lr_\mu^\la$
is an irreducible $\SAzero\ts$-module by Proposition \ref{proplus}.
By using Proposition \ref{proposition3.9}~we~get

\begin{corollary*}
\label{corollary3.11plus}
Suppose that $\la+\rho\in\h^*$ is nonsingular, and that\/ 
$\ES_\la=\Sym_\la\ts$. Then\/ $\Lr_\mu^\la$
is an irreducible $\Dwo$-module.
\end{corollary*}

Let us now regard $\Lr_\mu^\la$ 
as the quotient of the vector subspace $\Mr_\mu^\lambda$ by
its subspace 
\begin{equation}
\label{nrlm}
\Nr_\mu^\la=(N_\mu\ot\P\ts(U))_{\tts\n}^\la\ts.
\end{equation}
We will prove that under certain conditions on
$\la$ and $\mu\ts$, the subspace 
$\Nr_\mu^\la$ of $\Mr_\mu^\lambda$
coincides with the kernel $\Ker\,S^\lambda_\mu$
of the Shapovalov form on  $\Mr_\mu^\lambda\ts$.
Then the quotient $\Lr_\mu^\la$ of $\Mr_\mu^\lambda$
coincides with the quotient
$\Mr_\mu^\lambda\,/\ts\Ker\,S^\lambda_\mu\ts$.
In particular, then these two quotients are
the same as $\Dwo\ts$-modules.
Observe that for any vector $u\in\P\ts(U)$ of the weight $\la-\mu\ts$, 
the denominators of $\zeta(\p[\mu+\rho])\ts u$ get evaluated
as products of factors of the form
$$
(\ts\la-\mu\ts)(H_\al)+(\ts\mu+\rho\ts)(H_\al)+s=(\lambda+\rho)(H_\al)+s
$$
where $\al\in\De^+$ and $s$ is a positive integer.
These factors do not depend on $\mu\ts$. 
And they do not vanish, when the weight $\la+\rho$ is nonsingular.

Note that here the weight $\la-\mu$ of the vector 
$u\in\P\ts(U)$ is taken relative to the action
of $\g$ on $\P\ts(U)$ via the homomorphism
$\zeta:\U(\g)\to\Heist\ts$, by using the definition \eqref{Sh8a}.
If we used the adjoint action of $\g$ on $\P\ts(U)$
as on a subalgebra of $\Heist\subset\Ar\ts$, like we did in the previous
subsection, then the weight of the same vector $u$ 
would be equal to~\eqref{lamuka}.

\begin{proposition*}
\label{lmeq}
Suppose that $\la+\rho\in\h^*$ is nonsingular, and that
$\zeta(\p[\mu+\rho])\ts u\neq0$ for some vector 
$u\in\P\ts(U)$ of weight $\la-\mu\ts$. Then

\begin{equation}
\label{nkeq}
\Nr_\mu^\la=\Ker\,S^\lambda_\mu\ts.
\end{equation}
\end{proposition*}

\begin{proof}
Consider the vector
$
1_\mu\ot\ts u\in M_\mu\ot\P\ts(U)\ts.
$
This vector is of weight $\la\ts$.
We will show that the image $\io_\mu(u)$ of this vector
in the quotient space
$$
\Mr_\mu^\lambda=(M_\mu\ot\P\ts(U))^\la_{\tts\n}
$$
does not belong to the sum $\Nr_\mu^\la+\Ker\,S^\lambda_\mu\ts$.
This implies the equality \eqref{nkeq}.
Indeed, then
\begin{equation}
\label{nkeq1}
\Mr_\mu^\lambda\ts/\ts(\ts\Nr_\mu^\la+\Ker\,S^\lambda_\mu\ts)
\end{equation}
is a non-zero quotient of
\begin{equation}
\label{nkeq2}
\Mr_\mu^\lambda\ts/\ts\Nr_\mu^\la=\Lr_\mu^\la\ts,
\end{equation}
which is an irreducible $\SAzero\ts$-module by Proposition \ref{proplus}.
So the quotients \eqref{nkeq1} and \eqref{nkeq2}
of $\Mr_\mu^\lambda$ must coincide, and we obtain
the inclusion $\Ker\,S^\lambda_\mu\ts\subset\ts\Nr_\mu^\la\ts$.
The opposite inclusion is obtained from the 
irreducibility of $\Mr_\mu^\lambda\,/\ts\Ker\,S^\lambda_\mu$
as an $\SAzero\ts$-module, %in a similar way, 
see Proposition \ref{proposition2}.

Let us now assume that 
$\io_\mu(u)\in\Nr_\mu^\la+\Ker\,S^\lambda_\mu\ts$.
We shall bring this to a contradiction.
By our assumption, the element $\io_\mu(u)\in\Mr_\mu^\la$
equals the sum of an element of \eqref{nrlm} and of an element
$\io_\mu(v)\in\Ker\,S^\lambda_\mu$ where $v\in\P\ts(U)\ts$. 
Let us regard $\Mr_\mu^\la$ as a $\g\ts$-module, and
apply the extremal projector $\p=\p[\rho]$ to that equality.
Then we get a certain equality in the vector space 
$M_\mu\ot\P\ts(U)\ts$. Here we use the second property in \eqref{M4}.

By applying $\p$ to the element $\io_\mu(u)\in\Mr_\mu^\la$
we get the vector $\p\ts(1_\mu\ot\ts u)\in M_\mu\ot\P\ts(U)\ts$,
which belongs to
$$
1_\mu\ot\zeta(\p[\mu+\rho])\ts u+(\ts\n\,M_\mu)\ot\P\ts(U)\ts.
$$
Here $\zeta(\p[\mu+\rho])\ts u\neq0\ts$.
By applying $\p$ to any element of \eqref{nrlm} 
we get a vector from the subspace
$$
\p\ts(N_\mu\ot\P\ts(U))\subset N_\mu\ot\P\ts(U)
\subset(\ts\n\,M_\mu)\ot\P\ts(U)\ts.
$$
By applying $\p$ to the element $\io_\mu(v)\in\Mr_\mu^\la$
we get the vector $\p\ts(1_\mu\ot\ts v)\in M_\mu\ot\P\ts(U)\ts$,
which belongs to
$$
1_\mu\ot \zeta(\p[\mu+\rho])\ts v
+(\ts\n\,M_\mu)\ot\P\ts(U)\ts.
$$
But since $\io_\mu(v)\in\Ker\,S^\lambda_\mu\ts$, here we have
$\zeta(\p[\mu+\rho])\ts v=0$ by Proposition \ref{propositionSh14}.
We use the non-degeneracy of the restriction
of the form $\langle\ ,\,\rangle$
to any weight subspace of $\P\ts(U)\ts$,
proved in the very end of Subsection 2.3.
Thus we get a contradiction.
\qed
\end{proof}

Later we will produce vectors $u$ 
satisfying the conditions of Proposition \ref{lmeq},
for~each example of a Howe system and for certain weights $\mu\ts$;
see Propositions \ref{proposition3.7}, \ref{norm1} and \ref{norm2}.

%==============================================================================

\section*{\bf\normalsize 3. Intertwining operators}
\setcounter{section}{3}
\setcounter{equation}{0}
\setcounter{theorem*}{0}

%------------------------------------------------------------------------------

\subsection*{\it\normalsize 3.1.\ 
Zhelobenko operators and left ideals}

The vector space $W=U\op U^*$ has been equipped with
a non-degenerate bilinear form $B$, symmetric or alternating.
For each isotropic subspace $V\subset W$ 
denote by $\I_{\tts V}$ the left ideal of the algebra $\Ar$ generated by the
elements of $1\ot V\ts$. Also consider the left ideal
$$
\Ib_{\tts V}=\Uh\,\I_{\tts V}
$$ 
of algebra $\Ab\ts$. It is generated by the elements of the subspace
$1\ot V\subset\Ab\ts$. In particular,
$$
\I_{\ts U^*}=\I
\quad\text{and}\quad
\Ib_{\ts U^*}=\Ib\,.
$$ 
%If the subspace $V$ was not isotropic, then
%then the corresponding left ideal of $\Ar$ would be equal to $\Ar\ts$.
In Subsection 2.1 we assumed that
the group $\G$ acts on the vector space $W$
and preserves the bilinear form $B$. 
Thus $\sih\ts(V)$ is an isotropic subspace of $W$ for any $\si\in\Sym\ts$. 

\begin{lemma*}
\label{lemma2.12}
For any element $\si\in\Sym$ the operator $\xic_{\ts\si} $ 
on $%\Jb\backslash
\Ab\ts/\ts\Jb$ maps
$$
%\Jb\,\backslash\,
(\,\Ib%_{\ts U^*}
+\Jb\,)\ts/\ts\Jb
\,\to\, 
%\Jb\,\backslash\,
(\,\Ib_{\ts\sih\ts(U^*)}+\Jb\,)\ts/\ts\Jb\,.
$$
\end{lemma*}

\begin{proof}
We use the induction on the length $\ell(\si)$ of the element 
$\si\in\Sym\ts$. For $\ell(\si)=0$ the statement to prove is tautological.
Take any index $c\in\{1\lcd r\}$ with
$\ell(\si_c\,\si)=\ell(\si)+1\ts$. Then 
$\xic_{\ts\si_c\ts\si}=\xic_c\,\xic_{\ts\si}\ts$.
Note that here we have $\si^{-1}(\al_c)\in\De^+$. 
Take any element $X\in\Ib%_{\tts U^*}
\ts$. By the induction assumption,
the operator $\xic_{\ts\si}$ maps the coset of $X$ in 
$%\Jb\backslash
\Ab\ts/\ts\Jb$ to that of %the sum
$$
\sum_{k=1}^n\,Y_k\cdot1\ot\sih\ts(u'_k)
$$
where $u'_1\lcd u'_n$ are basis elements of $U^*$ and
$Y_1\lcd Y_n$ are certain elements of $\Ab\ts$.
By definition, the operator $\xic_c$ maps the latter coset to
$$
\sum_{k=1}^n\,\xib_c\ts(\ts\sih_c(Y_k)\cdot1\ot\sih_c\ts\sih\ts(u'_k))\ts.
$$
Due to the definition \eqref{M10}, 
for making the induction step it now suffices to show that
$$
[F_c\com1\ot\sih_c\,\sih\ts(U^*)]\,\subset\,1\ot\sih_c\,\sih\ts(U^*)\ts.
$$
Applying to this relation the automorphism $(\sih_c\,\sih)^{-1}$ 
of the algebra $\Ab\ts$, we~get~the~relation
$$
[E_\al\com 1\ot U^*]\subset 1\ot U^*
\quad\text{for}\quad
\al=\si^{-1}(\al_c)\in\De^+\ts,
$$
which holds by the condition (8) on a Howe system.
Here $(\sih_c\,\sih)^{-1}(F_c)$ is a multiple of $E_\al$,
because the action of $\Norm_{\,0}\TT$ on the subspace 
$\h\subset\g$ factors through that of $\Sym$.
\qed
\end{proof}

In Subsection 2.3 for arbitrary $\mu\in\h^*$
we introduced the left ideals $\Jp_\mu$ and $\Jpb_\mu$
of the algebras $\Ar$ and $\Ab$ respectively.
Generalizing \eqref{nmummu}, 
introduce the quotient vector spaces
\begin{equation}
\label{nmummuv}
\Ar_{\ts\mu,V}=\Ar\ts/\ts(\Jp_\mu+\I_{\tts V}\ts)
\quad\text{and}\quad 
\Mr_{\tts\mu,V}=%\J\,\backslash
\Ar_{\ts\mu,V}\ts/\ts\J
=\Ar\ts/\ts(\Jp_\mu+\I_{\tts V}\ts+\J)
\end{equation}
of $\Ar\ts$. Then $\Ar_{\ts\mu,V}$ is a left $\Ar\ts$-module
with the space of $\n\ts$-coinvariants $\Mr_{\tts\mu,V}$.
Note that for any element $\tau\in\DA$ the action
of its representative $\ttau\in\G$ on the algebra $\Ar$
determines a linear map
\begin{equation}
\label{taumap}
\Mr_{\tts\mu,V}\,\to\,\Mr_{\tts\tau(\mu),\ttau(V)}\ts.
\end{equation}
Indeed, the adjoint action of $\ttau$ on $\g$ preserves
the subalgebras $\n$ and $\np$.
The action of $\ttau$ 
on the algebra $\Ar$ preserves the right ideal $\J$, and 
maps the left ideals $\Jp_\mu$ and $\I_{\tts V}$ 
to $\Jp_{\tts\tau(\mu)}$ and $\I_{\,\tts\ttau(V)}$ respectively. 
The linear map \eqref{taumap} determined in this way
is clearly invertible.

If the weight $\mu+\ka$ is generic, 
then the spaces of double cosets 
$$
%\Jb\ts\backslash\ts
\Ab\ts/\ts(\Jpb_\mu+\Ib_{\tts V}+\Jb\ts)
\quad\text{and}\quad
%\J\ts\backslash\ts
\Ar\ts/\ts(\Jp_\mu+\I_{\tts V}+\J\ts)=\Mr_{\tts\mu,V}
$$
are naturally isomorphic. Here we use the equality \eqref{ijp}.
Furthermore, then $\Mr_{\tts\mu,V}$ is isomorphic to
the left $\ZA\ts$-module, defined as
the quotient of the algebra $\ZA$ by the left ideal generated by 
the cosets of the elements of $1\ot V$ and the cosets
of the elements $H-(\mu+\ka)(H)$ where $H\in \h\ts$.
Now denote respectively by $\Jp$ and $\Jpb$ the left ideals of 
$\Ar$ and $\Ab\ts$, generated by the elements $X\ot1$ with $X\in\nplus\ts$.

\begin{lemma*}
\label{lemma2.13}
For any element $\si\in\Sym$ the operator $\xic_{\ts\si}$ 
on $%\Jb\backslash
\Ab\ts/\ts\Jb$ maps
$$
%\Jb\,\backslash\,
(\,\Jpb+\Ib%_{\ts U^*}
+\Jb\,)\ts/\ts\Jb
\,\to\, 
%\Jb\,\backslash\,
(\,\Jpb+\Ib_{\ts\sih\ts(U^*)}+\Jb\,)\ts/\ts\Jb\,.
$$
\end{lemma*}

\begin{proof}
For each $c=1\lcd r$ denote by $\n'_c$ the vector subspace of 
$\n'$ spanned by all root vectors $E_\al$ except $E_c\ts$. 
Denote by $\Jp_c$ the left ideal of $\Ar$ generated by all elements 
$X\ot1$ with $X\in\n'_c\ts$, and by the element $E_c\ts$. 
Let $\Jpb_c=\Uh\,\J_c$ be the corresponding left ideal of $\Ab$.
We will first prove that the Zhelobenko operator
$\xic_{\ts c}$ on the vector space $%\Jb\backslash
\Ab\ts/\ts\Jb$ maps
\begin{equation}
\label{Int1}
%\Jb\,\backslash\,
(\,\Jpb_c+\Jb\ts\,)\ts/\ts\Jb
\,\to\, 
%\Jb\,\backslash\,
(\,\Jpb\ns+\Jb\,)\ts/\ts\Jb\,.
\end{equation}

The left ideal $\sih_c\,(\ts\Jpb_{\ts c}\ts)$ is generated by the
element $F_c$, and by the subspace of $\Ab$ formed by
all elements $X\ot1$ where $X\in\np_{\ts c}\ts$.
Observe that the latter subspace is preserved by the
adjoint action of the element $F_c\ts$.
Hence for any element $Y\in\Ab$ and any element $Z$ from that subspace,
$\bar{\xi}_{\ts c}\ts(\ts YZ)\in\ts\Jpb\ns+\ts\Jb$ by the definition
of the operator $\bar{\xi}_{\ts c}\ts$. The property \eqref{Int1} 
of the operator $\xic_{\ts c}$ on $%\Jb\backslash
\Ab\ts/\ts\Jb$ now follows from
the inclusion
$\bar{\xi}_{\ts c}\ts(\ts YF_c)\in\ts\Jb$
for every %\quad\text{for all}\quad 
$Y\in\Ab\ts$.
For a detailed proof of that inclusion see 
\cite[Section 7.2]{KO} or \cite[Section 3]{KN1}.

We will now prove Lemma \ref{lemma2.13}
by induction on the length $\ell(\si)$ of the element 
$\si\in\Sym\ts$. For $\ell(\si)=0$ the statement to prove is tautological.
Take any index $c\in\{1\lcd r\}$ with
$\ell(\si_c\,\si)=\ell(\si)+1\ts$. Then 
$\xic_{\ts\si_c\ts\si}=\xic_c\,\xic_{\ts\si}\ts$.
Note that here we have $\si^{-1}(\al_c)\in\De^+$. 
Take the statement of the lemma as the induction assumption.
Observe an equality of left ideals
\begin{equation}
\label{Int2}
\Jpb+\Ib_{\ts\sih\ts(U^*)}\,=\,\Jpb_c+\Ib_{\ts\sih\ts(U^*)}
\end{equation}
of the algebra $\Ab\ts$.
Indeed, by the definitions of left ideals $\Jpb$ and $\Jpb_c$ 
the equality \eqref{Int2} follows from the inclusion
$$
1\ot\zeta(E_c)\in\ts\Ib_{\ts\sih\ts(U^*)}\ts.
$$
Using the automorphism $\sih^{-1}$ of the algebra $\Ab\ts$, 
the latter inclusion is equivalent to 
$$
1\ot\zeta(E_\al)\in\Ib%_{\ts U^*}
\quad\text{for}\quad
\al=\si^{-1}(\al_c)\in\De^+,
$$
which holds by the condition (7) on a Howe system. To make
the induction step it now suffices to prove that the operator
$\xic_{\ts c}$ on $\Jb\backslash \Ab$ maps
$$
%\Jb\,\backslash\,
(\,\Jpb_c+\Ib_{\ts\sih\ts(U^*)}+\Jb\,)\ts/\ts\Jb
\,\to\, 
%\Jb\,\backslash\,
(\,\Jpb+\Ib_{\ts\sih_c\sih\ts(U^*)}+\Jb\,)\ts/\ts\Jb\,.
$$
But that follows from the property \eqref{Int1} of the operator 
$\xic_{\ts c}\ts$, and from the induction step we made when
proving Lemma \ref{lemma2.12}.
\qed
\end{proof}

\begin{proposition*}
\label{proposition2.13} 
Let $\mu+\ka\in\h^*$ be generic. Then %for any $\si\in\Sym$ 
the operator $\xic_{\ts\si}$ on $%\Jb\backslash
\Ab\ts/\ts\Jb$ maps
$$
%\Jb\,\backslash\,
(\,\Jpb_\mu+\Ib%_{\ts U^*}
+\Jb\,)\ts/\ts\Jb
\,\to\,
%\Jb\,\backslash\,
(\,\Jpb_{\si\circ\mu}+\Ib_{\ts\sih\ts(U^*)}+\Jb\,)\ts/\ts\Jb\ts.
$$
\end{proposition*}

\begin{proof} 
The left ideal $\Jpb_\mu$ of the algebra $\Ab$ 
is spanned by $\Jpb$ and by all the subspaces
\begin{equation}
\label{xhot}
\Ab\,(H\ot1-\mu(H))%=\Ab\,(H-1\ot\zeta(H)-\mu(H))
\end{equation} 
where $H\in\h\ts$.
Due to Lemma \ref{lemma2.13}, to prove Proposition \ref{proposition2.13}
it suffices to consider the action of $\xic_{\ts\si}$ on the
images  of the subspaces \eqref{xhot}
in the quotient $\Jb\backslash\Ab\ts$.
Rewrite~\eqref{xhot}~as
\begin{equation}
\label{xhotkak}
\Ab\,(H-(\mu+\ka)(H)-1\ot\zeta(H)+\ka(H))\ts.
\end{equation}
By \eqref{M21} the operator $\xic_{\ts\si}$ maps the image of
$\Ab\,H$ in $\Jb\backslash\Ab$ to that of $\Ab\,(\si\circ H)\ts$.
By \eqref{deka}
\begin{equation}
\label{zekaiu}
1\ot\zeta(H)-\ka(H)\in\Ib%_{\ts U^*}
\ts.
\end{equation}
Hence by using Lemma \ref{lemma2.12},
the image of \eqref{xhotkak} in $\Jb\backslash\Ab$
is mapped by $\xic_{\ts\si}$ 
to a subspace~of the image in $\Jb\backslash\Ab$ of 
$$
\Ab\,(\si\circ H-(\mu+\ka)(H))+\Ib_{\ts\sih\ts(U^*)}\ts.
$$
We will complete the proof of Proposition \ref{proposition2.13}
by showing that for every $H\in\h^*$,
\begin{equation}
\label{simuka}
\si\circ H-(\mu+\ka)(H)\in\Jpb_{\si\circ\mu}+\Ib_{\ts\sih\ts(U^*)}\ts.
\end{equation}

Replacing $H$ by $\si^{-1}\circ H$ 
at the left hand side of \eqref{simuka}, we get the element
$$
H-(\si\circ(\mu+\ka))(H)=H\ot1-(\si\circ\mu)(H)+1\ot\zeta(H)-\si(\ka)(H)\ts.
$$
By definition, here
$$
H\ot1-(\si\circ\mu)(H)\in\Jpb_{\si\circ\mu}\ts.
$$
But by applying to the relation \eqref{zekaiu}
the automorphism $\sih$ of $\Ab\ts$, 
using the condition (1) on a Howe system, and then replacing 
$H$ by $\si^{-1}(H)$ in the resulting relation, we get
$$
1\ot\zeta(H)-\si(\ka)(H)\in\I_{\ts\sih\ts(U^*)}\ts.
\eqno\qed
$$
\end{proof}

\begin{corollary*} 
\label{corolary2.13}
Let $\mu+\ka\in\h^*$ be generic. 
Then for any $\si\in\Sym$ the operator $\xic_{\ts\si}$
defines a $\ZA\ts$-equivariant linear map
\begin{equation}
\label{xisimu}
\Mr_{\mu}%,U^*}
\,\to\,
\Mr_{\ts\si\circ\mu\ts,\ts\sih\ts(U^*)}\ts.
\end{equation}
\end{corollary*}

Let us denote this linear map by the same symbol $\xic_{\tts\si}\ts$.
Then its $\ZA\ts$-equivariance means
$$
\xic_{\tts\si}(Z\mult f)=\xic_{\ts\si}(Z)\mult\xic_{\ts\si}(f)
\quad\text{for}\quad
Z\in\ZA
\quad\text{and}\quad
f\in\Mr_\mu%,U^*}
\ts.
$$
This property follows from Proposition \ref{proposition1.3}\ts(iii).
Note that by \eqref{M12}, the above defined
operator \eqref{xisimu} maps
the weight subspace $\Mr_{\tts\mu}^\la\subset\Mr_\mu$ to the
weight subspace
$$
\Mr_{\ts\si\circ\mu\ts,\ts\sih\ts(U^*)}^{\ts\si\circ\la}
\,\subset\,\Mr_{\ts\si\circ\mu\ts,\ts\sih\ts(U^*)}\ts.
$$

%------------------------------------------------------------------------------

\subsection*{\it\normalsize 3.2.\ Action of Zhelobenko operators}

For any element $\si\in\Sym$ and any weight $\mu\in\h^*$
consider the $\Ar\ts$-module 
\begin{equation}
\label{nsim}
\Ar_{\,\si\circ\mu\ts,\ts\sih\ts(U^*)}=
\Ar\,/\ts(\,\Jp_{\ts\si\circ\mu}+\I_{\ts\sih\ts(U^*)}\ts)\ts.
\end{equation}
Its vector space can be identified with that
of the tensor product of the Verma module
$M_{\ts\si\circ\mu}$ over $\g\ts$, multiplied by the $\Heist$-module 
$\Heist/\Heist\,\sih\ts(U^*)\ts$. The latter quotient space 
can be idenitified with
$\P\ts(\sih\ts(U))\ts$. The vector space of \eqref{nsim}
is then identified with $M_{\ts\si\circ\mu}\ot\P\ts(\sih\ts(U))\ts$. Further,
the space $\Mr_{\ts\si\circ\mu,\sih\ts(U^*)}$
of $\n\ts$-coinvariants of \eqref{nsim}
can be identified with 
$\P\ts(\sih\ts(U))\ts$, by assigning to any $w\in\P\ts(\sih\ts(U))$
the coset of $1\ot w\in\Ar$ in 
\begin{equation}
\label{msim}
\Mr_{\ts\si\circ\mu\ts,\ts\sih\ts(U^*)}=
\Ar\,/\ts(\,\Jp_{\ts\si\circ\mu}+\I_{\ts\sih\ts(U^*)}+\J\ts)\ts.
\end{equation}
In particular, when $\si\in\Sym$ is the identity element,
the space $\Mr_{\tts\mu,U^*}=\Mr_\mu$ gets identified 
with $\P\ts(U)\ts$, as in Subsection 2.3.

Suppose that the weight $\mu+\ka$ is generic. 
Using Corollary \ref{corolary2.13} together with
the above identifications of the source and target
vector spaces in \eqref{xisimu}, the %Zhelobenko 
operator $\xic_{\ts\si}$ on $%\Jb\backslash
\Ab\ts/\ts\Jb$
determines a linear map $\P\ts(U)\to\P\ts(\sih\ts(U))\ts$.
%To avoid confusion, 
The latter map will be denoted by $I_{\tts\si,\mu}\ts$.

\begin{proposition*} 
\label{proposition3.7}
Let $\mu+\ka\in\h^*$ be generic. Then for any $\si\in\Sym$
and any $u\in\P\ts(U)$
\begin{equation}
\label{Sh20}
I_{\tts\si,\mu}(u)=\sih\ts(\ts\zeta(\p_\si[\mu+\rho])\ts u)
\end{equation}
where $\zeta(\p_\si[\mu+\rho])\ts u$ is regarded as an element
of\/ $\P\ts(U)$ by using the definition \eqref{Sh8a}.
\end{proposition*}

\begin{proof}
We prove Proposition \ref{proposition3.7} 
by induction on the length $\ell(\si)$ of $\si\in\Sym$. When
$\si$ is the identity element of $\Sym$, the statement to prove
is trivial. Let us now use the statement as the induction assumption.
Take any index $c\in\{1\lcd r\}$ with
$\ell(\si_c\,\si)=\ell(\si)+1\ts$. Then 
$\si^{-1}(\al_c)$ is a positive root, let us denote it by $\al\ts$.
For short, denote by $v$ the~element 
$\zeta(\p_\si[\mu+\rho])\ts u\in\P\ts(U)\ts$.
%We will indicate by the symbol $\equiv$ the equalities in $\Ab$ modulo
We have
\begin{equation}
\label{absubspace}
\Mr_{\ts\si_c\si\circ\mu\ts,\ts\sih_c\sih\ts(U^*)}=
\Ab\,/\ts(\,
\Jpb_{\si_c\si\circ\mu}+\Ib_{\ts\si_c\si(U^*)}+\Jb\ts)\ts.
\end{equation}
Using the definition \eqref{M10} along with the 
induction assumption, the element $I_{\tts\si_c\si,\mu}(u)$
of $\P\ts(\sih\ts(U))$ can be identified with the coset in 
\eqref{absubspace}
of the sum
$$
\sum_{s=0}^\infty\,\, 
\prod_{t=1}^s\,\,
\bigl(\ts t\,(H_c-t+1)\bigr)^{-1}
\cdot 
E_c^{\ts s}\,
\ad^{\ts s}_{\ts F_c}(1\ot\sih_c\ts\sih\ts(v))\ts.
$$ 
Without changing the sum,
we can replace the operator $\ad_{\ts F_c}^{\ts s}$
by $1\ot\ad_{\ts\zeta(F_c)}^{\ts s}$ here,
because the elements $F_c\ot1$ and 
$1\ot\sih_c\ts\sih\ts(v)$ of $\Ab$ commute.
Without changing the coset,
we can then replace %every factor 
$E_c^{\ts s}$ by $1\ot\zeta(E_c^{\ts s})$ in resulting sum,
because $E_c\ot1$ commutes with 
$$
1\ot\ad_{\ts\zeta(F_c)}^{\ts s}(\ts\sih_c\ts\sih\ts(v))
$$
and belongs to the left ideal $\Jp\subset\Jp_{\si_c\si\circ\mu}\ts$.
So we get the sum 
\begin{align*}
&\sum_{s=0}^\infty\,\, 
\prod_{t=1}^s\,\,
\bigl(\ts t\,(H_c-t+1)\bigr)^{-1}
\cdot 1\ot\zeta(E_c^{\ts s})\,
\ad^{\ts s}_{\ts\zeta(F_c)}(\sih_c\ts\sih\ts(v))\,=
\\
&\sum_{s=0}^\infty\,\, 
\prod_{t=1}^s\,\,
\bigl(\ts t\,(H_c-t+1)\bigr)^{-1}
\cdot 1\ot\sih_c\ts\sih\ts(\ts
\zeta(F_\al^{\ts s})\,\ad^{\ts s}_{\ts\zeta(E_\al)}(v))\ts.
\end{align*}
By the property (7) of a Howe system, the element
$1\ot\zeta(E_\al)$ belongs to the left ideal $\I=\I_{\ts U^*}$ of $\Ar\ts$.
Hence
$$
1\ot\sih_c\ts\sih\ts(\zeta(E_\al))\in\I_{\ts\si_c\si(U^*)}\ts.
$$
Therefore in the last displayed sum, the element
$\ad^{\ts s}_{\ts\zeta(E_\al)}(v)\in\Heist$ can be replaced 
by the element $\zeta(E_\al^{\ts s})\,v\in\Heist\ts$, without changing
the coset of the sum in \eqref{absubspace}. We~get
\begin{equation}
\label{thatsum}
\sum_{s=0}^\infty\,\, 
\prod_{t=1}^s\,\,
\bigl(\ts t\,(H_c-t+1)\bigr)^{-1}
\cdot 1\ot\sih_c\ts\sih\ts(\ts\zeta(F_\al^{\ts s})\,
\zeta(E_\al^{\ts s})\,v\ts)\ts.
\end{equation}
Here $H_c$ stands for $H_c\ot1+1\ot\zeta(H_c)\ts$.
Modulo the left ideal $\Jp_{\si_c\si\circ\mu}$ of $\Ar\ts$,
the element $H_c\ot1$ equals
$$
(\si_c\si\circ\mu)(H_c)=
(\si_c\si\ts(\mu+\rho)-\rho)(H_c)=
$$
$$
(\mu+\rho)((\si_c\si)^{-1}(H_c))-\rho\ts(H_c)=
-\ts(\mu+\rho)\ts(H_\al)-1\ts.
$$
Hence the coset of the sum \eqref{thatsum} in \eqref{absubspace}
coincides with that of
\begin{align*}
&\sum_{s=0}^\infty\,\,
\prod_{t=1}^s\,\,
1\ot\bigl(\ts t\,(\ts\zeta(H_c)-(\mu+\rho)\ts(H_\al)-t\ts)\bigr)^{-1} 
\cdot1\ot\sih_c\ts\sih\ts(\ts\zeta(F_\al^{\ts s})\,
\zeta(E_\al^{\ts s})\,v\ts)\,=
\\
&\sum_{s=0}^\infty\,\,1\ot
\sih_c\ts\sih\ts\Bigl(\,\,
\prod_{t=1}^s\,\,
\bigl(\ts t\,(-\ts\zeta(H_\al)-(\mu+\rho)\ts(H_\al)-t\ts)\bigr)^{-1} 
\cdot\zeta(F_\al^{\ts s})\,
\zeta(E_\al^{\ts s})\,v\Bigr)\,=
\end{align*}

\vspace{0pt}
$$
1\ot\sih_c\ts\sih\ts(\ts\zeta(P_\al[\mu+\rho])\,v\ts)\,=\,
1\ot\sih_c\ts\sih\ts(\ts\zeta(P_{\ts\si_c\si}[\mu+\rho])\,u\ts)\ts.
\eqno\qed
$$
\end{proof}

\vspace{8pt}
The operator \eqref{xisimu} has been defined only when the
weight $\mu+\ka$ is generic. However, 
Proposition \ref{proposition3.7} yields the following result,
which is valid for any $\mu\in\h^*$.
Here we also use the remark made at the very end of
Subsection 3.1.

\begin{corollary*}
\label{corollary3.7}
If $\lambda+\rho$ is nonsingular, the operator $\xic_{\ts\si}$ on
$%\Jb\backslash
\Ab\ts/\ts\Jb$ defines a linear map
\begin{equation}
\label{zhelarb}
\Mr_\mu^\la\,\to\,
\Mr_{\ts\si\circ\mu\ts,\ts\sih\ts(U^*)}^{\ts\si\circ\la}\ts.
\end{equation}
\end{corollary*}

\begin{proof}
The source and target spaces in \eqref{xisimu}
can be identified with $\P\ts(U)$ and $\P\ts(\sih\ts(U))$
respectively. The first of these identifications
uses the bijection $\io_\mu$ as in Subsection 2.3.
Let $u\in\P\ts(U)$ be any element of weight $\la-\mu$ relative
to the action of $\h\ts$. Here any element $H\in\h$ acts on 
$\P\ts(U)$ via the left multiplication by $\zeta(H)\ts$, using 
the definition \eqref{Sh8a}. Then $\io_\mu(u)\in\Mr_\mu^\la\ts$.
When applying the operator $\zeta(\p_\si[\mu+\rho])$ to $u\ts$,
the denominators become, up to non-zero scalar multipliers,
certain products of the factors of the form 
$$
(\la-\mu)(H_\al)+(\mu+\rho)(H_\al)+s=
(\lambda+\rho)(H_\al)+s
$$
where $\al$ is a positive root and $s$ is a positive integer.
These factors do not depend on~$\mu\ts$, and do not vanish
if the weight $\lambda+\rho$ is nonsingular.
Thus by mapping $u$ to the element of $\ts\P\ts(\sih\ts(U))$
at the right hand side of \eqref{Sh20},
we get the required linear map \eqref{zhelarb}.
\qed
\end{proof}

%------------------------------------------------------------------------------

\subsection*{\it\normalsize 3.3.\ Irreducibility theorem}

Now take the longest element $\si_0$ of the Weyl group $\Sym\ts$. 
For any weight $\mu\in\h^*$ consider the $\R\ts$-module 
$\Mr_{\ts\si_0\circ\mu\tts,\tts\sih_0(U^*)}\ts$.
The next Proposition \ref{corollary3.8}
establishes a connection between 
the Shapovalov form $S_\mu^\la$~on~$\Mr_\mu^\la\ts$, 
and the linear map
\begin{equation}
\label{zhelong}
\Mr_\mu^\la\,\to\,
\Mr_{\ts\si_0\circ\mu\ts,\ts\sih_0(U^*)}^{\ts\si_0\circ\la}
\end{equation}
defined by the Zhelobenko operator $\xic_{\tts0}=\xic_{\ts\si_0}$ on
$%\Jb\backslash
\Ab\ts/\ts\Jb\ts$, when $\la+\rho\in\h^*$ is nonsingular.
The latter map will be denoted by $\xic_{\tts0}\ts|\ts\Mr_\mu^\la\,$.
Note that when $\mu+\ka\in\h^*$ fails to be generic,
the operator $\xic_{\tts0}$ on $%\Jb\backslash
\Ab\ts/\ts\Jb$
does not necessarily define any map from the whole
space~$\Mr_{\mu}\ts$.
By combining Proposition \ref {corollary3.8}
with Corollary \ref{corollary3.11}, 
we will obtain Theorem~\ref{theorem2}. 

\begin{proposition*}
\label{corollary3.8}
Let $\lambda+\rho\in\h^*$ be nonsingular. Then for any $\mu\in\h^*$
\begin{equation*}
\label{shapker}
\Ker\,S_\mu^\lambda 
=
\Ker\ts(\ts\xic_{\tts0}\ts|\ts\Mr_\mu^\la\ts)\ts.
\end{equation*}
\end{proposition*}

\begin{proof}
According to Subsection 3.2, the vector space of the
the $\R\ts$-module $\Mr_{\ts\si_0\circ\mu\tts,\tts\sih_0(U^*)}$
can be identified with $\P\ts(\sih_0(U))\ts$.
Denote by $\langle\ ,\,\rangle_{\si_0}$ the non-degenerate
$\CC\ts$-bilinear pairing
$$
\P\ts(U)\times\P\ts(\sih_0(U))\to\CC
$$
defined by the equality
\begin{equation}
\label{sinop}
\langle\ts u\com w\ts\rangle_{\si_0}=
\langle\ts u\com\sih_0^{\ts-1}(w)\ts\rangle\ts.
\end{equation}
Here $\langle\ ,\,\rangle$ is the non-degenerate symmetric bilinear form 
on the vector space $U$ selected in Subsection 2.1. Hence the pairing 
$\langle\ ,\,\rangle_{\si_0}$ is non-degenerate too.

Now for any $u,v\in\P\ts(U)$ put $f=\io_\mu(u)$ and $g=\io_\mu(v)\ts$. 
If the weight $\mu+\ka$ is generic, then by
Propositions \ref{propositionSh14} and \ref{proposition3.7} we have
\begin{equation}
\label{maineq}
S_\mu^\la\ts(f,g)=
\langle\ts u\,,\zeta(\p[\mu+\rho])\ts v\ts\rangle
=
\langle\ts u\,,\sih_0^{\ts-1}(I_{\tts\si_0,\mu}(v))\ts\rangle
=
\langle\ts u\,,I_{\tts\si_0,\mu}(v)\ts\rangle_{\si_0}\ts.
\end{equation}
When the difference $\la-\mu$ 
and the elements $u\com v\in\P\ts(U)$ are fixed,
the left hand side of the above equalities
becomes a rational functions of $\mu\in\h^*\ts$. 
If $\la+\rho$ is nonsingular,
this rational function has finite values for all $\mu$
by Proposition \ref{proposition3.6}.
Now Proposition~\ref{corollary3.8}
follows from the definition
of the operator $\xic_{\tts0}\ts|\ts\Mr_\mu^\la\ts$,
see our proof of Corollary~\ref{corollary3.7}.
\qed
\end{proof}

\begin{theorem*}
\label{theorem2}
Suppose that $\la+\rho\in\h^*$ is nonsingular, and that\/ 
$\ES_\la=\Sym_\la\ts$.
Then for any $\mu\in\h^*$ the quotient\/
$\Mr_\mu^\la\,/\,\Ker(\ts\xic_{\tts0}\ts|\ts\Mr_\mu^\la\ts)$
is an irreducible\/ $\Dwo\ts$-module.
\end{theorem*}

Recall that the algebra $\Ar$ is specified
as the tensor product \eqref{H2}, where
the group $\G$ acts diagonally. Consider the corresponding 
subalgebra $\Ar^\G\subset\Ar$ of $\G\ts$-invariants. 
This subalgebra
acts on $\Ar_{\ts\mu}$ by restricting the action of $\Ar\ts$.
Since $\Ar^\G\subset\Norm\,\J\ts$, the subalgebra
$\Ar^\G$ then
acts on the space $\Mr_\mu$ of $\n\ts$-coinvariants of $\Ar_{\ts\mu}\ts$.
The latter action preserves the subspace
$\Mr_\mu^\la\subset\Mr_\mu\ts$.
This subspace is also a $\Dwo\ts$-module, if
the weight $\la+\rho$ is nonsingular.
The above action of $\Ar^\G$ on $\Mr_\mu^\la$
can also be obtained by pulling the action of $\Dwo$ 
back through the isomorphism $\ga:\Ar^\G\to\Dwo\ts$, see Subsection 1.5.
Note that for any $\si\in\Sym$
the actions of $\Ar^\G$ on the source
and target vector spaces of the map \eqref{zhelarb}
defined by $\xic_{\ts\si}\ts$,
are intertwined by this map by its definition.
Using this observation when $\si=\si_0$,
we obtain a corollary to Theorem \ref{theorem2}. 

\begin{corollary*}
\label{ircor}
Suppose that $\la+\rho\in\h^*$ is nonsingular, and that\/ 
$\ES_\la=\Sym_\la\ts$.
Then for any $\mu\in\h^*$ the quotient\/
$\Mr_\mu^\la\,/\,\Ker(\ts\xic_{\tts0}\ts|\ts\Mr_\mu^\la\ts)$
is an irreducible\/ $\Ar^\G\ts$-module.
\end{corollary*}

%------------------------------------------------------------------------------

\subsection*{\it\normalsize 3.4.\ Contravariant pairing}

For any $\mu\in\h^*$ we can define a non-degenerate 
$\CC\ts$-bilinear pairing
$$
Q_\mu:
\Mr_\mu\times
\Mr_{\ts\si_0\circ\mu\ts,\ts\sih_0(U^*)}
\to\CC
$$
as follows. For any two elements $u\in\P\ts(U)$ and 
$w\in\P\ts(\sih_0(U))$ consider the cosets
in $\Mr_\mu$ and $\Mr_{\ts\si_0\circ\mu\ts,\ts\sih_0(U^*)}$
of $1\ot u$ and $1\ot w$ respectively; see \eqref{msim}.
By definition, the value of $Q_\mu$ on this pair of cosets is 
$\langle\ts u\com w\ts\rangle_{\si_0}$; see \eqref{sinop}.
By restricting $Q_\mu$ to weight subspaces,
we define a pairing 
\begin{equation}
\label{qlm}
Q_\mu^{\ts\la}:
\Mr_\mu^\la\times
\Mr_{\ts\si_0\circ\mu\ts,\ts\sih_0(U^*)}^{\ts\si_0\circ\la}
\to\CC
\end{equation}
for any $\la\in\h^*\ts$. The latter pairing is also non-degenerate,
see the end of Subsection 2.3. 

Now suppose that the weight $\la+\rho$ is nonsingular. 
Let us prove that then 
\begin{equation}
\label{qlme}
Q_\mu^{\ts\la}\ts(\ts \ep(Z)\mult f\com h\ts)=
Q_\mu^{\ts\la}\ts(\ts f\com\xic_{\tts0}(Z)\mult h\ts)
\quad\text{when}\quad 
Z\in\SAzero\ts.
\end{equation}
In particular, here
$\xic_{\tts0}(Z)\mult h\in
\Mr_{\ts\si_0\circ\mu\ts,\ts\sih_0(U^*)}^{\ts\si_0\circ\la}$
is defined, even though the element $\xic_{\tts0}(Z)$
of the algebra $\ZA$ may be 
not contained in the subspace $\SA\subset\ZA\ts$.

First let $\mu+\ka$ be generic. 
Then the map \eqref{zhelong}
defined by the Zhelobenko operator $\xic_{\tts0}$ on
$\Ab\ts/\ts\Jb\ts$, is invertible. Let $g\in\Mr_\mu^\la$ 
be the image of $h$ under the inverse map. By~\eqref{maineq}
$$
Q_\mu^{\ts\la}\ts(\ts\ep(Z)\mult f\com h\ts)=
S_\mu^{\ts\la}\ts(\ts\ep(Z)\mult f\com g\ts)=
S_\mu^{\ts\la}\ts(\ts f\com Z\mult g\ts)=
Q_\mu^{\ts\la}\ts(\ts f\com \xic_{\tts0}(Z)\mult h\ts)\ts.
$$
Here we used \eqref{shalamu} and Corollary \ref{corolary2.13}.
Thus we get \eqref{qlme} for generic $\mu+\ka\ts$. 
Now note that when $\la+\rho$ is nonsingular, the left hand
side of \eqref{qlme} is defined for any weight $\mu\ts$,
as a finite number. When $f,h,Z$ and the difference $\la-\mu$ are all
fixed while $\mu$ varies, this number becomes a rational function of $\mu\ts$.
Hence the right hand side of \eqref{qlme} can also be defined 
as a rational function of the weight $\mu\ts$, with only finite values.
Since the pairing \eqref{qlm} is non-degenerate,
we can then determine the vector $\xic_{\tts0}(Z)\mult h\ts$.

If $Z\in\Dwo\ts$, then $\xic_{\tts0}(Z)=Z$ by definition.
Hence the equality \eqref{qlme} implies that
for any $\mu$ and nonsingular $\la+\rho$ the pairing \eqref{qlm}
is $\Dwo\ts$-contravariant\tts:
$$ 
Q_\mu^{\ts\la}\ts(\ts \ep(Z)\mult f\com h\ts)=
Q_\mu^{\ts\la}\ts(\ts f\com Z\mult h\ts)
\quad\text{when}\quad 
Z\in\Dwo\ts.
$$
In particular, we have an action $\mult$
of the algebra $\Dwo$ on the vector space
$\Mr_{\ts\si_0\circ\mu\ts,\ts\sih_0(U^*)}^{\ts\si_0\circ\la}\ts$.

%==============================================================================

\section*{\bf\normalsize 4. Yangians and reductive dual pairs}

\setcounter{section}{4}
\setcounter{equation}{0}
\setcounter{theorem*}{0}

%------------------------------------------------------------------------------

\subsection*{\it\normalsize 4.1.\ Reductive dual pairs}

From now on we will work with examples of Howe systems,
as defined in Subsection~2.1. For each of the corresponding algebras
\eqref{H2} we will describe explicitly its subalgebra of $\G\tts$-invariants.
Then we will apply Corollary \ref{ircor} to each of these examples. 
The group $\G$ will be one of the classical complex
Lie groups $\GL_{\tts m}\ts$, $\SP_{\tts2m}$ and $\SO_{\tts2m}$ 
with any positive integer $m\ts$.
Respectively, $\g$ will be one of the 
Lie algebras $\gl_m\ts$, $\sp_{2m}$ and $\so_{2m}\ts$.

First consider $\G=\GL_m\ts$. 
Let the indices $a$ and $b$ run through %the sequence
$1\lcd m\ts$. Then $e_a$ will denote a vector of the standard basis 
in the vector space $\CC^m\ts$, while $E_{ab}\in\gl_m$ will be
a standard matrix unit. Choose the standard triangular
decomposition \eqref{hc0} of $\g=\gl_m$ where 
the subalgebras $\n\com\h\com\np$ are spanned
by the matrix units $E_{ab}$ with $a>b$, $a=b$, $a<b$ respectively.
The elements $E_{aa}$ form a basis in $\h\ts$,
and we will denote by $\eta_{\ts a}$ the vector
of the dual basis in $\h^*$ corresponding to $E_{aa}\ts$. 
The positive and negative roots are %the weights 
$\eta_{\ts a}-\eta_{\ts b}$ with $a<b$ and $a>b$ respectively.
The semisimple rank $r$ of $\gl_m$ is $m-1\ts$. For each $c=1\lcd m-1$
we will choose $\eta_{\ts c}-\eta_{\ts c+1}$ 
as the simple root $\al_c\ts$, and
$$
H_c=E_{cc}-E_{c+1,c+1}\ts,
\quad
E_c=E_{c,c+1}
\quad\text{and}\quad
F_c=E_{c+1,c}
$$
as the basis elements of the 
$\mathfrak{sl}_2\ts$-subalgebra $\g_c\subset\g\ts$.
Define the Chevalley anti-involution $\ep$ on $\gl_m$ by
setting $\ep(E_{ab})=E_{ba}\ts$. Choose the trivial
Cartan decomposition $\g=\g_{\tts+}\op\g_{\tts-}$
of $\gl_m$ so that $\g_{\tts+}=\gl_m$ and $\g_{\tts-}=\{0\}\ts$.

The algebraic group $\GL_m$ is connected, and the
maximal torus $\TT\subset\GL_m$ with the Lie algebra
$\h$ consists of all elements acting on each basis vector $e_a\in\CC^m$
by a  scalar multiplication. 
The action of 
$\si_c\in\Sym$ on $\h^*$ exchanges the dual basis vectors
$\eta_{\ts c}$ and $\eta_{\ts c+1}\ts$, 
leaving all other basis vectors fixed.
Hence the Weyl group $\Sym$ 
can be identified with the symmetric group $\Sym_m$.
Choose the representative $\sih_c\in\Norm\,\TT$ so that
its action on $\CC^m$ exchanges the basis vectors $e_c$ and $e_{c+1}\ts$,
leaving all other basis vectors of $\CC^m$ fixed.
Note that here the group $\DA$ is trivial.

Now consider $\G=\SP_{\tts2m}\com\SO_{\ts2m}\ts$. Let
$a$ and $b$ run through $-\ts m\lcd-1\com1\lcd m\ts$.
Then $e_a$ will denote a vector of the standard basis 
in $\CC^{\tts2m}\ts$, and $E_{ab}$ will be
a standard matrix unit in $\gl_{\tts2m}\ts$.
We will regard %this
$\G$ as the subgroup in $\GL_{\tts2m}$ 
preserving the bilinear form on $\CC^{\tts2m}$ whose value on
any pair $(e_a\com e_b)$ of the basis vectors is respectively
$$
\de_{a,-b}\cdot\sgn a 
\quad\text{or}\quad
\de_{a,-b}
$$
when $\G$ is $\SP_{\tts2m}$ or $\SO_{\tts2m}\ts$.
Then $\g$ is the Lie 
subalgebra of $\gl_{\tts2m}$ spanned by the elements
$$
F_{ab}=E_{ab}-\sgn  a\ts b\,\cdot E_{-b,-a}  
\quad\text{or}\quad
F_{ab}=E_{ab}-E_{-b,-a}\ts.  
$$

Choose the triangular
decomposition \eqref{hc0} of $\g=\sp_{2m}\com\so_{2m}$
where the subalgebras $\n\com\h\com\np$ are spanned
by the elements $F_{ab}$ with $a>b$, $a=b$, $a<b$ respectively.
The elements $F_{-a,-a}$ with $a>0$ 
form a basis in $\h\ts$. Here for any $a>0$
we will denote by $\eta_a$ the vector of the dual basis in $\h^*$ 
corresponding to the basis vector
$F_{a-m-1\tts,\ts a-m-1}\in\h\ts$.
The positive roots of $\sp_{2m}$ are %the weights 
$\eta_{\ts a}-\eta_{\ts b}$
and $\eta_{\ts a}+\eta_{\ts b}$ where $1\le a<b\le m\ts$, 
together with %the weights 
$2\ts\eta_{\ts a}$ where $1\le a\le m\ts$.
The positive roots of $\so_{2m}$ are only %the weights 
$\eta_{\ts a}-\eta_{\ts b}$
and $\eta_{\ts a}+\eta_{\ts b}$ where $1\le a<b\le m\ts$.
The semisimple rank $r$ of $\so_{2m}$ with $m>1$
and of $\sp_{2m}$ is $m\ts$. But the
the semisimple rank of $\so_2$ is zero, while the root system is empty.

For $\g=\sp_{2m}\com\so_{2m}$ and $c=1\lcd m-1$ choose 
$\al_c=\eta_{\ts c}-\eta_{\tts c+1}$ and
$$
H_c=F_{c-m-1\tts,\ts c-m-1}-F_{c-m,c-m}\ts,
\quad
E_c=F_{c-m-1\tts,\ts c-m}
\quad\text{and}\quad
F_c=F_{c-m\tts,\ts c-m-1}
$$
Further, if $\g=\sp_{2m}$ then choose $\al_m=2\ts\eta_m$ and
$$
H_m=F_{-1\tts,\tts-1}\,,
\quad
E_m=F_{-1\tts,\tts1}/2%=E_{-1\tts,\tts1}
\quad\text{and}\quad
F_m=F_{1\tts,\tts-1}/2%=E_{1\tts,\tts-1}
\,.
$$
If $\g=\so_{2m}$ and $m>1\ts$, then choose $\al_m=\eta_{m-1}+\eta_m$ and
$$
H_m=F_{-2\tts,\tts-2}+F_{-1\tts,\tts-1}
\,,
\quad
E_m=F_{-2\tts,\tts1}
\quad\text{and}\quad
F_m=F_{1\tts,\tts-2}
\,.
$$
For $\g=\sp_{2m}\com\so_{2m}$ choose the 
Cartan decomposition $\g=\g_{\tts+}\op\g_{\tts-}$
where $\g_{\tts+}$ is spanned by the elements $F_{ab}$ with
$a\com b>0$ while $\g_{\tts-}$ is spanned by
$F_{a,-b}$ and $F_{-a,b}$ with $a\com b>0\ts$.
The Chevalley anti-involution $\ep$ %on $\g$
will be defined by setting 
$\ep(F_{ab})=\sgn  a\ts b\,\cdot F_{ba}$ if $\th=1\ts$,
or $\ep(F_{ab})=F_{ba}$ if $\th=-1\ts$. Here the parameter 
$\th$ is the same as in Subsection~2.1.
This choice of $\ep$ is prescribed by the
condition (3) on our particular Howe systems, see below.
%which will be described below.

The algebraic group $\SP_{\tts2m}$ is connected, 
but $\SO_{\tts2m}$ has two connected components.
For $\G=\SP_{\tts2m}\com\SO_{\tts2m}$ the
maximal torus $\TT\subset\G_{\ts0}$ with the Lie algebra
$\h$ consists of 
all the elements of $\GL_{\tts2m}$
which multiply any two basis vectors $e_a\com e_{-a}\in\CC^{\tts2m}$
by scalars inverse to each other.
The Weyl group of $\sp_{2m}$ is isomorphic to the 
semidirect product $\Sym_m\ltimes\ZZ_2^m$ where
the symmetric group $\Sym_m$ acts by permutations of
the $m$ copies of $\ZZ_2\ts$. 
The Weyl group of $\so_{2m}$ is isomorphic to a
subgroup of $\Sym_m\ltimes\ZZ_2^m$ of index two.
For $c=1\lcd m-1$ and $\g=\sp_{2m}\com\so_{2m}$
the action of $\si_c\in\Sym$ on $\h^*$
exchanges the basis vectors
$\eta_{\ts c}$ and $\eta_{\ts c+1}\ts$, leaving 
other basis vectors fixed. Then
choose the representative $\sih_c\in\Norm_{\,0}\TT$ so that
its action on $\CC^{\tts2m}$ exchanges 
$e_{\ts c-m-1}$ and $e_{\ts c-m}\ts$, also exchanges  
$e_{\ts m-c+1}$ and $e_{\ts m-c}\ts$, leaving all other basis vectors 
of $\CC^{\tts2m}$ fixed.
For $\g=\sp_{2m}$ we have $\si_m(\eta_{\ts m})=-\ts\eta_{\ts m}$ and
$\si_m(\eta_{\ts a})=\eta_{\ts a}$ for $1\le a<m\ts$.
Choose the representative 
$\sih_m\in\SP_{\tts2m}$
so that $\sih_m(e_{-1})=e_1$ and $\sih_m(e_1)=-\ts e_{-1}$
while $\sih_m(e_a)=e_a$ for $|a|>1\ts$. Note that for $\g=\sp_{2m}$ 
the group $\DA$ is trivial, like it was in the case $\g=\gl_m\ts$.

Now take $\g=\so_{2m}\ts$.
Here the group $\DA$ is not trivial, but is isomorphic to $\ZZ_2\ts$.
Let $\tau_m$ be the generator of this group. Then
$\tau_m\ts(\eta_{\ts m})=-\ts\eta_{\ts m}$ and
$\tau_m\ts(\eta_{\ts a})=\eta_{\ts a}$ for $1\le a<m\ts$,
so that $\si_m=\tau_m\,\si_{m-1}\,\tau_m\ts$.
Choose the representative $\ttau_m\in\SO_{\tts2m}$
so that $\ttau_m\ts(e_{-1})=e_1$ and $\ttau_m\ts(e_1)=e_{-1}$
while $\ttau_m\ts(e_a)=e_a$ for $|a|>1\ts$.
Choose $\sih_m\in\SO_{\tts2m}$  
to be $\,\ttau_m\,\ts\sih_{m-1}\,\ttau_m\ts$.
We will need a representative in $\SO_{\tts2m}$
for every element $\om\in\ES\ts$. If
$\om\notin\Sym\ts$, then $\om=\tau_m\,\si$
for some $\si\in\Sym\ts$. In this case, the representative
of $\om$ in $\SO_{\tts2m}$ will be $\omh=\ttau_m\,\sih\ts$. 

We will now associate to any $\G$ another classical complex Lie group,
to be denoted by $\Gd$. Let $n$ be any positive integer.
If $\th=1$ then for $\G=\GL_m\com\SP_{\tts2m}\com\SO_{\tts2m}$
put $\Gd=\GL_n\com\SO_n\com\SP_n$ respectively. If $\th=-1$ then
put $\Gd=\GL_n\com\SP_n\com\SO_n$ respectively.
Here for $\Gd=\SP_n$ the integer $n$ is to be even. 
Then $(\G\com\Gd)$ is a \textit{reductive dual pair}~\cite{H1}.

Let $\gd$ be the Lie algebra of $\Gd$.
%We will always regard
%$\Gd=\SO_n\com\SP_n$ as subgroups in $\GL_n\ts$,
%and $\gd=\so_n\com\sp_n$ as Lie subalgebras in $\gl_n\ts$.
Let the indices $i$ and $j$ run through the sequence
$1\lcd n\ts$. Then $f_i$ will denote a vector of the standard basis 
in the vector space $\CC^n\ts$, while $E_{\ts ij}\in\gl_n$ will be
a standard matrix unit. If $i$ is even, put
$\bi=i-1\ts$. If $i$ is odd and $i<n$, put $\bi=i+1\ts$.
Finally, if $i=n$ and $n$ is odd, put $\bi=i\ts$. 
We will regard $\Gd=\SO_n$ or $\Gd=\SP_n$ as the subgroup in $\GL_n\ts$,
preserving the bilinear form on $\CC^n$ whose value on
any pair $(f_i\com f_j)$ of the basis vectors is 
$\th_i\,\de_{\ts\bi j}$
where $\th_i=1$ or $\th_i=(-1)^{\ts i-1}$ respectively.
Then $\gd$ is the Lie subalgebra of $\gl_n$ spanned by the elements
$
E_{ij}-\th_i\ts\th_j\ts E_{\ts\bj\ts\bi}\,.
$

Let $U$ be the tensor product of vector spaces $\CC^m\ot\CC^n\ts$,
and let $U^*$ be the dual vector space. For $a=1\lcd m$ and $i=1\lcd n$
let $x_{ai}$ denote the basis vector $e_a\ot f_i$ of $U\ts$.
Then let $\d_{ai}$ denote the corresponding vector of the dual basis
in $U^*\ts$. 
If $\th=1$ then $\P(U)$ is the algebra of polynomials
in $x_{ai}$ while $\Heist$ can be identified with the algebra
of differential operators on $\P(U)\ts$,
so that $\d_{ai}$ is the partial derivation corresponding to 
$x_{ai}\ts$. If $\th=-1$ then $\P(U)$ is the \textit{Grassmann algebra}
with $m\ts n$ anticommuting generators $x_{ai}\ts$.
The definition \eqref{Sh8a} then implies that the element
$x_{ai}\in\Heist$ acts on $\P(U)$ via left multiplication,
while $\d_{ai}\in\Heist$ acts as the left derivation relative to
$x_{ai}\ts$. The latter operator on the Grassmann algebra
is also called the \textit{inner multiplication} by $x_{ai}\ts$.
For any $\th$ the %non-degenerate symmetric bilinear 
form $\langle\ ,\,\rangle$ on $U$
will be chosen so that the basis of the %vectors
$x_{ai}$ is orthonormal. Then the involutive
anti-automorphism $\ep$ of the algebra $\Heist$ 
exchanges $x_{ai}$ with $\d_{ai}\ts$. 

Let us now consider the vector space $W=U\op U^*\ts$.
The groups $\GL_m$ and $\GL_n$ act on the vector space $U$, and their
actions commute which other. Hence we get the mutually commuting
actions of $\GL_m$ and $\GL_n$ on $W\ts$. Clearly,
both actions preserve the bilinear form $B$ in $W$
as introduced in Subsection 2.1.
If $\G=\SP_{\tts2m}\com\SO_{\tts2m}$ then we can identify
the vector space $W$ with the tensor product $\CC^{\tts2m}\ot\CC^n$
so that for $a=1\lcd m$ we have $x_{ai}=e_a\ot f_i$ as above,
and $\d_{ai}=e_{-a}\ot (\ts\th_i\,f_{\bi}\ts)\ts$. 
Then the bilinear form $B$ on $W$
gets identified with the tensor product of the forms on $\CC^{\tts2m}$
and $\CC^n\ts$ chosen above.
If $\th=1$ then the form $B$ is alternating and  
$\Gd=\SO_n\com\SP_n$ respectively.
If $\th=-1$ then the form $B$ is symmetric and  
$\Gd=\SP_n\com\SO_n$ respectively. 
Hence %for $\G=\SP_{\tts2m}\com\SO_{\tts2m}$ 
we always obtain mutually commuting
actions of $\G$ and $\Gd$ on the vector space $W\ts$,
preserving the bilinear form $B\ts$.

To complete the description of a Howe system on $U$ 
one requires a homomorphism $\zeta:\U(\g)\to\Heist$ obeying the conditions
(1) to (6) from Subsection 2.1. It is well~known 
and can be verified directly that $\zeta$ can be chosen so that
for $\g=\gl_m$ and $a\com b=1\lcd m$
\begin{equation}
\label{zetagl}
\zeta(E_{ab})\,=\,\sum_{i=1}^n\,
x_{ai}\,\d_{\ts bi}\,;
\end{equation}
while for $\g=\sp_{2m}\com\so_{2m}$ and 
the same indices $a\com b=1\lcd m$
\begin{gather}
\nonumber
\zeta(F_{ab})\,=\,
\th\,\de_{ab}\,n\ts/\ts2\,\,+\,
\sum_{i=1}^n\,x_{ai}\,\d_{\ts bi}\,,
\\
\label{zetag}
\zeta(F_{a,-b})\,=\,-\,\sum_{i=1}^n\,
\th\,\th_i\,x_{a\bi}\,x_{\ts bi}\,,
\ \quad
\zeta(F_{-a,b})\,=\,\sum_{i=1}^n\,
\th_i\,\d_{ai}\,\d_{\ts b\bi}\,.
\end{gather}
For $\g=\gl_m\com\sp_{2m}\com\so_{2m}$ and $\th=1,-1$
by \cite[Sections 2.3,\,3.5,\,3.8,\,4.2,\,4.3]{H2}
the image of the homomorphism\/ $\zeta$ %in the algebra $\Heist$ 
coincides with the subalgebra of\/ $\Gd$-invariant elements in\/ $\Heist$.

%\begin{proposition*}
%The image of the homomorphism\/ $\zeta$ %in the algebra $\Heist$ 
%coincides with the subalgebra of\/ $\Gd$-invariant elements in\/ $\Heist$.
%\end{proposition*}

Note that for any $i\com j=1\lcd n$ the element $E_{ij}\in\gl_n$
acts on $\P(U)$ as the %differential 
operator 
\begin{equation}
\label{gloper}
\sum_{a=1}^m\,
x_{ai}\,\d_{\ts aj}\,.
\end{equation}
Hence for $\gd=\so_n\com\sp_n$ the element
$
E_{ij}-\th_i\ts\th_j\ts E_{\ts\bj\ts\bi}\in\gd
$
acts on $\P(U)$ as the %differential 
operator 
\begin{equation}
\label{gdoper}
\sum_{a=1}^m\,
(\ts x_{ai}\,\d_{\ts aj}-\th_i\,\th_j\,x_{a\bj}\,\d_{\ts a\bi}\ts)\,.
\end{equation}
It is well known \cite{H1} that the subalgebra of
$\G\ts$-invariant elements in\/ $\Heist$ is generated by
the elements \eqref{gloper} for $\G=\GL_m$
or by the elements \eqref{gdoper} for $\G=\SP_{\tts2m}\com\SO_{\tts2m}\ts$. 
In Subsection 4.3 we give an analogue of this
result for the algebra \eqref{H2} instead~of~$\Heist\ts$.
 
%------------------------------------------------------------------------------

\vspace{-6pt}%%%%%%%%%%%%%%%%%%%%%%%%%%%%%%%%%%%%%%%%%%%%%%%%%%%%%%%%%%%%%%%%%%

\subsection*{\it\normalsize 4.2.\ Yangians}

First take the \textit{Yangian\/} $\Y(\gl_n)$ corresponding
to the Lie algebra $\gd=\gl_n\ts$.
This Yangian is a complex unital associative algebra 
with a family of generators
$
T_{ij}^{\ts(1)},T_{ij}^{\ts(2)},\ts\ldots
$
where $i\com j=1\lcd n\ts$.
Defining relations for these generators
can be written using the series
\begin{equation*}
T_{ij}(x)=
\de_{ij}+T_{ij}^{\ts(1)}x^{-\ns1}+T_{ij}^{\ts(2)}x^{-\ns2}+\,\ldots
\end{equation*}
where $x$ is a formal parameter. Let $y$ be another formal parameter.
Then the defining relations in the associative algebra $\Y(\gl_n)$
can be written as
\begin{equation}
\label{yrel}
(x-y)\,[\ts T_{ij}(x)\ts,T_{kl}(y)\ts]\ts=\;
T_{kj}(x)\ts T_{il}(y)-T_{kj}(y)\ts T_{il}(x)\,.
\end{equation}

The algebra $\Y(\gl_n)$ is commutative if $n=1\ts$.
By \eqref{yrel}, for any $z\in\CC$
the assignments
\begin{equation}
\label{tauz}
T_{ij}(x)\ts\mapsto\,T_{ij}(x+z)
\end{equation}
define an automorphism of the algebra $\Y(\gl_n)\ts$.
Here each of the formal
power series $T_{ij}(x+z)$ in $(x+z)^{-1}$ should be re-expanded
in $x^{-1}$, and every assignment \eqref{tauz} is a correspondence
between the respective coefficients of series in $x^{-1}$.
Relations \eqref{yrel} also show that for any
formal power series $g(x)$ in $x^{-1}$ with coefficients from
$\CC$ and leading term $1$, the assignments
\begin{equation}
\label{fut}
T_{ij}(x)\ts\mapsto\,g(x)\,T_{ij}(x)
\end{equation}
define an automorphism of the algebra $\Y(\gl_n)\ts$.
The subalgebra in $\Y(\gl_n)$ consisting
of all elements which are invariant under every automorphism 
of the form \eqref{fut}, is called the 
\textit{special Yangian\/} of $\gl_n\ts$, and will be denoted by
$\SY(\gl_n)\ts$. Two representations of the 
algebra $\Y(\gl_n)$ are called $\textit{similar\/}$ if they differ
by an automorphism of the form \eqref{fut}.
Similar representations of $\Y(\gl_n)$ have the same restriction
to the subalgebra $\SY(\gl_n)\ts$.

Using \eqref{yrel}, one can directly check that the assignments
\begin{equation}
\label{eval}
T_{ij}(x)\ts\mapsto\,\de_{ij}+E_{ij}\ts x^{-1}
\end{equation}
define a homomorphism of unital associative algebras
$\Y(\gl_n)\to\U(\gl_n)\ts$.
There is also an embedding $\U(\gl_n)\to\Y(\gl_n)\ts$,
defined by mapping $E_{ij}\mapsto T_{ij}^{\ts(1)}$. So
$\Y(\gl_n)$ contains the universal enveloping
algebra $\U(\gl_n)$ as a subalgebra. The homomorphism \eqref{eval} is
identical on the subalgebra $\U(\gl_n)\subset\Y(\gl_n)\ts$.

Let $T(x)$ be the $n\times n$ matrix
whose $i\com j$ entry is the series 
$T_{ij}(x)\ts$.
One can show that the assignment
\begin{equation}
\label{tin}
T(x)\mapsto T(-x)^{-1}
\end{equation}
defines an involutive automorphism of the algebra $\Y(\gl_n)\ts$.
Here each entry of the inverse matrix $T(-x)^{-1}$
is a formal power series in $x^{-1}$ with coefficients
from the algebra $\Y(\gl_n)\ts$, and the assignment \eqref{tin}
as a correspondence between the respective matrix entries.

The Yangian $\Y(\gl_n)$ is a Hopf algebra over the complex field $\CC\ts$.
The comultiplication $\Upde:\Y(\gl_n)\to\Y(\gl_n)\ot\Y(\gl_n)$ is defined by
the assignment
\begin{equation}
\label{1.33}
\Upde:\,T_{ij}(x)\ts\mapsto\ts\sum_{k=1}^n\ T_{ik}(x)\ot T_{kj}(x)\,.
\end{equation}
When taking tensor products of $\Y(\gl_n)\ts$-modules,
we use the comultiplication \eqref{1.33}.
The counit homomorphism
$\Y(\gl_n)\to\CC$ is defined by the assignment
$
T_{ij}(x)\ts\mapsto\ts\de_{ij}\ts.
$
The antipodal map $\Y(\gl_n)\to\Y(\gl_n)$ is defined by mapping
$
T(x)\mapsto T(x)^{-1}.
$
Note that the assignments
\begin{equation}
\label{tijtji}
T_{ij}(x)\mapsto T_{ji}(x)
\end{equation}
define an involutive anti-automorphism of
the associative algebra $\Y(\gl_n)$.
Moreover, they define a bialgebra anti-automorphism of $\Y(\gl_n)\ts$.

Note that the special Yangian $\SY(\gl_n)$ is a Hopf subalgebra
of $\Y(\gl_n)\ts$. It is isomorphic to
the Yangian $\Y(\mathfrak{sl}_n)$ of the special linear Lie
algebra $\mathfrak{sl}_n\subset\gl_n$ considered in \cite{D1,D2}.
For the proofs of the latter two assertions see 
\cite[Subsection 1.8]{M3}.

Now let $\gd$ be one of the two Lie algebras $\so_n\com\sp_n\ts$.
When considering these two cases simultaneously,
we will use the following convention. Whenever
the double sign $\ts\pm\ts$ or $\ts\mp\ts$ appears, 
the upper sign will correspond to the case
of a symmetric form on $\CC^{\ts n}$ so that $\gd=\so_n\ts$,
while the lower sign will corresponds to the case
of an alternating form on $\CC^{\ts n}$ so that $\gd=\sp_n\ts$.
Let $\Tt(x)$ be the transpose to the matrix $T(x)$
relative to that form on $\CC^{\ts n}\ts$.
The $i\com j$ entry of the matrix $\Tt(x)$ is
$\ts\th_i\ts\th_j\ts T_{\ts\bj\ts\bi\ts}(x)\ts$,
see Subsection 4.1.
An involutive automorphism of the algebra $\Y(\gl_n)$
can be then defined by the assignment
\begin{equation}
\label{transauto}
T(x)\mapsto\Tt(-x)\ts.
\end{equation}
This assignment is
understood as a correspondence between respective matrix entries.
Note that \eqref{transauto} defines an anti-automorphism
of the coalgebra $\Y(\gl_n)\ts$, like \eqref{tijtji} does.

Consider the product $\Tt(-x)\,T(x)\ts$.
The $i\com j$ entry of this %$n\times n$
matrix is the series
\begin{equation}
\label{yser}
\sum_{k=1}^n\,\th_i\ts\th_k\,T_{\,\bk\ts\bi\ts}(-x)\,T_{kj}(x)\,.
\end{equation}
The \textit{twisted Yangian\/}
corresponding to the Lie algebra $\gd$
is the subalgebra of $\Y(\gl_n)$ generated by
coefficients of all series \eqref{yser}.
We denote this subalgebra by $\Y(\gd)\ts$.
The subalgebra $\Y(\gd)\ts\cap\,\SY(\gl_n)$ of $\Y(\gl_n)$
is called the \textit{special twisted Yangian\/}
corresponding to $\gd$.
This subalgebra will be denoted by $\SY(\gd)\ts$.
The automorphism \eqref{fut} of $\Y(\gl_n)$
determines an automorphism of $\Y(\gd)$ which
multiplies the series \eqref{yser} by $g(x)\,g(-x)\ts$.
The subalgebra $\SY(\gd)$ of $\Y(\gd)$ consists of
the elements fixed by all such automorphisms. 
Two representations of %the algebra 
$\Y(\gd)$ are called $\textit{similar\/}$ if they differ
by such an automorphism.
Similar representations of $\Y(\gd)$ have the same restriction
to the subalgebra $\SY(\gd)\ts$.

To give defining relations for these generators of
$\Y(\gd)\ts$, let us introduce the \textit{extended twisted Yangian}
$\X(\ts\gd)\ts$. The complex unital associative algebra $\X(\ts\gd)$ 
has a family of generators
$
S_{ij}^{\ts(1)},S_{ij}^{\ts(2)},\ts\ldots
$
where $i\com j=1\lcd n\ts$.
Put
\begin{equation*}
S_{ij}(x)=
\de_{ij}+S_{ij}^{\ts(1)}x^{-\ns1}+S_{ij}^{\ts(2)}x^{-\ns2}+\,\ldots
\end{equation*}
and let $S(x)$ be the $n\times n$ matrix
whose $i\com j$ entry is the series $S_{ij}(x)\ts$.
Defining relations in the algebra $\X(\gd)$ can then be written as
\begin{gather}
\nonumber
(x^2-y^2)\,[\ts S_{ij}(x)\ts,S_{kl}(y)\ts]\ts=\ts
(x+y)(\ts S_{kj}(x)\,S_{il}(y)-S_{kj}(y)\,S_{il}(x))
\\[4pt]
\nonumber
\mp\,(x-y)\,(\,
\th_k\ts\th_j\,S_{i\ts\bk}(x)\,S_{\ts\bj\ts l}(y)-
\th_i\ts\th_l\,S_{k\ts\bi}(y)\,S_{\ts\bl\ts j}(x))
\\[2pt]
\label{xrel}
\pm\,\ts\th_i\ts\th_j\,
(\ts S_{k\ts\bi\ts}(x)\,S_{\ts\bj\ts l}(y)-
S_{k\ts\bi\ts}(y)\,S_{\ts\bj\ts l}(x))\,.
\end{gather}
These relations shows that for any
formal power series $f(x)$ in $x^{-1}$ with the coefficients from
$\CC$ and leading term $1$, an automorphism of the algebra
$\X(\gd)$ is defined by~mapping
\begin{equation}
\label{fus}
S_{ij}(x)\ts\mapsto\,f(x)\,S_{ij}(x)\ts.
\end{equation}

Let $\St(x)$ be the transpose to the matrix $S(x)$
relative to our form on $\CC^{\ts n}\ts$,
so that the $i\com j$ entry of $\St(x)$ is
$\ts\th_i\ts\th_j\ts S_{\ts\bj\ts\bi\ts}(x)\ts$.
By \cite[Theorem 2.3.13]{M3} there is a formal power series
$O(x)$ in $x^{-1}$ with the coefficients in the centre of
$\X(\gd)$ and leading term $1$,
such that
\begin{equation}
\label{5.2}
S(x)\mp2\ts x\,\St(x)=(1\mp2\ts x)\,O(x)\,S(-x)\ts.
\end{equation}
Moreover, then $O(x)\ts O(-x)=1\ts$.
Note that \eqref{5.2} yields an explicit formula for $O(x)\ts$.

One can define a homomorphism $\X(\gd)\to\Y(\gd)$ by mapping
\begin{equation}
\label{xy}
S(x)\,\mapsto\,\Tt(-x)\,T(x)\ts.
\end{equation}
Moreover, the kernel of the homomorphism \eqref{xy}
is generated by the coefficients of the series $O(x)-1\ts$.
For the proof of the last two statements see \cite[Section 2.13]{M3}.
By setting $O(x)=1$ in the matrix relation \eqref{5.2}
and then considering the $i\com j$ entry we get 
\begin{equation}
\label{srel}
S_{ij}(x)\mp2\ts x\,\th_i\ts\th_j\,S_{\ts\bj\ts\bi\ts}(x)
=(1\mp2\ts x)\,S_{ij}(-x)\ts.
\end{equation}
The homomorphism \eqref{xy} is surjective by its definition.
Thus the twisted Yangian $\Y(\gd)$ can also be defined as the
associative unital algebra with the generators
$
S_{ij}^{\ts(1)},S_{ij}^{\ts(2)},\ts\ldots
$
subject to the relations \eqref{xrel} and \eqref{srel}.
%For more details of the definition of the algebra $\Y(\gd)$
%see \cite[Chapter~2]{M3}.

Note that the anti-automorphism \eqref{tijtji}
and the automorphism \eqref{transauto} of the algebra $\Y(\gl_n)$
commute with each other. Their composition is an
involutive anti-automorphism of $\Y(\gl_n)$ which maps
\begin{equation}
\label{composit}
T_{ij}(x)\mapsto\ts\th_i\ts\th_j\ts T_{\ts\bi\ts\bj\ts}(-x)\ts.
\end{equation}
The composition \eqref{composit} preserves the
the subalgebra $\Y(\gd)\subset\Y(\gl_n)\ts$. The resulting
anti-automorphism of $\Y(\gd)$ can also be obtained as follows.
By using %the defining relations 
\eqref{xrel} only, one shows that the assignments
\begin{equation}
\label{sijsji}
S_{ij}(x)\mapsto S_{ji}(x)
\end{equation}
define an involutive anti-automorphism of the algebra $\X(\gd)\ts$;
see \cite[Proposition 2.3.4]{M3}. Moreover, due to
\eqref{srel} it factors to anti-automorphism of the algebra
$\Y(\gd)\ts$. The latter coincides with the restriction
of \eqref{composit} to the subalgebra $\Y(\gd)\subset\Y(\gl_n)\ts$.

Suppose we are given a representation of the algebra $\X(\gd)$
such that every coefficient of the series $O(x)$ is represented by
a scalar. Denote by $o\ts(x)$ the corresponding series
with scalar coefficients. Then $o\ts(x)\,o\ts(-x)=1$
and the leading term of $o\ts(x)$ is $1\ts$. Hence we can find
another formal power series $f(x)$ in $x^{-1}$ with scalar
coefficients and the leading term $1\ts$, such that 
$o\ts(x)=f(-x)/f(x)\ts$. By pulling back the given representation
of $\X(\gd)$ through the automorphism \eqref{fus} we then get
another representation of $\X(\gd)\ts$, which factors through
the homomorphism $\X(\gd)\to\Y(\gd)\ts$.
The series $f(x)$ and hence the resulting representation
of $\Y(\gd)$ are not unique. However, here we can only
replace $f(x)$ by $f(x)\,h(x)$ where $h(x)$ is a formal 
power series in $x^{-1}$ with scalar
coefficients and the leading term $1\ts$, such that
$h(x)=h(-x)\ts$. Then $h(x)=g(x)\,g(-x)$ for some series $g(x)$
as in \eqref{fut}. Hence all the representations of $\Y(\gd)$
corresponding to the given representation of $\X(\gd)$ are 
similar to each other, and have the same restriction to 
$\SY(\gd)\ts$.
  
The twisted Yangian $\Y(\gd)$ has
an analogue of the
homomorphism $\Y(\gl_n)\to\U(\gl_n)$ defined by \eqref{eval}.
Namely, one can define a homomorphism $\X(\gd)\to\U(\gd)$ by mapping
\begin{equation}
\label{pin}
S_{ij}(x)\,\mapsto\,\de_{ij}+
\frac{E_{ij}-\th_i\ts\th_j\ts E_{\ts\bj\ts\bi}}
{\textstyle x\pm\frac12}
\end{equation}
This can be proved by using the defining relations \eqref{xrel},
see \cite[Proposition~2.1.2]{M3}.
Moreover, the homomorphism \eqref{pin}
factors through the homomorphism $\X(\gd)\to\Y(\gd)$
defined by \eqref{xy}.
Further, there is an embedding $\U(\gd)\to\Y(\gd)$ defined by mapping
each element $E_{ij}-\,\th_i\ts\th_j\ts E_{\ts\bj\ts\bi}\in\gd$
%of the Lie algebra $\gd$
to the coefficient at $x^{-1}$ of the series \eqref{yser}. Hence
%the twisted Yangian
$\Y(\gd)$ contains the universal enveloping
algebra $\U(\gd)$ as a subalgebra.
The homomorphism $\Y(\gd)\to\U(\gd)$ corresponding to \eqref{pin}
is identical on the subalgebra $\U(\gd)\subset\Y(\gd)\ts$.

The twisted Yangian $\Y(\gd)$ is not only a subalgebra of
$\Y(\gl_n)\ts$, it is also a right coideal of the coalgebra
$\Y(\gl_n)$ relative to the comultiplication \eqref{1.33}. Indeed,
let us apply this comultiplication to the
$i\com j$ entry of the $n\times n$ matrix $\Tt(-x)\,T(x)\ts$.
We get the sum
\begin{gather*}
\sum_{k=1}^n\,\th_i\ts\th_k\,
\Upde\ts(\,T_{\,\bk\ts\bi\ts}(-x)\,T_{kj}(x))\,=
\\
\sum_{g,h,k=1}^n\,\th_i\ts\th_j\,
(\,T_{\ts\bk\ts\bg\ts}(-x)\ot T_{\ts\bg\ts\bi\ts}(-x))\,
(\,T_{kh}(x)\ot T_{hj}(x))\,=
\end{gather*}
%%%PAGEBREAK
$$
\sum_{g,h,k=1}^n\,\th_g\ts\th_k\,
T_{\ts\bk\ts\bg\ts}(-x)\,T_{kh}(x)
\ts\ot\ts
\th_i\ts\th_g\,T_{\ts\bg\ts\bi\ts}(-x)\,T_{hj}(x)\ts.\
$$
In the last displayed line, by performing the summation
over $k=1\lcd n$ in the first tensor factor we get
the $g\com h$ entry of the matrix $\Tt(-x)\,T(x)\ts$. Therefore
\begin{equation}
\label{rci}
\Upde\ts(\ts\Y(\gd))\subset\Y(\gd)\ot\Y(\gl_n)\ts.
\end{equation}

For the extended twisted Yangian $\X(\gd)\ts$,
one defines a homomorphism of associative algebras
$$
\X(\gd)\to\X(\gd)\ot\Y(\gl_n)
$$
by mapping
\begin{equation}
\label{comod}
S_{ij}(x)\,\mapsto
\sum_{g,h=1}^n\,S_{gh}(x)
\ts\ot\ts
\th_i\ts\th_g\,T_{\ts\bg\ts\bi\ts}(-x)\,T_{hj}(x)\ts.
\vspace{2pt}
\end{equation}
The homomorphism property can be verified directly, see \cite[Section 3]{KN3}.
Using the homomorphism \eqref{comod}, the tensor product
of any modules over the algebras $\X(\gd)$ and $\Y(\gl_n)$
becomes another module over $\X(\gd)\ts$.

Moreover, the homomorphism \eqref{comod} is a \textit{right coaction}
of the Hopf algebra $\Y(\gl_n)$ on the algebra $\X(\gd)\ts$.
Formally, one can define a homomorphism of associative algebras
$$
\X(\gd)\to\X(\gd)\ot\Y(\gl_n)\ot\Y(\gl_n)
$$
in two different ways: either by using the assignment \eqref{comod} twice,
or by using \eqref{comod} and then \eqref{1.33}. Both ways however
lead to the same result, see again \cite[Section 3]{KN3}.

%------------------------------------------------------------------------------

\vspace{-6pt}%%%%%%%%%%%%%%%%%%%%%%%%%%%%%%%%%%%%%%%%%%%%%%%%%%%%%%%%%%%%%%%%%%

\subsection*{\it\normalsize 4.3.\ Olshanski homomorphisms}

%Let $(\G\com\Gd)$ be any of the pair
%from Subsection 4.1, with $\th=1$ or $\th=-1\ts$. 
First consider the reductive dual pair 
$(\G\com\Gd)=(\GL_{\tts m}\com \GL_{\tts n})\ts$.
We will treat the cases of $\th=1$ and $\th=-1$ simultaneously. 
We shall use the Yangian $\Y(\gl_n)$ to describe the 
subalgebra of $\GL_m\ts$-invariant elements in the 
corresponding algebra \eqref{H2}.
Denote by $E$ the $m\times m$ matrix whose $a\com b$
entry is $E_{\ts ba}\in\gl_m\ts$. Note the transposition 
of the indices $a$ and $b$ here. 
The inverse matrix $(x+\th\ts E)^{\ts-1}$ will be regarded
as a formal power series in $x^{-1}$ whose coefficients
are certain $m\times m$ matrices with entries from the
algebra $\U(\gl_m)\ts$. Let $(x+\th\ts E)^{\ts-1}_{\ts ab}$
be the $a\com b$ entry of the inverse matrix. This entry equals
$$
\ts\de_{ab}\,x^{-1}-\th\ts E_{\ts ba}\ts x^{-2}
\,+\,\sum_{s=0}^\infty\,
\sum_{c_1,\ldots,c_s=1}^m
(-\th)^{\ts s+1}\,
E_{c_1a}\ts E_{c_2c_1}\ldots\ts E_{c_sc_{s-1}}\ts E_{\ts b\ts c_s}\,
x^{-s-2}\ts.
$$
The following result first appeared in \cite{O1}, 
although in another guise. For connections to our setting 
see \cite[Section 4]{KN1} or \cite[Section 4]{KN2} when
$\th=1$ or $\th=-1$ respectively.
Let $\U(\gl_m)^{\ts\GL_m}$ be
the subalgebra of $\GL_m\ts$-invariants in $\U(\gl_m)\ts$;
it coincides with the centre of the algebra $\U(\gl_m)\ts$.

\begin{proposition*}
\label{propo1}
{\rm(i)}
One can define a homomorphism $\Y(\gl_n)\to\U(\gl_m)\ot\Heist$ by %mapping
\begin{equation}
\label{homo1}
T_{ij}(x)\,\mapsto\,\de_{ij}\,+\sum_{a,b=1}^m
(\ts x+\th\,m/2+\th\ts E\ts)^{\ts-1}_{\ts ab}\ot x_{ai}\,\d_{\ts bj}\ts.
\end{equation}
{\rm(ii)}
The subalgebra of\/ $\GL_{\tts m}$-invariant elements in\/
$\Ar=\U(\gl_m)\ot\Heist$ coincides with the subalgebra, 
generated by\/ $\U(\gl_m)^{\ts\GL_m}\ot1$ and by
the image of homomorphism~\eqref{homo1}.
\end{proposition*}

\begin{proof}
Part (i) of the proposition was proved in
\cite[Section 1]{KN1} and \cite[Section~1]{KN2} for
$\th=1$ and $\th=-1$ respectively. The claim that
the image of the homomorphism~\eqref{homo1} belongs to 
the subalgebra $\Ar^{\GL_m}\subset\Ar\ts$,
was also proved therein. It remains to prove that the elements
\begin{equation}
\label{xd}
\sum_{a=1}^m\ 
1\ot x_{ai}\,\d_{\ts aj}
\end{equation}
and
\begin{equation}
\label{eexd}
\sum_{a,b,c_1,\ldots,c_s=1}^m
E_{c_1a}\ts E_{c_2c_1}\ldots\ts E_{c_sc_{s-1}}\ts E_{\ts b\ts c_s}
\ot\,x_{ai}\,\d_{\ts bj}\,,
\end{equation}
together with the elements of $\U(\gl_m)^{\ts\GL_m}\ot1\ts$, 
generate the whole subalgebra $\Ar^{\GL_m}\subset\Ar\ts$. Here
$s=0\com1\com2\com\ts\ldots$ and $i\com j=1\lcd n\ts$.
For $s=0$ the first tensor factor of the summand in 
\eqref{eexd} should be understood as $E_{\ts ba}\ts$. 

Let $\NN$ denote the additive semigroup of non-negative integers.
Take the standard $\NN\ts$-filtration on the algebra $\U(\gl_m)\ts$,
where any element of $\gl_m$ has degree $1\ts$.
The adjoint action of the group $\GL_m$ on $\U(\gl_m)$ 
preserves this filtration, and the corresponding graded algebra
is identified with the symmetric algebra $\S(\gl_m)\ts$.
The algebra $\Heist$ has an $\NN\times\NN\ts$-filtration, such
that the elements of $U$ and $U^\ast$ have degrees $(1\com0)$ and
$(0\com1)$ respectively.
This filtration is preserved by the action of the group
$\GL_m$ on $\Heist\ts$. The corresponding graded algebra
is identified with the symmetric algebra
of $U\op U^\ast$ in the case $\th=1\ts$,
or with the exterior algebra of $U\op U^*$ in the case $\th=-1\ts$.
In both cases, the graded algebra is denoted by
$\P\ts(\ts U\op U^\ast)\ts$.
It suffices to prove that the elements of the algebra 
$\S(\gl_m)\ot\P\ts(\ts U\op U^\ast)$
corresponding to \eqref{xd} and \eqref{eexd}, taken together with 
\begin{equation}
\label{sinv1}
\S(\gl_m)^{\ts\GL_m}\ot1
\,\subset\,
\S(\gl_m)\ot\P\ts(\ts U\op U^\ast)\,,
\end{equation}
generate the whole subalgebra of $\GL_m$-invariants. 
We will prove this for $n=1$ only.
The generalization of our proof to any $n\ge1$ will be obvious\ts;
cf.\ \cite[Subsection 2.9]{MO}.

So let us suppose that $U=\CC^m\ot\CC^1=\CC^m$. 
Let $e_1^{\ts\prime}\lcd e_m^{\ts\prime}$ be the basis
of the vector space $U^*$ dual to the standard basis $e_1\lcd e_m$ of $U\ts$.
First consider the $\GL_m$-invariants in the subspace
\begin{equation}
\label{pinv1}
1\ot\P\ts(\ts U\op U^\ast)
\,\subset\,
\S(\gl_m)\ot\P\ts(\ts U\op U^\ast)\ts.
\end{equation}
The subspace of $\P\ts(\ts U\op U^\ast)$ of degree
$(p\com q)\in\NN\times\NN$ is a subspace
in the tensor product $U^{\ts\ot\tts p}\ot(U^\ast){}^{\tts\ot\ts q}\ts$.
The latter tensor product contains non-zero 
$\GL_m$-invariant vectors, only if $p=q\ts$.
The $\GL_m$-invariants in $U^{\tts\ot\tts p}\ot(U^\ast){}^{\tts\ot\tts p}$
are any linear combinations of the~sums
\begin{equation}
\label{xxdd}
\sum_{c_1,\ldots,c_p=1}^m
e_{\ts c_{\varpi(1)}}\ot\ldots\ot e_{\ts c_{\varpi(p)}}
\ot 
e^{\ts\prime}_{c_1}\ot\ldots\ot e^{\ts\prime}_{c_p}
\end{equation}
where $\varpi$ ranges over all permutations of the indices $1\lcd p\ts$.
Here we used the classical invariant theory for $\GL_m\ts$;
see \cite[Theorem 1A]{H1}. 
For $\th=1$ or $\th=-1$ respectively,
by applying to the sum \eqref{xxdd}
the symmetrization or antisymmetrization 
in the first $p$ and also in the last $p$
tensor factors of the summand, we get an element of 
$$
\P\ts(U)\ot\P(U^\ast)=\P\ts(\ts U\op U^\ast)
$$
corresponding to %the element 
$$
\sum_{c_1,\ldots,c_p=1}^m
x_{c_11}\ldots x_{c_p1}
\,
\d_{c_11}\ldots\d_{c_p1}
\in\Heist
$$ 
multiplied by $\th^{\ts\ts\ell(\varpi)}\ts$.
This observation shows that for $n=1\ts$, the 
$\GL_m$-invariants in the subspace \eqref{pinv1}
are generated by the elements %of this subspace 
corresponding to \eqref{xd} with $i=j=1\ts$.

Now consider the subspace \eqref{sinv1}. Generators of the subalgebra
$\S(\gl_m)^{\ts\GL_m}\subset\S(\gl_m)$
are well known. We will reproduce a set of generators here,
because they will be used~in the subsequent argument.
Identify $\gl_m$ with $U\ot U^\ast$ as
a module of the group $\GL_m\ts$,
so that the matrix unit $E_{ab}\in\gl_m$ is identified with 
the vector $e_a\ot e_{\ts b}^{\ts\prime}\in U\ot U^*$.
Then the subspace in $\S(\gl_m)$
of degree $t$ becomes a subspace in 
$(\ts U\ot U^\ast){}^{\tts\ot\ts t}$.
The $\GL_m$-invariants in the latter tensor product
are any linear combinations of the sums
\begin{equation}
\label{ees}
\sum_{c_1,\ldots,c_t=1}^m
e_{\ts c_{\varpi(1)}}\ot e^{\ts\prime}_{c_1}
\ot\ldots\ot 
e_{\ts c_{\varpi(t)}}\ot e^{\ts\prime}_{c_t}
\end{equation}
where $\varpi$ ranges over all permutations of the indices $1\lcd t\ts$.
Let $s_1\lcd s_k$ be the cycle lengths of the permutation $\varpi\ts$,
so that $t=s_1+\ldots+s_k\ts$. By applying to the sum 
\eqref{ees} the symmetrization in the $t$ pairs of tensor factors 
of the summand, we get an element of $\S(\gl_m)\ts$, 
which corresponds to the product over $s=s_1\lcd s_k$ of the elements
\begin{equation}
\label{eeee}
\sum_{c_1,\ldots,c_s=1}^m
E_{c_2c_1}\ts E_{c_3c_2}\ldots\ts E_{c_sc_{s-1}}\ts E_{c_1c_s}
\in\U(\gl_m)^{\ts\GL_m}\,.
\end{equation}

\vbox{
Now take the subspace in $\S(\gl_m)\ot\P\ts(\ts U\op U^\ast)$
of degree $t+1$ in the first tensor factor and of degree $(p,q)$
in the second tensor factor, for any $t\ge0$ and $(p\com q)\neq(0\com0)\ts$.
Regard it as a subspace in the tensor product
$$
(\ts U\ot U^\ast){}^{\tts\ot\ts(t+1)}
\ot
U^{\ts\ot\tts p}\ot(U^\ast){}^{\tts\ot\ts q}\ts.
$$
The tensor product contains non-zero 
$\GL_m$-invariant vectors, only if $p=q\ts$.
Suppose this is the case.
Then the $\GL_m$-invariants in the tensor product
are any linear combinations of the sums over the indices 
$c_0\com c_1\lcd c_t\com c_{t+1}\lcd c_{t+p}=1\lcd m$ of
the vectors
\begin{equation}
\label{eeeeeeee}
e_{\ts c_{\varpi(0)}}\!\ot e^{\ts\prime}_{c_0}\!
\ot\ldots\ot
e_{\ts c_{\varpi(t)}}\!\ot e^{\ts\prime}_{c_t}\!
\ot
e_{\ts c_{\varpi(t+1)}}\!\ot\ldots\ot e_{\ts c_{\varpi(t+p)}}\!
\ot
e^{\ts\prime}_{c_{t+1}}\!\ot\ldots\ot e^{\ts\prime}_{c_{t+p}}\!
\end{equation}
where each of the sums corresponds to a permutation $\varpi$
of $0\com1\lcd t\com t+1\lcd t+p\ts$.
The $\GL_m$-invariants in $\S(\gl_m)\ot\P\ts(\ts U\op U^\ast)$
are obtained by applying to these linear combinations
the symmetrization of the first $t+1$ pairs of factors of 
the tensor product \eqref{eeeeeeee}, and %also
the symmetrization or antisymmetrization 
in the next $p$ and in the last $p$ factors,
for $\th=1$ or $\th=-1$ respectively.
Let $Q_\varpi$ be the element of 
$\S(\gl_m)\ot\P\ts(\ts U\op U^\ast)$
obtained by these operations from the sum  %of \eqref{eeeeeeee} 
corresponding to the permutation $\varpi\ts$.
}

If the cycle of an element of the set 
$\{\ts0\com1\lcd t\ts\}$ under
the action of the powers of $\varpi$
is contained in the set, then the element
$Q_\varpi\in\S(\gl_m)\ot\P\ts(\ts U\op U^\ast)$
is divisible in the first tensor factor
by the element of $\S(\gl_m)$
corresponding to \eqref{eeee},
where $s$ is the length of the cycle.
Then our argument reduces to a similar one with $t$
replaced by $t-s\ts$. Suppose that
the $\varpi\ts$-cycle of each element of
$\{\ts0\com1\lcd t\ts\}$ is not contained in this set.

Now consider the index $\varpi({t+1})\ts$. If this index 
belongs to the set $\{t+1\lcd t+p\ts\}$ then 
$Q_\varpi$ is divisible by the element of 
$\S(\gl_m)\ot\P\ts(\ts U\op U^\ast)$
corresponding to \eqref{xd} with $i=j=1\ts$,
and our argument reduces to a similar one with $p$
replaced by~$p-1\ts$. Suppose that
$\varpi({t+1})\in\{\ts0\com1\lcd t\ts\}\ts$.
Due to the symmetrization described above, 
without loss of generality we may assume that $\varpi({t+1})=0\ts$. 
Let $s\geqslant0$ be the minimal number such that 
$\varpi^{\ts s+1}(0)\notin\{\ts0\com1\lcd t\ts\}\ts$.
Again due to (anti)symmetrization, we may assume that 
$$
\varpi(0)=1\ts,\,\ 
\varpi(1)=2\ts,\,\ 
\ldots\ \ts,\ 
\varpi(s-1)=s
\ \quad\text{and}\ \quad
\varpi(s)=t+1\ts.
$$
Now we see that $Q_\varpi$ is divisible by the element of 
$\S(\gl_m)\ot\P\ts(\ts U\op U^\ast)$
corresponding to \eqref{eexd} with $i=j=1\ts$,
$a=c_0$ and $b=c_{t+1}\ts$. This observation
reduces our argument to a similar one where $t+1$ and $p$
are replaced by $t-s$ and $p-1$ respectively.
\qed
\end{proof}

Now let $(\G\com\Gd)$ be any of the reductive dual pairs
from Subsection 4.1 with $\G=\SP_{\tts2m}$ or $\G=\SO_{\tts2m}\ts$.
Again, we will treat the cases of $\th=1$ and $\th=-1$ simultaneously. 
Denote by $F$ the $2m\times 2m$ matrix whose $a\com b$
entry is the matrix unit $F_{ab}\in\g\ts$. 
Here the indices $a$ and $b$ run through %the sequence 
$-\ts m\lcd-1\com1\lcd m\ts$. Regard the
inverse matrix $(x-\th\ts F)^{\ts-1}$ 
as a formal power series in $x^{-1}$ whose coefficients
are $2m\times2m$ matrices with entries from the %universal enveloping 
algebra $\U(\g\ts)\ts$. Let $(x-\th\ts F)^{\ts-1}_{\ts ab}$
be the $a\com b$ entry of the inverse matrix. It equals
$$
\ts\de_{ab}\,x^{-1}+\th\ts F_{ab}\ts x^{-2}\,+\,
\sum_{s=1}^\infty\,
\sum_{|c_1|,\ldots,|c_s|=1}^m
\th^{\,s+1}\,
F_{a\ts c_1}\ts F_{c_1c_2}\ldots\ts F_{c_{s-1}c_s}\ts F_{\ts c_sb}\,
x^{-s-2}\ts.
$$
Like the indices $a$ and $b$, here the indices $c_1\lcd c_s$
run through $-\ts m\lcd-1\com1\lcd m\ts$.
In another guise, the next result appeared first in \cite{O2}.
For connections to our present setting 
see \cite[Section 6]{KN3} or \cite[Section 6]{KN4} when
$\th=1$ or $\th=-1$ respectively.

To state this result like Proposition \ref{propo1},
for $a<0$ put $p_{\ts ai}=x_{\ts-a,i}$ and $q_{\ts ai}=\d_{\ts-a,i}\ts$.
For $a>0$ put $p_{\ts ai}=-\,\th\,\th_i\,\d_{\ts a\ts\bi}$
and $q_{\ts ai}=\th_i\,x_{\ts a\ts\bi}\,$.
Note that then our definition of the homomorphism 
$\zeta:\U(\g)\to\Heist$ for $\g=\sp_{2m}\com\so_{2m}$ can be written as
$$
\zeta(F_{ab})\,=\,
\th\,\de_{ab}\,n\ts/\ts2\,\,-\,
\sum_{i=1}^n\,\th\,q_{\ts ai}\,p_{\,bi}
\ \quad\text{for}\ \quad
a\com b=-\ts m\lcd-1\com1\lcd m\ts.
$$
Let $\U(\g)^{\ts\G}$ be
the subalgebra of $\G\ts$-invariant elements in $\U(\g)\ts$.
This subalgebra coincides with the centre of $\U(\g)$
if $\g=\sp_{2m}\ts$, but is strictly contained in the centre if
$\g=\so_{2m}\ts$. 

\begin{proposition*}
\label{propo2}
{\rm(i)}
One can define a homomorphism\/ $\X(\gd)\to\U(\g)\ot\Heist$ by %mapping
\begin{equation}
\label{homo2}
S_{ij}(x)\,\mapsto\,
\de_{ij}\,+
\sum_{|a|,|b|=1}^m
(\ts x\pm{\textstyle\frac12}+\th\,m-\th\ts F\ts)^{\ts-1}_{\ts ab}\ot 
p_{\ts ai}\,q_{\ts bj}\,.
\end{equation}
{\rm(ii)}
The subalgebra of all\/ $\G$-invariant elements in\/
$\Ar=\U(\g)\ot\Heist\ts$ coincides with the subalgebra, generated by
$\U(\g)^{\ts\G}\ot1$ and by
the image of the homomorphism~\eqref{homo2}.
\end{proposition*}

\begin{proof}
Part (i) of the proposition was proved in
\cite[Section 6]{KN3} and \cite[Section~6]{KN4} for
$\th=1$ and $\th=-1$ respectively. The claim that
the image of the homomorphism~\eqref{homo2} belongs to 
the subalgebra $\Ar^{\G}\subset\Ar\ts$,
was also proved therein. It remains to prove that the elements
\begin{equation}
\label{pq}
\sum_{|a|=1}^m\ 
1\ot p_{\ts ai}\,q_{\ts aj}
\end{equation}
and
\begin{equation}
\label{ffpq}
\sum_{|a|,|b|,|c_1|,\ldots,|c_s|=1}^m
F_{ac_1}\ts F_{c_1c_2}\ldots\ts F_{c_{s-1}c_s}\ts E_{c_sb}
\ot p_{\ts ai}\,q_{\ts bj}\,,
\end{equation}
taken together with the elements of $\U(\g)^{\ts\G}\ot1\ts$, 
generate the whole subalgebra $\Ar^{\G}\subset\Ar\ts$. Here
$s=0\com1\com2\com\ts\ldots$ and $i\com j=1\lcd n\ts$.
For $s=0$ the first tensor factor of the summand in 
\eqref{eexd} should be understood as $F_{ab}\ts$. 

Consider the standard $\NN\ts$-filtration on the algebra $\U(\g)\ts$,
where any element of $\g$ has degree $1\ts$.
The adjoint action of the group $\G$ on $\U(\g)$ 
preserves this filtration, and the corresponding graded algebra
is identified with the symmetric algebra $\S(\g)\ts$.
The algebra $\Heist$ is generated by its subspace
$W=U\op U^\ast$, which is identified with
the tensor product $\CC^{\tts2m}\ot\CC^n\ts$.
By presenting $W$ as a direct sum of $n$ copies of 
$\CC^{\tts2m}$ we get an $\NN^{\ts\times n}$ filtration on
the algebra $\Heist\ts$, such that the $i\,$th direct summand 
$\CC^{\tts2m}$ of $W$ has degree $1$ 
in the $i\,$th factor $\NN$ of $\NN^{\ts\times n}$.
This filtration is preserved by the action of the group
$\G$ on $\Heist\ts$. The corresponding graded algebra
is identified with $\P(\CC^{\tts2m})^{\ts\ot\ts n}\ts$.
Here $\P(\CC^{\tts2m})$ denotes the symmetric algebra of 
$\CC^{\tts2m}$ if $\th=1\ts$, or the exterior algebra of
$\CC^{\tts2m}$ if $\th=-1\ts$.
It suffices to prove that the elements of the algebra 
$\S(\g)\ot\P(\CC^{\tts2m})^{\ts\ot\ts n}$
corresponding to \eqref{pq} and \eqref{ffpq}, 
taken together with the elements of the subalgebra
\begin{equation}
\label{sinv2}
\S(\g)^{\ts\G}\ot1
\,\subset\,
\S(\g)\ot\P(\CC^{\tts2m})^{\ts\ot\ts n}\ts,
\end{equation}
generate the whole subalgebra of $\G$-invariants. 
We will prove this by using 
the classical invariant theory for the group $\G\,$;
cf.\ \cite[Subsection 4.9]{MO}. 

Now $e_a$ with $a=-\ts m\lcd-1\com1\lcd m$
denotes a basis vector of $\CC^{\tts2m}$. Consider
the bilinear form on $\CC^{\tts2m}$
preserved by  the action of the group $\G\ts$.
Denote by $e^{\ts\prime}_a$ the vector
$$
\sgn a\cdot e_{-a}
\quad\text{or}\quad
e_{-a}
$$
when $\G\ts$ is $\SP_{\tts2m}$ or $\SO_{\tts2m}$ respectively.
The value of the %bilinear 
form on a pair $(e_a\com e^{\ts\prime}_b\ts)\ts$ equals $\de_{ab}$
for any indices $a$ and $b\ts$.
The odd tensor powers of $\CC^{\tts2m}$ do not contain
any non-zero $\G\ts$-invariant vectors.
The $\G\ts$-invariant vectors in $(\CC^{\tts2m})^{\tts\ot\tts2\tts t}$
are all linear combinations~of
\begin{equation}
\label{ppqq}
\sum_{|c_1|,\ldots,|c_t|=1}^m
e_{c_1}\ot e^{\ts\prime}_{c_1}
\ot\ldots\ot 
e_{c_p}\ot e^{\ts\prime}_{c_t}
\end{equation}
and of the sums, obtained from \eqref{ppqq} by any
simultaneous permutation of the $2\ts t$ tensor factors of every summand. 
For instance, see \cite[Theorem 1B]{H1}.

First let us consider the $\G\ts$-invariants in the subspace
\begin{equation}
\label{pinv2}
1\ot\P(\CC^{\tts2m})^{\ts\ot\ts n}
\,\subset\,
\S(\g)\ot\P(\CC^{\tts2m})^{\ts\ot\ts n}\ts.
\end{equation}
Regard the subspace of $\P(\CC^{\tts2m})^{\ts\ot\ts n}$
of degree $(\ts s_1\lcd s_n)\in\NN^{\ts\times n}$ as a subspace
in the tensor product
\begin{equation}
\label{cpcp}
(\CC^{\tts2m})^{\ts\ot\ts s_1}
\ot\,\ldots\,\ot
(\CC^{\tts2m})^{\ts\ot\ts s_n}\ts.
\end{equation}
Our subspace is obtained by applying to each
of the $n$ groups of the tensor factors $\CC^{\tts2m}$
of \eqref{cpcp} the symmetrization or antisymmetrization,
for $\th=1$ or $\th=-1$ respectively.
By using the explicit description of $\G\ts$-invariants 
in \eqref{cpcp} and arguing like in the %case $\G=\GL_m\ts$, 
proof of Proposition \ref{propo1}, 
one shows that the $\G\ts$-invariants in the subspace \eqref{pinv2} 
are generated by elements of this subspace
corresponding to \eqref{pq} with $i\com j=1\lcd n\ts$.
For example, if $n=2\ts t$ then %the sum 
\eqref{ppqq} may be regarded as a element of the graded algebra
$\P(\CC^{\tts2m})^{\ts\ot\ts n}$ of degree 
$(\ts s_1\lcd s_n)=(\ts1\lcd1)\ts$.
Then \eqref{ppqq} corresponds to the sum 
$$
\sum_{|c_1|,\ldots,|c_t|=1}^m
p_{\ts c_11}\,q_{c_11}
\ts\ldots\, 
p_{\ts c_tt}\,q_{c_tt}
\,\in\,\Heist\ts.
$$

Now consider the subspace \eqref{sinv2}. Generators of the subalgebra
$\S(\g)^{\ts\G}\subset\S(\g)$
are well known. We will reproduce a set of generators here,
because they will be used~in the subsequent argument.
If $\G=\SP_{\tts2m}$ then $\g$ can be identified with
the symmetric square of the $\G\ts$-module $\CC^{\tts2m}\ts$,
so that $F_{ab}\in\g$ is identified with the element
$$
(\ts e_a\ot e^{\ts\prime}_b+e^{\ts\prime}_b\ot e_a)/2
$$
of the symmetric square. If $\G=\SO_{\tts2m}$ then $\g$ 
can be identified with
the exterior square of the $\G\ts$-module $\CC^{\tts2m}\ts$,
so that $F_{ab}\in\g$ is identified with the element
$$
(\ts e_a\ot e^{\ts\prime}_b-e^{\ts\prime}_b\ot e_a)/2
$$
of the exterior square. Now for $\G=\SP_{\tts2m}\com\SO_{\tts2m}$
the subspace in $\S(\g)$ of degree $t$ becomes a subspace in 
the space of tensors $(\CC^{\tts2m})^{\tts\ot\tts2\tts t}$,
and we may use the explicit description of $\G\ts$-invariants 
in the latter space. This description implies that 
the elements of the graded algebra $\S(\g)$ corresponding to the elements 
$$
\sum_{|c_1|,\ldots,|c_s|=1}^m
F_{c_1c_2}\ts F_{c_2c_3}\ldots\ts F_{c_{s-1}c_s}\ts F_{c_sc_1}
\in\U(\g)^{\ts\G}
$$
with $s=1\com2\com\ts\ldots$
generate the subalgebra $\S(\g)^{\ts\G}\subset\S(\g)\ts$.
A similar argument was given
in the proof of Proposition \ref{propo1}\ts;
here we omit the details.

Now take the subspace in $\S(\g)\ot\P(\CC^{\tts2m})^{\ts\ot\ts n}$
of degree $t\ge1$ in the first tensor factor, and of any non-zero degree 
$(\ts s_1\lcd s_n)$ in the remaining $n$ tensor factors.
Regard it as a subspace in %the tensor product
\begin{equation}
\label{cpcpcp}
(\CC^{\tts2m})^{\tts\ot\tts2t}
\ot
(\CC^{\tts2m})^{\ts\ot\ts s_1}
\ot\,\ldots\,\ot
(\CC^{\tts2m})^{\ts\ot\ts s_n}\ts,
\end{equation}
by using the embedding of $\g$ into
the tensor square of $\CC^{\tts2m}$.
Our subspace is obtained by applying 
%the symmetrization or antisymmetrization
to each of the first $t$ pairs of tensor factors $\CC^{\tts2m}$ 
of \eqref{cpcpcp} the symmetrization or antisymmetrization
for $\g=\sp_{2m}$ or $\g=\so_{2m}$ respectively,
by symmetrizing these $t$ pairs of tensor factors $\CC^{\tts2m}$,
and by applying to each
of the $n$ groups of the remaining tensor factors $\CC^{\tts2m}$
of \eqref{cpcpcp} the symmetrization or antisymmetrization
for $\th=1$ or $\th=-1$ respectively.
Using the explicit description of $\G\ts$-invariants 
in \eqref{cpcpcp} and arguing like we did in the case $\G=\GL_m$
of Proposition \ref{propo1}, we complete the proof.
\qed
\end{proof}

%------------------------------------------------------------------------------

\vspace{-16pt}%%%%%%%%%%%%%%%%%%%%%%%%%%%%%%%%%%%%%%%%%%%%%%%%%%%%%%%%%%%%%%%%%

\subsection*{\it\normalsize 4.4.\ Tensor products}
 
Recall the following general definition. Let $\Fr$ be any
algebraic group over $\CC$ with a Lie algebra $\f\ts$. 
Let $\Br$ be any associative algebra where the group $\Fr$
acts by automorphisms. Consider the crossed product algebra
$\Fr\ltimes\Br\ts$. It is generated by the elements of 
$\Br$~and~$\Fr$ subject
to the relations in $\Br$ and $\Fr$ together with the cross relations
$$
\ups\,Z\,\ups^{-1}=\ups\ts(Z)
\quad\,\text{for}\,\quad
Z\in\Br
\quad\text{and}\quad
\ups\in\Fr\ts.
$$
Suppose there is also a homomorphism $\U(\f)\to\Br\ts$.
A module $K$ over the algebra $\Br$ and over the group $\Fr$
will be called a $(\Br\com\Fr)\ts$-\textit{module,} if
the joint actions of $\Br$ and $\Fr$ on $K$ make it a module over
the algebra $\Fr\ltimes\Br\ts$, and if
the action of $\f$ on $K$ corresponding to that of $\Fr$
coincides with the action of $\f$ obtained by pulling
the action of $\Br$ on $K$ 
back through the homomorphism $\U(\f)\to\Br\ts$.
%$$
%\renewcommand{\arraystretch}{1.25}
%\begin{array}{ccccc}
%\U(\f)\,&\xrightarrow{\ \ }&\!\!\Br\quad
%\\
%&\searrow&\!\!\big\downarrow\quad
%\\
%&&\!\!\!\!\!\!\!\!\!\!\End M\quad
%\end{array}
%$$
%where the diagonal arrow represents the action 
%of $\f$ on $M$ defined by that~of~$\Fr\ts$.
We suppose that $\Fr$ acts on $K$ locally finitely.
%so that the former action is well~defined.

Take any reductive dual pair
$(\G\com\Gd)$ from Subsection 4.1. We have $U=\CC^m\ot\CC^n\ts$.
The Lie algebra $\g$ acts on the vector space $\P\ts(U)$
via the homomorphism $\zeta:\U(\g)\to\Heist\ts$. 
The action of the algebra $\Heist$ on 
$\P\ts(U)$ is determined by \eqref{Sh8a}.
The group $\Gd$ acts on $\Heist$ via linear transformations of
the vector space $\CC^n\ts$. This action of $\Gd$
preserves the subspaces $U$ and $U^*$ of $\Heist\ts$. 
It also leaves invariant any element in the image of the
homomorphism $\zeta\ts$. Hence the actions of $\g$ and $\Gd$
on $\P\ts(U)$ commute with each other.

Now for any pair of weights $\la\com\mu\in\h^*$ 
consider the weight subspace
\begin{equation}
\label{lacoinv}
(M_\mu\ot\P\ts(U))_{\ts\n}^{\ts\la}
\subset 
(M_\mu\ot\P\ts(U))_{\ts\n}
\end{equation}
of the space of $\n\ts$-coinvariants of tensor product
of $\g\ts$-modules $M_\mu\ot\P\ts(U)\ts$.
We extend the action of the group $\Gd$ from $\P\ts(U)$
to the latter tensor product, so that $\Gd$ acts on the 
tensor factor $M_\mu$ trivially. Then $\Gd$ also acts on 
the pair of vector spaces \eqref{lacoinv}.

A bijective linear map
\begin{equation}
\label{iomu}
\P\ts(U)\to(M_\mu\ot\P\ts(U))_{\ts\n}
\end{equation}
can be defined by mapping any $u\in\P\ts(U)$ 
to the coset of $1_\mu\ot u$ in $(M_\mu\ot\P\ts(U))_{\ts\n}\ts$.
The algebra $\Ar$ acts on the vector space $M_\mu\ot\P\ts(U)$ by definition. 
By restricting the latter action to the subalgebra $\U(\g)\subset\Ar$
we get the diagonal action of $\g$ on $M_\mu\ot\P\ts(U)\ts$, as above.
When the $\Ar\ts$-module $M_\mu\ot\P\ts(U)$ is identified
with $\Ar_{\ts\mu}$ as in Subsection~2.3, then the pair of vector spaces
\eqref{lacoinv} gets identified with the pair $\Mr_\mu^\la\subset\Mr_\mu\ts$.

The group $\Gd$ acts by automorphisms of the algebra
$\Ar=\U(\g)\ot\Heist$ via its action on $\Heist\ts$.
The action of $\Gd$ on $\Ar$ preserves each of the ideals
$\I\com\J$ and $\Jp_\mu\ts$. Hence
$\Gd$ also acts on the pair $\Mr_\mu^\la\subset\Mr_\mu\ts$.
Its identification with  
\eqref{lacoinv} is that of pairs of $\Gd$-modules.

Note that the action of the Lie group $\Gd$ on the subspace
$\P\ts(U)\subset\Heist$ provides a homomorphism of algebras
$\U(\gd)\to\Heist\ts$. Any element of this image is
$\G\ts$-invariant. 
For $\Gd=\GL_n$ the image of the element
$E_{ij}\in\gl_n$ under this homomorphism equals \eqref{gloper}.
For $\Gd=\SO_n$ or $\Gd=\SP_n$ the image of the element
$
E_{ij}-\th_i\ts\th_j\ts E_{\ts\bj\ts\bi}\in\gd
$
equals \eqref{gdoper}. By identifying $\Heist$
with the subalgebra $1\ot\Heist\subset\Ar$ we get a homomorphism
$\U(\gd)\to\Ar^\G\ts$. Then \eqref{lacoinv} is a pair
of $(\Ar^\G\com\Gd)\ts$-modules, and so is the pair 
$\Mr_\mu^\la\subset\Mr_\mu\ts$.

More generally, for any element $\om$ of the group $\ES=\Sym\ltimes\DA$
consider the subspace 
\begin{equation}
\label{silacoinv}
(M_{\ts\om\tts\circ\mu}\ot\P\ts(\ts\omh\ts(U)))_{\ts\n}^{\,\om\tts\circ\la}
\subset 
(M_{\ts\om\tts\circ\mu}\ot\P\ts(\ts\omh\ts(U)))_{\ts\n}
\end{equation}
of the space of $\n\ts$-coinvariants of the tensor product
of $\g\ts$-modules $M_{\ts\om\tts\circ\mu}\ot\P\ts(\ts\omh\ts(U))\ts$.
A bijective linear map
\begin{equation}
\label{iomuom}
\P\ts(\ts\omh\ts(U))
\to
(M_{\ts\om\tts\circ\mu}\ot\P\ts(\ts\omh\ts(U)))_{\ts\n}
\end{equation}
can be defined by mapping any $w\in\P\ts(\ts\omh\ts(U))$
to the coset of $1_{\ts\om\tts\circ\mu}\ot w$ 
in the space of coinvariants.
Similarly to \eqref{Sh8a} the action of the algebra $\Heist$ on 
$\P\ts(\ts\omh\ts(U))$ is defined~by
$$
X(u)=v
\quad\text{if}\quad 
X\ts u-v\in\Heist\cdot\omh\ts(U^*)
$$
where $u\com v\in\P\ts(\ts\omh\ts(U))\ts$.
Then %the Lie algebra 
$\g$ acts on the vector space 
$\P\ts(\ts\omh\ts(U))$
via the homomorphism $\zeta:\U(\g)\to\Heist\ts$.
The pair of vector spaces \eqref{silacoinv} gets identified with the pair
\begin{equation}
\label{morgenpair}
\Mr_{\ts\om\tts\circ\mu\ts,\ts\omh\ts(U^*)}^{\,\om\tts\circ\la}
\subset 
\Mr_{\ts\om\tts\circ\mu\ts,\ts\omh\ts(U^*)}\ts.
\end{equation}
%see the definitions \eqref{nmummuv}.

The action of $\Gd$ on $\Heist$ commutes
with that of $\G\ts$. In particular, the action of
$\Gd$ preserves the subspaces $\omh\ts(U)$ and $\omh\ts(U^*)$
of $\Heist\ts$. Hence the group $\Gd$ also acts on the two pairs 
\eqref{silacoinv} and \eqref{morgenpair}. 
Our identification of them is that of pairs of $\Gd$-modules.
Note that both \eqref{silacoinv} and \eqref{morgenpair}
are also $(\Ar^\G\com\Gd)\ts$-modules,
like they were in the particular case above, when $\om$
was the identity element of $\ES\ts$.

The subalgebra $\Ar^\G$ of $\G\ts$-invariant elements of
$\Ar$ acts on the weight subspace \eqref{silacoinv}
by definition.  
First suppose that $\G=\GL_m\ts$. Like in the previous subsection, 
we will treat the cases of $\th=1$ and $\th=-1$ simultaneously.
Via the homomorphism $\Y(\gl_n)\to\Ar^{\GL_m}$ 
defined by \eqref{homo1}, the subspace in \eqref{silacoinv}
becomes a module over the Yangian $\Y(\gl_n)\ts$.
We shall now describe this module explicitly,
by using the comultiplication \eqref{1.33}
and the homomorphism \eqref{eval}. For $\g=\gl_m$
the group $\DA$ is trivial, while any $\omh$ preserves $U$ and $U^*$. 
Hence it suffices to consider only the case when
$\om$ is the identity element. %of $\Sym\ts$.

Every weight $\mu$ of $\gl_m$ is determined by the sequence 
$\mu_1\lcd\mu_m$ of its labels where $\mu_a=\mu(E_{aa})$
for each $a=1\lcd m\ts$. Note that here for the half-sum $\rho$ of the
positive roots $\rho_a=m/2-a+{\textstyle\frac12}\,$. 
For $\g=\gl_m$ we get $\ka=0$ by the definitions 
\eqref{deka},\eqref{zetagl}. Then $\nu=\la-\mu$ by \eqref{lamuka}.
Take the sequence $\nu_1\lcd\nu_m$ of labels
of $\nu\ts$. If $\th=1$ then suppose that each label
$\nu_a\in\{0\com1\com2\com\ts\ldots\ts\}\ts$.
If $\th=-1$ then suppose that each
$\nu_a\in\{0\com1\lcd n\ts\}\ts$.
Otherwise the weight subspace in \eqref{lacoinv} would contain zero only.

If $\th=1$ then denote by $\Ph^{\ts k}$
the $k\,$th symmetric power of the defining $\gl_n$-module $\CC^n\ts$.
If $\th=-1$ then denote by the same symbol $\Ph^{\ts k}$
the $k\,$th exterior power of the $\gl_n$-module $\CC^n\ts$.
The group $\Gd=\GL_n$ also acts on $\Ph^{\ts k}\ts$. 
Using the homomorphism \eqref{eval}, regard $\Ph^{\ts k}$ as
a module over the Yangian $\Y(\gl_n)\ts$. For $t\in\CC$
denote by $\Ph^{\ts k}_{\tts t}$ the $\Y(\gl_n)\ts$-module obtained by
pulling the $\Y(\gl_n)\ts$-module $\Ph^{\ts k}$ back through the
automorphism \eqref{tauz} where $z=\th\,t\ts$.

For $\th=1$ denote by
$\P\ts(\CC^n)$ the symmetric algebra of the vector space $\CC^n\ts$. 
For $\th=-1$ denote by $\P\ts(\CC^n)$ the exterior algebra of $\CC^n\ts$. 
For any $\th$ and $t$ the underlying vector space of the 
$\Y(\gl_n)\ts$-module $\Ph^{\ts k}_{\tts t}$
consists of all homogeneous elements of $\P(\CC^n)$ of degree $k\ts$.
Using the standard basis $e_1\lcd e_m$ of $\CC^m\ts$,
decompose the vector space 
$U=\CC^m\ot\CC^n$ into a direct sum of $m$ copies of $\CC^n\ts$.
Then the vector space $\P\ts(U)$ gets identified with the
tensor product of $m$ copies of $\P\ts(\CC^n)\ts$.
For the proof of the following proposition see
\cite[Section 2]{KN1} or \cite[Section~2]{KN2} when
$\th=1$ or $\th=-1$ respectively.

\begin{proposition*}
\label{iso1}
Under the above assumptions on the weight\/ $\nu\ts$,
the subspace in \eqref{lacoinv} is equivalent as a $\Y(\gl_n)\ts$-module
to the tensor product
\begin{equation}
\label{baspro1}
\Ph^{\,\nu_1}_{\mu_1+\rho_1+\frac12}
\ot\ts\ldots\ts\ot\ts
\Ph^{\,\nu_m}_{\mu_m+\rho_m+\frac12}\,.
\end{equation}
An equivalence map from the latter\/ $\Y(\gl_n)$-module
to the former is defined by regarding \eqref{baspro1}
as a subspace of\/ $\P\ts(U)$, and then applying
the map \eqref{iomu} to that subspace.
\end{proposition*}

The equivalence map here is also that 
of $\GL_n$-modules. This follows from the $\GL_n\ts$-equivariance
of the map \eqref{iomu}. 
The group $\GL_n$ also acts by automorphisms of the Hopf algebra
$\Y(\gl_n)\ts$. When an element $\ups\in\GL_n$ is regarded as
a $n\times n$ matrix by using the standard basis $f_1\lcd f_n$ of $\CC^n$,
the corresponding automorphism of
$\Y(\gl_n)$ is defined~by
%the assignment
$$
T(x)\mapsto\ups^{-1}\,T(x)\,\ups\ts.
$$
The Olshanski homomorphism $\Y(\gl_n)\to\Ar^{\GL_m}$ 
is $\GL_n\ts$-equivariant. This can be verified directly by
using the definition \eqref{homo1}. Furthermore,
we can consider the crossed product %algebra 
$\GL_n\ltimes\Y(\gl_n)\ts$. Moreover, we have
an embedding $\U(\gl_n)\to\Y(\gl_n)\ts$.
Then both modules in
Proposition \ref{iso1} become $(\ts\Y(\gl_n)\com\GL_n)\ts$-modules.
For the latter module this can be 
verified directly, by using the definition of $\Ph^{\ts k}_{\tts t}\ts$.
For the former this follows from %the 
$\GL_n\ts$-equivariance~of~\eqref{homo1},
by regarding $M_\mu\ot\P\ts(U)\ts$ as a module over the algebra 
$\GL_n\ltimes\Ar\ts$.

Now take $\om=\si$ with arbitrary $\si\in\Sym\ts$.
By replacing in Proposition~\ref{iso1} 
the weights $\la$ and $\mu$ respectively by 
$\si\circ\la$ and $\si\circ\mu\ts$, 
we get a description of the subspace in \eqref{silacoinv}
with $\om=\si\ts$
as a module over the Yangian $\Y(\gl_n)\ts$. 
Then $\nu$ is replaced by $\si(\nu)\ts$.
Let us identify
the Weyl group $\Sym$ with the symmetric group $\Sym_m$
so that any permutation of the dual basis vectors $\eta_1\lcd\eta_m\in\h^*$
by an element $\si\in\Sym$ 
corresponds to a permutation of the numbers $1\lcd m\ts$. 
The latter permutation will be denoted by the same symbol $\si\ts$.
Put
$$
\mut_a=\mu_{\ts\si^{-1}(a)}\ts,
\quad 
\nut_a=\nu_{\ts\si^{-1}(a)}\ts,
\quad
\rhot_a=\rho_{\ts\si^{-1}(a)}\ts.
$$
These are the $a\,$th labels of the weights 
$\si(\mu)\com\si(\nu)\com\si(\rho)$
respectively. In this notation, the subspace in \eqref{silacoinv} 
with $\om=\si$ is equivalent as a module
over $\Y(\gl_n)$ and $\GL_n$ to the tensor product
\begin{equation}
\label{bremod}
\Ph^{\,\nut_1}_{\ts\mut_1+\rhot_1+\frac12}
\ot\ts\ldots\ts\ot\ts
\Ph^{\,\nut_m}_{\ts\mut_m+\rhot_m+\frac12}\,.
\end{equation}

Now suppose that the weight $\la+\rho$ of $\gl_m$ is nonsingular.
By Corollary~\ref{corollary3.7}, then the Zhelobenko operator 
$\xic_{\ts\si}$ on $%\Jb\,\backslash\ts
\Ab\ts/\ts\Jb$ defines a linear map 
\eqref{zhelarb}. This map commutes with the actions of the 
algebra $\Ar^{\GL_m}$ on the source and the target vector spaces 
of \eqref{zhelarb}.
Hence it is an intertwining operator of $\Y(\gl_n)\ts$-modules.
By replacing the modules by their equivalents, we now get
an intertwining operator
\begin{equation}
\label{inter1}
\Ph^{\,\nu_1}_{\mu_1+\rho_1+\frac12}
\ot\ts\ldots\ts\ot\ts
\Ph^{\,\nu_m}_{\mu_m+\rho_m+\frac12}
\,\to\,
\Ph^{\,\nut_1}_{\ts\mut_1+\rhot_1+\frac12}
\ot\ts\ldots\ts\ot\ts
\Ph^{\,\nut_m}_{\ts\mut_m+\rhot_m+\frac12}
\end{equation}
of two tensor products of $\Y(\gl_n)\ts$-modules. It is well known 
that both tensor products are irreducible and equivalent to each other
if (but not only if) the weight $\mu$ is generic,
that is if $\mu_a-\mu_b\notin\ZZ$ whenever $a\neq b\,$;
see \cite[Theorem 4.8]{NT} for a more general result.
Hence for generic $\mu\ts$, an intertwining operator
between the two tensor products is unique up to scalar factor.
For our particular intertwining operator, this factor is
determined by the next proposition. This proposition is valid
for any weight $\mu\ts$, not  necessarily generic.
%For the proof see \cite[Section 3]{KN1} or \cite[Section~3]{KN2} when
%$\th=1$ or $\th=-1$ respectively.

Choose $\ph_k\in\Ph^{\ts k}$ as follows.
If $\th=1$ then $\Ph^{\ts k}$ is 
the $k\,$th symmetric power of $\CC^n\ts$,
and we put $\ph_k=f_1^{\tts k}\ts$.
If $\th=-1$ then $\Ph^{\ts k}$ is 
the $k\,$th exterior power of $\CC^n\ts$,
and we put $\ph_k=f_1\wedge\ldots\wedge f_k\ts$.
The vector $\ph_k$ is annihilated by
the action of the elements $E_{ij}\in\gl_n$ with $i<j\ts$. 
Take any positive root $\al=\eta_{\ts a}-\eta_{\ts b}$ of $\gl_m$
with~$a<b\ts$. For $\th=1$ define
$$
z_{\ts\al}\,=\,\prod_{s=1}^{\nu_b}\ 
\frac{\,\mu_a-\mu_b+\rho_a-\rho_b-s}{\,\la_a-\la_b+\rho_a-\rho_b+s}\ .
$$
Here the denominator corresponding to the running index $s$ 
equals $(\la+\rho)(H_\al)+s$\ts.
Hence the denominator does not vanish for any nonsingular
$\la+\rho\ts$. For $\th=-1\ts$ define 
$$
z_{\ts\al}\,=\,
(-1)^{\ts\nu_a\nu_b}\,
\left\{
\begin{array}{ll}
\displaystyle
\frac{\,\la_a-\la_b+\rho_a-\rho_b}{\,\mu_a-\mu_b+\rho_a-\rho_b}
&\quad\textrm{if}\quad\ \nu_a<\nu_b\,;
\\[12pt]
\hspace{39pt}1
&\quad\textrm{if}\quad\ \nu_a\ge\nu_b\,.
\end{array}
\right.   
$$
Here $\mu_a=\la_a-\nu_a$ for any index $a\ts$. Therefore in the first
of the last two cases,
the denominator equals $(\la+\rho)(H_\al)-\nu_a+\nu_b\ts$.
Hence it does not vanish for any nonsingular weight $\la+\rho\ts$,
under the condition $\nu_a<\nu_b$ when this denominator occurs.

\begin{proposition*}
\label{norm1}
Let $\la+\rho$ be nonsingular.
%Under the above assumptions on $\nu\ts$, 
The operator \eqref{inter1} determined by\/ $\xic_{\ts\si}$ 
maps the vector 
$\ph_{\ts\nu_1}\ot\ldots\ot\ph_{\ts\nu_m}$ of \eqref{baspro1}
to the vector
$\ph_{\ts\nut_1}\ot\ldots\ot\ph_{\ts\nut_m}$ of \eqref{bremod}
multiplied by
all those $z_{\ts\al}$ where $\al\in\De^+$ but $\si(\al)\notin\De^+\ts$.
%$$
%\prod\limits_{\substack{\al\in\De^+
%\\ 
%\si(\al)\notin\De^+}}\!\!\!z_{\ts\al}\,.
%$$
\end{proposition*}

\begin{proof}
For a generic weight $\mu$, this proposition has been proved
in \cite[Section 3]{KN1} and \cite[Section~3]{KN2} when
$\th=1$ or $\th=-1$ respectively. But when the weight $\nu$ is fixed,
our operator \eqref{inter1} depends on $\mu\in\h^*$
continiously, see the proof of Corollary \ref{corollary3.7}.
\qed
\end{proof}

For any $\mu$ and nonsingular $\la+\rho\ts$, 
Corollary~\ref{ircor} and Proposition~\ref{propo1} give
that for $\si=\si_0$ the quotient by
the kernel, or equivalently the image of our intertwining operator
\begin{equation}
\label{ak1}
\Ph^{\,\nu_1}_{\mu_1+\rho_1+\frac12}
\ot\ts\ldots\ts\ot\ts
\Ph^{\,\nu_m}_{\mu_m+\rho_m+\frac12}
\,\to\,
\Ph^{\,\nu_m}_{\mu_m+\rho_m+\frac12}
\ot\ts\ldots\ts\ot\ts
\Ph^{\,\nu_1}_{\mu_1+\rho_1+\frac12}
\end{equation}
is an irreducible $\Y(\gl_n)\ts$-module.
Here we also use the observation that 
any element of the subalgebra 
$\U(\gl_m)^{\tts\GL_m}\ot1\subset\Ar^{\tts\GL_m}$
acts on the subspace \eqref{lacoinv} by scalar multiplication.
In the next section we will show that up to equivalence
and similarity, every irreducible finite-dimensional 
$\Y(\gl_n)\ts$-module arises as the image of \eqref{ak1}
for $\th=-1$ and some~$\la\com\mu\ts$.
%For similar results, see \cite[Section 4]{AK}.

Now let $(\G\com\Gd)$ be any of the reductive dual pairs
from Subsection 4.1 with $\G=\SP_{\tts2m}$ or $\G=\SO_{\tts2m}\ts$.
We will continue treating the cases of $\th=1$ and $\th=-1$ simultaneously.
Via the homomorphism $\X(\gd)\to\Ar^{\G}$ 
defined by \eqref{homo2}, the subspace in \eqref{silacoinv}
becomes a module over the extended twisted Yangian $\X(\gd)\ts$.
We shall now describe this module explicitly.
For $\g\neq\gl_m$ we do not have $\ts\omh\ts(U)=U$ in general.
We will give an analogue of Proposition \ref{iso1} for arbitrary 
$\om\in\ES\ts$. For the proof of this analogue see
\cite[Section~5]{KN3} and \cite[Section~5]{KN4} in the cases
$\th=1$ and $\th=-1$ respectively.
We will keep using the $\Y(\gl_n)\ts$-modules $\Ph^{\ts k}_{\tts t}$
for integers $k\geqslant0\ts$. %as introduced above. 
Further, denote by $\Ph^{\ts-k}_{\ts t}$
the $\Y(\gl_n)\ts$-module obtained by
pulling the $\Y(\gl_n)\ts$-module $\Ph^{\ts k}_{\tts t}$ 
back through the automorphism \eqref{transauto}.

We regard $\gd$ is a Lie subalgebra of $\gl_n\ts$.
Any $\Y(\gl_n)\ts$-module can also be regarded as a
$\X\ts(\gd)\ts$-module, first by restricting from $\Y(\gl_n)$
to the subalgebra $\Y(\gd)\ts$, and then by pulling back through
the homomorphism $\X(\gd)\to\Y(\gd)$ defined by \eqref{xy}.
The resulting $\X(\gd)\ts$-module can also be defined 
as the tensor product of the initial $\Y(\gl_n)\ts$-module
by the one-dimensional trivial $\X(\gd)\ts$-module,
defined by the assignment $S_{ij}(x)\mapsto\de_{ij}$.
Here we use the coaction of $\Y(\gl_n)$ on $\X(\gd)$
defined by \eqref{comod}. Either way, any tensor product 
of $\Y(\gl_n)\ts$-modules
of the form $\Ph^{\ts k}_{\tts t}$ or $\Ph^{\ts-k}_{\ts t}$
will be regarded as a $\X(\gd)\ts$-module.
The group $\Gd$ acts on $\Ph^{\ts k}_{\tts t}$ 
and on $\Ph^{\ts-k}_{\ts t}$ 
by restricting the natural action of $\GL_n$ on $\Ph^{\ts k}$.

For any weight $\mu$ of $\g=\sp_{2m}$ or $\g=\so_{2m}$ define the sequence 
$\mu_1\lcd\mu_m$ of its labels by setting
$\mu_a=\mu\ts(F_{aa})$ for $a=1\lcd m\ts$. Note 
that unlike in the case of $\g=\gl_m\ts$, here %we have
\begin{equation}
\label{mulab}
\mu=-\,\mu_m\,\eta_1-\ldots-\mu_1\,\eta_m\,.
\end{equation}
For the half-sum of positive roots we get
$\rho_a=-\ts a$ if $\g=\sp_{2m}\ts$, or 
$\rho_a=1-a$ if $\g=\so_{2m}\ts$. By the definitions 
\eqref{deka},\eqref{zetag} here $\ka_a=\th\,n/2\ts$.
Consider the sequence $\nu_1\lcd\nu_m$ of the labels
the weight \eqref{lamuka}. Like for $\g=\gl_m\ts$,
suppose that each label
$\nu_a\in\{0\com1\com2\com\ts\ldots\ts\}$ if $\th=1\ts$.
If $\th=-1$ then suppose that each
$\nu_a\in\{0\com1\lcd n\ts\}\ts$.
Otherwise the weight subspace in \eqref{silacoinv} would contain zero only. 

The Weyl group of $\sp_{2m}$ is isomorphic to 
the semidirect product $\Sym_m\ltimes\ZZ_2^m\ts$,
and the group $\DA$ is trivial in this case.
The Weyl group of $\so_{2m}$ is isomorphic to
a subgroup of $\Sym_m\ltimes\ZZ_2^m$ of index two,
and $\DA$ is isomorphic to $\ZZ_2$.
The extended Weyl group $\ES=\Sym\ltimes\DA$ of $\so_{2m}$
is isomorphic to $\Sym_m\ltimes\ZZ_2^m$.
Fix the isomorphisms as follows.
Regard $\Sym_m\ltimes\ZZ_2^m$ as the group
of permutations of the indices $-\ts m\lcd-1\com1\lcd m$ 
such that if $a\mapsto b$ under a permutation, then 
$-\ts a\mapsto-\ts b$ under the same permutation.
The image of $\om\in\ES$ in $\Sym_m\ltimes\ZZ_2^m$
will be denoted by $\omb\ts$.
Then for any $c=1\lcd m-1$ the permutation $\sib_c$
only exchanges $c-m-1$ with $c-m\ts$, and
$m-c$ with $m-c+1\ts$.
For $\g=\sp_{2m}$ the permutation $\sib_m$ exchanges only $-1$ and $1\ts$.
For $\g=\so_{2m}$ the transposition of $-1$ and $1$ will be
denoted by $\bar\tau_m\,$; then 
$\bar\si_m=\bar\tau_m\,\sib_{m-1}\,\bar\tau_m\ts$.
The transposition $\bar\tau_m$ is the image in
$\Sym_m\ltimes\ZZ_2^m$ of the generator $\tau_m$ of $\DA$
for $\g=\so_{2m}\ts$. Now put 
$$
\mut_a=\mu_{\,|\ts\omb^{\ts-1}(a)|}\ts,
\quad 
\nut_a=\nu_{\,|\ts\omb^{\ts-1}(a)|}\ts,
\quad
\rhot_a=\rho_{\,|\ts\omb^{\ts-1}(a)|}
\quad\text{and}\quad
\de_a=\sgn\,\omb^{\,-1}(a)
$$
for $a=1\lcd m\ts$. In this notation,
the $a\,$th label of the weight $\om\ts(\mu)$
is equal to $\de_a\,\mut_a\ts$.

For any $a=1\lcd m$ let $\omh_{\tts a}$ be the element of the group
$\G$ such that $\omh_{\tts a}(e_{-a})=e_{a}\ts$, and
$\omh_{\tts a}(e_{a})=-\ts e_{-a}$
or
$\omh_{\tts a}(e_{a})=e_{-a}\ts$
depending on whether 
$\G=\SP_{\tts2m}$
or 
$\G=\SO_{\tts2m}\ts$.
By definition, 
all other basis vectors of $\CC^{\tts2m}$ are invariant
under the action of $\omh_{\tts a}\ts$. 
We have $\omh_{\tts a}\in\Norm\,\TT\ts$.
The image of $\omh_{\tts a}$ in the group $\ES$ will be denoted by $\om_a\ts$.
In particular,
$\om_{\tts 1}=\si_m$ in the case $\G=\SP_{\tts2m}\ts$, but
$\om_{\tts 1}=\tau_m$ in the case $\G=\SO_{\tts2m}\ts$.
The action of the element $\omh_{\tts a}\in\G$ on the vector space 
$W=\CC^{\tts2m}\ot\CC^n$ 
determines an automorphism of the algebra $\Heist$ such that 
for all $i=1\lcd n$
$$
x_{\ts ai}\mapsto-\,\th\,\th_i\,\d_{\ts a\bi}
\quad\ \text{and}\ \quad
\d_{\ts ai}\mapsto\th_i\,x_{\ts a\bi}\,,
$$
while the elements $x_{bi}\com\d_{bi}\in\Heist$ with $b\neq a$
are invariant under this automorphism.
This automorphism of $\Heist$ will
denoted by the same symbol $\omh_{\tts a}\ts$.

\begin{proposition*}
\label{iso2}
Under the above assumptions on the weight\/ $\nu\ts$,
the subspace in \eqref{silacoinv} is equivalent as module
over\/ $\X(\gd)$ to the tensor product
\begin{equation}
\label{pro2}
\Ph^{\,\ts\de_1\nut_1}_{\ts\mut_1+\rhot_1+\frac12}
\ot\ts\ldots\ts\ot\ts
\Ph^{\,\ts\de_m\nut_m}_{\ts\mut_m+\rhot_m+\frac12}
\end{equation}
pulled back through the automorphism \eqref{fus} of\/ $\X(\gd)$ where
$$
f(x)\,=\,\prod_{a=1}^m\,\ts
\frac{\,\th\ts x-\mu_a-\rho_a-\frac12}
{\,\th\ts x-\mu_a-\rho_a+\frac12}\ .
$$
An equivalence map from the latter $\X(\gd)\ts$-module 
to the former is defined by regarding \eqref{pro2}
as a subspace of\/ $\P\ts(U)\subset\Heist$, then applying to this subspace
all automorphisms $\omh_{\tts a}^{\ts-1}$ with\/ $\de_a=-1$, 
and then applying
the map \eqref{iomuom} to the resulting subspace of\/ $\Heist$.
\end{proposition*}

The equivalence map here is also that 
of $\Gd$-modules. This fact follows from the~$\Gd$-equivariance
of the map \eqref{iomuom}. 
The group $\Gd$ also acts by automorphisms of the 
right coideal subalgebra $\Y(\gd)\subset\Y(\gl_n)\ts$,
and by automorphisms of the right $\Y(\gl_n)\ts$-comodule 
algebra $\X(\gd)\ts$. When an element $\ups\in\Gd$ is regarded as
a $n\times n$ matrix  
by using the standard basis $f_1\lcd f_n$ of $\CC^n$,
the corresponding automorphism of
$\X(\gd)$ is defined by %the assignment
\begin{equation}
\label{upsaus}
S(x)\mapsto\ups^{-1}\,S(x)\,\ups\ts.
\end{equation}
It factors to an automorphism of the quotient $\Y(\gd)$ of 
the algebra $\X(\gd)\ts$.
The Olshanski homomorphism $\X(\gd)\to\Ar^{\G}$ 
is $\Gd$-equivariant. This can be verified directly by
using the definition \eqref{homo2}. 
We can consider the crossed product algebra 
$\Gd\ltimes\Y(\gd)\ts$. We also have an embedding 
$\U(\gd)\to\Y(\gd)\ts$,
defined by mapping
each element $E_{ij}-\,\th_i\ts\th_j\ts E_{\ts\bj\ts\bi}\in\gd$
to the coefficient at $x^{-1}$ of the series \eqref{yser}. 
Then \eqref{pro2} becomes a $(\ts\Y(\gd)\com\Gd)\ts$-module.
This statement can be verified directly, by using the definitions 
of $\Ph^{\ts k}_{\tts t}$ and $\Ph^{\ts-k}_{\tts t}\ts$.
%For the former module this follows from the 
%$\Gd$-equivariance of \eqref{homo2},
%by regarding $M_\mu\ot\P\ts(U)\ts$ as a module over the algebra 
%$\Gd\ltimes\Ar\ts$.

Note that the series $f(x)$ 
in Proposition \ref{iso2} does not depend on the choice
of $\om\in\ES\ts$. By using this proposition in the basic case
when $\om$ is the identity element, %of $\Sym\ts$,
the subspace in \eqref{lacoinv} is equivalent as an $\X(\gd)\ts$-module
to the tensor product of the form \eqref{baspro1}
pulled back through the automorphism \eqref{fus} of\/ $\X(\gd)$ where
the series $f(x)$ is as above. But here the labels
$\mu_a\com\nu_a\com\rho_a$ for $a=1\lcd m$ correspond
to the Lie algebra $\g=\sp_{2m}$ or $\g=\so_{2m}\ts$, not to
$\g=\gl_m\ts$.
In particular, here the labels $\rho_a$ 
of the half-sum $\rho$ of the positive roots 
are different from those for $\g=\gl_m\ts$.
Our choice \eqref{mulab} of the labels
of $\mu$ for  $\g=\sp_{2m}$ or $\g=\so_{2m}$
has been made so that here the description of the subspace
\eqref{lacoinv} becomes similar to that in the case $\g=\gl_m\ts$.
For the same purpose, in the definition
\eqref{homo1} we employed the matrix
$(\ts x+\th\,m/2+\th\ts E\ts)^{\ts-1}$ rather than 
the matrix $(x+\th\ts E)^{\ts-1}$.

Now we suppose that the weight $\la+\rho$ of $\g=\sp_{2m}$ or $\g=\so_{2m}$
is nonsingular. Then by Corollary~\ref{corollary3.7}, the Zhelobenko operator 
$\xic_{\ts\si}$ on $%\Jb\,\backslash\ts
\Ab\ts/\ts\Jb$ determines a linear map 
\eqref{zhelarb}. This map commutes with the actions of the algebra $\Ar^{\G}$
on the source and target vector spaces of \eqref{zhelarb}.
Further, for $\g=\so_{2m}$ and $\om=\tau_m\,\si$ the representative
$\ttau_m\in\G$ of $\tau\in\DA$ yields a linear map
\begin{equation}
\label{taint}
\Mr_{\ts\si\circ\mu\ts,\ts\sih\ts(U^*)}^{\ts\si\circ\la}
\,\to\,
\Mr_{\ts\om\tts\circ\mu\ts,\ts\omh\ts(U^*)}^{\,\om\tts\circ\la}\ts,
\end{equation}
see \eqref{taumap}.
This map also commutes with the actions of $\Ar^{\G}$ 
on the source and the target vector spaces. 
Note that $\ttau_m\in\SO_{2m}$ is an involution, and
the corresponding map \eqref{taint} is invertible.
Take the composition of the latter map with
the map \eqref{zhelarb} determined~by~$\xic_{\ts\si}\ts$.

Thus for $\g=\sp_{2m}$ or $\g=\so_{2m}$ and
for any $\om\in\ES$ we get a linear map 
\begin{equation}
\label{omint}
\Mr_\mu^\la\,\to\,
\Mr_{\ts\om\tts\circ\mu\ts,\ts\omh\ts(U^*)}^{\,\om\tts\circ\la}\ts,
\end{equation} 
commuting with the actions of the algebra $\Ar^{\G}$
on the source and target vector spaces.
Hence this is an intertwining operator of $\X(\gd)\ts$-modules.
By replacing the modules in \eqref{omint} by their equivalents, and
by using the observation that the series $f(x)$ for both modules 
is the same, we get an intertwining operator of $\X(\gd)\ts$-modules
\begin{equation}
\label{inter2}
\Ph^{\,\nu_1}_{\mu_1+\rho_1+\frac12}
\ot\ts\ldots\ts\ot\ts
\Ph^{\,\nu_m}_{\mu_m+\rho_m+\frac12}
\,\to\,
\Ph^{\,\ts\de_1\nut_1}_{\ts\mut_1+\rhot_1+\frac12}
\ot\ts\ldots\ts\ot\ts
\Ph^{\,\ts\de_m\nut_m}_{\ts\mut_m+\rhot_m+\frac12}
\,.
\end{equation}
The last two can also be regarded as $\Y(\gd)\ts$-modules,
and our operator intertwines~them. In this case we first take
the tensor products of $\Y(\gl_n)\ts$-modules, and then
restrict both tensor products to the subalgebra $\Y(\gd)\subset\Y(\gl_n)\ts$.
%These two $\Y(\gd)\ts$-modules are irreducible and equivalent 
%to each other if (but not only~if)
%the weight $\mu$ is generic, see \cite{MN}.

We will now give an analogue of Proposition \ref{norm1} for
$\G=\SP_{\tts2m}\com\SO_{\tts2m}\ts$. 
Let us arrange the indices $1\lcd n$ into the sequence
\begin{equation}
\label{corder}
1\com3\lcd n-1\com n\lcd4\com 2
\ \quad\text{or}\ \quad
1\com3\lcd n-2\com n\com n-1\lcd4\com 2
\end{equation}
when %the number
$n$ is even or odd respectively.
The mapping $i\mapsto\bi$ reverses the sequence \eqref{corder}.
We will write $i\prec j$ when $i$ precedes $j\/$ in this sequence. 
The elements
$
E_{\ts ij}-\th_i\ts\th_j\tts E_{\ts\bj\ts\bi}\in\gl_n
$
with $i\prec j$ or $i=j$ span a Borel subalgebra of $\gd\subset\gl_n\ts$,
while the elements $E_{\ts ii}-E_{\ts\bi\ts\bi}$
span the corresponding Cartan subalgebra of $\gd$.
Choose a vector $\psi_k\in\Phi^{\ts k}$ as follows.
For $\th=1$ put $\psi_k=f_1^{\tts k}\ts$,
so that $\ph_k=\psi_k$ in this case.
However, for $\th=-1$ let $\psi_k$ be the exterior product
of the vectors $f_i$ taken over the first $k$ indices in the sequence
\eqref{corder}. For instance, $\psi_2=f_1\wedge f_3$ if $n\ge3\ts$.
Note that the vector $\psi_k$ is always annihilated by
the action of the elements 
$E_{\ts ij}%%%-\th_i\ts\th_j\tts E_{\ts\bj\ts\bi}
\in\gl_n$ with $i\prec j\ts$.

For each positive root $\al$ of $\g=\sp_{2m}$ or $\g=\so_{2m}$
define $z_{\ts\al}\in\CC$ as follows. 
%For short, denote $\bar a=m-b+1$ for each $a=1\lcd m\ts$. 
If $\th=1\ts$,
$$
\hspace{-47pt}
z_\al\,=\,
\left\{
\begin{array}{ll}
\hspace{3.5pt}
\displaystyle
\prod_{s=1}^{\nu_{\ts b}}\,\hspace{5pt}
\frac{\mu_{\ts a}-\mu_{\ts b}+\rho_{\ts a}-\rho_{\ts b}-s}
{\la_{\ts a}-\la_{\ts b}+\rho_{\ts a}-\rho_{\ts b}+s}
&\quad\text{if}\quad
\al=\eta_{\ts m-b+1}-\eta_{\ts m-a+1}\,,
\\[16pt]
\hspace{3.5pt}
\displaystyle
\prod_{s=1}^{\nu_{\ts b}}\hspace{5pt}
\frac{\mu_{\ts a}+\mu_{\ts b}+\rho_{\ts a}+\rho_{\ts b}+s}
{\la_{\ts a}+\la_{\ts b}+\rho_{\ts a}+\rho_{\ts b}-s}
&\quad\text{if}\quad
\al=\eta_{\ts m-b+1}+\eta_{\ts m-a+1}\,,
\\[12pt]
\displaystyle
\prod_{s=1}^{[\nu_{\ts a}/2]}\,
\frac{\mu_{\ts a}+\rho_{\ts a}+s}{\la_{\ts a}+\rho_{\ts a}-s}
&\quad\text{if}\quad
\al=2\,\eta_{\ts m-a+1}\,.
\end{array}
\right.
$$
In the first two cases we have $1\le a<b\le m\ts$,
while in the third case $1\le a\le m$ and $\g=\sp_{2m}\ts$.
In the first case, the denominator
corresponding to $s$ equals $(\la+\rho)(H_\al)+s\ts$,
while in each of the last two cases the denominator equals 
$-\,(\la+\rho)(H_\al)-s\ts$. 
All these denominators do not vanish for a nonsingular $\la+\rho\ts$. 
If $\th=-1$, then
$z_{\ts\al}=z_{\ts\al}^{\ts\prime}\,z_{\ts\al}^{\ts\prime\prime}$ where 
\begin{align*}
z_{\ts\al}^{\ts\prime}\,&=\,
\left\{
\begin{array}{cll}
(-1)^{\ts\nu_a\nu_b}
&\quad\text{if}\quad
\al=\eta_{\ts m-b+1}-\eta_{\ts m-a+1}
&\quad\text{or}\quad
\al=\eta_{\ts m-b+1}+\eta_{\ts m-a+1}\ts,
\\[5pt]
\,1
&\quad\text{otherwise\ts;}&
\end{array}
\right.
\\[8pt]
z_\al^{\ts\prime\prime}\,&=\,
\left\{
\begin{array}{cll}
\,\dfrac{\la_{\ts a}-\la_{\ts b}+\rho_{\ts a}-\rho_{\ts b}}
{\mu_{\ts a}-\mu_{\ts b}+\rho_{\ts a}-\rho_{\ts b}}
&\quad\text{if}\quad
\al=\eta_{\ts m-b+1}-\eta_{\ts m-a+1}
&\quad\text{and}\quad\nu_a<\nu_b\,,
\\[12pt]
\,\dfrac{\la_{\ts a}+\la_{\ts b}+\rho_{\ts a}+\rho_{\ts b}}
{\mu_{\ts a}+\mu_{\ts b}+\rho_{\ts a}+\rho_{\ts b}}
&\quad\text{if}\quad
\al=\eta_{\ts m-b+1}+\eta_{\ts m-a+1}
&\quad\text{and}\quad\nu_a+\nu_b>n\,,
\\[12pt]
\,\dfrac{\la_{\ts a}+\rho_{\ts a}}
{\mu_{\ts a}+\rho_{\ts a}}
&\quad\text{if}\quad
\al=2\,\eta_{\ts m-a+1}
&\quad\text{and}\quad2\ts\nu_a>n\,,
\\[12pt]
\,1
&\quad\text{otherwise\ts.}&
\end{array}
\right.
\end{align*}
Here for $\th=-1$ we have $\mu_a=\la_a-\nu_a+n/2$ for each index $a\ts$.
Hence in the first of the last four cases, 
the denominator equals $(\la+\rho)(H_\al)-\nu_a+\nu_b\ts$.
In the second of these four cases,
the denominator equals $-\,(\la+\rho)(H_\al)-\nu_a-\nu_b+n\ts$.
In the third case, the denominator equals
$-\,(\la+\rho)(H_\al)-\nu_a+n/2\ts$.
These denominators do not vanish for any
nonsingular $\la+\rho\ts$, under the conditions
they occur with.
Denote by $\De^{++}$ 
the set of \textit{compact\/} positive roots of $\g\ts$,
these are the weights $\eta_a-\eta_{\tts b}$
where $1\le a<b\le m\ts$.

\begin{proposition*}
\label{norm2}
Let $\la+\rho$ be nonsingular.
The operator \eqref{inter2} determined by\/ $\xic_{\ts\si}$ 
maps the vector
%Under the above assumptions on $\nu\ts$, 
\begin{equation}
\label{ourvec2}
\psi_{\ts\nu_1}\ot\ldots\ot\psi_{\ts\nu_m}\in\,
\Ph^{\,\nu_1}_{\mu_1+\rho_1+\frac12}
\ot\ts\ldots\ts\ot\ts
\Ph^{\,\nu_m}_{\mu_m+\rho_m+\frac12}
\end{equation}
to the vector $\psi_{\ts\nut_1}\ot\ldots\ot\psi_{\ts\nut_m}$
of \eqref{pro2} multiplied by those $z_\al$ 
where $\om(\al)\notin\De^+\ts$ while
$$
\al\,\in\,
\left\{
\begin{array}{ll}
\De^{++}
&\quad\text{if}\quad
\th=1
\quad\text{and}\quad
n>1\,,
\\[2pt]
\De^+
&\quad\text{otherwise\ts.}
\end{array}
\right.
$$
\end{proposition*}

\begin{proof}
For a generic weight $\mu$, this proposition has been proved
in \cite[Section 5]{KN3} and \cite[Section~5]{KN4} when
$\th=1$ or $\th=-1$ respectively. But when the weight $\nu$ is fixed,
our operator \eqref{inter2} depends on $\mu\in\h^*$
continiously, see the proof of Corollary \ref{corollary3.7}.
\qed
\end{proof}

Our intertwining operator of $\Y(\gd)\ts$-modules
\eqref{inter2} has been defined for any element $\om\in\ES\ts$. 
By definition, the corresponding element $\omb\in\Sym_m\ltimes\ZZ_2^m$ 
is a certain permuation of the indices $-\ts m\lcd-1\com1\lcd m\ts$.  
Now consider the special case when the permutation
$\omb$ only changes the sign of each of these indices, 
so that $\om=\om_1\ldots\ts\om_m$
in the notation introduced just before stating
Proposition \ref{iso2}.
Thus we get an intertwining operator
\begin{equation}
\label{ak2}
\Ph^{\,\nu_1}_{\mu_1+\rho_1+\frac12}
\ot\ts\ldots\ts\ot\ts
\Ph^{\,\nu_m}_{\mu_m+\rho_m+\frac12}
\,\to\,
\Ph^{\,-\nu_1}_{\mu_1+\rho_1+\frac12}
\ot\ts\ldots\ts\ot\ts
\Ph^{\,-\nu_m}_{\mu_m+\rho_m+\frac12}
\end{equation}
of $\Y(\gd)\ts$-modules.
It corresponds to the longest
element $\om=\si_0$ of the Weyl group $\Sym\ts$,
if $\g=\sp_{2m}$ or if $\g=\so_{2m}$ and $m$ is even.
If $\g=\so_{2m}$ and $m$ is odd, then 
$\si_0=\om_2\ldots\ts\om_m$ so that %the intertwining operator 
\eqref{ak2} corresponds to
the element $\om=\tau_m\,\si_0$ of the extended Weyl group $\ES\ts$. 
%From now to the end of this section,
%we will consider $\g=\sp_{2m}$ and $\g=\so_{2m}$ separately.

For any $\g=\sp_{\tts2m}\com\so_{\tts2m}$ %and any $m$ 
the element $\om_2\ldots\ts\om_m\in\ES$ gives the
intertwining operator
\begin{equation}
\label{4.56}
\Ph^{\,\nu_1}_{\mu_1+\rho_1+\frac12}
\ot\ts\ldots\ts\ot\ts
\Ph^{\,\nu_m}_{\mu_m+\rho_m+\frac12}
\,\to\,
\Ph^{\,\nu_1}_{\mu_1+\rho_1+\frac12}
\ot
\Ph^{\,-\nu_2}_{\mu_2+\rho_2+\frac12}
\ot\ts\ldots\ts\ot\ts
\Ph^{\,-\nu_m}_{\mu_m+\rho_m+\frac12}
\end{equation}
of $\Y(\gd)\ts$-modules. The underlying 
vector spaces of the two $\Y(\gl_n)\ts$-modules
appearing as the first tensor factors on the right hand sides
of \eqref{ak2} and \eqref{4.56} are the same by definition.
Moreover, for $\g=\so_{\tts2m}$
our operators \eqref{ak2} and \eqref{4.56} 
are the same due to Proposition \ref{iso2}.
In particular, for $\g=\so_2$ the intertwining operator
$$
\Ph^{\,\nu_1}_{\mu_1+\rho_1+\frac12}
\,\to\,
\Ph^{\,-\nu_1}_{\mu_1+\rho_1+\frac12}
$$
corresponding to $\om=\om_1$ is the identity map.
This explains the following fact from \cite{N}.
If $\gd=\sp_n$ and $\th=1\ts$, or if
$\gd=\so_n$ and $\th=-1\ts$, that is if 
$\g=\so_{\tts2m}$ in the context of the present article,
then for any $t\in\CC$
the restriction of 
the $\Y(\gl_n)\ts$-module $\Ph^{\ts-k}_{\ts t}$
to the subalgebra $\Y(\gd)\subset\Y(\gl_n)$
coincides with the restriction of the
$\Y(\gl_n)\ts$-module $\Ph^{\ts k}_{\tts t}\ts$.
 
Now suppose that $\g=\sp_{2m}\ts$. Then we have $\ES=\Sym\ts$,
so that $\ES_\la=\Sym_\la$ for any $\la\in\h^*$ automatically.
Thus for any $\mu$ and nonsingular $\la+\rho\ts$, 
Corollary~\ref{ircor} and Proposition~\ref{propo2} imply that
the quotient by the kernel %(or equivalently the image) 
of our intertwining operator \eqref{ak2}
is an irreducible $\Y(\gd)\ts$-module.
Here we also use the observation that 
any element of the subalgebra 
$\U(\sp_{\tts2m})^{\tts\SP_{\tts2m}}\ot1\subset\Ar^{\tts\SP_{\tts2m}}$
acts on the subspace \eqref{lacoinv} via scalar multiplication.
Note that here 
$\gd=\so_n$ or $\gd=\sp_n$ respectively for $\th=1$ or $\th=-1\ts$.
In the next section we will show that up to equivalence
and similarity, any finite-dimensional irreducible
$\Y(\sp_n)\ts$-module arises as the image of
\eqref{ak2} for $\th=-1$ and some $\la\com\mu\ts$.

Next suppose that $\g=\so_{2m}\ts$. Then $\ES\neq\Sym\ts$.
Here for any $\mu$ and nonsingular $\la+\rho\ts$, 
Corollary~\ref{ircor} and Proposition~\ref{propo2} imply that
the quotient by the kernel, or equivalently the image 
of our intertwining operator \eqref{ak2}
is an irreducible $\Y(\gd)\ts$-module, under the
extra condition that $\ES_\la=\Sym_\la\ts$.
We also use the fact that any element of the subalgebra 
$\U(\so_{\tts2m})^{\tts\SO_{\tts2m}}\ot1\subset\Ar^{\tts\SO_{\tts2m}}$
acts on the subspace \eqref{lacoinv} via scalar multiplication.
In the next subsection we %will 
study the quotient by the kernel
of \eqref{ak2} without imposing that extra condition, 
but for $\th=-1$ only. We will show that then
the quotient is either an irreducible $\Y(\so_n)\ts$-module, 
or splits into a direct sum of two non-equivalent irreducible
$\Y(\so_n)\ts$-modules. In Section 5 we will 
explain which irreducible $\Y(\so_n)\ts$-modules arise 
in this particular way.

%------------------------------------------------------------------------------

\vspace{-6pt}%%%%%%%%%%%%%%%%%%%%%%%%%%%%%%%%%%%%%%%%%%%%%%%%%%%%%%%%%%%%%%%%%%

\subsection*{\it\normalsize 4.5.\ Crossed product algebras}

%We will begin this subsection with remarks valid 
For any pair $(\G\com\Gd)$
consider the crossed product algebra $\Gd\ltimes\Ar\ts$.
The action of $\Gd$ on $\Ar$ commutes with that of the group $\G\ts$,
and leaves invariant any element of the subalgebra $\U(\g)\subset\Ar\ts$.
Therefore the group $\Gd$ 
acts by automorphisms of the double coset algebra $\ZA\ts$, 
and we can also consider the crossed product algebra $\Gd\ltimes\ZA\ts$.
%The action of $\Gd$ on $\ZA$ preserves the 
%subalgebra $\SAzero\subset\ZA\ts$. Hence we get
%the crossed product algebra $\Gd\ltimes\SAzero\ts$.

The action of $\Gd$ on $\ZA$ commutes with 
the Zhelobenko automorphisms $\xic_1\lcd\xic_r\ts$.
Since it also commutes with the action of $\ttau\in\G$ on $\ZA$
for any $\tau\in\DA\,$,
the action of $\Gd$ preserves the subalgebra $\Dwo\subset\ZA\ts$.
So we get the crossed product algebra $\Gd\ltimes\Dwo\ts$.
The homomorphism $\U(\gd)\to\Ar^\G$ used in Subsection 4.4
yields a homomorphism $\U(\gd)\to\Dwo\ts$.

For any $\la\com\mu\in\h^*$ 
the subspace $\Mr_\mu^\la\subset\Mr_\mu$
is a module over the~subalgebra 
$$
\Gd\ltimes\Ar^\G\subset\Gd\ltimes\Ar\ts,
$$
and moreover an $(\Ar^\G\com\Gd)\ts$-module.
%see Subsection 4.4.
If $\la+\rho$ is nonsingular, the subspace $\Mr_\mu^\la\subset\Mr_\mu$
is also an $\Dwo\ts$-module, and moreover
an $(\Dwo\com\Gd)\ts$-module.
Then the action of $\Ar^\G$ on $\Mr_\mu^\la$
can also be obtained by pulling the action of %the algebra
$\Dwo$ on $\Mr_\mu^\la$
back through the isomorphism $\ga:\Ar^\G\to\Dwo\ts$, see 
Subsection~3.3. Moreover, the isomorphism $\ga$
is $\Gd$-equivariant. Hence by using $\ga$ we get the same 
structure of an $(\Ar^\G\com\Gd)\ts$-module on $\Mr_\mu^\la$ as above.

The Zhelobenko operator $\xic_{\tts0}$ on
$%\Jb\backslash
\Ab\ts/\ts\Jb$ is $\Gd$-equivariant.
So is the corresponding linear map \eqref{zhelong}
for any $\mu$ and nonsingular $\la+\rho\ts$.
Hence the action of the group $\Gd$
preserves the subspace 
$\Ker\ts(\ts\xic_{\tts0}\ts|\ts\Mr_\mu^\la\ts)$ of
$\Mr_\mu^\la\ts$.
Another way to see this is to use 
Proposition \ref{corollary3.8} and
$\Gd$-contravariance of
the Shapovalov form $S_\mu^\la$ on $\Mr_\mu^\la\ts$. 
The latter property means~that 
$$
S_\mu^\la(\ts\ups\ts(f)\com g)=S_\mu^\la(f\com\ups^{\ts\prime}(g))
\quad\text{for}\quad f\com g\in\Mr_\mu^\la
$$
where $\ups\mapsto\ups^{\ts\prime}$ is the anti-involution on $\Gd$
defined by the matrix transposition.
Here the elements of $\Gd$ are regarded as matrices
by using the standard basis $f_1\lcd f_n$ of $\CC^n\ts$.
Either way, the quotient space
of $\Mr_\mu^\la$ by $\Ker\ts(\ts\xic_{\tts0}\ts|\ts\Mr_\mu^\la\ts)$
becomes an $(\Dwo\com\Gd)\ts$-module.

When the Lie group $\Gd$ is connected, that is 
when $\Gd=\GL_n$ or $\Gd=\SP_n\ts$, we will not need
to use the action of $\Gd\ltimes\Dwo$ on $\Mr_\mu^\la\ts$.
It will suffice to use only the action of $\Dwo\ts$.
But the Lie group $\Gd=\SO_n$ is not connected.
Until the end of this subsection we will be considering
only the case when $(\G\com\Gd)=(\SO_{\tts2m}\com\SO_n)$
so that $\th=-1\ts$.
The space $U$ has been identified with 
the tensor product $\CC^m\ot\CC^n\ts$. In our case
$\Heist$ is the Clifford algebra generated  
by the elements of the vector space $W=\CC^{\tts2m}\ot\CC^n$
subject to the relations \eqref{H00}, where
$B$ is the tensor product of the symmetric forms  
on $\CC^{\tts2m}$ and $\CC^n$ preserved by the actions of the groups
$\SO_{\tts2m}$ and $\SO_n$ respectively.
Here $x_{ai}=e_a\ot f_i$ and $\d_{ai}=e_{-a}\ot f_{\bi}\,$
for $a=1\lcd m$ and for all $i=1\lcd n\ts$. 
Note that $\th_i=1$ in this case.

Choose any vector $f_0\in\CC^n$ of length $\sqrt{2}$ 
with respect to the form preserved by $\SO_n\ts$. 
The corresponding orthogonal reflection is
\begin{equation*}
\label{reflect}
\ups_{\ts0}:\,\CC^n\to\CC^n:\,u\,\mapsto\,u-z\ts f_0
\end{equation*}
where $z$ is the value of the symmetric form on $\CC^n$ taken on the pair 
of vectors $(f_0\com u)\ts$. This reflection is
an element of the group $\SO_n\ts$.
It determines an automorphism of the algebra $\Heist\ts$,
which preserves the subalgebra $\P\ts(U)$
generated by all the elements $x_{ai}\ts$.

For $a=1\lcd m$ define
the vectors $x_a=e_a\ot f_0$ and $\d_a=e_{-a}\ot f_0$
of $\CC^{\tts2m}\ot\CC^n\ts$. By \eqref{H00} we get 
\begin{equation}
\label{xddx}
x_a\ts\d_a+\d_a\ts x_a=B\ts(x_a\com\d_a)=2\ts.
\end{equation}
Consider the product
$$
A_{\ts0}\,=\,\prod_{a=1}^m\,(\ts1-x_a\ts\d_a)\,\in\,\Heist\ts.
$$
%in $\Heist\ts$. 
The $m$ factors of this product pairwise
commute. We will use the next properties of~$A_{\ts0}\ts$.

\begin{lemma*}
\label{iiiiii}
{\rm(i)}
There is an equality $\ttau_m\ts(A_{\ts0})=-\ts A_{\ts0}\ts$.
\\
{\rm(ii)}
The action of the element $A_{\ts0}$ on $\P\ts(U)$
coincides with that of the reflection $\ups_{\ts0}\ts$.
\\
{\rm(iii)}
The element $A_{\ts0}\in\Heist$ commutes with 
$\zeta(X)$ for every $X\in\so_{2m}\ts$.
\end{lemma*}

\begin{proof}
In our case, 
the action of the element $\ttau_m\in\SO_{\tts2m}$ 
on $\CC^{\tts2m}$ exchanges
$e_1$ with~$e_{-1}\ts$, and leaves other basis vectors of $\CC^{\tts2m}$
fixed. So the action of $\ttau_m$ on $\Heist$ exchanges the element 
$x_1$ with $\d_1\ts$, and leaves fixed any element $x_a$ or $\d_a$ 
with $a>1\ts$. Now (i) follows from the relation \eqref{xddx} 
with $a=1\ts$. Part (ii) can be obtained by a direct calculation, which
reduces to the case $m=1\ts$. Part (iii) follows from (ii),
because every element $\zeta(X)\in\Heist$ with $X\in\so_{2m}$ is
$\SO_n\ts$-invariant.
\qed
\end{proof}

\begin{proposition*}
\label{proposition4.8}
Suppose that $\la+\rho$ is nonsingular and\/
$(\G\com\Gd)=(\SO_{\tts2m}\com\SO_n)\ts$.
Then for any weight\/ $\mu\in\h^*$ the quotient\/
$\Mr_\mu^\la\,/\,\Ker(\ts\xic_{\tts0}\ts|\ts\Mr_\mu^\la\ts)$
is an irreducible\/ $(\ts\Dwo\com\SO_n)\ts$-module.
\end{proposition*}

\begin{proof}
By Propositions \ref{proposition2} and \ref{corollary3.8} 
the quotient
$N=\Mr_\mu^\la\,/\,\Ker(\ts\xic_{\tts0}\ts|\ts\Mr_\mu^\la\ts)$
is an irreducible $\SAzero\ts$-module.
By applying Proposition \ref{proposition3.9}
to this $\SAzero\ts$-module, we get Proposition \ref{proposition4.8} 
when $\ES_\la=\Sym_\la\ts$.
In this case Proposition \ref{proposition4.8}  
follows from Theorem \ref{theorem2}.
We will now modify~our proof of Proposition \ref{proposition3.9},
to get Proposition \ref{proposition4.8} without assuming that
$\ES_\la=\Sym_\la\ts$.

Instead of the element \eqref{essym} of the algebra $\Dwo$
let us first consider an element of $\SAzero\ts$,
$$
Z'=|\Sym_\lambda|^{-1}\sum_{\si\in\Sym}\xic_{\ts\si}(X'\,Y)\ts.
$$
Here $X'$ and $Y$ are the same as in \eqref{essym},
but the sum is taken over the Weyl group $\Sym\ts$.
%not the extended Weyl group $\ES\ts$.
The arguments given just after \eqref{essym} show that
$\chi(Y)=\chi(Z')\ts$. We also get the equality
$\xic_{\ts\si}(Z')=Z'$ for every $\si\in\Sym\ts$.
But we might have $\ttau_m\ts(Z')\neq Z'$,
so that $Z'\notin\Dwo$ then.

Then put $Z=Z'+\ttau_m\ts(Z')\ts$.
Since $\ttau_m^{\,2}=1\ts$,
we get the equality $\ttau_m\ts(Z)=Z\ts$.
The action of the element $\ttau_m\in\SO_{2m}$ on $\SAzero$
exchanges the operator $\xic_{m-1}$ with $\xic_{m}\ts$,
and commutes with the operators $\xic_1\lcd\xic_{m-2}\ts$.
Hence $\xic_{\ts\si}(Z)=Z$ for any $\si\in\Sym\ts$. 
Thus we have $Z\in\Dwo\ts$, like we had for the element \eqref{essym}.  

Now put $Z''=Z'-Z\ts$. Then 
$\xic_{\ts\si}(Z'')=Z''$ for any $\si\in\Sym\ts$,
while $\ttau_m\ts(Z'')=-\ts Z''\ts$.
But here we have $Z''A_{\ts0}\in\Dwo\ts$.
Indeed, by Part (iii) of Lemma \ref{iiiiii}
we have $\xic_{\ts\si}(Z''A_{\ts0})=Z''A_{\ts0}$ for any $\si\in\Sym\ts$,
while $\ttau_m\ts(Z''A_{\ts0})=Z''A_{\ts0}$ by Part (i). 
Now consider the element 
$$
Z''A_{\ts0}\,\ups_{\ts0}^{-1}\in\SO_n\ltimes\Dwo\ts.
$$
Its action on $N$ coincides with that of the element $Z''$
due to Part (ii) of Lemma \ref{iiiiii}.
Hence the action of
$$
Z+Z''A_{\ts0}\,\ups_{\ts0}^{-1}\in\SO_n\ltimes\Dwo
\hspace{24pt}
$$
coincides with that of the element $Z'\ts$.
But the latter action coincides with that of $Y$.
\qed
\end{proof}

\begin{corollary*}
\label{iror}
Suppose $\la+\rho$ is nonsingular and\/
$(\G\com\Gd)=(\SO_{\tts2m}\com\SO_n)\ts$,
so that $\th=-1\ts$.
\\[2pt]
{\rm(i)}
If\/ $n$ is odd then the quotient\/
$\Mr_\mu^\la\,/\,\Ker(\ts\xic_{\tts0}\ts|\ts\Mr_\mu^\la\ts)$
is an irreducible\/ $\Ar^{\SO_{\tts2m}}$-module.
\\
{\rm(ii)}
If\/ $n$ is even then\/ %the quotient\/
$\Mr_\mu^\la\,/\,\Ker(\ts\xic_{\tts0}\ts|\ts\Mr_\mu^\la\ts)$
is either an irreducible\/ $\Ar^{\SO_{\tts2m}}$-module,
or splits into a direct sum of two irreducible non-equivalent\/
$\Ar^{\SO_{\tts2m}}$-modules.
\end{corollary*}

\begin{proof}
If $n$ is odd then the group $\SO_n$ splits as a direct product
$\ZZ_2\times\mathrm{SO}_n$ where $\ZZ_2$ is the subgroup of $\SO_n$
generated by the minus identity element. This element acts
on $\Mr_\mu^\la$ as the multiplication by
$(-1)^{\ts\nu_1+\ts\ldots\ts+\nu_m}$.
Instead of the action of the connected Lie group $\mathrm{SO}_n$ on 
$\Mr_\mu^\la$ it suffices to consider the action of the Lie algebra
$\so_n\ts$. But the latter action can
also be obtained by pulling the action of $\Ar^{\SO_{2m}}$ on
$\Mr_\mu^\la$ back through the homomorphism
$\U(\so_n)\to\Ar^{\SO_{2m}}$, as mentioned earlier 
in this subsection. Proposition~\ref{knv} now yields~(i).

When $n$ is even, we can use the action of the Lie algebra 
$\so_n$ on $\Mr_\mu^\la$ 
instead of the action of the subgroup $\mathrm{SO}_n\subset\SO_n\ts$.
This is a normal subgroup of index two.
We can also use the homomorphism $\U(\so_n)\to\Ar^{\SO_{2m}}$,
as we did for an odd $n\ts$. Proposition \ref{knv} now implies that
the quotient $\Mr_\mu^\la\,/\,\Ker(\ts\xic_{\tts0}\ts|\ts\Mr_\mu^\la\ts)$
is irreducible over the joint action of %the algebra
$\Ar^{\SO_{2m}}$ and of any element $\ups$ from the complement
to $\mathrm{SO}_n$ in $\SO_n\ts$.
Now general arguments from \cite[Section V.8]{W}
yield (ii). Moreover, if the quotient
is a direct sum of two irreducible
$\Ar^{\SO_{\tts2m}}$-modules and $N$ is one of them, %then 
the space of the other equals $\ups\ts N$ for any~$\ups\ts$.
\qed
\end{proof}

%------------------------------------------------------------------------------

\vspace{-14pt}%%%%%%%%%%%%%%%%%%%%%%%%%%%%%%%%%%%%%%%%%%%%%%%%%%%%%%%%%%%%%%%%%

\subsection*{\it\normalsize 4.6.\ Equivalent modules}

Let $(\G\com\Gd)$ be any dual pair from Subsection 4.1.
For any weights $\la$ and $\mu$ of $\g\ts$, the
subspace \eqref{2.4444} can be regarded as an $\Ar^\G\ts$-module.
We have denoted the subspace~by~$\Lr_\mu^\la\ts$.
First consider $(\G\com\Gd)=(\GL_m\com\GL_n)\ts$.
By using the homomorphism $\Y(\gl_n)\to\Ar^{\GL_m}$
we can also regard $\Lr_\mu^\la$ as an $\Y(\gl_n)\ts$-module.
By Corollary \ref{corollary3.11plus}
the $\Y(\gl_n)\ts$-module $\Lr_\mu^\la$ is irreducible
for any $\mu$ and nonsingular $\la+\rho\ts$.
Here we also use Propositions \ref{knv} and \ref{propo1}.

Further suppose that the weight
$\nu=\la-\mu$ satisfies the conditions of
Proposition~\ref{iso1}. Using that proposition, 
the $\Y(\gl_n)\ts$-module $\Lr_\mu^\la$ is equivalent
to a certain quotient of the tensor product \eqref{baspro1}. 
The latter tensor product can be regarded as a subspace
in $\P\ts(U)\ts$. Denote by  
$u$ the vector of $\P\ts(U)$ corresponding to
the vector $\ph_{\ts\nu_1}\ot\ldots\ot\ph_{\ts\nu_m}$ 
of \eqref{baspro1}. If the weight $\la+\rho$ is nonsingular, then
%we have the equality 
$\zeta(\p[\mu+\rho])\ts u=z\ts u$
where $z$ stands for the product of all
$z_\al$ with $\al\in\De^+\ts$.
This equality
follows from Propositions \ref{proposition3.7} and \ref{norm1}.
It can also be obtained directly from
\cite[Section 3]{KN1} or \cite[Section~3]{KN2}
when $\th=1$ or $\th=-1$ respectively.
If $z\neq0$ then by Propositions \ref{lmeq} and \ref{corollary3.8}
the $\Y(\gl_n)\ts$-module $\Lr_\mu^\la$
is equivalent to the image of our intertwining operator \eqref{ak1}. 
Then $\Lr_\mu^\la$ is not zero.

Now take any pair $(\G\com\Gd)$ from Subsection 4.1
other than $(\GL_m\com\GL_n)\ts$.
By using the homomorphism $\X(\gd)\to\Ar^\G$ we can regard
$\Lr_\mu^\la$ as an $\X(\gd)\ts$-module.
Moreover, $\Lr_\mu^\la$ can be then regarded as
a $(\ts\X(\gd),\Gd)\ts$-module, see the beginning of Subsection 4.5.
Take any $\mu$ and nonsingular $\la+\rho\ts$.  
By Corollary \ref{corollary3.11plus}
the $\X(\gd)\ts$-module $\Lr_\mu^\la$ is irreducible,
if $\ES_\la=\Sym_\la\ts$. Here we also
use Propositions~\ref{knv}~and~\ref{propo2}.
Next take $(\G\com\Gd)=(\SO_{\tts2m}\com\SO_n)\ts$,
so that $\th=-1\ts$.
Then our proof of Proposition \ref{proposition4.8}
demonstrates that $\Lr_\mu^\la$ is an irreducible
$(\ts\X(\so_n)\com\SO_n)\ts$-module.
Our proof of Corollary \ref{iror}
demonstrates that $\Lr_\mu^\la$ is an irreducible
$\X(\so_n)\ts$-module, if $n$ is odd.
If $n$ is even then $\Lr_\mu^\la$
is either an irreducible $\X(\so_n)\ts$-module,
or a direct sum of two irreducible non-equivalent
$\X(\so_n)\ts$-modules.

Now consider again any pair $(\G\com\Gd)$ from Subsection 4.1
other than $(\GL_m\com\GL_n)\ts$.
%but not only $(\SO_{\tts2m}\com\SO_n)\ts$.
Suppose that the weight
$\nu=\la-\mu-\ka$ satisfies  the conditions of
Proposition~\ref{iso2}. Using that proposition in the case
when $\om$ is the identity element of the group $\ES\ts$,
the $\X(\gd)\ts$-module $\Lr_\mu^\la$ is equivalent
to a certain quotient of the tensor product of the form \eqref{baspro1}.
This tensor product can be regarded as a subspace
in $\P\ts(U)\ts$. 
Let $u$ be the vector of $\P\ts(U)$ corresponding to
\eqref{ourvec2}. If $\la+\rho$ is nonsingular, then
we again have the equality 
$\zeta(\p[\mu+\rho])\ts u=z\ts u\ts$.
But here $z$ stands for the product of all
$z_\al$ with $\al\in\De^{++}$ if
$\th=1$ and $n>1\ts$. Otherwise 
$z$ stands for the product of all
$z_\al$ with $\al\in\De^+\ts$.
The equality follows from Propositions 
\ref{proposition3.7} and \ref{norm2}.
It can also be obtained directly 
from \cite[Erratum]{KN3} or from \cite[Section~5]{KN4}
when $\th=1$ or $\th=-1$ respectively.
If $z\neq0$ then by Propositions \ref{lmeq} and \ref{corollary3.8}
the $\X(\gd)\ts$-module $\Lr_\mu^\la$
is equivalent to the image of our intertwining operator \eqref{ak2}. 
Moreover, then $\Lr_\mu^\la$ is
equivalent to the image of \eqref{ak2}
as an $(\ts\X(\gd)\com\Gd)\ts$-module.
Note that then the quotient $\Lr_\mu^\la$ of $\Mr_\mu^\la$
is not zero, see Subsection 2.5.

%------------------------------------------------------------------------------

\subsection*{\it\normalsize 4.7.\ Dual modules}

In Subsection 3.4 for any $\mu$ and nonsingular $\la+\rho$
we defined a $\Dwo\ts$-contravariant pairing \eqref{qlm}. 
Let us now describe that pairing in terms of representations of Yangians.
First take $(\G\com\Gd)=(\GL_m\com\GL_n)\ts$. The Chevalley
anti-involution $\ep$ on the Lie algebra $\g=\gl_m$ is defined
by $\ep\ts(E_{ab})=E_{\ts ba}\ts$. It extends to an involutive
anti-automorphism of the algebra $\Ar\ts$, also denoted by $\ep\tts$,
so that the extension exchanges the generators $x_{ai}$ and $\d_{ai}$
of the subalgebra $\Heist\subset\Ar\ts$. First
applying the homomorphism $\Y(\gl_n)\to\Ar$ defined by \eqref{homo1},
and then applying the anti-automorphism $\ep$ of $\Ar\ts$, amounts to
first applying the anti-automorphism \eqref{tijtji} of $\Y(\gl_n)\ts$, and then
applying the homomorphism \eqref{homo1}. Therefore,
when we regard the two vector spaces in the pairing \eqref{qlm}
as $\Y(\gl_n)\ts$-modules, the pairing becomes contravariant relative to the 
%involutive 
anti-automorphism \eqref{tijtji} of $\Y(\gl_n)\ts$.

By replacing these two $\Y(\gl_n)\ts$-modules by their equivalents, we get
a non-degenerate contravariant pairing of the source and target
$\Y(\gl_n)\ts$-modules in \eqref{ak1}.
In particular, the $\Y(\gl_n)\ts$-module 
dual to the source in \eqref{ak1}, is equivalent to the target.
%$\Y(\gl_n)\ts$-module. 
The latter equivalence can be proved directly, for all $\la$ and $\mu\ts$.
Indeed, because \eqref{tijtji} is a coalgebra anti-automorphism,
it is enough to consider the case $m=1$ only. But it is 
well known that the $\Y(\gl_n)\ts$-modules $\Ph^{\ts k}_{\tts t}$ 
are equivalent to their duals relative to \eqref{tijtji}, 
whenever $\th=1$ or $\th=-1\ts$.
See \cite[Proposition 1.7]{NT1} for a more general result.

Now let $(\G\com\Gd)$ be any pair from \eqref{holist}
other than $(\GL_m\com\GL_n)\ts$. The Chevalley
anti-involution $\ep$ on the Lie algebra $\g$ is defined by setting
for $a\com b=-\ts m\lcd-1\com1\lcd m$
$$ 
\ep\ts(F_{ab})=\sgn\ts a\ts b\ts\cdot F_{\ts ba}
\quad\text{or}\quad
\ep\ts(F_{ab})=F_{\ts ba}
$$
when $\th=1$ or $\th=-1$ respectively.
It extends to an involutive
anti-automorphism of the algebra $\Ar$, also denoted by $\ep\ts$, so that 
the extension exchanges the element $p_{\ts ai}\in\Heist$ with the element
$$ 
-\ts\sgn\ts a\cdot q_{\ts ai}
\quad\text{or}\quad
q_{\ts ai}
$$
when $\th=1$ or $\th=-1$ respectively.
Here we use the notation from Proposition \ref{propo2}.

First applying the homomorphism $\X(\gd)\to\Ar$ defined by \eqref{homo2},
and then applying the anti-automorphism $\ep$ of $\Ar\ts$, amounts to
first applying the anti-automorphism \eqref{sijsji} of $\X(\gd)\ts$, and then
applying the homomorphism \eqref{homo2}. Therefore,
when we regard the two vector spaces in the pairing \eqref{qlm}
as $\X(\gd)\ts$-modules, the pairing becomes contravariant relative to 
the involutive anti-automorphism \eqref{sijsji} of $\X(\gd)\ts$.

By replacing these two $\X(\gd)\ts$-modules by their equivalents, we get
a non-degenerate contravariant pairing of the source and target
$\X(\gd)\ts$-modules in \eqref{ak2},
if $\g=\sp_{2m}$ or if $\g=\so_{2m}$ and $m$ is even.
If $\g=\so_{2m}$ and $m$ is odd, then we get a 
%non-degenerate contravariant 
pairing of the source and target
$\X(\gd)\ts$-modules in \eqref{4.56}. But if
$\g=\so_{2m}\ts$, then the target $\X(\gd)\ts$-modules in 
\eqref{ak2} and \eqref{4.56} are equivalent, see
the end of Subsection 4.4. Thus~for~~both 
$\g=\sp_{2m}\com\so_{2m}$ and for any $m$ we get
a non-degenerate contravariant pairing of the source and target
$\X(\gd)\ts$-modules in \eqref{ak2}.
Let us now regard them as $\Y(\gd)\ts$-modules,
like we did in Subsection 4.4. Then we~get
a non-degenerate pairing of them,
which is contravariant relative to the restriction
of the involutive anti-automorphism \eqref{composit}
of $\Y(\gl_n)$ to the subalgebra $\Y(\gd)\subset\Y(\gl_n)\ts$.
We also use the fact that any automorphism
of $\X(\gd)$ defined by \eqref{fus} commutes with
the anti-automorphism defined by \eqref{sijsji}.
In particular, the $\Y(\gd)\ts$-module 
dual to the source in \eqref{ak2}, is equivalent to the target.

Like in the case $\gd=\gl_n\tts$,
the equivalence here can be proved for all $\la$ and $\mu$ directly.
Indeed, the $\Y(\gd)\ts$-module dual to the source in \eqref{ak2} 
can be defined by first considering the source as a 
$\Y(\gl_n)\ts$-module, then taking its dual relative
to the anti-automorphism \eqref{composit} of $\Y(\gl_n)\ts$,
and then restricting the resulting $\Y(\gl_n)\ts$-module to the 
subalgebra $\Y(\gd)\subset\Y(\gl_n)\ts$. But \eqref{composit}
is the composition of \eqref{tijtji} and \eqref{transauto}.
Hence taking the dual relative to \eqref{composit}
amounts to first taking the dual relative to \eqref{tijtji},
and then pulling the result back through the
automorphism \eqref{transauto}. 

We have already proved that relative to \eqref{tijtji},
the $\Y(\gl_n)\ts$-module dual to the source tensor product in \eqref{ak2}
is equivalent to the tensor product of the same factors,
but taken in the reversed order. 
Pulling the latter tensor product back through \eqref{transauto}
replaces each factor $\Ph^{\ts k}_{\tts t}$ by $\Ph^{\ts-k}_{\tts t}\ts$,
and also reverses the order of the factors once again.
Thus we get the target in \eqref{ak2} as $\Y(\gl_n)\ts$-module.
Here we used the definition of the $\Y(\gl_n)\ts$-module
$\Ph^{\ts-k}_{\tts t}\ts$, and the fact that \eqref{transauto}
defines an anti-automorphism of the coalgebra $\Y(\gl_n)\ts$.

%==============================================================================

\section*{\bf\normalsize 5. Irreducible representations of Yangians}
\setcounter{section}{5}
\setcounter{equation}{0}
\setcounter{theorem*}{0}

%------------------------------------------------------------------------------

\smallskip\noindent
{\it 5.1.\ Irreducible representations of\/ $\Y(\gl_n)$}

\vspace{10pt}\noindent
Let $\Ph$ be a non-zero finite-dimensional $\Y(\gl_n)\ts$-module. 
A non-zero vector of $\ts\Ph$ is called {\it highest} if it is annihilated
by all the coefficients of the series $T_{ij}(x)$ with $i<j\ts$.
If $\Ph$ is irreducible then a highest vector
$\ph\in\Ph$ is unique up to a scalar multiplier.
Moreover, then $\ph$ is an eigenvector for the coefficients of all
series $T_{ii}(x)\ts$, and for $i=1\lcd n-1$
$$
\textstyle
T_{ii}(x)\,T_{i+1,i+1}(x)^{\tts-1}\,\ph
=
P_i(x+\frac12)\,P_i(x-\frac12)^{-1}\,\ph
$$
where $P_i(x)$ is a monic polynomial in $x$ with coefficients in $\CC\ts$.
Then $P_1(x)\lcd P_{n-1}(x)$ are called the 
\textit{Drinfeld polynomials\/} of $\Ph\ts$.
Any sequence of $n-1$ monic polynomials with complex coefficients 
arises in this way.
Furthermore, two irreducible finite-dimensional $\Y(\gl_n)\ts$-modules 
have the same Drinfeld polynomials if and only if
their restrictions to the subalgebra %special Yangian
$\SY(\gl_n)\subset\Y(\gl_n)$ are equivalent. 
Thus up to equivalence and similarity,
all the non-zero irreducible finite-dimensional $\Y(\gl_n)\ts$-modules
are parametrized by their Drinfeld polynomials 
\cite[Theorem 2]{D2}. 
For example, consider the \textit{trivial\/} $\Y(\gl_n)\ts$-module. 
It is one-dimensional, and is defined by the counit
homomorphism $\Y(\gl_n)\to\CC\ts$. Then the corresponding Drinfeld 
polynomials are also trivial: $P_1(x)=\ldots=P_{n-1}(x)=1\ts$.

In this subsection we will assume that
$\th=-1\ts$ and $\g=\gl_m\ts$, so that $\gd=\gl_n\ts$.
For any $k\in\{0\com1\lcd n\ts\}$ the vector space
$\Ph^{\ts k}$ is irreducible under the action of $\gl_n\ts$.
Hence for any $t\in\CC$ 
the $\Y(\gl_n)\ts$-module $\Ph^{\ts k}_{\tts t}$ 
is irreducible. The vector 
$\ph_k=f_1\wedge\ldots\wedge f_k$ of this module
is highest, see the definitions \eqref{tauz} and \eqref{eval}.
Moreover, by these definitions
$$
T_{ii}(x)\,\ph_k=
\left\{
\begin{array}{ll}
(x-t+1)\,(x-t)^{-1}\,\ph_k
&\quad\textrm{if}\quad\ 1\le i\le k\,;
\\[2pt]
\hspace{44pt}
\ph_k
&\quad\textrm{if}\quad\ k<i\le n\ts.
\end{array}
\right. 
$$
Hence for $i=1\lcd n-1$ the Drinfeld polynomial $P_i(x)$
of the $\Y(\gl_n)\ts$-module $\Ph^{\ts k}_{\tts t}$ is 
$$
P_{i}(x)=
\left\{
\begin{array}{ll}
x-t+\frac12
&\quad\textrm{if}\quad\ i=k\,;
\\[2pt]
\hspace{17pt}
1
&\quad\textrm{otherwise}.
\end{array}
\right. 
$$

For $\g=\gl_m$ we have
%for each index $a=1\lcd m$ we have $\rho_a=m/2-a+\frac12\,$. 
$\nu=\la-\mu\ts$.
Suppose that each label $\nu_a\in\{0\com1\lcd n\}\ts$.
Further suppose that the weight $\la+\rho$ is nonsingular.
The latter condition means here that 
\begin{equation}
\label{domcon1}
\la_b-\la_a+\rho_b-\rho_a\,\neq\,1,2,\,\ldots
\quad\text{for all}\quad
1\le a<b\le m\ts.
\end{equation}
By Corollary~\ref{ircor} and Proposition~\ref{propo1}, 
the quotient by the kernel %(or equivalently the image) 
of our %intertwining 
operator \eqref{ak1} is then an irreducible $\Y(\gl_n)\ts$-module.
But the definition \eqref{1.33} of the comultiplication on $\Y(\gl_n)$
implies that the vector 
$\ph_{\ts\nu_1}\ot\ldots\ot\ph_{\ts\nu_m}$
of \eqref{baspro1} is highest. Suppose that
\begin{equation}
\label{chercon}
\nu_a\ge\nu_b
\quad\text{whenever}\quad
\la_a+\rho_a=\la_b+\rho_b
\quad\text{and}\quad
a<b\ts.
\end{equation}
Due to the nonsingularity of the weight $\la+\rho\ts$, then
Proposition \ref{norm1} implies that the image
of the vector $\ph_{\ts\nu_1}\ot\ldots\ot\ph_{\ts\nu_m}$ 
in the quotient is not zero. Hence this image is highest
relative to the action of the Yangian $\Y(\gl_n)$
on the quotient. Since the vectors
$\ph_{\ts\nu_1}\lcd\ph_{\ts\nu_m}$
are highest in their $\Y(\gl_n)\ts$-modules, their tensor product
is an eigenvector for all coefficients
of the series $T_{ii}(x)\ts$, %for each $i=1\lcd n\ts$;
see again %the definition 
\eqref{1.33}. Moreover,
$$
T_{ii}(x)\,(\ts\ph_{\ts\nu_1}\ot\ldots\ot\ph_{\ts\nu_m})=
(\ts T_{ii}(x)\,\ph_{\ts\nu_1})\ot\ldots\ot(\ts T_{ii}(x)\,\ph_{\ts\nu_m})\,.
$$
Therefore for any index $i=1\lcd n-1$
the Drinfeld polynomial $P_i(x)$ of the
quotient is equal to the product of the Drinfeld
polynomials with the same index $i$
of the $m$ tensor factors of \eqref{baspro1}.
Thus we get the following theorem.
Recall that here $\rho_a=m/2-a+\frac12\,$.

\begin{theorem*}
\label{drin1}
Let the labels $\la_1\lcd\la_m$ satisfy the condition
\eqref{domcon1}, while the labels 
$$
\nu_1=\la_1-\mu_1
\ts\,,\,\ldots\,,\,
\nu_m=\la_m-\mu_m
$$
belong to the set\/
$\{0\com1\lcd n\}$ and satisfy the condition \eqref{chercon}.
Then the quotient by the kernel of our intertwining operator
\eqref{ak1} is a non-zero irreducible $\Y(\gl_n)\ts$-module.
For any\/ $i=1\lcd n-1$ the Drinfeld polynomial $P_i(x)$ of this module 
is the product of the differences\/ $x-\mu_a-\rho_a$ taken over 
all indices\/ $a$ such that $\nu_a=i\ts$.
\end{theorem*}

Note that if $\nu_1=\ldots=\nu_m=0\ts$, then
both the source and the target
$\Y(\gl_n)\ts$-modules in \eqref{ak1} are trivial
for any $\mu\ts$. If moreover
$\la+\rho$ is nonsingular, then by Proposition \ref{norm1} 
our operator \eqref{ak1} is the identity map $\CC\to\CC\ts$.

Now take any sequence of $n-1$ monic polynomials 
$P_1(x)\lcd P_{n-1}(x)$ with complex coefficients. 
Denote by $\mathcal{P}$ the
collection of pairs $(\ts i\com z)$
where $i=1\lcd n-1$ and $z$ ranges over all roots
of the polynomial $P_i(x)\ts$; 
the roots are taken with their multiplicities.
Note that the collection $\mathcal{P}$ is unordered.
Let $m$ be the total number of elements in $\mathcal{P}\ts$, it is equal
to the sum of the degrees of %the polynomials 
$P_1(x)\lcd P_{n-1}(x)\ts$.
Suppose that at least one of the polynomials is not trivial,
so that $m>0\ts$.
Let $\la$ and $\mu$ be any weights of $\gl_m$ such that for $\nu=\la-\mu$
the collection of pairs 
$(\ts\nu_a\com \mu_a+\rho_a\ts)$
for $a=1\lcd m$ coincides with $\mathcal{P}\ts$.
In particular $\nu_a\in\{1\lcd n-1\}\ts$.
We do not yet impose any other conditions on $\la$ and $\mu\ts$.

Then $\la$ and $\mu$ are determined up to any permutation 
of the $m$ pairs $(\ts\nu_a\com \mu_a+\rho_a\ts)\ts$.
Equivalenly, they are determined up to a
permutation of the $m$ pairs $(\ts\la_a+\rho_a\com\mu_a+\rho_a\ts)\ts$.
In other words, $\la$ and $\mu$
are determined up to
the (simultaneous) shifted action of the Weyl group $\Sym$
of $\gl_m$ on them, as on elements of $\h^*$. 
Therefore we can choose $\la$ to satisfy the conditions \eqref{domcon1},
so that the weight $\la+\rho$ is nonsingular. For $a<b$
the equality on the right hand
side of \eqref{chercon} means that %the weight 
$\la$ is invariant
under the shifted action of that element of $\Sym$ which
exchanges $\eta_a$ with $\eta_b\ts$, and leaves all other
basis vectors of $\h^*$ fixed. 
This action amounts to exchanging the pair
$(\ts\nu_a\com \mu_a+\rho_a\ts)$ with 
$(\ts\nu_b\com \mu_b+\rho_b\ts)\ts$.
By using this action and keeping $\la$ fixed, 
we can choose $\mu$ so that the condition \eqref{chercon} is satisfied.

Now consider our %intertwining 
operator \eqref{ak1}
corresponding to the weights $\la$ and $\mu$ chosen above.
Due to Theorem \ref{drin1}, the quotient by the kernel of this 
operator is a non-zero irreducible $\Y(\gl_n)\ts$-module,
and has the given Drinfeld polynomials $P_1(x)\lcd P_{n-1}(x)\ts$.
Thus up to equivalence and similarity, 
every non-zero irreducible finite-dimensional 
$\Y(\gl_n)\ts$-module arises as such a quotient.
Note that the choice $\th=-1$ here is essential.

%------------------------------------------------------------------------------

\vspace{-6pt}%%%%%%%%%%%%%%%%%%%%%%%%%%%%%%%%%%%%%%%%%%%%%%%%%%%%%%%%%%%%%%%%%%

\subsection*{\it\normalsize 5.2.\ Representations of twisted Yangians}

\noindent
From now on we will regard the coefficients of all the series $S_{ij}(x)$
as generators of the algebra $\Y(\gd)\ts$, with the relations
\eqref{xrel} and \eqref{srel} imposed on them.
In this subsection we collect some general results on
finite-dimensional $\Y(\gd)\ts$-modules for both $\gd=\sp_n$ and
$\gd=\so_n\ts$. Let $\ts\Psi$ be any of these modules.
If $\Psi$ is obtained by 
restricting the trivial %one-diminsional
$\Y(\gl_n)\ts$-module to the subalgebra
$\Y(\gd)\subset\Y(\gl_n)\ts$, 
then $\Psi$ will be also called \textit{trivial}. 

Now assume that $\Psi$ is non-zero.
We will use the ordering \eqref{corder}
of the indices $1\lcd n\ts$. We will keep writing $i\prec j$
when $i$ precedes $j\ts$ in the sequence \eqref{corder}. 
A non-zero vector of $\ts\Psi$ is called {\it highest} if it is annihilated
by all the coefficients of the series $S_{ij}(x)$ with $i\prec j\ts$.
If $\ts\Psi$ is irreducible then a highest vector
$\psi\in\Psi$ is unique up to a scalar multiplier.
Then $\psi$ is an eigenvector for the coefficients of all
series $S_{ii}(x)\ts$; see \cite[Theorem~4.2.6 and Corollary 4.2.7]{M3}.
Then $\ts\Psi$ is determined by
the corresponding eigenvalues up to equivalence. 
In the next %three 
subsections we will use the description
from \cite[Chapter~4]{M3} of all possible eigenvalues.
Note that the relation \eqref{srel} with $i=j$ takes the form
$$
S_{ii}(x)\mp2\ts x\,S_{\ts\bi\ts\bi\ts}(x)
=(1\mp2\ts x)\,S_{ii}(-x)\ts.
$$
Hence it suffices to describe the eigenvalues of the coefficients
of only one of every two series $S_{ii}(x)$ and 
$S_{\ts\bi\ts\bi\ts}(x)\ts$. We choose the series 
$S_{ii}(x)$ with $i\succcurlyeq\bi\ts$,
or equivalently with $i\succcurlyeq n\ts$.
This choice is explained by the next lemma which 
valid for any non-zero $\Y(\gl_n)\ts$-module $\Ph\ts$,
not necessarily irreducible\ts; cf. \cite[Corollary 4.2.10]{M3}.
We will also regard $\Ph$ as module over the subalgebra 
$\Y(\gd)\subset\Y(\gl_n)\ts$, by restriction.

\begin{lemma*}
\label{higi}
Let $\ph\in\Ph$ be a non-zero vector annihilated 
by all coefficients of the series $T_{ij}(x)$ with $i\prec j\ts$.
Suppose\/
$
T_{ii}(x)\,\ph=h_i(x)\,\ph
$
for\/ $i=1\lcd n$
where $h_i(x)$ is a formal power series in $x^{-1}$ with 
coefficients from\/ $\CC\ts$. %and the leading term $1\ts$.
Then $\ph$ is a highest vector for\/ $\Y(\gd)$~and moreover\/
$
S_{ii}(x)\,\ph=g_i(x)\,\ph
$
for\/ $i=1\lcd n$ where
$$
g_i(x)=
\left\{
\begin{array}{ll}
(\ts1\mp(2\ts x)^{-1})\ts\,h_i(x)\,h_{\ts\bi}(-x)
\pm
(2\ts x)^{-1}\,h_i(-x)\,h_{\ts\bi}(x)
&\quad\textrm{if}\quad\ i\prec\bi\ts;
\\[2pt]
\hspace{80pt}
h_i(x)\,h_{\ts\bi}(-x)
&\quad\textrm{if}\quad\ i\succcurlyeq\bi\,.
\end{array}
\right. 
$$
\end{lemma*}

\begin{proof}
Take %any indices 
$i\com j\in\{1\lcd n\}$ such that
$i\preccurlyeq j\ts$. The coefficients of the series $S_{ij}(x)$ 
act on the vector $\ph\in\Ph$ 
as the corresponding coefficients of the series 
\eqref{yser}, giving the sum of
\begin{equation}
\label{suma}
\th_i\ts\th_k\,T_{\,\bk\ts\bi\ts}(-x)\,T_{kj}(x)\,\ph
\end{equation}
over the indices $k\succcurlyeq j\ts$. This is because 
the vector $\ph$ is highest relative to $\Y(\gl_n)\ts$. Using the commutation 
relations \eqref{yrel}, the summand \eqref{suma} is equal to the difference 
\begin{gather}
\label{sumb}
\th_i\ts\th_k\,T_{kj}(x)\,T_{\,\bk\ts\bi\ts}(-x)\,\ph
\\[2pt]
\label{sumc}
-\,\,\th_i\ts\th_k\,(2\ts x)^{-1}
(\ts T_{k\ts\bi\ts}(-x)\,T_{\,\bk j}(x)-
T_{k\ts\bi\ts}(-x)\,T_{\,\bk j}(x))\,\ph\ts.
\end{gather}
Here \eqref{sumb} may be non-zero only if $\ts\bk\succcurlyeq\bi\ts$,
that is only if $k\preccurlyeq i\ts$. The latter condition together with
$i\preccurlyeq j$ and $j\preccurlyeq k$ implies that $i=j=k\ts$. Then
\eqref{sumb} equals 
\begin{equation}
\label{sumf}
h_i(x)\,h_{\ts\bi}(-x)\,\ph\ts.
\end{equation}

Further, \eqref{sumc} may be non-zero only if 
$\ts\bk\succcurlyeq j\ts$. When presenting \eqref{suma}
as the sum of \eqref{sumb} and of \eqref{sumc} by using \eqref{yrel}, 
we could replace \eqref{sumc} by the expression 
\begin{equation}
\label{sumd}
-\,\,\th_i\ts\th_k\,(2\ts x)^{-1}
(\ts T_{\,\bk j}(x)\,T_{k\ts\bi\ts}(-x)-
T_{\,\bk j}(x)\,T_{k\ts\bi\ts}(-x))\,\ph\ts.
\end{equation}
But \eqref{sumd} may be non-zero only if 
$k\succcurlyeq\bi\ts$, that is only if $\bk\preccurlyeq i\ts$.
The latter condition together with $i\preccurlyeq j$ and 
$j\preccurlyeq\bk\ts$ implies that $i=j=\bk\ts$.
Then \eqref{sumc} and \eqref{sumd} are equal to
\begin{equation}
\label{sume}
\mp\ts(2\ts x)^{-1}\,
(\ts h_i(x)\,h_{\ts\bi}(-x)-h_i(-x)\,h_{\ts\bi}(x))\,\ph\ts.
\end{equation}
We have also assumed that $k\succcurlyeq j\ts$, which 
together with the equalities $i=j=\bk$ implies that $i\preccurlyeq\bi\ts$.
Moreover, if $i=\bi$ then \eqref{sume} obviously vanishes.

Thus under the assumption $i\preccurlyeq j$ 
we have proved that $S_{ij}(x)\,\ph$ may difer from zero
only of $i=j\ts$. If $i\prec\bi$ then $S_{ii}(x)\,\ph$
is equal to the sum of \eqref{sumf} and \eqref {sume}. 
But if $i\succcurlyeq\bi$ then $S_{ii}(x)\,\ph$ is equal
to the expression \eqref{sumf} alone.
\qed
\end{proof}

We will keep assuming that $\th=-1\ts$. Hence
for $\g=\sp_{2m}$ or $\g=\so_{2m}$ we will have
$\gd=\sp_n$ or $\gd=\so_n$ respectively. 
In these next subsections, we will consider 
the two cases separately. Moreover, we will
separate the cases of $\gd=\so_n$ with $n$ even and odd. 
If $\gd=\sp_n$ then $n$ has to be even. 
We will write $n=2\tts l$ if
$n$ is even, or $n=2l+1$ if $n$ is odd.
Then the condition $i\succcurlyeq n$ will mean respectively that $i$
is one of the indices
$$
2\tts l\com2\tts l-2\lcd 2
\ \quad\text{or}\ \quad
2\tts l+1\com
2\tts l\com2\tts l-2\lcd 2\ts.
$$

As the first application of Lemma \ref{higi},
consider the case when $\Ph=\Ph^{\ts k}_{\tts t}$ and 
$\ph=\psi_k$ is the vector of $\Ph^{\ts k}$ defined in 
Subsection 4.4. We assume that $k\in\{0\com1\lcd n\ts\}\ts$. 
All the conditions of Lemma \ref{higi}
are then satisfied. Here
$h_i(x)=1+(x-t)^{-1}$ if $i$ is one of the first $k$
indices in the sequence \eqref{corder}, otherwise $h_i(x)=1\ts$.
Hence $h_{\ts\bi}(-x)=1-(x+t)^{-1}$ if $i$ is one of 
the last $k$ indices in the sequence \eqref{corder},
otherwise $h_{\ts\bi}(-x)=1\ts$. If
$i\succcurlyeq n\ts$, then
$g_i(x)=h_i(x)\,h_{\ts\bi}(-x)\ts$.
This implies the following relations,
to be used later on. 

Take any $i\succ n\ts$. Let $j$ be
the index occuring just before $i$ in the sequence
\eqref{corder}. Here $j=i+2$ unless $n=2\tts l+1$ and $i=2\tts l\ts$,
in which case $j=i+1\ts$.
Note that $j\succcurlyeq n\ts$. Then %we have the relation
$$
g_j(x)\,g_i(x)^{-1}=
\left\{
\begin{array}{ll}
(x+t)\,(x+t-1)^{-1}
&\quad\textrm{if}\quad\ 2\tts k=i\,;
\\[2pt]
(x-t+1)\,(x-t)^{-1}
&\quad\textrm{if}\quad\ 2\tts (n-k)=i\,;
\\[2pt]
\hspace{40pt}
1
&\quad\textrm{otherwise}.
\end{array}
\right. 
$$
Indeed, for $i\succ n$ we have $h_i(x)\neq h_j(x)$
only if $2(n-k)=i\ts$. In the latter case we have
$h_i(x)=1$ while $h_j(x)=1+(x-t)^{-1}\ts$.
Further, for any $i\succ n$ we have $h_{\ts\bi}(x)\neq h_{\ts\bj}(x)$
only if $2\tts k=i\ts$.
In the latter case we have $h_{\ts\bi}(x)=1+(x-t)^{-1}\ts$
while $h_{\ts\bj}(x)=1\ts$. Thus we obtain the relation displayed above.

If $n=2\tts l\ts$, then we also have the relation
$$
g_n(-x)\,g_n(x)^{-1}=
\left\{
\begin{array}{ll}
(x-t+1)\,(x+t)\,(x-t)^{-1}\,(x+t-1)^{-1}
&\quad\textrm{if}\quad\ k=l\,;
\\[2pt]
\hspace{88pt}
1
&\quad\textrm{otherwise}.
\end{array}
\right. 
$$
Indeed, if $n=2\tts l$ then $\tilde n=2\tts l-1\ts$. 
Then for $k\le l$ we have $h_n(x)=1\ts$, while for $k>l$ we have
$h_n(x)=1+(x-t)^{-1}$. Further, for $k<l$ we have 
$h_{\ts\tilde n}(x)=1\ts$, while for $k\ge l$ we have
$h_{\ts\tilde n}(x)=1+(x-t)^{-1}\ts$.
Therefore $g_n(x)=g_n(-x)$ for any $k\neq l\ts$.
But if $k=l$ then $g_n(x)=1+(x-t)^{-1}\ts$.
Thus we obtain the last displayed relation.

We will end this section by introducing a notion
applicable to any finite-dimensional module $\Psi$ over the twisted
Yangian $\Y(\gd)$ where $\gd=\so_n$ and $n$ is any positive integer,
even or odd. Consider the restriction of $\Psi$ to the subalgebra
$\U(\so_n)\subset\Y(\so_n)\ts$. Suppose that the restriction
integrates to a module of the complex special orthogonal group 
$\mathrm{SO}_n\ts$. Thus we exclude the spinor representations of~$\so_n\ts$.  
Then we call the $\Y(\so_n)\ts$-module $\Psi$
\textit{integrable\ts}. 
Equivalently, then $\Psi$ is 
an $(\ts\Y(\so_n)\com\mathrm{SO}_n)\ts$-module\tts; see the beginning
of Subsection 4.4.
For instance, any $\Y(\so_n)\ts$-module
of the form \eqref{pro2} is integrable, and so is any quotient of 
such a module. 
We will not treat non-integrable $\Y(\so_n)\ts$-modules~here.

%------------------------------------------------------------------------------

\bigskip\noindent
{\it 5.3.\ Irreducible representations of\/ $\Y(\sp_n)$}

\vspace{10pt}\noindent
In this subsection, we consider the case
when $\th=-1$ and $\g=\sp_{2m}\ts$, so that $\gd=\sp_n$ 
where $n=2\ts l\ts$. Let $\Psi$ be a
non-zero irreducible finite-dimensional $\Y(\sp_n)\ts$-module. 
Let $\psi\in\Psi$ be a highest vector.
By \cite[Theorem 4.3.8]{M3} for any $k=1\lcd l-1$ we have
$$
\textstyle 
S_{2k+2,2k+2}(x)\,S_{2k,2k}(x)^{\tts-1}\,\psi
\,=\,
Q_k(x+\frac12)\,Q_k(x-\frac12)^{-1}\,\psi
$$
where $Q_k(x)$ is a monic polynomial in $x$ with coefficients in 
$\CC\ts$. Further, we have
$$
\hspace{26pt}
\textstyle 
S_{nn}(-\ts x)\,S_{nn}(x)^{\tts-1}\,\psi
\,=\,
Q_l(x+\frac12)\,Q_l(x-\frac12)^{-1}\,\psi
$$
where $Q_l(x)$ is an even monic polynomial in $x$ with coefficients in
$\CC\ts$.
Any sequence of $l$ monic polynomials with complex coefficients 
arises in this way, provided that the last polynomial in the sequence is even.
Furthermore, two irreducible finite-dimensional $\Y(\sp_n)\ts$-modules 
have the same sequence of polynomials $Q_1(x)\lcd Q_l(x)$ 
if and only if their restrictions to the subalgebra
$\SY(\sp_n)\subset\Y(\sp_n)$ are equivalent.
Thus the non-zero irreducible finite-dimensional $\Y(\sp_n)\ts$-modules
are parametrized by their polynomials $Q_1(x)\lcd Q_l(x)$ up to
equivalence and similarity \cite[Corollary 4.3.11]{M3}. 
For example, if $\Psi$ is the trivial $\Y(\sp_n)\ts$-module
then $Q_1(x)=\ldots=Q_l(x)=1\ts$.

In this subsection, $\rho_a=-\ts a$ and $\ka_a=-\ts l$ 
for each index $a=1\lcd m\ts$. Recall the definition
\eqref{lamuka} of the weight $\nu\ts$. 
Suppose that each label $\nu_a\in\{0\com1\lcd n\}\ts$.
Further suppose that the weight $\la+\rho$ is nonsingular.
The latter condition means here that 
\begin{align}
\label{domcon31}
\la_b-\la_a+\rho_b-\rho_a\,\neq\,1,2,\,\ldots
&\quad\text{for all}\quad
1\le a<b\le m\ts;
\\
\label{domcon32}
\la_a+\la_b+\rho_a+\rho_b\,\neq\,1,2,\,\ldots
&\quad\text{for all}\quad
1\le a<b\le m\ts;
\\
\label{domcon33}
\la_a+\rho_a\,\neq\,1,2,\,\ldots
&\quad\text{for all}\quad
1\le a\le m\ts.
\end{align}
By Corollary~\ref{ircor} and Proposition~\ref{propo2}, 
the quotient by the kernel %(or equivalently the image) 
of our intertwining operator
\eqref{ak2} is then an irreducible $\Y(\sp_n)\ts$-module.
But the definition \eqref{1.33} of the comultiplication on $\Y(\gl_n)$
implies that the vector \eqref{ourvec2}
is highest relative to the twisted Yangian $\Y(\sp_n)\ts$. 
Here we also use Lemma \ref{higi}.
Further suppose that 
\begin{align}
\label{chercon31}
\nu_a\ge\nu_b
&\quad\text{whenever}\quad
\la_a-\la_b+\rho_a-\rho_b=0
\quad\text{and}\quad
a<b\ts;
\\
\label{chercon32}
\nu_a+\nu_b\le n\ 
&\quad\text{whenever}\quad
\la_a+\la_b+\rho_a+\rho_b=0
\quad\text{and}\quad
a<b\ts;
\\
\label{chercon33}
\nu_a\le l\ \ 
&\quad\text{whenever}\quad
\la_a+\rho_a=0\ts.
\end{align}

Due to the nonsingularity of the weight $\la+\rho\ts$,
Proposition \ref{norm2} then implies that the image
of the vector \eqref{ourvec2} in the quotient
is not zero. Hence this image is highest
relative to the action of the twisted Yangian $\Y(\sp_n)$
on the quotient. Lemma \ref{higi} implies that each of
the polynomials $Q_1(x)\lcd Q_l(x)$ of the quotient
is multiplicative with respect to the $m$ tensor factors
of the vector \eqref{ourvec2}. But in the
case $m=1$ these polynomials are transparent
from the relations given in the very end of Subsection 5.2.
Thus we get

\begin{theorem*}
\label{drin3}
Put $n=2l$.
Let $\la_1\lcd\la_m$ satisfy %the conditions
\eqref{domcon31},\eqref{domcon32},\eqref{domcon33} while the labels
\begin{equation}
\label{nu3}
\nu_1=\la_1-\mu_1+l
\ts\,,\,\ldots\,,\,
\nu_m=\la_m-\mu_m+l
\end{equation}
belong to the set\/
$\{0\com1\lcd n\}$ and satisfy
\eqref{chercon31},\eqref{chercon32},\eqref{chercon33}.
Then the quotient by the kernel of our intertwining operator
\eqref{ak2} is a non-zero irreducible $\Y(\sp_n)\ts$-module.
For any\/ $k=1\lcd l-1$ the polynomial $Q_k(x)$ of this module 
is the product of the sums $x+\mu_a+\rho_a$ taken over 
all indices\/ $a$ such that $\nu_a=k\ts$, and
of the differences\/ $x-\mu_a-\rho_a$ taken over 
all indices\/ $a$ such that $\nu_a=n-k\ts$.
The polynomial $Q_l(x)$ of this module 
is the product of the differences\/ 
$x^{\ts2}-(\mu_a+\rho_a){}^{\tts2}$ taken over 
all indices\/ $a$ such that $\nu_a=l\ts$.
\end{theorem*}

Note that if $\nu_1=\ldots=\nu_m=0\ts$, then
both the source and the target $\Y(\sp_n)\ts$-modules
in \eqref{ak2} are trivial
for any $\mu\ts$. If moreover
$\la+\rho$ is nonsingular, then by Proposition \ref{norm2} 
our operator \eqref{ak2} is the identity map $\CC\to\CC\ts$.

Now let us take any sequence of $l$ monic polynomials 
$Q_1(x)\lcd Q_l(x)$ with complex coefficients,
such that the polynomial $Q_l(x)$ is even.
Let $m$ be the sum of the degrees of $Q_1(x)\lcd Q_{l-1}(x)$
plus half of the degree of $Q_l(x)\ts$.
Suppose that at least one of all the $l$ polynomials is not trivial,
so that $m>0\ts$. Let $\la$ and $\mu$ be any weights of $\sp_{2m}$
such that the corresponding labels \eqref{nu3}
belong to the set $\{1\lcd n-1\}$ and such that
the given polynomials $Q_1(x)\lcd Q_l(x)$
are obtained from $\la$ and $\mu$ as in Theorem \ref{drin3}.
But we do not yet impose any other conditions on $\la$ and $\mu\ts$,
such as nonsingularity of $\la+\rho\ts$.

Then $\la$ and $\mu$ are determined up to permuting 
the $m$ pairs $(\ts\nu_a\com \mu_a+\rho_a\ts)\ts$,
and up to replacing $(\ts\nu_a\com \mu_a+\rho_a\ts)$
by $(\ts n-\nu_a\com-\mu_a-\rho_a\ts)$ for any number
of indices $a\ts$.
Equivalenly, $\la$ and $\mu$ are determined up to a
permutation of the $m$ pairs $(\ts\la_a+\rho_a\com\mu_a+\rho_a\ts)\ts$,
and up to replacing $(\ts\la_a+\rho_a\com \mu_a+\rho_a\ts)$
by $(-\la_a-\rho_a\com-\mu_a-\rho_a\ts)$ for any number
of indices $a\ts$.
In other words, $\la$ and $\mu$
are determined up to
the (simultaneous) shifted action of the Weyl group 
$\Sym$ of $\sp_{2m}$ on them,
as on elements of $\h^*$. 
Therefore we can choose $\la$ to satisfy the conditions 
\eqref{domcon31},\eqref{domcon32},\eqref{domcon33}
so that the weight $\la+\rho$ is nonsingular. 

The equality on the right hand
side of \eqref{chercon31} for $a<b$
means that $\la$ is invariant
under the shifted action of that element of $\Sym$
which exchanges $\eta_{m-a+1}$ with $\eta_{m-b+1}\ts$, and
leaves all other basis vectors of $\h^*$ fixed.
This action amounts to exchanging the pair
$(\ts\nu_a\com \mu_a+\rho_a\ts)$ with
$(\ts\nu_b\com \mu_b+\rho_b\ts)\ts$.
By using this action and keeping $\la$ fixed, 
we can choose $\mu$ so that the condition \eqref{chercon31} is satisfied.

The equality on the right hand
side of \eqref{chercon33} means that the weight $\la$ is invariant
under the shifted action of that element of $\Sym$
which maps $\eta_{m-a+1}$ to $-\ts\eta_{m-a+1}\ts$,
and leaves all other basis vectors of $\h^*$ fixed.
This action amounts to replacing the pair
$(\ts\nu_a\com \mu_a+\rho_a\ts)$ by the pair 
$(\ts n-\nu_a\com-\mu_a-\rho_a\ts)\ts$.
By using this action and keeping $\la$ fixed, 
we can choose $\mu$ so that the condition \eqref{chercon33} is satisfied.

Finally, the equality on the right hand
side of \eqref{chercon32} for $a<b$ means that the weight 
$\la$ is invariant under the shifted action of that element of 
the group $\Sym$
which maps $\eta_{m-a+1}$ and $\eta_{m-b+1}$
respectively to $-\ts\eta_{m-b+1}$ and $-\ts\eta_{m-a+1}\ts$,
leaving all other basis vectors of $\h^*$ fixed.
This action amounts to replacing the pairs
$(\ts\nu_a\com \mu_a+\rho_a\ts)$ and 
$(\ts\nu_b\com \mu_b+\rho_b\ts)$ by the pairs
$(\ts n-\nu_b\com-\mu_b-\rho_b\ts)$ and 
$(\ts n-\nu_a\com-\mu_a-\rho_a\ts)$ respectively.
By using this action and keeping $\la$ fixed, 
we can choose the weight $\mu$ 
so that the condition \eqref{chercon32} is satisfied
when $\la_a+\rho_a\neq0\ts$, or equivalently when
$\la_b+\rho_b\neq0\ts$. When $\la_a+\rho_a=0$ and
$\la_b+\rho_b=0\ts$, the condition
\eqref{chercon32} is already satisfied, because then
$\nu_a\le l$ and $\eta_b\le l$ due to \eqref{chercon33}.

Now consider our %intertwining 
operator \eqref{ak2}
corresponding to $\la$ and $\mu\ts$.
Due to Theorem~\ref{drin3} the quotient by the kernel of this 
operator is a non-zero irreducible $\Y(\sp_n)\ts$-module,
and has the given polynomials $Q_1(x)\lcd Q_l(x)\ts$.
Thus up to equivalence and similarity, 
every non-zero irreducible finite-dimensional 
$\Y(\sp_n)\ts$-module arises as such a quotient.
%Note that the choice $\th=-1$ here is essential.

%------------------------------------------------------------------------------

\bigskip\noindent
{\it 5.4.\ Irreducible representations of\/ $\Y(\so_n)$ for odd $n$}

\vspace{10pt}\noindent
In this subsection, we consider the case of $\gd=\so_n$ 
where $n$ is odd. Hence $n=2\tts l+1$ where $l$ is a non-negative integer.
We assume that $\th=-1\ts$, so that $\g=\so_{2m}\ts$.
Let $\Psi$ be any integrable finite-dimensional $\Y(\so_n)\ts$-module,
see the end of Subsection 5.2.
Further suppose that the $\Y(\so_n)\ts$-module $\Psi$ is irreducible
and non-zero. Let $\psi\in\Psi$ be a highest vector.
By \cite[Theorem 4.5.9]{M3} there exist
monic polynomials $Q_1(x)\lcd Q_l(x)$
in $x$ with complex coefficients such that
for any $k=1\lcd l-1$ %we have
$$
\textstyle 
S_{2k+2,2k+2}(x)\,S_{2k,2k}(x)^{\tts-1}\,\psi
\,=\,
Q_k(x+\frac12)\,Q_k(x-\frac12)^{-1}\,\psi\ts,
$$
while
$$
\hspace{13pt}
\textstyle 
S_{nn}(x)\,S_{n-1,n-1}(x)^{\tts-1}\,\psi
\,=\,
Q_l(x+\frac12)\,Q_l(x-\frac12)^{-1}\,\psi\ts.
$$
Every sequence of $l$ monic polynomials with coefficients from $\CC$ 
arises in this way.
Two irreducible integrable finite-dimensional $\Y(\so_n)\ts$-modules 
have the same polynomials $Q_1(x)\lcd Q_l(x)$ 
if and only if their restrictions to the subalgebra
$\SY(\so_n)\subset\Y(\so_n)$ are equivalent \cite[Corollary~4.5.12]{M3}.
Thus all the non-zero irreducible integrable 
finite-dimensional $\Y(\so_n)\ts$-modules
are parametrized by their polynomials $Q_1(x)\lcd Q_l(x)$ 
up to equivalence and {similarity}.
For example, if $\Psi$ is the trivial $\Y(\so_n)\ts$-module
then we have $Q_1(x)=\ldots=Q_l(x)=1\ts$.

In this subsection,
$\rho_a=1-a$ and $\ka_a=-\,l+\frac12\ts$ 
for each index $a=1\lcd m\ts$. Recall the definition
\eqref{lamuka} of the weight $\nu\ts$. 
Suppose that each label $\nu_a\in\{0\com1\lcd n\}\ts$.
Further suppose that the weight $\la+\rho$ is nonsingular.
The latter condition means here that 
\begin{align}
\label{domcon41}
\la_b-\la_a+\rho_b-\rho_a\,\neq\,1,2,\,\ldots
&\quad\text{for all}\quad
1\le a<b\le m\ts;
\\
\label{domcon42}
\la_a+\la_b+\rho_a+\rho_b\,\neq\,1,2,\,\ldots
&\quad\text{for all}\quad
1\le a<b\le m\ts.
\end{align}
In the end of Subsection 4.4
we observed that for $\g=\so_{\tts2m}$
our intertwining operators \eqref{ak2} and \eqref{4.56}
are the same. Using Proposition \ref{propo2} and Corollary~\ref{iror},
the quotient by the kernel %(or equivalently the image) 
of this operator is an irreducible $\Y(\so_n)\ts$-module.
But the definition \eqref{1.33} of the comultiplication on $\Y(\gl_n)$
implies that the vector \eqref{ourvec2}
is highest relative to the twisted Yangian $\Y(\so_n)\ts$. 
Here we also use Lemma \ref{higi}.
Further suppose that 
\begin{align}
\label{chercon41}
\nu_a\ge\nu_b
&\quad\text{whenever}\quad
\la_a-\la_b+\rho_a-\rho_b=0
\quad\text{and}\quad
a<b\ts;
\\
\label{chercon42}
\nu_a+\nu_b\le n\ 
&\quad\text{whenever}\quad
\la_a+\la_b+\rho_a+\rho_b=0
\quad\text{and}\quad
a<b\ts.
\end{align}
We have already noted that for $\si=\si_0$ 
and $\om=\tau_m\,\si_0$ our operator \eqref{taint}
is invertible.
Due to the nonsingularity of the weight $\la+\rho\ts$,
Proposition \ref{norm2} now implies that the image
of the vector \eqref{ourvec2} in the quotient
is not zero. Hence this image is highest
relative to the action of the twisted Yangian $\Y(\so_n)$
on the quotient. Lemma \ref{higi} implies that each of
the polynomials $Q_1(x)\lcd Q_l(x)$ of the quotient
is multiplicative with respect to the $m$ tensor factors
of the vector \eqref{ourvec2}. But in the
case $m=1$ these polynomials are transparent
from the relations given in the very end of Subsection 5.2.
Thus we get

\begin{theorem*}
\label{drin4}
Put $n=2l+1\ts$.
Let $\la_1\lcd\la_m$ satisfy %the conditions
\eqref{domcon41},\eqref{domcon42} while the labels
\begin{equation}
\label{nu4}
\textstyle
\nu_1=\la_1-\mu_1+l+\frac12\,
\,,\,\ldots\,,\,
\nu_m=\la_m-\mu_m+l+\frac12
\end{equation}
belong to the set\/
$\{0\com1\lcd n\}$ and satisfy %the conditions
\eqref{chercon41},\eqref{chercon42}.
Then the quotient by the kernel of our intertwining operator
\eqref{ak2} is a non-zero irreducible $\Y(\so_n)\ts$-module.
For any\/ $k=1\lcd l$ the polynomial $Q_k(x)$ of this module 
is the product of the sums $x+\mu_a+\rho_a$ taken over 
all indices\/ $a$ such that $\nu_a=k\ts$, and
of the differences\/ $x-\mu_a-\rho_a$ taken over 
all indices\/ $a$ such that $\nu_a=n-k\ts$.
\end{theorem*}

Note that if $\nu_1=\ldots=\nu_m=0\ts$, then
both the source and the target $\Y(\so_n)\ts$-modules
in \eqref{ak2} are trivial
for any $\mu\ts$. If moreover
$\la+\rho$ is nonsingular, then by Proposition \ref{norm2} 
our operator \eqref{ak2} is the identity map $\CC\to\CC\ts$.

Now let us take any sequence of $l$ monic polynomials 
$Q_1(x)\lcd Q_l(x)$ with complex coefficients.
Let $m$ be the sum of the degrees of $Q_1(x)\lcd Q_l(x)\ts$.
Suppose that 
%at least one of all the $l$ polynomials is not trivial, so that 
$m>0\ts$. Let $\la$ and $\mu$ be any weights of $\so_{2m}$
such that the corresponding labels \eqref{nu4}
belong to the set $\{1\lcd n-1\}$ and such that
the given polynomials $Q_1(x)\lcd Q_l(x)$
are obtained from $\la$ and $\mu$ as in Theorem \ref{drin4}.
We do not yet impose any other conditions on $\la$~and~$\mu\ts$.
%such as the nonsingularity of $\la+\rho\ts$.

Then $\la$ and $\mu$ are determined up to permuting 
the $m$ pairs $(\ts\nu_a\com \mu_a+\rho_a\ts)\ts$,
and up to replacing $(\ts\nu_a\com \mu_a+\rho_a\ts)$
by $(\ts n-\nu_a\com-\mu_a-\rho_a\ts)$ for any number
of indices $a\ts$.
Equivalenly, $\la$ and $\mu$ are determined up to a
permutation of the $m$ pairs $(\ts\la_a+\rho_a\com\mu_a+\rho_a\ts)\ts$,
and up to replacing $(\ts\la_a+\rho_a\com \mu_a+\rho_a\ts)$
by $(-\la_a-\rho_a\com-\mu_a-\rho_a\ts)$ for any number
of indices $a\ts$.
Thus $\la$ and $\mu$ are determined up to
the (simultaneous) shifted action of the group 
$\ES$ of $\so_{2m}$ on them,
as on elements of $\h^*$. Using only the action
of the subgroup $\Sym\subset\ES\ts$,
we can choose $\la$ to satisfy the conditions 
\eqref{domcon41},\eqref{domcon42}
so that the weight $\la+\rho$ is nonsingular. 

The equality on the right hand
side of \eqref{chercon41} for $a<b$
means that $\la$ is invariant
under the shifted action of that element of $\Sym$
which exchanges $\eta_{m-a+1}$ with $\eta_{m-b+1}\ts$, and
leaves all other basis vectors of $\h^*$ fixed.
This action amounts to exchanging the~pair
$(\ts\nu_a\com \mu_a+\rho_a\ts)$ with
$(\ts\nu_b\com \mu_b+\rho_b\ts)\ts$.
By using this action and keeping $\la$ fixed, 
we can choose $\mu$ so that the condition \eqref{chercon41} is satisfied.

The equality on the right hand
side of \eqref{chercon42} for $a<b$ means that the weight 
$\la$ is invariant under the shifted action of that element of 
the group $\Sym$
which maps $\eta_{m-a+1}$ and $\eta_{m-b+1}$
respectively to $-\ts\eta_{m-b+1}$ and $-\ts\eta_{m-a+1}\ts$,
leaving all other basis vectors of $\h^*$ fixed.
This action amounts to replacing the pairs
$(\ts\nu_a\com \mu_a+\rho_a\ts)$ and 
$(\ts\nu_b\com \mu_b+\rho_b\ts)$ by the pairs
$(\ts n-\nu_b\com-\mu_b-\rho_b\ts)$ and 
$(\ts n-\nu_a\com-\mu_a-\rho_a\ts)$ respectively.
By using this action and keeping $\la$ fixed, 
we can choose the weight $\mu$ 
so that the condition \eqref{chercon42} is satisfied.

Now consider our intertwining 
operator \eqref{ak2}
corresponding to these $\la$ and $\mu\ts$.
Due to Theorem~\ref{drin4}, the quotient by the kernel of this 
operator is a non-zero irreducible integrable $\Y(\so_n)\ts$-module,
and has the given polynomials $Q_1(x)\lcd Q_l(x)\ts$.
Thus up to equivalence and similarity, 
every non-zero irreducible integrable finite-dimensional 
$\Y(\so_n)\ts$-module arises as such a quotient.
Note that the choice $\th=-1$ here is essential.

%------------------------------------------------------------------------------

\bigskip\noindent
{\it 5.5.\ Irreducible representations of\/ $\Y(\so_n)$ for even $n$}

\vspace{10pt}\noindent
In this subsection, we will consider the case of $\gd=\so_n$ 
where $n$ is even. Hence $n=2\tts l$ where $l$ is a positive integer.
We keep assuming that $\th=-1\ts$, so that $\g=\so_{2m}\ts$.
Let $\Psi$ be any $\Y(\so_n)\ts$-module.
Let $\ups\in\SO_n$ be the element which exchanges the basis vector
$f_{n-1}$ with $f_n\ts$, and leaves all other basis vectors of $\CC^n$ fixed.
Consider the corresponding automorphism \eqref{upsaus} of the algebra
$\Y(\so_n)\ts$. Our element $\ups$ is involutive, and so is
the corresponding automorphism of $\Y(\so_n)\ts$.
Denote by $\Psicon$ the $\Y(\so_n)\ts$-module
obtained by pulling the action of $\Y(\so_n)$ on $\Psi$ back through
this automorphism.
Any $\Y(\so_n)\ts$-module equivalent to
$\Psicon$ will be called \textit{conjugate} to $\Psi\ts$.

The $\Y(\so_n)\ts$-modules $\Psi$ and $\Psicon$ may be equivalent or not.
Suppose that the $\Y(\so_n)\ts$-module $\Psi$ is irreducible, 
finite-dimensional,
integrable and non-zero. So is the $\Y(\so_n)\ts$-module $\Psicon$ then.
%Let $\psi\in\Psi$ be a highest vector.
By \cite[Theorem 4.4.14]{M3} for a highest vector $\psi$
of at least one of the two modules $\Psi$ and $\Psicon$
we have for $k=1\lcd l-1$
\begin{equation}
\label{qk}
\textstyle 
S_{2k+2,2k+2}(x)\,S_{2k,2k}(x)^{\tts-1}\,\psi
\,=\,
Q_k(x+\frac12)\,Q_k(x-\frac12)^{-1}\,\psi
\end{equation}
where $Q_k(x)$ is a monic polynomial in $x$ with coefficients in 
$\CC\ts$. Further, we have
\begin{equation}
\label{ql}
\hspace{26pt}
\textstyle 
S_{nn}(-\ts x)\,S_{nn}(x)^{\tts-1}\,\psi
\,=\,
Q_l(x+\frac12)\,Q_l(x-\frac12)^{-1}\,\psi
\end{equation}
where $Q_l(x)$ is an even monic polynomial in $x$ with coefficients in
$\CC\ts$.
Any sequence of $l$ monic polynomials with complex coefficients 
arises in this way, provided that the last polynomial in the sequence is even.
Moreover, 
the two $\Y(\so_n)\ts$-modules $\Psi$ and $\Psicon$ are equivalent,
if and only if zero is not a root of the 
corresponding polynomial $Q_l(x)\ts$. 

Let us consider the case when $\Psi$ and $\Psicon$ are not equivalent,
so that zero is a root of the polynomial $Q_l(x)\ts$. Let $h$
be the positive integer such that $0\com1\com\ldots\com h-1$
are roots of $Q_l(x)\ts$, but $h$ is not. The above vector $\psi$
has been a highest vector of one of the two
$\Y(\so_n)\ts$-modules $\Psi$ and $\Psicon\ts$.
Let $\psicon$ be a highest vector of the other of the two.
For each index $i\succcurlyeq n$ we have
$
S_{ii}(x)\,\psi=g_i(x)\,\psi
$
where $g_i(x)$ is a formal power series in $x^{-1}$ 
with the coefficients from $\CC\ts$. The proof of
\cite[Theorem 4.4.14]{M3} demonstrates that then for each $i\succ n$
we also have the equality
$
S_{ii}(x)\,\psicon=g_i(x)\,\psicon\ts,
$
while
$$
\textstyle
S_{nn}(x)\,\psicon\ts=g(x)\,g_n(x)\,\psicon
$$
where 
$$
\textstyle
g(x)=(x+h+\frac12)\,(x-h+\frac12)^{\ts-1}\ts.
$$
By the definition of the polynomial $Q_l(x)$ here we have 
$$
\textstyle
g_n(-x)\,g_n(x)^{-1}=Q_l(x+\frac12)\,Q_l(x-\frac12)^{-1}\ts.
$$
However,
$$
\textstyle
g(-x)\,g(x)^{-1}\neq
Q(x+\frac12)\,Q(x-\frac12)^{-1}
$$
for any polynomial $Q(x)\ts$, because the integer $h$ is positive.
This implies that only one of the two non-equivalent $\Y(\so_n)\ts$-modules
$\Psi$ and $\Psicon$ gives rise
to $l$ polynomials by using its highest vector $\psi\ts$,
as in \eqref{qk} and \eqref{ql}.
If $\psi$ is a highest vector of $\Psi\ts$, 
assign to $\Psi$ a label $\de=1\ts$.
Otherwise, that is if $\psi$ is a highest vector of $\Psicon\ts$,
assign to $\Psi$ a label $\de=-1\ts$.
The polynomials $Q_1(x)\lcd Q_l(x)$ will be from now on
associated to both $\Psi$~and~$\Psicon\ts$. 

Thus to every non-zero irreducible integrable finite-dimensional
$\Y(\so_n)\ts$-module $\Psi$ we have associated a sequence
of monic polynomials $Q_1(x)\lcd Q_l(x)$ where the last 
polynomial is even. 
If zero is a root of the polynomial $Q_l(x)\ts$,
that is if $\Psi$ is not equivalent to $\Psicon\ts$,
then we also have associated to $\Psi$ a label
$\de\in\{+1\com-1\}\ts$. The modules $\Psi$ 
are parametrized by their polynomials $Q_1(x)\lcd Q_l(x)$ 
and by their labels $\de$ (where the latter exist) 
up to equivalence and similarity \cite[Corollary 4.4.17]{M3}. 
For example, if $\Psi$ is the trivial $\Y(\so_n)\ts$-module
then $Q_1(x)=\ldots=Q_l(x)=1\ts$, and there is no label $\de\ts$.

In this subsection $\g=\so_{2m}\ts$, as
in Subsection 5.4. In particular, here $\rho_a=1-a$ 
for each index $a=1\lcd m\ts$. 
Suppose that the weight $\la+\rho$ is nonsingular.
This condition can be written as the collection of
inequalities \eqref{domcon41} and \eqref{domcon42}.
However, now $\ka_a=-\,l\ts$. 
Using the definition
\eqref{lamuka} of the weight $\nu\ts$, 
suppose that each label $\nu_a\in\{0\com1\lcd n\}\ts$.

In the end of Subsection 4.4
we observed that for $\g=\so_{\tts2m}$
our intertwining operators \eqref{ak2} and \eqref{4.56} are the same.
Denote by $N$ the quotient by the kernel of this operator.
By Proposition \ref{propo2} and Corollary~\ref{iror},
$N$ is either an irreducible $\Y(\so_n)\ts$-module
or splits into a direct sum of two irreducible non-equivalent
$\Y(\so_n)\ts$-modules. In the latter case, our $N$
is irreducible as a module over the %crossed product 
algebra $\SO_n\ltimes\Y(\so_n)$ by Proposition~\ref{proposition4.8}.
%Here we also use the $\SO_n\ts$-invariance 
%of the series $O(x)-1$ in $x^{-1}$ whose coefficients
%generate the kernel of the homomorphism 
%$\X(\so_n)\to\Y(\so_n)\ts$; see \eqref{5.2} and \eqref{xy}.

Due to the $\SO_n\ts$-equivariance of the Olshanski homomorphism
$\X(\so_n)\to\Ar^{\SO_{2m}}$, for any $\ups\in\SO_n$
pulling the action of $\Y(\so_n)$ on $N$
back through the automorphism \eqref{upsaus}
amounts to pulling that action forward through
the automorphism of $X\mapsto \ups\,X\,\ups^{-1}$ of 
the algebra $\End N\ts$. Let us apply this observation
to the element $\ups$ used in the beginning of this subsection.
If the $\Y(\so_n)\ts$-module $N$ is irreducible,
then it is of the form $\Psi\ts$,
where $\Psi$ and $\Psicon$ are equivalent. 
If $N$ splits to a direct sum of two irreducible non-equivalent
submodules and one of them is denoted by $\Psi\ts$,
the other submodule is equivalent to $\Psicon\ts$. 

Like in Subsection 5.4, 
the definition \eqref{1.33} of the comultiplication on $\Y(\gl_n)$
implies that the vector \eqref{ourvec2}
is highest for the twisted Yangian $\Y(\so_n)\ts$. 
Here we again used Lemma \ref{higi}.
Suppose that the weight $\nu$ satisfies the inequalities
\eqref{chercon41} and \eqref{chercon42}.
Note that the weights $\ka$
used here and in Subsection 5.4
to determine $\nu$ are different.
It is the form of the inequalities
\eqref{chercon41} and \eqref{chercon42} for $\nu$ that is the same.

We have already noted that for $\si=\si_0$ 
and $\om=\tau_m\,\si_0$ our operator \eqref{taint}
is invertible.
Due to the nonsingularity of the weight $\la+\rho\ts$,
Proposition \ref{norm2} now implies that the image
of the vector \eqref{ourvec2} in the quotient $N$
is not zero. This image is going to be our vector $\psi\ts$.
In particular, it will satisfy \eqref{qk} and \eqref{ql}
for certain polynomials $Q_1(x)\lcd Q_l(x)\ts$.
Lemma \ref{higi} implies that each of
these polynomials
is multiplicative with respect to the $m$ tensor factors
of the vector \eqref{ourvec2}. But in the
case $m=1$ these polynomials are transparent
from the relations given in the very end of Subsection 5.2.
Thus we get

\begin{theorem*}
\label{drin5}
Put $n=2\ts l$.
Let $\la_1\lcd\la_m$ satisfy %the conditions
\eqref{domcon41},\eqref{domcon42} while the labels
\begin{equation}
\label{nu5}
\textstyle
\nu_1=\la_1-\mu_1+l\,
\,,\,\ldots\,,\,
\nu_m=\la_m-\mu_m+l
\end{equation}
belong to the set\/
$\{0\com1\lcd n\}$ and satisfy %the conditions
\eqref{chercon41},\eqref{chercon42}.
Then the quotient by the kernel of our intertwining operator
\eqref{ak2} is either a non-zero irreducible self-conjugate 
$\Y(\so_n)\ts$-module, 
or splits into a direct sum of two non-equivalent non-zero
irreducible $\Y(\so_n)\ts$-modules conjugate to each other.
For every\/ $k=1\lcd l-1$
polynomial\/ $Q_k(x)$ of any of the irreducible modules
is the product of the sums\/ $x+\mu_a+\rho_a$ taken over 
all indices\/ $a$ such that $\nu_a=k\ts$, and
of the differences\/ $x-\mu_a-\rho_a$ taken over 
all indices\/ $a$ such that $\nu_a=n-k\ts$.
The polynomial\/ $Q_l(x)$ of any of the irreducible modules
is the product of the differences\/ 
$x^{\ts2}-(\mu_a+\rho_a){}^{\tts2}$ taken over 
all indices\/ $a$ such that $\nu_a=l\ts$.
The splitting occurs if and only if\/ $\mu_a+\rho_a=0$ for 
at least one index $a$ such that $\nu_a=l\ts$.
\end{theorem*}

In view of \eqref{nu5},
the splitting in Theorem \ref{drin5}
occurs if and only if 
$\la_a+\rho_a=0$ and $\mu_a+\rho_a=0$ 
simultaneously for at least one index $a\ts$.  
If $\nu_1=\ldots=\nu_m=0$ then
both the source and the target $\Y(\so_n)\ts$-modules
in \eqref{ak2} are trivial
for any $\mu\ts$. If moreover
$\la+\rho$ is nonsingular, then by Proposition \ref{norm2} 
our operator \eqref{ak2} is the identity map $\CC\to\CC\ts$.

Now let us take any sequence of $\ts l$ monic polynomials 
$Q_1(x)\lcd Q_l(x)$ with complex coefficients,
such that the polynomial $Q_l(x)$ is even. 
Let $m$ be the sum of the degrees of $Q_1(x)\lcd Q_{l-1}(x)$
plus half of the degree of $Q_l(x)\ts$.
Suppose that 
$m>0\ts$. Let $\la$ and $\mu$ be any weights of $\so_{2m}$
such that the corresponding labels \eqref{nu5}
belong to the set $\{1\lcd n-1\}$ and such that
the given polynomials $Q_1(x)\lcd Q_l(x)$
are obtained from $\la$ and $\mu$ as in Theorem \ref{drin5}.
We do not yet impose any other conditions on $\la$~and~$\mu\ts$.

Then $\la$ and $\mu$ are determined up to permuting 
the $m$ pairs $(\ts\nu_a\com \mu_a+\rho_a\ts)\ts$,
and up to replacing $(\ts\nu_a\com \mu_a+\rho_a\ts)$
by $(\ts n-\nu_a\com-\mu_a-\rho_a\ts)$ for any number
of indices $a\ts$.
Equivalenly, $\la$ and $\mu$ are determined up to a
permutation of the $m$ pairs $(\ts\la_a+\rho_a\com\mu_a+\rho_a\ts)\ts$,
and up to replacing $(\ts\la_a+\rho_a\com \mu_a+\rho_a\ts)$
by $(-\la_a-\rho_a\com-\mu_a-\rho_a\ts)$ for any number
of indices $a\ts$.
Thus $\la$ and $\mu$ are determined up to
the (simultaneous) shifted action of the group 
$\ES$ of $\so_{2m}$ on them,
as on elements of $\h^*$. Using only the action
of the subgroup $\Sym\subset\ES\ts$,
we can choose $\la$ to satisfy the conditions 
\eqref{domcon41},\eqref{domcon42}
so that the weight $\la+\rho$ is nonsingular. 

By using only the shifted action on $\mu$ of those 
elements of the subgroup $\Sym\subset\ES$ which leave $\la$
invariant, we can choose $\mu$ 
so that the conditions \eqref{chercon41},\eqref{chercon42} are
satisfied. The arguments are the same as in the end of
Subsection~5.4, and we do not repeat~them~here.

Now consider our intertwining 
operator \eqref{ak2}
corresponding to these $\la$ and $\mu\ts$.
Due to Theorem~\ref{drin5}, the quotient by the kernel of this 
operator is either a non-zero irreducible integrable
$\Y(\so_n)\ts$-module, or splits into a direct sum of two 
non-equivalent irreducible integrable
$\Y(\so_n)\ts$-modules. To any of the irreducible modules 
we associate the given polynomials $Q_1(x)\lcd Q_l(x)\ts$.
It the quotient is irreducible,
there is no label $\de$ associated to it.
If the quotient splits into two irreducible modules,
they have the labels $\de=1$ and $\de=-1$ associated to them. 
Therefore up to equivalence and similarity, 
every non-zero irreducible integrable finite-dimensional 
$\Y(\so_n)\ts$-module arises either as such a quotient,
or as one of its two direct summands.
Note that the choice $\th=-1$ here is essential.

%==============================================================================

\newpage%%%%%%%%%%%%%%%%%%%%%%%%%%%%%%%%%%%%%%%%%%%%%%%%%%%%%%%%%%%%%%%%%%%%%%%

\section*{Acknowledgments}

We are grateful to Ernest Vinberg 
for collaborating with us on \cite{KNV}. 
The present work began when we visited MPIM,
and continued when we visited IHES. 
We are grateful to the staff of both institutes for
their kind help and generous hospitality.
The first named author was supported by the RFBR grant 
08-01-00392,
joint grant 09-01-93106,
interdisciplinary grant 09-01-12185-ofi-m,
and the grant for Support of Scientific Schools 3036-2008-2. 
The second named author was supported by the EPSRC grant C511166.

%==============================================================================

%==============================================================================

\end{document}